\documentclass[11pt,twoside,letterpaper]{book}

\usepackage{amsmath,amsthm}
\usepackage{amsfonts}
\usepackage{bbm}
\usepackage{enumitem}
\usepackage{appendix}
\usepackage{verbatim}
\usepackage{latexsym}
\usepackage{graphicx}
\usepackage{makeidx}
\usepackage{times,fancyhdr} 
\usepackage{amsbsy}
\usepackage{amssymb}
\usepackage{longtable}

\usepackage{tikz}

\makeindex

\newcounter{ExampCount}[section]
\numberwithin{ExampCount}{section} \numberwithin{figure}{section}

\newcommand{\R}{\mathbb R}

\newcommand{\F}{{\cal F}}

\newtheorem{theorem}{Theorem}[chapter]
\newtheorem{corollary}{Corollary}[chapter]
\newtheorem{lemma}{Lemma}[chapter]

\newtheorem{definition}{Definition}[chapter]
\newtheorem{example}{Example}[chapter]
\newtheorem{remark}{Remark}[chapter]

\newtheorem*{dov}{Proof}

\makeatletter
\renewcommand{\section}{\@startsection{section}{1}{\parindent}{3.5ex plus 1ex minus .2ex}{2.3ex plus .2ex}{\normalfont\large\bfseries\S\;\,.\!\!\!\!\!\!\!\!\!\!\!}}
\renewcommand{\subsection}{\@startsection{subsection}{2}{\parindent}{-3.25ex plus -1ex minus-.2ex}{1.5ex plus .2ex}{\normalfont\bfseries}}

\makeatother

\theoremstyle{remark}

\theoremstyle{definition}

\begin{document}

\pagestyle{plain}

\begin{titlepage}
  \begin{center}
  \vspace*{5cm}
  \huge
  \bigskip
  \textbf{\scshape\bfseries Quasi-Banach spaces \\ of \\ random variables \\ and \\ stochastic processes}
  \bigskip
  \vfill

\vspace{.2in}

{\Large\scshape\bfseries
Yuriy {Kozachenko}
\\
\smallskip
\smallskip
{Yuriy {Mlavets}
\\
\smallskip
{Oleksandr {Mokliachuk}}}
\smallskip
\smallskip
}

\end{center}

\end{titlepage}

\noindent \textbf{\large{Abstract}} \hspace{2pt}

Recently, the theory of random variables and stochastic processes has been developed in Orlicz spaces, $Sub_\varphi(\Omega)$ spaces, $V(\varphi,\psi)$  spaces.
This book is devoted to the theory of quasi-Banach $K_\sigma$-spaces $\mathbf{F}_\psi(\Omega)$, $\mathbf{F}_\psi^*(\Omega)$ and $D_{V,W}(\Omega)$ spaces of random variables and stochastic processes.
The book consists of ten chapters.
The first chapter examines stochastic processes from the quasi-Banach $K_{\sigma}$-spaces of random variables.
The second chapter is dedicated to studying the fundamental properties of the spaces $\mathbf{F}_\psi(\Omega)$.
The third chapter provides estimates of the distribution of suprema of stochastic processes from the space $\mathbf{F}_\psi^*(\Omega)$.
The fourth chapter studies approximation of stochastic processes from the space  $SF_\psi(\Omega)$.
The fifth  chapter explores the main properties of Orlicz spaces and examines their connection with the spaces $\mathbf{F}_\psi(\Omega)$.
The sixth  chapter addresses the pre-Banach $K_\sigma$-spaces $D_{V,W}(\Omega)$ and investigates the essential properties of these spaces.
The seventh chapter evaluates the accuracy and reliability of constructing models of stochastic processes from $D_{V,W}(\Omega)$ within this space.
In the eighth chapter, estimates for the norm distributions in $L_p(T)$ of stochastic processes from the spaces $\mathbf{F}_\psi(\Omega)$ are provided.
The ninth chapter examines the Monte Carlo method for calculating multiple integrals defined on $\mathbb R^n$ with specified reliability and accuracy.
The tenth chapter evaluates the modeling of stochastic processes from the spaces $Sub_\varphi(\Omega)$,
which are subclasses of $K_\sigma$-spaces, with given reliability and accuracy in the spaces $L_p(T)$.
The eleventh chapter assesses the modeling of stochastic processes from the spaces $Sub_\varphi(\Omega)$, which are subclasses of $K_\sigma$-spaces, with specified reliability and accuracy in the spaces $C(T)$.
The material is substantially based on the results obtained by the authors and their co-authors.
\\

\noindent\textbf{{Keywords}} \hspace{2pt}
 Models of stochastic processes, quasi-Banach spaces, distribution of suprema, Karhunen-Lo{\`e}ve model, Monte Carlo method, multiple integrals, reliability and accuracy, sub-Gaussian random processes, Orlicz spaces.
\\


\tableofcontents

\chapter*{Introduction}
\addcontentsline{toc}{chapter}{Preface}

\medskip

One of the current challenges in the theory of stochastic processes is the study of various classes of stochastic processes,
the investigation of their general properties, and the derivation of distribution estimates for functionals of processes from different classes of random variables.

Recently, the theory of stochastic processes has been developed in Orlicz spaces, $Sub_\varphi(\Omega)$ spaces, and $D_{V,W}(\Omega)$ spaces.
As a rule, it is difficult or impossible to find norms in these spaces. Therefore, instead of norms, equivalent moment norms are used in these spaces.
When estimating distributions of various functionals of stochastic processes from these spaces,
the estimates slightly deteriorated, which led to the need to consider Banach spaces with moment norms.
This approach makes it possible, in many cases, to achieve better estimates than those obtained using, for example, the norms from Orlicz spaces.

The first result dedicated to the study of local properties of stochastic processes belongs to A. N. Kolmogorov \cite{Kolm}.
His theorem on sample continuity with probability one of stochastic processes was published in the work of E. E. Slutsky (1937) \cite{Slut37}.
This theorem laid the foundation for an entire direction in the theory of stochastic processes.
Other works in this field concerning general conditions of sample continuity and membership in Lipschitz classes of random fields
include the work by M. I. Yadrenko (1967) \cite{jadr67} and his monograph (1980) \cite{jadr80}.

In the works of Yu. V. Kozachenko and M. I. Yadrenko (1976) \cite{koz76,koz761},
local properties of sample functions from
various classes of random fields were investigated.
Conditions for the continuity of random functions on compact sets in Hilbert spaces were studied in the works of A. V. Skorokhod (1973) \cite{skor73} and M.~I.~Yadrenko (1968) \cite{jadr68}.

Many researchers studied properties of distributions of the suprema of random processes, explored problems related to the existence of moments and exponential moments of the distributions of process suprema. Significant attention has been given to the problem of estimating the probability $P\{\sup_{t \in T}|X(t)|>\varepsilon\}$.
Results of investigation of these problems were published in the books by
H.~Cram\'er and M. R. Leadbetter (1967) \cite{cram67}, M. B. Marcus and J. Pisier (1981) \cite{marc81}, M.~R.~Leadbetter, G. Lindgren, and H. Rootz\'en (1983) \cite{lead83}, R. D. Adler (1990) \cite{adle90}, M.~Ledoux and M. Talagrand (1991) \cite{ledo91}, as well as V. V. Buldygin and Yu.~V.~Kozachenko (1998) \cite{byld98}.
The study of properties of stochastic processes from certain classes of stochastic processes has been carried out by
D. M. P. Albin (1990, 1998) \cite{albi90, albi98}, M. Ledoux (1990) \cite{ledo90}, M. B. Marcus (1988) \cite{marc88}, and Ye.~I.~Ostrovskij (1990) \cite{ostr90}.

In the 1960s, studies of the local properties of Gaussian stochastic processes began.
Specifically, Yu. K. Belyaev obtained conditions for the continuity of Gaussian stationary processes in terms of spectral functions
and the well-known ``Belyaev's alternative'' (1960, 1961) \cite{bely60, bely61}.
Various methods were used to derive continuity conditions for Gaussian processes by
R. M. Dudley (1965) \cite{dudl65} and D. Delport (1964) \cite{delp64}.
Fundamental results concerning the properties of Gaussian stochastic processes were obtained by
D. Lindgren (1971) \cite{lind71}, R. M. Dudley (1973) \cite{dudl73}, C. Borell (1978) \cite{bore78}, M. Talagrand (1987) \cite{tala87}, R. D. Adler (1990) \cite{adle90}, M. Ledoux (1990) \cite{ledo90}, and V. I. Piterbarg (1996) \cite{pite96}.
The work by I. V. Bondarenko and O. V. Ivanov (1992) \cite{bond92} studied properties of sample functions of random fields with stable increments.

The evaluation of exponential moments and the distributions of the suprema of Gaussian stochastic processes has been the focus of works by
A. V. Skorokhod (1970)\cite{skor70}, G. Y. Landau and L. A. Shepp (1970) \cite{land70}, H. Fernique (1970) \cite{fern70},
N. N. Vakhania, V. I. Tarieladze, and S. A. Chobanyan (1985) \cite{vaha85}, M. Ledoux and M. Talagrand (1991) \cite{ledo91},
M. A. Lifshits (1995) \cite{lifs95}, and V. Yurinsky (1995) \cite{yuri95}.
The distribution of the suprema of Gaussian processes was investigated in the works of
Yu. K. Belyaev (1968) \cite{bely68}, S. Berman (1973) \cite{berm73}, V. I. Piterbarg (1982) \cite{pite82}, and D. I. Pickands (1969) \cite{pick69}.

In the 1960s, studies emerged that examined broader classes of random variables and stochastic processes beyond Gaussian ones.
The concept of a sub-Gaussian random variable was introduced by J. P. Kahane (1960) \cite{kaha60}.
Later, in the work of V.~V.~Buldygin and Yu. V. Kozachenko (1980) \cite{byld80},
it was proved that the space of sub-Gaussian random variables is a Banach space with respect to the sub-Gaussian standard.
The properties and various applications of sub-Gaussian and strictly sub-Gaussian random variables were studied in the works of
V.~V.~Buldygin (1977, 1980) \cite{byld77, byld801},
V.~V.~Buldygin and Yu. V. Kozachenko (1974, 1980, 1987, 1992, 1994, 1995, 1997, 1998) \cite{byld74, byld80, byld87, byld92, byld93, byld95, byld97, byld98},
V.~V.~Buldygin and S.~A.~Solntsev (1997) \cite{byld97},
A. R. Giuliano, Yu. V. Kozachenko and A.~M.~Tegza (2002) \cite{dzyl02, julantonini2002-2},
A. R. Giuliano, Yu. V. Kozachenko and T.~Nikitina (2003) \cite{giul03},
A. R. Giuliano et al (2013) \cite{giul2013},
E. I. Ostrovskij (1990) \cite{ostr90},
A. O. Pashko (1998) \cite{pash98},
N. K. Jain (1978) \cite{jain78},
J. P. Kahane (1960) \cite{kaha60},
M. Ledoux (1991) \cite{ledo91}, and
M. B. Marcus (1981) \cite{marc81}.

Yu. V. Kozachenko introduced the concept of sub-Gaussian stochastic processes in 1968 \cite{koz68}.
The properties of this class of stochastic processes were studied in the works of
V. V. Buldygin (1977) \cite{byld77},
as well as N. K. Jain and M. B. Marcus (1978) \cite{jain78}.

Exponential moments and estimates of the distributions of suprema for sub-Gaussian and related processes have been explored in the works of
J. P. Kahane (1968) \cite{kaha68}, N. K. Jain and M. B. Marcus (1975, 1978) \cite{jain75,jain78},
M. B. Marcus and J. Pisier (1981) \cite{marc81}, Ye. I. Ostrovskij (1990) \cite{ostr90}, D. Fukuda (1990) \cite{fuku90},
as well as M. Ledoux and M. Talagrand (1991) \cite{ledo91}.

Yu. V. Kozachenko and Ye. I. Ostrovskij (1985) \cite{koz85} introduced the concept of Banach spaces of the sub-Gaussian type,
specifically the $Sub_\varphi(\Omega)$ spaces of random variables and stochastic processes,
which generalize the spaces of sub-Gaussian random variables and stochastic processes.
The $\varphi$-sub-Gaussian spaces of random variables are spaces of centered random variables with specific growth of exponential moments.
The properties of such spaces, their estimates, and the conditions for the convergence of sums of independent random variables from these spaces,
where the process is defined on a space with a pseudometric induced by this process, were examined in the monograph by
V. V. Buldygin and Yu. V. Kozachenko (1998) \cite{byld98}.
The properties of $\varphi$-sub-Gaussian spaces were also studied in the work of
A. R. Giuliano, Yu. V. Kozachenko, and T. Nikitina (2003) \cite{giul03},
as well as in the monograph by
O. I. Vasylyk, Yu. V. Kozachenko, and R. Ye. Yamnenko (2008) \cite{vasyl2008}.
Estimates of the distributions of the suprema of $\varphi$-sub-Gaussian stochastic processes were investigated in the works of
Yu. V. Kozachenko, O. I. Vasylyk, and R. Ye. Yamnenko (2003) \cite{koz03}.

Pre-Gaussian stochastic processes were introduced in the works of V.~V.~Buldygin and Yu. V. Kozachenko (1974) \cite{byld74}.
The properties of pre-Gaussian stochastic processes and the distributions of their suprema were studied in the works of
V.~V.~Buldygin and Yu. V. Kozachenko (1992) \cite{byld92}, (1993) \cite{byld93}.

A broader class of random variables than Gaussian ones includes $\varphi$-sub-Gaussian and pre-Gaussian random variables from Orlicz spaces.
The theory of Orlicz functional spaces is thoroughly presented in the books by M. A. Krasnosel'ski{\u{\i}} and Ya. B. Rutitski{\u{\i}} (1958) \cite{kras58},
as well as M. M. Rao and Z. D. Ren (1991) \cite{raom91}, (2001) \cite{raom01}..
The theory of Orlicz spaces of random variables and the theory of stochastic processes from Orlicz spaces of random variables are detailed in the monograph by
V. V. Buldygin and Yu. V. Kozachenko (1998) \cite{byld98}.
Systematic research into the properties of Orlicz processes and processes with Orlicz increments began with the works of
Yu. V. Kozachenko (1983, 1984, 1984) \cite{koz83, koz84, koz841},
H.~Fernique (1983) \cite{fern83},
K. Nanopoulos and F. Nobelis (1978) \cite{nano78},
J. Pisier (1979, 1983) \cite{pisi79, pisi83},
and M. Weber (1983) \cite{webe81}.
Additionally, the properties of stochastic processes in certain Orlicz spaces were addressed in the works of
N. K\^{o}no (1980) \cite{kono80} and
Yu.~V.~Kozachenko (1985) \cite{koz851}.
The continuity conditions for the realizations of processes with Orlicz increments were discussed in the studies by
E.~A.~Abzhanov and Yu. V. Kozachenko (1985) \cite{abza85},
J. Pisier (1979, 1983) \cite{pisi79, pisi83},
and M.~Ledoux and M. Talagrand (1991) \cite{ledo91}.

Orlicz spaces of exponential type have been studied in the context of solving various problems in the theory of stochastic processes, notably in the works of
H.~Fernique (1975) \cite{fern75} and E. I. Ostrovskij (1982) \cite{ostr82}.
Estimates for the distributions of the norms of suprema of various Orlicz processes were examined in the studies
by H. Fernique (1983) \cite{fern83},
N. K\^{o}no (1980) \cite{kono80},
K. Nanopoulos and F.~Nobelis (1978) \cite{nano78},
as well as J. Pisier (1979, 1983) \cite{pisi79, pisi83}.

The space $\mathbf{F}_\psi(\Omega)$ was introduced by S. V. Ermakov and Ye. I. Ostrovskij in their 1986 work \cite{erma86},
where it was proved that this space is a Banach space with the norm
$\left\|\xi\right\|_\psi=\sup\limits_{u\geq1}\left(\frac{\left(\mathsf{E}\left|\xi\right|^{u}\right)^{1/u}}{\psi(u)}\right)$.

Conditions for the almost sure membership of stochastic processes from the Orlicz space of random variables to functional
Orlicz spaces and Sobolev-Orlicz spaces were established in the works of
Yu. V. Kozachenko and T. O. Yakovenko (2004, 2006) \cite{yako04, yako06}.

In the 1990s, a new direction in the theory of stochastic processes began to develop: wavelet analysis.
It turned out that wavelet analysis is a novel branch of mathematics with its own unique challenges and problems. Moreover,
this theory has proven to be highly effective in practical applications. For example, it is utilized for recording information,
sound, and images onto compact discs and computers. Recording and storing information using wavelets is significantly more efficient than other methods.
Wavelets are also widely used in information encoding.

Significant contributions to the development of the theory of wavelet analysis were made by scientists from Western Europe and North America, such as
S.~G.~Mallat (1998) \cite{mall98}, I. Meyer (1990) \cite{meye90}, I. Daubechies (1992) \cite{daub92}, and C.~K.~Chui (1992) \cite{chui92}.
Research on uniform convergence with probability one over a bounded interval of wavelet expansions of random processes from Orlicz spaces of random variables was conducted by
Yu. V. Kozachenko and M.~M.~Perestyuk (2007, 2008) \cite{pere06, pere07, pere08}.
See  the book by I. V. Darijchuk, Yu.~V.~Kozachenko and M. M. Perestyuk \cite{dari11} for more details.

In the works of Yu. V. Kozachenko and O. Ye. Kamenshchikova (2007-- 2014) \cite{kame07, kame071, koz-kam, kozKam-2014}, estimates of the deviation process are provided
for approximating stochastic processes using broken lines, Bernstein polynomials, and splines, along with the construction of
approximations of $SSub_{\varphi}(\Omega)$ random processes in the spaces $C[a,b]$ and $L_{p}(T)$ with specified accuracy and reliability.

As early as the 1930s, Enrico Fermi, and later John von Neumann and Stanislaw Ulam in the 1940s,
anticipated that the connection between stochastic processes and differential equations could be utilized "in reverse."
They proposed using a stochastic approach to approximate multidimensional integrals in transport equations that arose
in relation to problems of neutron movement in isotropic media. The idea was further developed by Stanislaw Ulam,
who suggested that instead of employing traditional combinatorial reasoning, one could simply conduct ``experiments'' numerous times and,
by counting the number of successful outcomes, estimate the probability of events. He also proposed using computers to perform calculations via the Monte Carlo method.

The advent of the first electronic computers, capable of rapidly generating pseudo-random numbers, greatly expanded the range of problems
for which the stochastic approach proved to be more effective than other mathematical methods.
The Monte Carlo method began to be applied to many tasks; however, its use was not always justified due to the large number
of computations required to obtain answers with the desired accuracy.

The year of birth of the Monte Carlo method is considered to be 1949, when the article by N. Metropolis and S. Ulam was published \cite{metr49}.

In the 1950s, the Monte Carlo method was employed for calculations during the development of the hydrogen bomb. The primary contributions to the advancement of the Monte Carlo method at that time were made by staff members of the U.S.~Air Force laboratories and the RAND Corporation.

In the 1960s, within the emerging field of computational mathematics, a class of problems was identified whose complexity (the number of computations needed to obtain an accurate solution) grows exponentially with the dimensionality of the problem. In some cases, sacrificing accuracy can allow the development of algorithms with slower growth in complexity; however, there exists a significant number of problems for which this is not feasible (e.g., determining the volume of a convex body in $n$-dimensional Euclidean space). For such problems, the Monte Carlo method remains the only viable option for obtaining sufficiently accurate solutions within an acceptable time frame.

The Monte Carlo method for calculating multidimensional integrals has been studied in numerous works, notably in the
 books by S. M. Ermakov (1975) \cite{erma75}
and S. M. Ermakov and G. A. Mikhailov (1976) \cite{erma76}.

The rapid advancement of computational technology has significantly driven the development of stochastic modeling methods, including numerical modeling of random processes.
These methods are widely applied in various fields of natural and social sciences, such as financial mathematics, meteorology, radio engineering, sociology, physics, and more.

Many specialists have contributed to the development of the theory of modeling random processes and fields.
This topic is covered in numerous works, such as those by
S. M. Ermakov and G. A. Mikhailov (1982) \cite{ermmikh1982},
A. S. Shalygin (1986) \cite{shaligin1985},
N. P. Buslenko (1978) \cite{busenko1978},
and B. D. Ripley (1987) \cite{ripley1987}.

The most thoroughly studied section of modeling is dedicated to methods for simulating Gaussian random processes and fields.
M. Y. Yadrenko and G. Rakhimov \cite{rakhimovyadremko1993}
describe the modeling of isotropic and homogeneous random fields on a plane.
M. Y. Yadrenko and Z. O. Vyzhva \cite{vyzhvayadrenko2000, vyzhva2003, vyzhva2014, vyzhva2019}
investigated the statistical modeling of isotropic random fields on a sphere.

Traditional methods for modeling Gaussian random processes and fields--such as linear transformation methods,
autoregressive and moving average techniques, spectrum randomization, moving summation, canonical representations,
and others--are explored in the works of
V. V. Bykov (1978) \cite{bykov1978},
S. M. Ermakov and G. A. Mikhailov (1982) \cite{ermmikh1982},
A. S. Shalygin and Y. I. Palagin (1986) \cite{shaligin1985},
V. A. Ogorodnikov and S. M. Prigarin (1990) \cite{ogorodnik1996},
T. Sana and M. Chaiki (1997) \cite{sun1997}.
However, these methods often overlook the questions of accuracy and reliability of the constructed model.

Certain methods for constructing models that approximate a Gaussian random process with a given accuracy and reliability
in specific functional spaces were explored in the works of
I. N. Zelepugina and Yu. V. Kozachenko (1982, 1988) \cite{zelepugkoz1982,zelepugkoz1988},
as well as Yu. V. Kozachenko and L. F. Kozachenko (1991, 1992) \cite{kozkoz1991,kozkoz1992}.
The model of a Gaussian random process introduced in the works of G. O. Mikhailov was further examined in detail by
Yu. V. Kozachenko, A.~M.~Tegza and A.~R.~Giuliano \cite{dzyl02,julantonini2002-2,koztegza2002}.
This research analyzed the accuracy and reliability of the Gaussian process model in certain functional spaces.
 Yu.~V.~Kozachenko and A. O. Pashko studied sub-Gaussian random processes and methods for constructing models of
 these processes with specified accuracy and reliability in various Banach spaces in their works
 \cite{kozpash1998,kozpash1999,kozpash1999-2,kozpash2000,pashko2001}
 and their book \cite{koz-pash}.
 A. O. Pashko also investigated accuracy estimates in the uniform metric for modeling sub-Gaussian random
 fields and Gaussian random fields on a sphere \cite{pashko2000,pashko2001-2}.
 The modeling of stationary random processes from $Sub_\varphi(\Omega)$ was studied by
 Yu. V. Kozachenko and I. V. Rozora \cite{rozora2004, koz-roz2019}.
 Yu. V. Kozachenko, A. O. Pashko and I. V. Rozora \cite{koz-pash-roz2007}

One of the relatively new branches of the theory of random processes is the theory of random variables and processes in generalized Banach spaces.
The properties of random processes in pre-Banach and quasi-Banach spaces were studied in the works of
V. V. Buldygin and Yu. V. Kozachenko \cite{byld98}. Examples of quasi-Banach spaces were examined in the works of
Yu. V. Kozachenko and E. I. Ostrovskij \cite{koz85},
including Orlicz spaces, $D(\Omega)$-spaces, and $\Gamma_\psi(\Omega)$-spaces.
Functional $K_\sigma$-spaces were considered by
L. V. Kantorovich and G. P. Akilov \cite{kant84}.
Quasi-$K_\sigma$ spaces and $K_\sigma$ spaces of random variables were studied by
Ye. A. Abzhanov and Yu.~V.~Kozachenko \cite{abza85, abza86},
Yu.~V.~Kozachenko and Yu.~Yu.~Mlavets \cite{mlav111, mlav11, mlav12, mlav143, mlav2015, mlav2016, mlav2018, mlav121, mlav122, mlav141, mlav142, mlav144}.

The modeling of random processes is widely applied in geological sciences. The development and creation of
new approaches to mathematical modeling in geology are the focus of the works by
S. A. Vyzhva and H. T. Prodayvoda \cite{prodayvoda1999, prodayvoda2009}.

Expansions of random processes into series with independent components play a crucial role in the theory of modeling random processes.
One such expansion is the Karhunen-Lo\'eve decomposition, which represents a random process as an infinite linear
combination of orthogonal functions defined by the process's covariance function \cite{koz-pash}.
This decomposition is an extension of the Karhunen-Lo\'eve transform, or principal component method, proposed by
K. Pearson in 1901 \cite{pearson1901}.
Yu.~V.~Kozachenko and I. V. Rozora \cite{koz-roz2019},
Yu.~V.~Kozachenko, I. V. Rozora, and Ye. V. Turchyn  \cite{koz-roz-turch},
O.~M.~Moklyachuk \cite{mokl2007, mokl2012, mokl2014, mokl2022}
studied the expansion of a random process into a series of orthogonal components based on functions determined by the process correlation function.
Problems related to modeling  of stochastic processes with given accuracy and reliability in various spaces ($C[0,T]$, ${L_p}([0,{T}])$ etc)
were investigated in the papers by I.~V.~Rozora \cite{rozora2004, rozora2009, rozora2015, rozora2018, rozora2018b, rozora2018c, rozora2020, rozora2022, rozora2024a, rozora2024b}.
Properties of generalized sub-{Gaussian} fractional {Brownian} motion queueing models were investigated by
Yu.~V.~Kozachenko, R. Ye. Yamnenko and D. Bushmitch \cite{yamn2014}.
Problems related to estimation of probability of exceeding a fixed level function by sub-{Gaussian} random processes were investigated by
Yu.~V.~Kozachenko, O. I. Vasylyk and R.~Ye.~Yamnenko \cite{koz03, vasyl2005},
R. Ye.Yamnenko and N. Yurchenko \cite{yamn2020}

Thus, the study of the accuracy and reliability of models of random processes and fields in various functional spaces is
highly relevant due to the widespread application of stochastic modeling of random processes and fields in various fields of natural and social sciences.

The book consists of ten chapters.

The first chapter examines random processes from the quasi-Banach $K_{\sigma}$-spaces of random variables.

The second chapter is dedicated to studying the fundamental properties of the spaces $\mathbf{F}_\psi(\Omega)$.

The third chapter provides estimates of the distribution of suprema of random processes  from the space $\mathbf{F}_\psi^*(\Omega)$.

The fourth chapter explores the main properties of Orlicz spaces and examines their connection with the spaces $\mathbf{F}_\psi(\Omega)$.

The fifth chapter addresses the pre-Banach $K_\sigma$-spaces $D_{V,W}(\Omega)$ and investigates the essential properties of these spaces.

The sixth chapter evaluates the accuracy and reliability of constructing models of processes from $D_{V,W}(\Omega)$ within this space.

In the seventh chapter, estimates for the norms' distributions in $L_p(T)$ of random processes from the spaces $\mathbf{F}_\psi(\Omega)$ are provided.

The eighth chapter examines the Monte Carlo method for calculating multiple integrals defined on $\R^n$ with specified reliability and accuracy.

The ninth chapter evaluates the modeling of random processes from the spaces $Sub_\varphi(\Omega)$,
which are subclasses of $K_\sigma$-spaces, with given reliability and accuracy in the spaces $L_p(T)$.

The tenth chapter assesses the modeling of random processes from the spaces $Sub_\varphi(\Omega)$,
which are subclasses of $K_\sigma$-spaces, with specified reliability and accuracy in the spaces $C(T)$.

\newpage

\chapter{Random processes from quasi-Banach $K_{\sigma}$-spaces of random variables}

\medskip

In the first chapter we consider random processes from quasi-Banach $K_{\sigma}$-spaces of random variables. The conditions of sample continuity with probability one and the conditions that the suprema of processes belong to the same space as the processes themselves are investigated. Estimates of the distributions of these processes are found.

\section{Quasi-Banach $K_{\sigma}$-spaces of random variables}

Let $\left\{\Omega,\F,P\right\}$ be a standard probability space, and let $K(\Omega)$ be a space of random variables $\xi=\xi(\omega), \omega\in\Omega$.

\begin{definition}
The function
$\Theta=(\Theta(\xi),\xi\in\mathcal{M})$ is called a prenorm
if, for all $\xi\in\mathcal{M}$,

1. $\Theta(\xi)\in[0,\infty)$;

2. $\Theta(0)=0$;

3. $\Theta(-\xi)=\Theta(\xi).$

\end{definition}

\begin{definition}
A space $\mathcal{M}$ that is complete with respect to the prenorm $\Theta$
will be called a pre-Banach space.
\end{definition}

\begin{definition}
\label{om-prebanach}
A pre-Banach space $\mathcal{M}$ is called a pre-$K_\sigma$-space if the following properties hold true:
\begin{itemize}
\item[$a_{1}$)] if $\xi,\eta$ $\in$ $\mathcal{M}$, then $max(\xi,\eta)\in \mathcal{M}$ and $min(\xi,\eta)\in \mathcal{M}$. Hence, and $|\xi|\in \mathcal{M}$.
\item[$a_{2}$)] if $\eta$ $\in$ $\mathcal{M}$ and $|\xi|\leq |\eta|$ almost
everywhere, then $|\xi|\in \mathcal{M}$.
\item[$a_{3}$] if for a sequence $\{\xi_n,n\geq1\}$ of random variables from $\mathcal{M}$ there exists a random {\nolinebreak variable $\eta \in \mathcal{M}$ such that $\sup\limits_{n\geq1}|\xi_n|\leq \eta$, then $\sup\limits_{n\geq1}|\xi_n|\in \mathcal{M}$.}
\end{itemize}
\end{definition}

\begin{definition}
\label{premetric}
A premetric is a function
$\rho(t,s)$, $t,s\in T$, such that $\rho(t,s)\in [0,\infty)$,
$\rho(t,t)=0$, $\rho(t,s)=\rho(s,t)$
\end{definition}

\begin{definition} A quasinorm on the space $K(\Omega)$ is called a functional $||\cdot||$ that associates a non-negative number $||\xi||$ to each random variable $\xi\in K(\Omega)$ such that the following conditions are satisfied:
\begin{itemize}
\item[1] $||\xi||=0\Leftrightarrow\xi=0$ with probability one;
\item[2] $||\xi+\eta||\leq||\xi||+||\eta||$;
\item[3] for $|\lambda|\leq 1$ $||\lambda\xi||\leq||\xi||$.
\end{itemize}
\end{definition}

\begin{remark} When instead of condition 3 we have the equality $||\lambda\xi||=|\lambda|||\xi||$, then $||\cdot||$ is the usual norm on $K(\Omega)$. The equality $\rho(\xi,\eta)=||\xi-\eta||$ defines a metric on $K(\Omega)$, so $K(\Omega)$ can be considered as a metric space.
\end{remark}

\begin{definition} Let $j=\left\{j(\lambda),\lambda\in\left[0,1\right]\right\}$ be a monotonically non-decreasing function such that $j(\lambda)\geq0$, $j(\lambda)\rightarrow0$, if $\lambda\rightarrow0$. If for a quasinorm $||\cdot||$ on $K(\Omega)$ the inequality $||\lambda\xi||\leq j(|\lambda|)||\xi||$ holds, if $|\lambda|\leq1$, then we will call this quasinorm subordinate to the function $j$.
\end{definition}

\begin{remark}
In the case where $||\cdot||$ is a norm, it obeys the function $j(|\lambda|)=|\lambda|$.
\end{remark}

\begin{definition}
A space $K(\Omega)$ that is complete with respect to the quasinorm $||\cdot||$ will be called a quasi-Banach space.
\end{definition}

\begin{definition}
A quasi-Banach space $K(\Omega)$ is called a quasi-$K_{\sigma}$-space if the following properties hold true:
\begin{itemize}
\item[$a_{1}$] if $\xi,\eta\in K(\Omega)$, then $\max(\xi,\eta)\in K(\Omega)$ and $\min(\xi,\eta)\in K(\Omega)$, hence $|\xi|\in K(\Omega)$;
\item[$a_{2}$] if $\xi,\eta\in K(\Omega)$ and $|\xi|\leq|\eta|$ almost everywhere, then $||\xi||<||\eta||$;
\item[$a_{3}$] if for a sequence $(\xi_{n}, n\geq1)$ of random variables from $K(\Omega)$ there exists a random variable $\eta\in K(\Omega)$ such that $\sup\limits_{n\geq1}|\xi_{n}|\leq \eta$, then $\sup\limits_{n\geq1}|\xi_{n}|\in K(\Omega)$.
\end{itemize}
If in this case $K(\Omega)$ is a Banach space, then we will call it a $K_{\sigma}$-space.
\end{definition}

\begin{remark}
If $K(\Omega)$ is a quasi-$K_{\sigma}$-space, then the third condition of the quasinorm is fulfilled automatically. Indeed, if $|\lambda|\leq1$, then $|\lambda\xi|\leq|\xi|$ and by property $a_{2})$ $||\lambda\xi||\leq||\xi||$.
\end{remark}

This remark is significant in the case when the quasinorm $||\cdot||$ is generated by some metric.

\begin{remark}
Functional $K_{\sigma}$-spaces were considered, for example, in the book \cite{kant84}, $K_{\sigma}$-spaces of random variables in the paper \cite{abza85}, and quasi-$K_{\sigma}$-spaces in the paper \cite{abza86}.
\end{remark}

\begin{definition} \label{de01}
A non-decreasing numerical sequence $\left(\varkappa(n),n\geq1\right)$ is called a majorizing characteristic of a quasi-$K_{\sigma}$-space $K(\Omega)$ if for any random variables $\xi_i\in K(\Omega),\ i=1,2,\ldots,n$ the following inequality holds:
\[
\left\|\max\limits_{1\leq i\leq n}\left|\xi_i\right|\right\|\leq\varkappa(n)\max\limits_{1\leq i\leq n}\left\|\xi_i\right\|.
\]
\end{definition}

The concept of a majorizing characteristic for Orlicz spaces was first introduced in the work \cite{koz84}, for $K_{\sigma}$-spaces -- in \cite{abza85}, and for quasi-$K_{\sigma}$-spaces -- in \cite{abza86}.

For quasi-$K_{\sigma}$-spaces, the following lemma holds.

\begin{lemma} \label{le01}
Let $\xi_{n},n=1,2,...$ be a sequence of random variables such that $\xi_{n}\in K(\Omega)$, where $K(\Omega)$ is a quasi-$K_{\sigma}$-space. If there exists a random variable $\xi_{n}\in K(\Omega)$ such that $||\xi_{n}-\xi||\rightarrow0$ as $n\rightarrow\infty$, then $\xi_{n}\rightarrow\xi$ by probability.
\end{lemma}

For $K_{\sigma}$-spaces this lemma was proved in the work \cite{abza85}. For quasi-$K_{\sigma}$-spaces the proof is similar.

\section{Orlicz spaces of random variables}

\begin{definition} \cite{byld98}
A continuous even convex function $U=\left\{U(x)\right., \linebreak \left. x\in \R\right\}$ is called
a $C$-function if $U(x)$ monotonically increases for $x>0$ and $U(0)=0$.
\end{definition}

\begin{example} The following functions are examples of Orlicz $C$-functions:
\begin{enumerate}
\item[1.] $U(x)=A|x|^\alpha,\; x\in \R,\; A>0,\alpha\geq 1$;
\item[2.] $U(x)=C(\exp\{B|x|^\beta\}-1),\; x\in \R,\;
C>0,B>0,\beta\ge 1$;
\item[3.] $U(x)=C(\exp\{\varphi(x)\}-1),\;
x\in \mathbb{R}, C>0$ and $\varphi(x)$, $x\in\mathbb{R}$, is an arbitrary $C$-function;
\item[4.] $U(x)=\left\{%
\begin{array}{ll}
 \left(\frac{e\alpha}{2}\right)^{\frac{2}{\alpha}}x^2, & \mbox{when}
 \quad
 |x|\le\left(\frac{2}{\alpha}\right)^{\frac{1}{\alpha}}; \\
 \exp\{|x|^\alpha\}, & \mbox{when} \quad
 |x|>\left(\frac{2}{\alpha}\right)^{\frac{1}{\alpha}};
\end{array}%
\right.$ $ 0<\alpha<1$;
\item[5.] $U(x)=D|x|^{\alpha}\ln(|x|+1),\;
x\in \R,\; D>0,\alpha\geq 1$.
\end{enumerate}
\end{example}

\begin{definition} \cite{byld98}
Let $U$ be an arbitrary $C$-function. An Orlicz space of random
variables $L_U(\Omega)$ is a family of random variables such that for each
$\xi\in L_U(\Omega)$ there exists a constant $r_\xi>0$ such that
\[
EU\left(\frac{\xi}{r_\xi}\right)<\infty.
\]
\end{definition}
An Orlicz space is a Banach space with the norm
\[
\|\xi\|_{U}=\inf\biggl\{r>0:
EU\left(\frac{\xi}{r}\right)\leq1\biggr\},
\]
which is called the Luxembourg norm.

\begin{lemma} \cite{byld98} \label{le:o4chapter-1}
Let $\xi \in L_U(\Omega)$ and $\left\|\xi\right\|_{U}>0$. Then for all $x>0$ the inequality holds
\begin{equation}\label{de02}
P\left\{\left|\xi\right|\geq x\right\}\leq\frac{1}{U\left(x/\left\|\xi\right\|_{U}\right)}.
\end{equation}
\end{lemma}

\section{Spaces $M(\Omega)$ of random variables}

Let $U(x), V(x), W(x)$ be $C$-functions such that for each $D>0$ the function $U\left(W(x)/D\right)$ is a $C$-function, and the function $W(x/D)/(1+V(x))$ is monotonically non-decreasing for $x>0$. These conditions are fulfilled, for example, when $U(x)=\left|x\right|^p$, $W(x)=\left|x\right|^q$, $V(x)=\left|x\right|^s$, $p\geq1$, $q\geq1$, $s\geq1$, $q\geq s$, or when $U(x)=\exp \left\{\left|x\right|^\alpha\right\}-1$, $\alpha\geq1$, $W(x)=\left|x\right|^q$, $q\geq1$, $V(x)=\left|x\right|^s$, $s\geq1$, $q\geq s$.

\begin{definition}
The space $M(\Omega)$ is the space of random variables $\xi$ such that for each of them there exists a constant $C_{\xi}$ such that
$$EU(W(\xi/C_{\xi})/(1+V(\xi)))<\infty.$$
\end{definition}

The spaces $M(\Omega)$ were introduced as a special case in the work \cite{abza86}.

\begin{theorem} The space $M(\Omega)$ is a quasi-Banach $K_{\sigma}$-space with respect to the quasi-norm
\begin{equation}\label{eq01}
||\xi||_{M}=\inf\left\{r>0: EU(W(\xi/r)/(1+V(\xi)))<1\right\},
\end{equation}
the major characteristic of this space is any of the functions $(x_{0}>0)$
\[
\mu(n)=(1+U(W(x_{0}))){S}^{W}_{x_{0}}({S}^{U}_{W(x_{0})}(n)),
\]
where ${S}^{g}_{n_{0}}(n)$ for any $z_{0}>0$ and any $C$-function is defined Yes:
\[
{S}^{g}_{z_{0}}(n)=\sup_{x\geq z_{0}}\frac{1}{x}{g}^{(-1)}(ng(x)).
\]
\end{theorem}

\begin{remark}
In particular cases, Theorem 1 was proved in the paper \cite{abza86}.
\end{remark}

\begin{dov} Let us first prove that $||\xi||_{M}=||\xi||$ is a quasinorm. Conditions 1 and 3 are obvious, as is condition 2, when $||\xi||=0$ and $||\eta||\neq0$.
If $||\xi||\neq0$, $||\eta||\neq0$, then condition 2 follows from the relations
\[
EU\left(\frac{1}{1+V(\xi+\eta)}W\left(\frac{\xi+\eta}{||\xi||+||\eta||}\right)\right)\leq
\]
\[
\leq EU\left(\frac{1}{1+V(|\xi|+|\eta|)}W\left(\frac{||\xi||}{||\xi||+||\eta||}\frac{|\xi|}{||\xi||}+\frac{|\eta|}{||\eta||}\frac{||\eta||}{||\xi||+||\eta||}\right)\right)\leq
\]
\[
\leq\frac{||\xi||}{||\xi||+||\eta||}EU\left(\frac{1}{1+V(|\xi|)}W\left(\frac{|\xi|}{||\xi||}\right)\right)+
\]
\[
+\frac{||\eta||}{||\xi||+||\eta||}EU\left(\frac{1}{1+V(|\eta|)}W\left(\frac{|\eta|}{||\eta||}\right)\right)\leq 1.
\]

The proof of the completeness of the space is no different from the case of Orlicz spaces. The fact that $M(\Omega)$ is a $K_{\sigma}$-space is obvious.

Let us prove the equality \eqref{eq01}. Let $\xi_{i}$, $i=1,...,n$ be random variables such that $\xi_{i}\in M(\Omega)$. Let $\theta=\max\limits_{i=1,...,n}|\xi_{i}|$, $a=\max\limits_{i=1,...,n}||\xi_{i}||$, $\chi(A)$ be the indicator of the set $A$. For any $r>0$, $x_{0}>0$, the relations hold
\[
EU\left(W(\theta{r}^{-1})\frac{1}{1+V(\theta)}\right)\leq EU\left(W(\theta{r}^{-1})\frac{1}{1+V(\theta)}\right)\chi(|\theta|{r}^{-1}\leq x_{0})+
\]
\[
+EU\left(W(\theta{r}^{-1})\frac{1}{1+V(\theta)}\right)\chi(|\theta|{r}^{-1}>x_{0})\leq EU\left(W(x_{0})\frac{1}{1+V(\theta)}\right)+
\]
\[
+\max_{1\leq i\leq n}nE\chi(|\xi_{i}|{r}^{-1}>x_{0})U\left(W(|\xi_{i}|{r}^{-1})\frac{1}{1+V(|\xi_{i}|)}\right)\leq
\]
\[
\leq U(W(x_{0}))+\max_{1\leq i\leq n}E\chi(|\xi_{i}|{r}^{-1}>x_{0})U\left({S}^{U}_{W(x_{0})}(n)\frac{W(|\xi_{i}|{r}^{-1})}{1+V(|\xi_{i}|)}\right)\leq
\]
\[
\leq U(W(x_{0}))+\max_{1\leq i\leq n}E\chi(|\xi_{i}|{r}^{-1}>x_{0})\times
\]
\[
\times U\left(W\left(|\xi_{i}|{r}^{-1}{S}^{W}_{x_{0}}\left({S}^{U}_{W(x_{0})}(n)\right)\right)\frac{1}{1+V(|\xi_{i}|)}\right)\leq
\]
\[
\leq U(W(x_{0}))+\max_{1\leq i\leq n}EU\left(W\left(|\xi_{i}|{r}^{-1}{S}^{W}_{x_{0}}\left({S}^{U}_{W(x_{0})}(n)\right)\right)\frac{1}{1+V(|\xi_{i}|)}\right).
\]

If we put $\theta=a{S}^{W}_{x_{0}}\left({S}^{U}_{W(x_{0})}(n)\right)$, then we get

$$EU\left(W(\theta{r}^{-1})\frac{1}{1+V(\theta)}\right)\leq U(W(x_{0}))+1.$$
That is, since the function $U(W(x){D}^{-1})$ is a $C$-function, then
\[
EU\left(W\left(\theta{r}^{-1}\left(U{(W(x_{0})+1)}^{-1}\right)\right)\frac{1}{1+V(\theta)}\right)\leq
\]
\[
\leq\frac{1}{U(W(x_{0}))+1}EU\left(W(\theta{r}^{-1})\frac{1}{1+V(\theta)}\right)\leq 1.
\]
Therefore, $||\theta||\leq a {S}^{W}_{x_{0}}\left({S}^{U}_{W(x_{0})}(n)\right)(1+U(W(x_{0}))).$
\end{dov}

\begin{lemma}
Let the random variable $\xi$ belong to the space $M(\Omega)$. Then for any $\epsilon>0$ the inequality holds
\[
P\left\{|\xi|>\varepsilon\right\}\leq{\left(U\left(W\left(\frac{\varepsilon}{||\xi||}\right)\frac{1}{1+V(\varepsilon)}\right)\right)}^{-1}.
\]
\end{lemma}

\begin{dov}
From the Chebyshev inequality and the properties of the functions $U(x)$, $V(x)$, $W(x)$ the following inequalities follow
\[
P\left\{|\xi|>\varepsilon\right\}\leq EU\left(W\left(\frac{\varepsilon}{||\xi||}\right)\frac{1}{1+V(\xi)}\right){\left(U\left(W\left(\frac{\varepsilon}{||\xi||}\right)\frac{1}{1+V(\varepsilon)}\right)\right)}^{-1}\leq
\]
\[
\leq{\left(U\left(W\left(\frac{\varepsilon}{||\xi||}\right)\frac{1}{1+V(\varepsilon)}\right)\right)}^{-1}.
\]
\end{dov}

\section{Random processes from quasi-Banach $K_{\sigma}$-spaces}

Let $T$ be a parametric set on which a pseudometric $\rho$ is given. Recall that a pseudometric differs from a metric only in that for a pseudometric the relation $\rho(t,s)=0$ does not necessarily imply that $t=s$.

We consider spaces with pseudometrics, not just metric spaces, because random processes from Banach or quasi-Banach spaces generate pseudometrics on parametric sets, not metrics. Note also that we will use only those properties of pseudometrics that do not differ from the properties of metrics.

\begin{definition}
A random process $X=(X(t),t\in T)$ belongs to a quasi-Banach $K_{\sigma}$-space $K(\Omega)$, $X\in K(\Omega)$, if for each $t\in T$ the random variable $X(t)$ belongs to $K(\Omega)$.
\end{definition}

\begin{definition}
A pseudometric $\rho(t,s)$ is generated by a random \linebreak process $X$ if $\rho(t,s)=||\xi(t)-\xi(s)||_{K}$.
\end{definition}

Let us make some assumptions. We will assume that the space $(T,\rho)$ is separable. If $B$ is a compact set with $(T,\rho)$, then denote by $N_{B}(\varepsilon)$, $\varepsilon>0$ the number of elements of the minimal cover of the set $B$ by open balls of radius $\varepsilon$.

Consider a random process $X=(X(t),t\in T)\in K(\Omega)$, where $K(\Omega)$ is a quasi-$K_{\sigma}$-space with quasi-norm $||\cdot||_{K}$, subject to the function $j=\left\{j(\lambda),|\lambda|<1\right\}$ and majorizing characteristic $\varkappa(n)$. Let there be a continuous monotonically increasing function $\sigma=\left\{\sigma(h),h\geq0\right\}$, $\sigma(0)=0$ such that $\sup\limits_{\rho(t,s)\leq h}||x(t)-x(s)||_{K}\leq \sigma(h)$. Note that in the case where $\rho(t,s)=||x(t)-x(s)||_{K}$, $\sigma(h)=h$. For any bounded function $f(t)$ on the set $B\subset T$, we denote $||f(t)||_{C(B)}=\sup\limits_{t\in B}|f(t)|$.

\begin{theorem} \label{th01}
Let $c(t)$ be a function on $(T,\rho)$ such that $|c(t)|<1$, and let the random process $X=(X(t),t\in T)$ be separable, $X\in K(\Omega)$, $T=\bigcup\limits_{k=1}^{\infty}B_{k}$, where $B_{k}$ are compact sets. Then the equality holds
\begin{equation}\label{eq02}
\left\|||c(t)X(t)||_{C(T)}\right\|_{K}\leq \sum\limits_{k=1}^{\infty}j(\gamma_{k})(\delta_{k}+\sum\limits_{l=0}^{\infty}\mu\left(N_{B_{k}}(\varepsilon_{k,l})\right)\sigma(\varepsilon_{k,l})),
\end{equation}
where $\varepsilon_{k,l}$, $l=1,2,...,\infty$ -- any monotonically decreasing sequence, $\varepsilon_{k,l}\rightarrow 0$, $l\rightarrow\infty$, where $\varepsilon_{k,0}$ -- is the smallest number such that for all $t\in B_{k}$ there exists $t_{0}\in B_{k}$, such that $\rho(t,t_{0})\leq\varepsilon_{k,0}$, $\varepsilon_{k,0}\leq\sup\limits_{l,s\in B_{k}}\rho(t,s)$, $\delta_{k}=\sup\limits_{t\in B_{k}}||X(t)||_{K}=\left\|||X(t)||_{K}\right\|_{C(B_{k})}, \gamma_{k}=||c(t)||_{C(B_{k})}$.
\end{theorem}

\begin{dov}
If the series on the right-hand side of \eqref{eq02} diverge, then continuity is trivial. Therefore, we will assume that
\[
\sum\limits_{k=1}^{\infty}j(\gamma_{k})\left(\delta_{k}+\sum\limits_{l=0}^{\infty}\mu(N_{B_{k}}(\varepsilon_{k,l+1}))\sigma(\varepsilon_{k,l})\right)<\infty.
\]
The inequality holds true
\[
\left\|c(t)X(t)\right\|_{C(T)}\leq\sup\limits_{k=1,...,\infty}\left\|c(t)X(t)\right\|_{C(B_{k})}\leq
\]
\[
\leq\sum^{\infty}_{k=1}\left\|c(t)X(t)\right\|_{C(B_{k})}\leq\sum^{\infty}_{k=1}\gamma_{k}\left\|X(t)\right\|_{C(B_{k})}.
\]
From this inequality it follows
\begin{equation}\label{eq03}
\left\|\left\|c(t)X(t)\right\|_{C(T)}\right\|_{K}\leq \sum\limits_{k=1}^{\infty}j(\gamma_{k})\left\|\left\|X(t)\right\|_{C(B_k)}\right\|_K.
\end{equation}

Let us estimate $\left\|\left\|X(t)\right\|_{C(B_k)}\right\|_K$. Let us consider the minimal covering of the compact set $B_k$ by balls of radius $\varepsilon_{k,l}$.

Let $V_{\varepsilon_{k,l}}$ be the set of ball centers of this covering. The number of points in $V_{\varepsilon_{k,l}}$ is equal to $N_{B_{k}}(\varepsilon_{k,l})$. The inequality is
\begin{equation}\label{eq04}
||X(t)||_{C(B_{k})}\leq \sum\limits_{l=0}^{\infty}\sup\limits_{t\in {V_{\varepsilon_{k,l+1}}}}|X(t)-X(\alpha_{1}(t))|+|X(t_{k})|,
\end{equation}
where $t_{k}$ is a fixed point in $B_k$, $\alpha_{1}(t)$ is some fixed point in $V_{\varepsilon_{k,l}}$ such that $\rho(t,\alpha_{1}(t))<\varepsilon_{k,l}$.
The inequality \eqref{eq04} can be obtained in exactly the same way as a similar inequality when proving Theorem 1 from the work \cite{abza85},
or from \cite{koz84}. From Definition~\ref{de01} we get the relation
\[
\left\|\sup\limits_{t\in V_{\varepsilon_{k,l+1}}}|X(t)-X(\alpha_{1}(t))|\right\|_{K}\leq
\]
\[
\leq\mu(N_{B_{k}}(\varepsilon_{k,l+1}))\sup\limits_{t\in V_{\varepsilon_{k,l+1}}}\left\|X(t)-X(\alpha_{1}(t))\right\|_{K}\leq\mu(N_{B_{k}}(\varepsilon_{k,l+1}))\sigma(\varepsilon_{k,l}).
\]
From this relation and inequality \eqref{eq04} we obtain
\begin{eqnarray}
\left\|\left\|X(t)\right\|_{C(B_k)}\right\|_K\leq\left\|X(t_{k})\right\|_{K}+\sum^{\infty}_{l=0}\mu(N_{B_{k}}(\varepsilon_{k,l+1}))\sigma(\varepsilon_{k,l}).
\nonumber
\end{eqnarray}
From the inequality \eqref{eq03} and the last inequality, the statement of the Theorem follows, if we take into account that $\left\|X(t_{k})\right\|\leq\delta_{k}$.
\end{dov}

\begin{corollary} \label{co01}
Let the conditions of Theorem \ref{th01} be satisfied and the series\\
$\sum\limits_{k=1}^{\infty}j(\gamma_{k})\left(\delta_{k}+\sum\limits_{l=0}^{\infty}\mu \left(N_{B_{k}}(\varepsilon_{k,l+1})\right)\sigma(\varepsilon_{k,l})\right)=A<\infty,
$
then with probability one $\left\|c(t)X(t)\right\|_{C(T)}\in K(\Omega)$ and $\left\|\left\|c(t)X(t)\right\|_{C(T)}\right\|_{K}\leq A$.
\end{corollary}

\begin{corollary}\label{co02}
Let the condition of Theorem \ref{th01} be satisfied and, in addition, for any $0<q<1$ the condition is satisfied
\[
\sum\limits_{k=1}^{\infty}j(\gamma_{k})\left(\delta_{k}\frac{1}{q(1-q)}\int^{q\sigma(\varepsilon_{k,0})}_{0}\mu\left(N_{B_{k}}(\sigma^{(-1)}(u))\right)du\right)=
\]
\[
=\sum\limits_{k=1}^{\infty}j(\gamma_{k})\left(\delta_{k}\frac{1}{q(1-q)}\int^{\gamma_{k,0}}_{0}\mu(N_{B_{k}}(t))d\sigma(t)\right)=W_{q}<\infty,
\]
where $\gamma_{k,0}=\sigma^{(-1)}(q(\sigma(\varepsilon_{k,0})))\leq\varepsilon_{k,0}$.
Then with probability one
\linebreak $\left\|c(t)X(t)\right\|_{C(T)}\in k(\Omega)$ and $\left\|\left\|c(t)X(t)\right\|_{C(T)}\right\|_{K}\leq W_{q}$.
\end{corollary}

\begin{dov}
Denote $t_{k,l}=\sigma(\varepsilon_{k,l})$, then
$\varepsilon_{k,l}=\sigma^{(-1)}(t_{k,l})$
\[
\sum\limits_{l=0}^{\infty}\mu(N_{B_{k}}(\varepsilon_{k,l+1}))\sigma(\varepsilon_{k,l})=\sum\limits_{l=0}^{\infty}\mu\left(N_{B_{k}}(\sigma^{(-1)}(t_{k,l+1}))\right)t_{k,l}.
\]
Let us choose a sequence $\varepsilon_{k,l}$ such that $t_{k,l}=t_{k,0}q^l$, $0<q<1$, i.e. $\varepsilon_{k,l}=\sigma^{(-1)}(t_{k,0}q^l)$. In this case
\[
t_{k,l}\mu\left(N_{B_{k}}(\sigma^{(-1)}(t_{k,l+1}))\right)\leq\frac{t_{k,l}}{t_{k,l+1}-t_{k,l+2}}\int^{t_{k,l+1}}_{t_{k,l+2}}\mu\left(N_{B_{k}}(\sigma^{(-1)}(u))\right)du=
\]
\[
=\frac{1}{q(1-q)}\int^{t_{k,l+1}}_{t_{k,l+2}}\mu\left(N_{B_{k}}(\sigma^{(-1)}(u))\right)du.
\]
Therefore,
\[
\sum\limits_{l=0}^{\infty}\mu\left(N_{B_{k}}(\varepsilon_{k,l+1})\right)\sigma(\varepsilon_{k,l})\leq\frac{1}{q(1-q)}\int^{t_{k,o}q}_{0}\mu\left(N_{B_{k}}(\sigma^{(-1)}(u))\right)du=
\]
\[
=\frac{1}{q(1-q)}\int^{\sigma^{(-1)}(q\sigma(\varepsilon_{k,0}))}_{0}\mu\left(N_{B_{k}}(t)\right)d\sigma(t).
\]

From the last inequalities we obtain that for such a sequence $A\leq W_{q}$.
Therefore, the statement of Corollary \ref{co02} follows from Corollary \ref{co01}.
\end{dov}

\begin{corollary} \label{co03}
Let $(T,\rho)$ be a compact set. If the integral converges
\[
\int^{q\sigma(\varepsilon_0)}_{0}\mu(N_{T}(\sigma^{(-1)}(u)))du<\infty,
\]
where $\varepsilon_{0}=\sup\limits_{t,s\in T}\rho(t,s)$, then $\left\|X(t)\right\|_{C(T)}\in K(\Omega)$ and
$\left\|\left\|X(t)\right\|_{C(T)}\right\|_{K}\leq V_{q}$, $0<q<1$,
where
\[
V_{q}=\delta+\frac{1}{q(1-q)}\int^{q\sigma(\varepsilon_{0})}_{0}\mu(N_{T}(\sigma^{(-1)}(u)))du=\delta+\frac{1}{q(1-q)}\int^{\gamma_{0}}_{0}\mu(N_{T}(t))d\sigma(t),
\]
\[
\delta=\left\|\left\|X(t)\right\|_{K}\right\|_{C(T)}, \gamma_{0}=\sigma^{(-1)}(q\sigma(\varepsilon_{0})).
\]
\end{corollary}

Corollary \ref{co03} follows from Corollary \ref{co02} if we put $B_{1}=T, B_{k}=0, k>1$.

\begin{remark}
If we put $\rho(t,s)=\left\|X(t)-X(s)\right\|_{K}$, then Corollary \ref{co03} implies the statement of Theorem 1 from the work \cite{abza85}. In this case
\[
V_{q}=\delta+\frac{1}{q(1-q)}\int^{\gamma_{0}}_{0}\mu(N_{T}(t))d\sigma(t).
\]
\end{remark}

\begin{remark}
In the case where $K(\Omega)$ is the Orlicz space $L_{U}(\Omega)$, the inequality \eqref{de02} implies that for any $\varepsilon>0$
\[
P\left\{\left\|c(t)X(t)\right\|_{C(T)}>\varepsilon\right\}\leq(U(\varepsilon/z))^{-1},
\]
where $z$ is equal to $A$ when the conditions of the Corollary \ref{co01} are satisfied, or $V_{q}$ if the conditions of the Corollary \ref{co03} are satisfied.
\end{remark}



\begin{theorem} \label{th02}
Let $X=(X(t),t\in T)$ be a random process, where $(T,\rho)$ is a compact set such that $X$ belongs to the quasi-$K_{\sigma}$-space $K(\Omega)$.
Let the condition hold
\begin{equation}\label{de03}
\sum^{\infty}_{l=0}\mu(N_{T}(\varepsilon_{l+1}))\sigma(\varepsilon_{l})<\infty,
\end{equation}
where $\varepsilon_{l}$ is any monotonically decreasing sequence such that \linebreak $\varepsilon_{0}=\sup\limits_{t,s\in T}\rho(t,s)$ and $\varepsilon_{l}\rightarrow 0$ when $l\rightarrow\infty$. If $X$ is a separable process, then it is sample continuous with probability one.
\end{theorem}

\begin{dov} Let $V_{\varepsilon_{l}}$ be the set of centers of the minimal covering of the compact set $T$ by balls of radius $\varepsilon_{l}$. If we repeat the proof of Theorem 1 from the work \cite{abza85} without changes, we can prove that for any $k>0$ there exists a positive number $d>0$ such that the inequality
\[
\sup\limits_{\rho(t,s)<d}|X(t)-X(s)|\leq 4\sum^{\infty}_{l=k}\sup\limits_{t\in V_{\varepsilon_{l+1}}}|X(t)-X(\alpha_{l}(t))|,
\]
where $\alpha_{l}(t)$ is a fixed point in $V_{\varepsilon_{l}}$ such that $\rho(t,\alpha_{l}(t))<\varepsilon_{l}$.

The last inequality yields the relations
\begin{multline}\label{de04}
\left\|\sup\limits_{\rho(t,s)<d}|X(t)-X(s)|\right\|_{K}\leq 4\sum^{\infty}_{l=k}\left\|\sup\limits_{t\in V_{\varepsilon_{l+1}}}|X(t)-X(\alpha_{l}(t))|\right\|_{K}\leq \\
\leq4\sum^{\infty}_{l=k}\mu(N_{T}(\varepsilon_{l+1}))\sigma(\varepsilon_{l}).
\end{multline}
From the condition \eqref{de03} we have that $\sum\limits_{l=k}^{\infty}\mu(N_{T}(\varepsilon_{l+1}))\sigma(\varepsilon_{l})\rightarrow 0$ as $k\rightarrow\infty$. Therefore, from \eqref{de04} it follows that $\left\|\sup\limits_{\rho(t,s)<\varepsilon}|X(t)-X(s)|\right\|_{K}\rightarrow 0$ as $\varepsilon\rightarrow 0$. From the lemma \ref{le01} it follows that $\sup\limits_{\rho(t,s)<\varepsilon}|X(t)-X(s)|\rightarrow 0$ as $\varepsilon\rightarrow 0$ by probability. Therefore, there exists a sequence $\varepsilon_{n}$ such that $\sup\limits_{\rho(t,s)<\varepsilon_{n}}|X(t)-X(s)|\rightarrow 0$ with probability one.
Thus, the process $X(t)$ is sample continuous with probability one.
\end{dov}

\begin{corollary} \label{co04}
Theorem \ref{th02} holds if instead of condition \eqref{de03} we require that the integral
\begin{equation}\label{de05}
\int^{\delta}_{0}\mu(N_{T}(t))d\sigma(t)<\infty,
\end{equation}
where $\delta>0$ is any constant.
\end{corollary}

The statement of Corollary \ref{co04} is proved in the same way as the statement of Corollary \ref{co02}.
\begin{remark} In the case when $\rho(t,s)=\left\|X(t)-X(s)\right\|_{K}$, that is, if we can put $\sigma(h)=h$, the condition \eqref{de05} has the form
\[
\int^{\delta}_{0}\mu(N_{T}(t))d(t)<\infty.
\]
This condition was obtained in the work \cite{abza86}, and in special cases in \cite{koz84, abza85}.
\end{remark}

\section{Random processes from quasi-Banach $K_{\sigma}$-spaces in the case where $T=\mathbb R^{1}$}

Let us consider applications of the previous results in the case where $T=\mathbb R^{1}$ with the usual metric. In the case where $T=\mathbb R^{n}$, the same results can be obtained.

\begin{theorem} Let $X=(X(t),t\in R^1)$ be a separable random process, $X\in K(\Omega)$ , where $K(\Omega)$ is a quasi-$K_\sigma$-space with quasinorm $||\cdot||_K$, subject to the function $j=\left\{j(\lambda),|\lambda|<1\right\}$, and with majorizing characteristic $\mu(n)$. Let there exist a continuous monotone increasing function $\sigma=\left\{\sigma(h),h>0\right\}$, where $\sigma(h)\rightarrow 0$, as $h\rightarrow 0$, such that
\[
\sup\limits_{|t-s|\leq h}\left\|X(t)-X(s)\right\|\leq\sigma(h).
\]
Let for any $v>0$ there exist an integral
\[
\int^{v}_{0}\mu\left(\left[\frac{v}{t}\right]\right)d\sigma(t)<\infty.
\]
Then on any closed interval I the process $X(t)$ is sample continuous with probability one. If, in addition, for some continuous function $c(t), |c(t)|<1$, the series
\[
\sum^{\infty}_{k=1}j(\gamma_k)\left(\delta_{k}+\frac{1}{q(1-q)}\int^{\alpha_{k,0}}_{0}\mu\left(\left[\frac{z_{k}}{t}\right]\right)d\sigma(t)\right)=\Phi_{q}<\infty,
\]
where $0<q<1$, $z_{k}$ is the length of the closed intervals $B_{k}$ such that $R^{1}\subset\bigcup\limits_{k=1}^{\infty}B_{k},\gamma_{k}=\sup\limits_{t\in B_{k}}|c(t)|,\delta_{k}=\sup\limits_{t\in B_{k}}\left\|X(t)\right\|_{K},\alpha_{k,0}=\sigma^{(-1)}(q(\sigma(z_{k}/2)))\leq z_{k}/2$, then with probability one
\[
\left\|c(t)X(t)\right\|^{1}_{C(R)}\leq\eta,
\]
where $\eta$ is a random variable from the space $K(\Omega)$ such that $\left\|\eta\right\|_{K}\leq\Phi_{q}$.
\end{theorem}

The Theorem follows from the Corollary \ref{co02} and Theorem \ref{th02}, if we note that in this case
\begin{eqnarray}
N_{B_{k}}(t)\leq\left[\frac{z_{k}}{2t}\right]+1\leq\left[\frac{z_{k}}{t}\right].
\nonumber
\end{eqnarray}

\begin{corollary} \label{co05}
Let $X=(X(t),t\in R^{1})$ be a separable random process, $X\in K(\Omega), \sup\limits_{t\in R^1}\left\|X(t)\right\|_{K}=\delta<\infty$. If the integral $(\varepsilon>0)$ converges
\[
\int^{\varepsilon}_{0}\mu\left(\left[\frac{I}{t}\right]\right)d\sigma(t)<\infty,
\]
then on any interval I with probability one the process $X$ is sample-continuous.
For any even continuous function $c(t),0<c(t)\leq 1$, which monotonically approaches zero as $n\rightarrow \infty$ and for which  $\sum\limits_{k=1}^{\infty}j(c(k))<\infty$ with probability one, the inequality holds
\[
\left\|c(t)X(t)\right\|^{1}_{C(R)}\leq\theta,
\]
where $\theta\in K(\Omega)$,
\[
\left\|\theta\right\|\leq 2\left[\delta+\frac{1}{q(1-q)}\int^{\alpha_{0}}_{0}\mu\left(\left[\frac{I}{t}\right]\right)d\sigma(t)\right]\sum^{\infty}_{k=1}j(c(k+1)),
\]
$\alpha_{0}=\sigma^{(-1)}(q(\sigma(1/2)))$.
\end{corollary}

To obtain the statement of the Corollary, let $B_{k}=[k,k+1],k=0,\pm1,\pm2,...$. In this case $z_{k}=1$, $\delta_{k}\leq\delta$, $\alpha_{k,0}=\sigma^{(-1)}(q(\sigma(1/2)))$. Therefore,
\[
\Phi_{q}=2\left[\delta+\frac{1}{q(1-q)}\int^{\alpha_{0}}_{0}\mu\left(\left[\frac{I}{t}\right]\right)d\sigma(t)\right]\sum^{\infty}_{k=0}j(c(k+1)).
\]

\begin{remark}
In the case where $K(\Omega)$ is a Banach $K_{\sigma}$-space, the conditions of the Corollary \ref{co05} are satisfied, for example, by the function $c(t)=(1+|t|_{\alpha})^{-1},\alpha>1$.
\end{remark}

It is clear that in special cases a more precise order of growth can be obtained.

\section*{Conclusions to the first chapter}\

In the first section, we consider random processes from quasi-Banach $K_{\sigma}$-spaces of random variables.
We investigate the conditions of sample continuity with probability one and the conditions that the suprema of the processes belong to the same space as the processes themselves.
We obtain estimates of the distributions of these processes. We study the behavior of random processes $X(t)$ on $\mathbb R^1$ as $t$ approaches infinity.

\chapter{Spaces $\mathbf{F}_\psi(\Omega)$ of random variables}
\label{ch:o2series}

In this chapter \ref{ch:o2series}, the basic properties of the spaces $\mathbf{F}_\psi(\Omega)$ are studied.
Examples of random variables from these spaces are given. Large deviation inequalities
and a majorizing characteristic for random variables from the space $\mathbf{F}_\psi(\Omega)$ are found.
The spaces $\mathbf{F}_{S_k,\psi,r}(\Omega)$ are considered and the conditions under which these spaces
are equivalent to the space $\mathbf{F}_\psi(\Omega)$ are studied. For the spaces $\mathbf{F}_\psi(\Omega)$, the conditions under which the condition $\mathbf{H}$ is satisfied are studied.

\section{Spaces $\mathbf{F}_\psi(\Omega)$ of random variables. Basic properties}

\begin{definition}
Let $\psi(u)>0$, $u\geq1$ be a monotonically increasing continuous function such that $\psi(u)\rightarrow\infty$ as $u\rightarrow\infty$.
A random variable $\xi$ belongs to the space $\mathbf{F}_\psi(\Omega)$ if the following condition holds:
\[
\sup\limits_{u\geq1}\frac{\left(E\left|\xi\right|^{u}\right)^{1/u}}{\psi(u)}<\infty.
\]
\end{definition}

A similar definition was formulated in the work of S.~V.~Ermakov and E.~I.~Ostrovskii \cite{erma86}.
But there it is required that $E\xi=0$ if $\xi \in \mathbf{F}_\psi(\Omega)$.
In addition, random variables were considered such that $E\left|\xi\right|^u=\infty$ for some $u>0$.

The following theorem holds true..

\begin{theorem} \cite{erma86}
The space $\mathbf{F}_\psi(\Omega)$ is a Banach space with the norm
\begin{equation} \label{eq:o1chapter_1}
\left\|\xi\right\|_\psi=\sup\limits_{u\geq1}\frac{\left(E\left|\xi\right|^{u}\right)^{1/u}}{\psi(u)}.
\end{equation}
\end{theorem}

\begin{proof}
Let us first prove that $\mathbf{F}_\psi(\Omega)$ is a linear normed space. It is obvious that $\left\|\xi\right\|_\psi=0$ if and only if $\xi=0$ with probability one. The equality is true
\[
\left\|\alpha\xi\right\|_\psi=\sup\limits_{u\geq1}\frac{\left(E\left|\alpha\xi\right|^{u}\right)^{1/u}}{\psi(u)}=\sup\limits_{u\geq1}\frac{\left|\alpha\right|\left(E\left|\xi\right|^{u}\right)^{1/u}}{\psi(u)}=\left|\alpha\right|\left\|\xi\right\|_\psi.
\]
The inequality of the triangle is also obvious. Indeed,
\[
\left\|\xi_1+\xi_2\right\|_\psi=\sup\limits_{u\geq1}\frac{\left(E\left|\xi_1+\xi_2\right|^{u}\right)^{1/u}}{\psi(u)}\leq \sup\limits_{u\geq1}\frac{\left(E\left|\xi_1\right|^{u}\right)^{1/u}+\left(E\left|\xi_2\right|^{u}\right)^{1/u}}{\psi(u)}\leq
\]
\[
\leq \sup\limits_{u\geq1}\frac{\left(E\left|\xi_1\right|^{u}\right)^{1/u}}{\psi(u)}+\sup\limits_{u\geq1}\frac{\left(E\left|\xi_2\right|^{u}\right)^{1/u}}{\psi(u)}=\left\|\xi_1\right\|_\psi+\left\|\xi_2\right\|_\psi.
\]
Let us now show that $\mathbf{F}_\psi(\Omega)$ is a complete space, i.e., if $\xi_n \in \mathbf{F}_\psi(\Omega)$ and $\left\|\xi_n-\xi_l\right\|_\psi\to 0$ for $n,l \to\infty$, then there exists a random variable $\xi$ such that $\xi \in \mathbf{F}_\psi(\Omega)$ and $\left\|\xi_n-\xi\right\|_\psi\to 0$ as $n \rightarrow\infty$.
From the definition of the norm it follows that for any $u\geq1$
\begin{equation} \label{eq:o1chapter_1111}
\left(E\left|\xi_n-\xi_l\right|^{u}\right)^{1/u}\leq \psi(u)\left\|\xi_n-\xi_l\right\|_\psi.
\end{equation}

Since $\left\|\xi_n-\xi_l\right\|_\psi\rightarrow 0$ for $n,l \rightarrow\infty$, then
$\left(E\left|\xi_n-\xi_l\right|^{u}\right)^{1/u}\rightarrow 0$ for $n,l \rightarrow\infty$. The space $L_u(\Omega)$, $u\geq1$ is complete because there exists a random variable $\xi \in L_u(\Omega)$ such that $\xi_n\rightarrow\xi$ for $n\rightarrow \infty$ is in the norm of this space. It is easy to see that there exists $\xi \in L_u(\Omega)$ for all $u\geq1$ such that $\left(E\left|\xi_n-\xi_l\right|^{u}\right)^{1/u}\rightarrow 0$ for $n\rightarrow \infty$. Indeed, if $\xi_n\rightarrow\xi$ is in the norm of the space $L_u(\Omega)$, then $\xi_n\rightarrow\xi$ is in the norm of the space $L_v(\Omega)$, where $v<u$. Let $\eta_s$, $s=1,2,...$ be random variables such that $\xi_n\rightarrow\eta_s$ is in the norm of the spaces $L_u(\Omega)$, where $s-1<u\leq s$. Then there exist subsequences $\xi_{n_s}$ that converge to $\eta_s$ with probability one. Let $A_s$ be a set $P\left(A_s\right)=1$ on which $\xi_{n_s}$ converges to $\eta_s$. Then on the set $\bigcap\limits_{s=1}^\infty A_s$ all $\eta_s$ are equal and $P\left\{\bigcap\limits_{s=1}^\infty A_s\right\}=1$.

Now let $\xi$ be a random variable equal to $\eta_s$ on $\bigcap\limits_{s=1}^\infty A_s$. It is clear that $P\left\{\eta_s\neq\xi\right\}=0$. Therefore
$\xi_n\rightarrow\xi$ in $L_u(\Omega)$, while $\xi_n\rightarrow\eta_s$ is in the same space.

So, for all $u\geq1$, the inequality \eqref{eq:o1chapter_1111} implies that
\[
\left(E\left|\xi_n-\xi_l\right|^{u}\right)^{1/u}\leq \psi(u) \sup\limits_{r>n} \left\|\xi_n-\xi_r\right\|_\psi<\infty.
\]
If in the last inequality we direct $l$ to infinity, then we obtain that for all $n\geq1$ and $u\geq1$
\begin{equation} \label{eq:o1chapter_11111}
\left(E\left|\xi_n-\xi\right|^{u}\right)^{1/u}\leq \psi(u) \sup\limits_{r>n} \left\|\xi_n-\xi_r\right\|_\psi<\infty.
\end{equation}
Therefore, $\sup\limits_{u\geq1}\frac{\left(E\left|\xi_n-\xi\right|^{u}\right)^{1/u}}{\psi(u)}<\infty$. That is, the random variables $\xi_n-\xi$ for $n\geq1$ belong to the space $\mathbf{F}_\psi(\Omega)$. Since $\left\|\xi\right\|_\psi\leq\left\|\xi-\xi_n\right\|_\psi+\left\|\xi_n\right\|_\psi<\infty$, then the random variable $\xi$ belongs to the space $\mathbf{F}_\psi(\Omega)$. Now from the inequality \eqref{eq:o1chapter_11111} it follows that $\left\|\xi_n-\xi\right\|_\psi\leq \sup\limits_{r>n} \left\|\xi_n-\xi_r\right\|_\psi\rightarrow 0$ for $n\rightarrow \infty$. That is, $\xi_n\rightarrow\xi$ is in the norm of the space $\mathbf{F}_\psi(\Omega)$.
\end{proof}

\begin{theorem}\label{th:o1chapter_1}
Let a random variable $\xi$ belong to the space $\mathbf{F}_\psi(\Omega)$, then for any $\varepsilon>0$ the following inequality holds:
\begin{equation} \label{eq:o1chapter_2}
P\left\{\left|\xi\right|>\varepsilon\right\}\leq \inf\limits_{u\geq1}\frac{\left\|\xi\right\|_\psi^u(\psi(u))^u}{\varepsilon^u}.
\end{equation}
\end{theorem}

\begin{proof}
From Chebyshev's inequality it follows that for $u>0$ the following inequality holds:
\[
P\left\{\left|\xi\right|>\varepsilon\right\}\leq \frac{E\left|\xi\right|^{u}}{\varepsilon^u} = \frac{E\left|\xi\right|^{u}(\psi(u))^u}{(\psi(u))^u \varepsilon^u} \leq\frac{\left\|\xi\right\|_\psi^u(\psi(u))^u}{\varepsilon^u}.
\]
\end{proof}

\begin{theorem}\label{th:o1chapter_2}
Let the random variable $\xi$ belong to the space $\mathbf{F}_\psi(\Omega)$ and $\psi(u)=u^\alpha$, where $\alpha>0$, then for any $\varepsilon\geq e^\alpha\left\|\xi\right\|_\psi$ the inequality holds:
\begin{equation} \label{eq:o1chapter_3}
P\left\{\left|\xi\right|>\varepsilon\right\}\leq \exp \left\{-\frac{\alpha}{e}\left(\frac{\varepsilon}{\left\|\xi\right\|_\psi}\right)^{1/\alpha}\right\}.
\end{equation}
\end{theorem}

\begin{proof}
Using Theorem \ref{th:o1chapter_1} we have that
\begin{equation} \label{eq:o1chapter_4}
P\left\{\left|\xi\right|>\varepsilon\right\}\leq \inf\limits_{u\geq1}\dfrac{\left\|\xi\right\|_\psi^u u^{\alpha u}}{\varepsilon^u}.
\end{equation}
Let us denote by $\frac{\left\|\xi\right\|_\psi}{\varepsilon}=b$, then from the equalities
\[
\left(\ln\left(b^u u^{\alpha u}\right)\right)^{'}= \left(u\ln b + \alpha u \ln u\right)^{'}
=\ln b+\alpha \ln u+\alpha=0;
\]
\[
\ln u=-\frac{\ln b+\alpha}{\alpha}
\]
it follows that the infinitude is reached at the point $u=\frac{1}{e} b^{-1/\alpha}$. Since $u\geq1$, then the inequality $\varepsilon\geq e^\alpha\left\|\xi\right\|_\psi$ must hold. Substituting this value for the quantity $u$ into the inequality \eqref{eq:o1chapter_4}, we obtain:
\[
P\left\{\left|\xi\right|>\varepsilon\right\}\leq b^{\frac{1}{e}b^{-\frac{1}{\alpha}}} \left(\frac{1}{e}b^{-\frac{1}{\alpha}}\right)^{\alpha\frac{1}{e}b^{-\frac{1}{\alpha}}}=\exp \left\{-\frac{\alpha}{e}\left(\frac{1}{b}\right)^{1/\alpha}\right\}.
\]
This leads to the statement of Theorem \ref{th:o1chapter_2}.
\end{proof}

\begin{theorem}\label{th:o1chapter_3}
Let the random variable $\xi$ belong to the space $\mathbf{F}_\psi(\Omega)$ and $\psi(u)=e^{au^\beta}$, where $a>0$, $\beta>0$, then for any $\varepsilon\geq e^{a(\beta+1)}\left\|\xi\right\|_\psi$ the inequality holds:
\begin{equation} \label{eq:o1chapter_5}
P\left\{\left|\xi\right|>\varepsilon\right\}\leq \exp \left\{-\frac{\beta}{a^{1/\beta}}\left(\frac{\ln \frac{\varepsilon}{\left\|\xi\right\|_\psi}}{\beta+1}\right)^{\frac{\beta+1}{\beta}}\right\}.
\end{equation}
\end{theorem}

\begin{proof}
From Theorem \ref{th:o1chapter_1} we obtain that
\begin{equation} \label{eq:o1chapter_6}
P\left\{\left|\xi\right|>\varepsilon\right\}\leq \inf\limits_{u\geq1}\frac{\left\|\xi\right\|_\psi^u e^{au^{\beta+1}}}{\varepsilon^u}.
\end{equation}
Let us denote $\frac{\left\|\xi\right\|_\psi}{\varepsilon}=b$. From the equalities
\[
\left(\ln\left(b^u e^{au^{\beta+1}}\right)\right)^{'}= \left(u\ln b + a u^{\beta+1}\right)^{'}=\ln b+a(\beta+1)u^{\beta}=0;
\]
it follows that the infinitude is reached at the point $u=\left(-\frac{\ln b}{a(\beta+1)}\right)^{1/\beta}$. Since $u\geq1$, then the inequality $\varepsilon\geq e^{a(\beta+1)}\left\|\xi\right\|_\psi$ must hold. Substituting this value of $u$ into the inequality \eqref{eq:o1chapter_6}, we obtain
\[
P\left\{\left|\xi\right|>\varepsilon\right\}\leq b^{\left(-\frac{\ln b}{a(\beta+1)}\right)^{1/\beta}}e^{a\left(-\frac{\ln b}{a(\beta+1)}\right)^{\frac{\beta+1}{\beta}}}=
\exp \left\{-\frac{\beta}{a^{1/\beta}}\left(\frac{\ln \frac{1}{b}}{\beta+1}\right)^{\frac{\beta+1}{\beta}}\right\},
\]
which had to be proven.
\end{proof}

\begin{theorem}\label{th:o1chapter_4}
Let the random variable $\xi$ belong to the space $\mathbf{F}_\psi(\Omega)$ and $\psi(u)=\left(\ln (u+1)\right)^\lambda$, where $\lambda>0$, then for any
$\varepsilon\geq\left(e \ln 2\right)^\lambda \left\|\xi\right\|_\psi$ the inequality holds:
\begin{equation} \label{eq:o1chapter_7}
P\left\{\left|\xi\right|>\varepsilon\right\}\leq e^\lambda \exp \left\{-\lambda  \exp\left\{\left(\frac{\varepsilon}{\left\|\xi\right\|_\psi}\right)^{1/\lambda}\frac{1}{e}\right\}\right\}.
\end{equation}
\end{theorem}

\begin{proof}
Since it follows from Theorem \ref{th:o1chapter_1} that
\begin{equation} \label{eq:o1chapter_8}
P\left\{\left|\xi\right|>\varepsilon\right\}\leq \inf\limits_{u\geq1}\frac{\left\|\xi\right\|_\psi^u (\ln(u+1))^{\lambda u}}{\varepsilon^u},
\end{equation}
then we put $u+1=\exp \left\{\left(\frac{\varepsilon}{\left\|\xi\right\|_\psi}\right)^{1/\lambda}\frac{1}{z}\right\}$, where $z>0$. Then, substituting this expression into the inequality \eqref{eq:o1chapter_8}, we obtain:
\begin{multline*}
\left(\frac{\left\|\xi\right\|_\psi(\ln(u+1))^{\lambda }}{\varepsilon}\right)^u=\frac{1}{z^{\lambda u}}=\exp \left\{-\lambda u \ln z\right\}= \\[1ex]
=\exp \left\{-\lambda  (\ln z)\left(\exp \left\{\left(\frac{\varepsilon}{\left\|\xi\right\|_\psi}\right)^{1/\lambda}\frac{1}{z}\right\}-1\right)\right\}= \\[1ex]
=z^\lambda \exp \left\{-\lambda  (\ln z)\exp \left\{\left(\frac{\varepsilon}{\left\|\xi\right\|_\psi}\right)^{1/\lambda}\frac{1}{z}\right\}\right\}.
\end{multline*}
Let us put $z=e$ in this equality, then we will obtain the statement of the Theorem:
\[
P\left\{\left|\xi\right|>\varepsilon\right\}\leq e^\lambda \exp \left\{-\lambda  \exp\left\{\left(\frac{\varepsilon}{\left\|\xi\right\|_\psi}\right)^{1/\lambda}\frac{1}{e}\right\}\right\}.
\]
\end{proof}

Let us give examples of random variables from the spaces $\mathbf{F}_\psi(\Omega)$.

\begin{example}\label{ex1}
A centered normal random variable $\xi \sim N(0,\sigma^2)$  belongs
to the space $\mathbf{F}_\psi(\Omega)$ with
$$\psi(u)=u^\frac12.$$
Indeed,
\begin{eqnarray}\label{eq-ex2}
E\vert\xi\vert^u&=&\int_{-\infty}^\infty \vert  x \vert^u \frac1{\sqrt{2\pi}\sigma}
e^{-\frac{x^2}{2\sigma^2}}dx=\frac{\sqrt{2}}{\sqrt{\pi}\sigma}\int_{0}^\infty
x^u e^{-\frac{x^2}{2\sigma^2}}dx
\nonumber\\
& = &\frac{2^{\frac{u}{2}} \sigma^{u}}{\sqrt{\pi}} \int_{0}^\infty
s^{\frac{u}{2}-\frac12} e^{-s} ds= \frac{2^{\frac{u}2}
\sigma^u}{\sqrt{\pi}} \Gamma \left(\frac{u}2+\frac12\right),
\end{eqnarray}
where $\Gamma (\alpha)=\int_{0}^\infty t^{\alpha-1} e^{-t} dt$  is
 the Gamma function or Euler integral of the se\-cond kind. It is shown
 in \cite{askey2010} that

\begin{equation}\label{eq-ex3}
\Gamma(z)\leq e^{-z}z^{z-\frac12}\sqrt{2\pi} e^{\frac1{12z}}.
\end{equation}
By \eqref{eq-ex3} we have
\begin{eqnarray}\label{eq-ex4}
\Gamma \left(\frac{u}2+\frac12\right)& \leq &
e^{-\frac{u}2-\frac12}\left(\frac{u}2+\frac12\right)^{\frac{u}2}\sqrt{2\pi}
e^{\frac1{12(\frac{u}2+\frac12)}}\nonumber\\
&=& \sqrt{2\pi} e^{-\frac{u}2}(u+1)^{\frac{u}2}
2^{-\frac{u}2}e^{\frac1{6(u+1)}-\frac12}.
\end{eqnarray}
Substituting \eqref{eq-ex4} in \eqref{eq-ex2} we obtain

$$
E\vert  \xi \vert^u \leq \sqrt{2}\sigma^u  e^{-\frac{u}2}(u+1)^{\frac{u}2}
e^{\frac1{6(u+1)}-\frac1{2}}.
$$
Since the value $u$ is always greater or equal to one, $u\geq1,$
then
\begin{eqnarray*}
\left(E\vert  \xi \vert^u\right)^{\frac1{u}} & \leq &  2^{\frac1{2u}}\sigma
e^{-\frac{1}2}(u+1)^{\frac{1}2} e^{\frac1{6u(u+1)}-\frac1{2u}}\\
&\leq &  \sqrt{ 2}\sigma e^{-\frac{1}2}(u+1)^{\frac{1}2}
e^{\frac1{12}}.
\end{eqnarray*}
For $u\geq 1$ we have $(u+1)^\frac12\leq (2u)^\frac12$. Finally, we obtain
that
$$
\left(E\vert  \xi \vert^u\right)^{\frac1{u}}  \leq  2\sigma u^{\frac{1}2}
e^{-\frac5{12}}.
$$
So, for a normal centered random variable $\xi$ the function
$\psi(u)$ can be taken as $\psi(u)=\sqrt{u}$ and the norm of $\xi$
in the space $\mathbf{F}_\psi(\Omega)$ is estimated in the following
way
\begin{equation}\label{eq-ex5}
\left\|\xi\right\|_\psi \leq  R\cdot \left(E \xi^2\right)^{1/2},
\end{equation}
where the constant $R=2e^{-\frac5{12}}\approx 1.32.$ It means that
Gaussian random variables belong to the space
$\mathbf{F}_\psi(\Omega)$ with $\psi(u)=\sqrt{u}$.
\end{example}

\begin{example}\label{ex2}
Suppose that a random variable $\xi$ has a symmetric exponential
distribution with the rate $\lambda$. Its pdf $f(x,\lambda)$ is
$f(x,\lambda)=\frac{\lambda}2 e^{-\lambda\vert  x \vert},\,
\lambda>0$. It is easy to show that $\xi$ belongs to
$\mathbf{F}_\psi(\Omega)$ with $\psi(u)=u.$ In fact, similarly to
the previous example apply inequality for Gamma function from
\cite{askey2010}

$$
E\vert  \xi \vert^u  = \frac{ \Gamma (u+1)}{\lambda^u}\leq
\frac{\sqrt{2\pi}}{\lambda^u}e^{-(u+1)}(u+1)^{u+\frac12}e^{\frac1{12(u+1)}}.
$$
Then
\begin{eqnarray*}
\left(E\vert  \xi \vert^u\right)^{\frac1{u}}& \leq &
\frac{(2\pi)^{\frac1{2u}}}{\lambda}e^{-(1+1/u)}(u+1)^{1+\frac1{2u}}e^{\frac1{12u(u+1)}}\\
&\leq& \frac{(2\pi)^{\frac1{2}} e^{-1}}{\lambda}(2u) (1+u)^{\frac1{2u}}e^{\frac1{12u(u+1)}-\frac1{u}},\\
 &\leq& \frac{4(\pi)^{\frac1{2}} e^{-1}}{\lambda}u e^{\frac1{12u(u+1)}-\frac1{u}},\quad \text{for}\quad u\geq 1.\\
\end{eqnarray*}
Since the function $\frac1{12u(u+1)}-\frac1{u}$ is increasing for
$u\geq 1$ and
$$\lim\limits_{u\to+\infty}\left(\frac1{12u(u+1)}-\frac1{u}\right)=0,$$
then $ e^{\frac1{12u(u+1)}-\frac1{u}}\leq 1. $ Hence,
$$
\left(E\vert  \xi \vert^u\right)^{\frac1{u}}\leq \frac{4\sqrt{\pi}
e^{-1}}{\lambda}u.
$$
So, for a symmetric exponential random variable $\xi$ the function
$\psi(u)$ can be taken as $\psi(u)=u $ and the norm of $\xi$ in the
space $\mathbf{F}_\psi(\Omega)$ is estimated in the following way
\begin{equation}\label{eq-ex5-exp}
\left\|\xi\right\|_\psi \leq  R\cdot \left(E \xi^2\right)^{1/2},
\end{equation}
where the constant $R=4\sqrt{\pi}e^{-1}\approx 2.61.$ It means that
exponential type random variables belong to the space
$\mathbf{F}_\psi(\Omega)$ with $\psi(u)=u$.
\end{example}

\begin{definition}
The space $\mathbf{F}_\psi(\Omega)$ will be called the space $\mathbf{\check{F}}_\psi(\Omega)$ if the condition is satisfied for the function $\psi(u)$:
\begin{equation} \label{eq:o1chapter_9}
\sup\limits_{u\geq1}\frac{\psi(u+v)}{\psi(u)}<\infty,
\end{equation}
where $v>0$ is any number.
\end{definition}

It is obvious that the condition \eqref{eq:o1chapter_9} is satisfied for the functions:

1) $\psi(u)=e^{au^\beta}$, where $0<\beta\leq1,\ a>0$;

2) $\psi(u)=Au^\alpha$, where $\alpha>0,\ A>0$.

\begin{definition}\label{eq:o2chapter-3}
A non-decreasing numerical sequence $\left(\varkappa(n),n\geq1\right)$ is called an $M$-characteristic (majorizing characteristic) of the space $\mathbf{F}_\psi(\Omega)$,
if for any random variables $\xi_i,\ i=1,2,\ldots,n$ from this space the inequality holds:
\begin{equation} \label{eq:o1chapter_10}
\left\|\max\limits_{1\leq i\leq n}\left|\xi_i\right|\right\|_\psi\leq\varkappa(n)\max\limits_{1\leq i\leq n}\left\|\xi_i\right\|_\psi.
\end{equation}
\end{definition}

Similarly, the majorizing characteristic is defined in all Banach spaces.

\begin{theorem} \label{th:o1chapter_5}
The sequence
\begin{equation} \label{eq:o1chapter_11}
\varkappa(n)=\sup\limits_{u\geq 1}\inf\limits_{v>0}n^{\frac{1}{u+v}}\frac{\psi(u+v)}{\psi(u)}
\end{equation}
is the majorizing characteristic of the space $\mathbf{F}_\psi(\Omega)$.
\end{theorem}

\begin{proof}
Let $\xi_1,\xi_2,\ldots,\xi_n$ be random variables from the space $\mathbf{F}_\psi(\Omega)$, then the chain of inequalities holds:
\begin{multline*}
\frac{\left(E\left(\max\limits_{1\leq i\leq n}\left|\xi_i\right|\right)^{u}\right)^{1/u}}{\psi(u)}\leq \frac{\left(E\left(\max\limits_{1\leq i\leq n}\left|\xi_i\right|\right)^{u+v}\right)^{\frac{1}{u+v}}}{\psi(u)}\leq \\[1ex]
\leq\max\limits_{1\leq i\leq n} n^{\frac{1}{u+v}}\frac{\left(E\left|\xi_i\right|^{u+v}\right)^{\frac{1}{u+v}}}{\psi(u+v)}\cdot \frac{\psi(u+v)}{\psi(u)}\leq \max\limits_{1\leq i\leq n}\left\|\xi_i\right\|_\psi  n^{\frac{1}{u+v}}\frac{\psi(u+v)}{\psi(u)}.
\end{multline*}
From the last inequality it follows that for all $u\geq1$
\[
\frac{\left(E\left(\max\limits_{1\leq i\leq n}\left|\xi_i\right|\right)^{u}\right)^{1/u}}{\psi(u)}\leq \max\limits_{1\leq i\leq n}\left\|\xi_i\right\|_\psi \inf\limits_{v>0}n^{\frac{1}{u+v}}\frac{\psi(u+v)}{\psi(u)}.
\]
Therefore,
\[
\left\|\max\limits_{1\leq i\leq n}\left|\xi_i\right|\right\|_\psi=\sup\limits_{u\geq1}\frac{\left(E\left(\max\limits_{1\leq i\leq n}\left|\xi_i\right|\right)^{u}\right)^{1/u}}{\psi(u)}\leq
\]
\[
\leq\max\limits_{1\leq i\leq n}\left\|\xi_i\right\|_\psi\sup\limits_{u\geq 1}\inf\limits_{v>0}n^{\frac{1}{u+v}}\frac{\psi(u+v)}{\psi(u)}.
\]
Since in this inequality $v$ is any number ($v>0$), then the inequality \eqref{eq:o1chapter_10} implies the statement of the Theorem.
\end{proof}

For the spaces $\mathbf{\check{F}}_\psi(\Omega)$, the formula for calculating the majorizing characteristic $\varkappa(n)$ can be simplified, i.e. the following Corollary is valid.

\begin{corollary} \label{co:o1chapter_1}
The sequence
\begin{equation} \label{eq:o1chapter_12}
\varkappa(n)=\inf\limits_{v>0}z(v)n^{\frac{1}{v+1}},
\end{equation}
where $z(v)=\sup\limits_{u\geq 1}\frac{\psi(u+v)}{\psi(u)}$, is a majorizing characteristic of the space $\mathbf{\check{F}}_\psi(\Omega)$.
\end{corollary}

\begin{proof}
Let $\xi_1,\xi_2,\ldots,\xi_n$ be random variables from the space $\mathbf{\check{F}}_\psi(\Omega)$, then the chain of inequalities holds:\\
$
\left\|\max\limits_{1\leq i\leq n}\left|\xi_i\right|\right\|_\psi=\sup\limits_{u\geq1}\frac{\left(E\left(\max\limits_{1\leq i\leq n}\left|\xi_i\right|\right)^{u}\right)^{1/u}}{\psi(u)}\leq
\sup\limits_{u\geq 1}\frac{\left(E\left(\max\limits_{1\leq i\leq n}\left|\xi_i\right|\right)^{u+v}\right)^{\frac{1}{u+v}}}{\psi(u)}\leq
$ \\
$
\leq\sup\limits_{u\geq 1}\max\limits_{1\leq i\leq n} n^{\frac{1}{u+v}}\frac{\left(E\left|\xi\right|^{u+v}\right)^{\frac{1}{u+v}}}{\psi(u+v)}\cdot
\frac{\psi(u+v)}{\psi(u)}\leq \max\limits_{1\leq i\leq n}\left\|\xi_i\right\|_\psi \sup\limits_{u\geq 1}n^{\frac{1}{u+v}}\frac{\psi(u+v)}{\psi(u)}\leq
$ \\
$
\leq\max\limits_{1\leq i\leq n}\left\|\xi_i\right\|_\psi n^{\frac{1}{v+1}}\sup\limits_{u\geq 1}\frac{\psi(u+v)}{\psi(u)} \leq \max\limits_{1\leq i\leq n}\left\|\xi_i\right\|_\psi \inf\limits_{v>0} n^{\frac{1}{v+1}}\sup\limits_{u\geq 1}\frac{\psi(u+v)}{\psi(u)}.
$

\noindent Since this inequality holds for any $v>0$, then the statement of the Corollary~\ref{co:o1chapter_1} follows from the inequality \eqref{eq:o1chapter_10}.
\end{proof}

\begin{theorem}
\label{th:o1chapter_6}
The sequence
\[
\varkappa(n)=\frac{1}{e^a}\exp \left\{S(a,\beta)(\ln n)^{\frac{\beta}{\beta+1}}\right\},
\]
where $S(a,\beta)=(\beta a)^{\frac{1}{\beta+1}}(\beta^{-1}+1)$ is a majorizing characteristic of the space $\mathbf{F}_\psi(\Omega)$, where $\psi(u)=e^{au^\beta}, a>0$, $\beta>0$, and $\varkappa(1)=1$.
\end{theorem}

\begin{proof}
Let us consider two cases. In the first case, when $0<\beta\leq 1$, it follows from the Corollary \ref{co:o1chapter_1} that
\[
z(v)=\sup\limits_{u\geq 1}\frac{\psi(u+v)}{\psi(u)}=\sup\limits_{u\geq 1} e^{a(u+v)^\beta-au^\beta}=e^{a\left((1+v)^\beta-1\right)},
\]
then, according to the equality \eqref{eq:o1chapter_12}, we have:
\begin{equation} \label{eq:o1chapter_13}
\varkappa(n)=\inf\limits_{v>0}e^{a\left((1+v)^\beta-1\right)}n^{\frac{1}{v+1}}.
\end{equation}
Since, the equality is true
\[
\left(\ln\left(e^{a\left((1+v)^\beta-1\right)}n^{\frac{1}{v+1}}\right)\right)^{'}= \left(\frac{\ln n}{v+1}+a(1+v)^\beta-a\right)^{'}=
\]
\[
=-\frac{\ln n}{(v+1)^2}+a\beta(1+v)^{\beta-1}=0,
\]
then the infinitude is reached at the point $v=\left(\frac{\ln n}{a\beta}\right)^{\frac{1}{\beta+1}}-1$ and, substituting this value of $v$ into the equality \eqref{eq:o1chapter_13}, we obtain:
\[
\varkappa(n)=n^{\left(\frac{a\beta}{\ln n}\right)^{\frac{1}{\beta+1}}} \exp \left\{a\left(\left(\frac{a\beta}{\ln n}\right)^{\frac{\beta}{\beta+1}}-1\right)\right\}=\frac{1}{e^a}\exp \left\{S(a,\beta)(\ln n)^{\frac{\beta}{\beta+1}}\right\}.
\]
In the second case, when $\beta>1$, from Theorem \ref{th:o1chapter_5} we have:
\begin{equation} \label{eq:o1chapter_13_1}
\varkappa(n)=\sup\limits_{u\geq 1}\inf\limits_{v>0}n^{\frac{1}{u+v}}\frac{\psi(u+v)}{\psi(u)}= \sup\limits_{u\geq 1}\inf\limits_{v>0}n^{\frac{1}{u+v}} e^{a(u+v)^\beta-au^\beta}.
\end{equation}
First, let's find the infinitude. Let's consider the equalities
\[
\left(\ln\left(e^{a\left((u+v)^\beta-1\right)}n^{\frac{1}{u+v}}\right)\right)^{'}= \left(\frac{\ln n}{u+v}+a(u+v)^\beta-au\right)^{'}=
\]
\[
=-\frac{\ln n}{(u+v)^2}+a\beta(u+v)^{\beta-1}=0.
\]
So, from the last equality it follows that the infinitude is reached at the point $v=\left(\frac{\ln n}{a\beta}\right)^{\frac{1}{\beta+1}}-u$. Thus, substituting this value of $v$ into the equality \eqref{eq:o1chapter_13_1}, we have:
\[
\varkappa(n)=\sup\limits_{u\geq 1}n^{\left(\frac{a\beta}{\ln n}\right)^{\frac{1}{\beta+1}}} \exp \left\{a\left(\left(\frac{a\beta}{\ln n}\right)^{\frac{\beta}{\beta+1}}-u\right)\right\}=
\]
\[
=n^{\left(\frac{a\beta}{\ln n}\right)^{\frac{1}{\beta+1}}} \exp \left\{a\left(\left(\frac{a\beta}{\ln n}\right)^{\frac{\beta}{\beta+1}}-1\right)\right\}=
\frac{1}{e^a}\exp \left\{S(a,\beta)(\ln n)^{\frac{\beta}{\beta+1}}\right\}.
\]
The statement of the Theorem follows from the two cases considered.
\end{proof}

\begin{remark}
Theorem \ref{th:o1chapter_6} illustrates that the equality \eqref{eq:o1chapter_12} cannot always be applied.
\end{remark}

\begin{theorem}
\label{th:o1chapter_7}
The sequence
\[
\varkappa(n)=\left(\frac{e}{\alpha}\right)^\alpha\left(\ln n\right)^\alpha
\]
is a majorizing characteristic of the space $\mathbf{F}_\psi(\Omega)$ for $n>1$, where $\psi(u)=u^\alpha, \alpha>0$, and $\varkappa(1)=1$.
\end{theorem}

\begin{proof}
From the Corollary \ref{co:o1chapter_1} it follows that
\[
z(v)=\sup\limits_{u\geq 1} \left(1+\frac{v}{u}\right)^\alpha=(1+v)^\alpha,
\]
then from the equality \eqref{eq:o1chapter_12} we have:
\begin{equation} \label{eq:o1chapter_14}
\varkappa(n)=\inf\limits_{v>0} (1+v)^\alpha n^{\frac{1}{v+1}}.
\end{equation}
From the equalities
\[
\left(\ln\left((1+v)^\alpha n^{\frac{1}{v+1}}\right)\right)^{'}= \left(\alpha \ln(1+v)+\frac{\ln n}{1+v}\right)^{'}=
\]
\[
=\frac{\alpha}{1+v}-\frac{\ln n}{(1+v)^2}=0
\]
we obtain that the infinitude is reached at the point $v=\frac{\ln n}{\alpha}-1$ and, substituting this value of $v$ into the equality \eqref{eq:o1chapter_14}, the statement of the Theorem becomes obvious:
\[
\varkappa(n)=\left(\frac{\ln n}{\alpha}\right)^\alpha n^{\frac{\alpha}{\ln n}}=\left(\ln n\right)^\alpha\left(\frac{e}{\alpha}\right)^\alpha.
\]
\end{proof}

\begin{theorem}\label{th:o1chapter-8}
The sequence
\[
\varkappa(n)=e \left(\frac{\ln(\ln n+2)}{\ln 2}\right)^\lambda
\]
is a majorizing characteristic of the space $\mathbf{F}_\psi(\Omega)$, where $\psi(u)=\left(\ln (u+1)\right)^\lambda$, $\lambda>0$ for $n>1$, and $\varkappa(1)=1$.
\end{theorem}

\begin{proof}
From the Corollary \ref{co:o1chapter_1} it follows that
\[
z(v)=\sup\limits_{u\geq 1} \left(\frac{\ln(u+v+1)}{\ln(u+1)}\right)^\lambda,
\]
and since $\left(\frac{\ln(u+v+1)}{\ln(u+1)}\right)'\leq 0$, then $z(v)=\left(\frac{\ln(v+2)}{\ln 2}\right)^\lambda$. Hence we have:
\begin{equation} \label{eq:o1chapter_15}
\varkappa(n)=\inf\limits_{v>0}\left(\frac{\ln(v+2)}{\ln 2}\right)^\lambda n^{\frac{1}{v+1}}.
\end{equation}

In the equality \eqref{eq:o1chapter_15} we put $v=\ln n$: \\
$
\varkappa(n)=\left(\frac{\ln(\ln n+2)}{\ln 2}\right)^\lambda n^{\frac{1}{\ln n+1}}\leq \left(\frac{\ln(\ln n+2)}{\ln 2}\right)^\lambda e,
$
which was to be proved.
\end{proof}

\section{Spaces $\mathbf{F}_{S_k,\psi,r}(\Omega)$ of random variables}

\begin{definition}
Let $S_k$ be an increasing numerical sequence ($S_k\geq 1$) and $S_k\to \infty$ as $k\to \infty$. Consider a monotonically increasing continuous function $\psi(u)>0$, $u\geq1$ such that $\psi(u)\to\infty$ as $u\to\infty$. A random variable $\xi$ belongs to the space $\mathbf{F}_{S_k,\psi,r}(\Omega)$ if the condition is satisfied:
\[
\sup\limits_{k\geq r}\frac{\left(E\left|\xi\right|^{S_k}\right)^{1/{S_k}}}{\psi(S_k)}<\infty,
\]
where the number $r$ is such that $S_r\geq1$.
\end{definition}

As in the previous case, it is easy to prove that the spaces $\mathbf{F}_{S_k,\psi,r}(\Omega)$ are Banach spaces with the norms
\begin{equation} \label{eq:o1chapter_16}
\left\|\xi\right\|_{S_k,\psi,r}=\sup\limits_{k\geq r}\frac{\left(E\left|\xi\right|^{S_k}\right)^{1/{S_k}}}{\psi(S_k)}.
\end{equation}

\begin{theorem} \label{th:o1chapter_9}
If the function $\psi$ satisfies the condition \eqref{eq:o1chapter_9} and there exists such $B_r>0$ that
\[
\frac{\psi\left(S_k\right)}{\psi\left(S_{k-1}\right)}\leq B_r,\ k\geq r,
\]
then the spaces $\mathbf{F}_{S_k,\psi,r}(\Omega)$ contain the same elements as the spaces $\mathbf{\check{F}}_\psi(\Omega)$, and the norms \eqref{eq:o1chapter_1} and \eqref{eq:o1chapter_16} are equivalent and the following inequalities hold:
\[
\left\|\xi\right\|_{S_k,\psi,r}\leq \left\|\xi\right\|_\psi,
\]
\[
\left\|\xi\right\|_\psi\leq \max \left(C_r, \widetilde{C}_r\right) \left\|\xi\right\|_{S_k,\psi,r},
\]
where $C_r= \sup\limits_{k\geq r}\frac{\psi\left(S_k\right)}{\psi\left(S_{k-1}\right)}$, $\widetilde{C}_r=\sup\limits_{1\leq u\leq S_r}\frac{\psi(S_r)}{\psi(u)}$.
\end{theorem}

\begin{proof}
It is obvious that
\[
\left\|\xi\right\|_{S_k,\psi,r}\leq \left\|\xi\right\|_\psi.
\]
On the other hand, from the Lyapunov inequality, for $S_{k-1}\leq u \leq S_k$, where $k-1\geq r$ it follows that
\begin{gather*}
\frac{\left(E\left|\xi\right|^{u}\right)^{1/u}}{\psi(u)}\leq \frac{\left(E\left|\xi\right|^{S_k}\right)^{1/{S_k}}}{\psi(u)}= \frac{\left(E\left|\xi\right|^{S_k}\right)^{1/{S_k}}}{\psi(S_k)}\cdot \frac{\psi(S_k)}{\psi(u)}\leq \\[1ex]
\leq\left\|\xi\right\|_{S_k,\psi,r}\cdot\frac{\psi(S_k)}{\psi(u)}\leq\left\|\xi\right\|_{S_k,\psi,r}\cdot\frac{\psi\left(S_k\right)}{\psi\left(S_{k-1}\right)}\leq C_r \left\|\xi\right\|_{S_k,\psi,r}.
\end{gather*}
For $1\leq u \leq S_r$ we have:
\[
\frac{\left(E\left|\xi\right|^{u}\right)^{1/u}}{\psi(u)}\leq \frac{\left(E\left|\xi\right|^{S_r}\right)^{1/{S_r}}}{\psi(S_r)}\cdot \frac{\psi(S_r)}{\psi(u)}\leq \widetilde{C}_r \left\|\xi\right\|_{S_k,\psi,r}.
\]
Then
\[
\left\|\xi\right\|_\psi\leq \max \left(\widetilde{C}_r,C_r\right) \left\|\xi\right\|_{S_k,\psi,r}.
\]
\end{proof}

Among the spaces $\mathbf{F}_{S_k,\psi,r}(\Omega)$, the most important for us is the space where $S_k=2k$. We denote the norm in this space by $\left\|\xi\right\|_{2k,\psi,r}=\sup\limits_{k\geq r}\frac{\left(E\left|\xi\right|^{2k}\right)^{1/2k}}{\psi(2k)}$. It is obvious that in cases where the condition \eqref{eq:o1chapter_9} is satisfied for $\psi$, then the spaces $\mathbf{\check{F}}_\psi(\Omega)$ and $\mathbf{F}_\psi(\Omega)$ coincide, and the norms in these spaces are equivalent.

Indeed, according to the previous Theorem, we have that
\begin{equation} \label{eq:o1chapter_17}
\left\|\xi\right\|_{2k,\psi,r}\leq \left\|\xi\right\|_\psi.
\end{equation}
Note that
\[
,
\]
that is, it takes place
\begin{equation} \label{eq:o1chapter_18}
\left\|\xi\right\|_\psi\leq \widehat{\psi}_r\left\|\xi\right\|_{2k,\psi,r},
\end{equation}
where $\widehat{\psi}_r=\max \left(\overline{\psi}_r,\widetilde{C}_r\right)$.

The following Theorem is proved similarly to the proof of Theorem \ref{th:o1chapter_1}.

\begin{theorem}\label{th:o1chapter_10}
Let the random variable $\xi$ belong to the space \\ $\mathbf{F}_{S_k,\psi,r}(\Omega)$, then for any $\varepsilon>0$ the inequality holds:
\begin{equation} \label{eq:o1chapter_19}
P\left\{\left|\xi\right|>\varepsilon\right\}\leq \inf\limits_{k\geq r}\frac{\left\|\xi\right\|_{S_k,\psi,r}^{S_k}(\psi(S_k))^{S_k}}{\varepsilon^{S_k}}.
\end{equation}

In particular, if $S_k=2k$, we get:
\[
P\left\{\left|\xi\right|>\varepsilon\right\}\leq \inf\limits_{k\geq r}\frac{\left\|\xi\right\|_{2k,\psi,r}^{2k}(\psi(2k))^{2k}}{\varepsilon^{2k}}.
\]
\end{theorem}

\begin{theorem}\label{th:o1chapter_11}
The sequence
\begin{equation} \label{eq:o1chapter_20}
\varkappa(n)=\sup\limits_{k\geq r}\inf\limits_{v>0}n^{\frac{1}{S_k+v}}\frac{\psi(S_k+v)}{\psi(S_k)}
\end{equation}
is the majorizing characteristic of the space $\mathbf{F}_{S_k,\psi,r}(\Omega)$.
\end{theorem}

\begin{proof}
Similarly to the proof of Theorem \ref{th:o1chapter_5}, consider the random variables \linebreak $\xi_1,\xi_2,\ldots,\xi_n$ from the space $\mathbf{F}_{S_k,\psi,r}(\Omega)$. Then the following inequalities hold:
\begin{gather*}
\frac{\left(E\left(\max\limits_{1\leq i\leq n}\left|\xi_i\right|\right)^{S_k}\right)^{1/S_k}}{\psi(S_k)}\leq \frac{\left(E\left(\max\limits_{1\leq i\leq n}\left|\xi_i\right|\right)^{S_k+v}\right)^{\frac{1}{S_k+v}}}{\psi(S_k)}\leq \\
\leq\max\limits_{1\leq i\leq n} n^{\frac{1}{S_k+v}}\frac{\left(E\left|\xi_i\right|^{S_k+v}\right)^{\frac{1}{S_k+v}}}{\psi(S_k+v)}\cdot \frac{\psi(S_k+v)}{\psi(S_k)}\leq \\
\leq\max\limits_{1\leq i\leq n}\left\|\xi_i\right\|_{S_k,\psi,r}  n^{\frac{1}{S_k+v}}\frac{\psi(S_k+v)}{\psi(S_k)}.
\end{gather*}
From the last inequality it follows that for all $S_k\geq1$
\[
\frac{\left(E\left(\max\limits_{1\leq i\leq n}\left|\xi_i\right|\right)^{S_k}\right)^{1/S_k}}{\psi(S_k)}\leq \max\limits_{1\leq i\leq n}\left\|\xi_i\right\|_{S_k,\psi,r} \inf\limits_{v>0}n^{\frac{1}{S_k+v}}\frac{\psi(S_k+v)}{\psi(S_k)}.
\]
Therefore,
\[
\left\|\max\limits_{1\leq i\leq n}\left|\xi_i\right|\right\|_{S_k,\psi,r}=\sup\limits_{k\geq r}\frac{\left(E\left(\max\limits_{1\leq i\leq n}\left|\xi_i\right|\right)^{S_k}\right)^{1/S_k}}{\psi(S_k)}\leq
\]
\[
\leq\max\limits_{1\leq i\leq n}\left\|\xi_i\right\|_{S_k,\psi,r}\sup\limits_{k\geq r}\inf\limits_{v>0}n^{\frac{1}{S_k+v}}\frac{\psi(S_k+v)}{\psi(S_k)}.
\]
Since in the last inequality $v$ is any number ($v>0$), then the inequality \eqref{eq:o1chapter_10} implies the statement of the Theorem.
\end{proof}

\section{Condition $\mathbf{H}$ for spaces $\mathbf{F}_\psi(\Omega)$}

\begin{definition} \label{de:o1chapter_5}
We say that for Banach spaces $B(\Omega)$ of random variables the condition $\mathbf{H}$ is satisfied if there exists an absolute constant $C_B$ such that for any centered independent random variables $\xi_1,\xi_2,\ldots,\xi_n$ from $B(\Omega)$ the inequality holds:
\begin{equation} \label{eq:o1chapter_21}
\left\|\sum_{i=1}^{n} \xi_i\right\| ^2\leq C_B \sum_{i=1}^{n}\left\|\xi_i\right\|^2.
\end{equation}
\end{definition}
The constant $C_B$ is called the scale constant of the space $B(\Omega)$. For the spaces $\mathbf{F}_\psi(\Omega)$, the constant $C_{\mathbf{F}_\psi(\Omega)}$ will be denoted by $C_\psi$.

Let us find the conditions under which the condition $\bf{H}$ is satisfied for the spaces $\mathbf{F}_\psi(\Omega)$, and also find the value of the constant $C_\psi$.

\begin{theorem} \label{th:o1chapter_12}
Let $\xi_1,\xi_2,\ldots,\xi_n$ be independent centered random variables from the space $\mathbf{F}_{2k,\psi,r}(\Omega)$. If $\xi_i$ are symmetric random variables and for $k\geq \max (r,2)$ the condition is satisfied
\begin{equation} \label{eq:o1chapter_22}
C_{2k}^{2l}\frac{\left(\psi(2l)\right)^{2l}\left(\psi(2k-2l)\right)^{2k-2l}}{\left(\psi(2k)\right)^{2k}}\leq C_{k}^{l}, \ l=\overline{1,k-1},
\end{equation}
then the inequality holds
\begin{equation} \label{eq:o1chapter_23}
\left\|\sum\limits_{i=1}^n\xi_i\right\|_{2k,\psi,r}^2\leq\sum\limits_{i=1}^n\left\|\xi_i\right\|_{2k,\psi,r}^2.
\end{equation}
That is, in this case, the condition $\mathbf{H}$ with the constant $C_\psi=1$ is satisfied for the space $\mathbf{F}_{2k,\psi,r}(\Omega)$.

If we abandon the symmetry condition, then from the condition \eqref{eq:o1chapter_22} we have the inequality:
\begin{equation} \label{eq:o1chapter_24}
\left\|\sum\limits_{i=1}^n\xi_i\right\|_{2k,\psi,r}^2\leq4\sum\limits_{i=1}^n\left\|\xi_i\right\|_{2k,\psi,r}^2.
\end{equation}
That is, in this case, for the space $\mathbf{F}_{2k,\psi,r}(\Omega)$, the condition $\mathbf{H}$ with the constant $C_\psi=4$ is satisfied.

If $\xi_i$ are not symmetric and the condition is fulfilled
\begin{equation} \label{eq:o1chapter_25}
C_{2k}^{2l}\left(1+\frac{k}{3}\right)\frac{\left(\psi(2l)\right)^{2l}\left(\psi(2k-2l)\right)^{2k-2l}}{\left(\psi(2k)\right)^{2k}}\leq C_{k}^{l}, \ l=\overline{1,k-1},
\end{equation}
then we have the inequality
\begin{equation} \label{eq:o1chapter_26}
\left\|\sum\limits_{i=1}^n\xi_i\right\|_{2k,\psi,r}^2\leq\sum\limits_{i=1}^n\left\|\xi_i\right\|_{2k,\psi,r}^2,
\end{equation}
i.e. in this case, for the space $\mathbf{F}_{2k,\psi,r}(\Omega)$, the condition $\mathbf{H}$ with the constant $C_\psi=1$ is satisfied.
\end{theorem}

In order to prove the Theorem \ref{th:o1chapter_12}, it is necessary to prove an additional statement.

\begin {lemma}
\label{le:o1chapter_1}
Let $\xi$ and $\eta$ be random variables from the space $\mathbf{F}_{2k, \psi, r}(\Omega)$. If $\xi$ and $\eta$ are independent and $E\eta=0$, then the inequality holds:
\begin{equation} \label{eq:o1chapter_27}
\left\|\xi\right\|_{2k, \psi, r}\leq\left\|\xi-\eta\right\|_{2k, \psi, r}.
\end{equation}
\end{lemma}

\begin{proof}
From Fubini's Theorem it follows that for $2k>1$
\begin{equation} \label{eq:o1chapter_28}
E\left|\xi-\eta\right|^{2k}=E_\xi\left(E_\eta\left|\xi-\eta\right|^{2k}\right),
\end{equation}
where $E_\xi$ is the mathematical expectation with respect to $\xi$, and $E_\eta$ is the mathematical expectation with respect to $\eta$. From the Lyapunov inequality it follows that for $2k\geq1$
\[
E_\eta\left|\xi-\eta\right|^{2k}\geq\left(E_\eta\left|\xi-\eta\right|\right)^{2k}\geq\left|E_\eta\left(\xi-\eta\right)\right|^{2k}=\left|\xi-E\eta\right|^{2k}=\left|\xi\right|^{2k}.
\]

Therefore, from the equality \eqref{eq:o1chapter_28} we obtain that
\[
E\left|\xi-\eta\right|^{2k}\geq E\left|\xi\right|^{2k}.
\]
From the last inequality it is obvious that the inequality \eqref{eq:o1chapter_27} follows.
\end{proof}

\begin{proof}[Proof of the Theorem]
If $\xi_1,\xi_2,\ldots,\xi_n$ are symmetric random variables, then all their odd moments are zero. Therefore,
\[
E\left(\xi_1+\xi_2\right)^{2k}=E\xi_1^{2k}+\sum\limits_{s=2}^{2k-2}C_{2k}^sE\xi_1^s\xi_2^{2k-s}+E\xi_2^{2k}=
\]
\[
=E\xi_1^{2k}+\sum\limits_{r=1}^{k-1}C_{2k}^{2r}E\xi_1^{2r}E\xi_2^{2k-2r}+E\xi_2^{2k}.
\]
Since $E\left|\xi_i\right|^{2k}\leq\left(\psi(2k)\right)^{2k}\left\|\xi_i\right\|_{2k,\psi,r}^{2k}$, then
\begin{gather*}
\frac{E\left(\xi_1+\xi_2\right)^{2k}}{\left(\psi(2k)\right)^{2k}}\leq\left\|\xi_1\right\|_{2k,\psi,r}^{2k}+\\[1ex]
+\sum\limits_{r=1}^{k-1}C_{2k}^{2r}\left(\psi(2r)\right)^{2r}\left(\psi(2k-2r)\right)^{2k-2r}\frac{\left\|\xi_1\right\|_{2k,\psi,r}^{2r}\left\|\xi_2\right\|_{2k,\psi,r}^{2k-2r}}{\left(\psi(2k)\right)^{2k}}+\left\|\xi_2\right\|_{2k,\psi,r}^{2k}\leq\\[1ex]
\leq\left\|\xi_1\right\|_{2k,\psi,r}^{2k}+\sum\limits_{r=1}^{k-1}C_{k}^{r}\left\|\xi_1\right\|_{2k,\psi,r}^{2r}\left\|\xi_2\right\|_{2k,\psi,r}^{2k-2r}+\left\|\xi_2\right\|_{2k,\psi,r}^{2k}=\\[1ex]
=\left(\left\|\xi_1\right\|_{2k,\psi,r}^2+\left\|\xi_2\right\|_{2k,\psi,r}^2\right)^k.
\end{gather*}
The last inequality for $n=2$ implies the inequality \eqref{eq:o1chapter_23}. Therefore, for any $n$ the inequality \eqref{eq:o1chapter_23} is also valid.

Let $\xi_1,\xi_2,\ldots,\xi_n$ be independent centered random variables from the space $\mathbf{F}_{2k,\psi,r}(\Omega)$, and $\xi_1^*,\xi_2^*,\ldots,\xi_n^*$ be independent random variables that have the same distribution as $\xi_i$ and do not depend on $\xi_i$ ($i=\overline{1,n}$). The random variables $\xi_i-\xi_i^*$ are symmetric. From the lemma \ref{le:o1chapter_1} it follows that
\begin{gather*}
\left\|\sum\limits_{i=1}^n\xi_i\right\|_{2k,\psi,r}^2\leq\left\|\sum\limits_{i=1}^n(\xi_i-\xi_i^*)\right\|_{2k,\psi,r}^2\leq
\sum\limits_{i=1}^n \left\|(\xi_i-\xi_i^*)\right\|_{2k,\psi,r}^2\leq \\[1ex]
\leq\sum\limits_{i=1}^n \left(\left\|\xi_i\right\|_{2k,\psi,r}+\left\|\xi_i^*\right\|_{2k,\psi,r}\right)^2=4\sum\limits_{i=1}^n \left\|\xi_i\right\|_{2k,\psi,r}^2,
\end{gather*}
since $\left\|\xi_i\right\|_{2k,\psi,r}=\left\|\xi_i^*\right\|_{2k,\psi,r}$. The inequality \eqref{eq:o1chapter_24} is proved.

The inequality \eqref{eq:o1chapter_26} is proved similarly to the inequality \eqref{eq:o1chapter_23}. Since
\[
E\left(\xi_1+\xi_2\right)^{2k}=E\xi_1^{2k}+\sum\limits_{s=2}^{2k-2}C_{2k}^sE\xi_1^s\xi_2^{2k-s}+E\xi_2^{2k}
\]
and for odd $s$

\[
\left|E\xi_1^s\xi_2^{2k-s}\right|\leq\frac{1}{2}\left(E\left|\xi_1\right|^{s+1}E\left|\xi_2\right|^{2k-s-1}+E\left|\xi_1\right|^{s-1}E\left|\xi_2\right|^{2k-s+1}\right),
\]
then
\[
E\left(\xi_1+\xi_2\right)^{2k}\leq E\left|\xi_1\right|^{2k}+\sum\limits_{l=1}^{k-1}R_{2k}^{2l}E\left|\xi_1\right|^{2l}E\left|\xi_2\right|^{2k-2l}+E\left|\xi_2\right|^{2k},
\]
where $R_{2k}^{2}=R_{2k}^{2k-2}=C_{2k}^{2}+0,5 C_{2k}^{3}$, $R_{2k}^{2l}=C_{2k}^{2l}+0,5\left(C_{2k}^{2l+1}+C_{2k}^{2l-1}\right)$, $l\neq1,\ l\neq k-1$.
It is obvious that $R_{2k}^{2l}\leq\left(1+\frac{k}{3}\right)C_{2k}^{2l}$. Indeed,
\begin{gather*}
R_{2k}^{2l}=C_{2k}^{2l}+\frac{1}{2}\left(C_{2k}^{2l+1}+C_{2k}^{2l-1}\right)= C_{2k}^{2l}\left(1+\frac{1}{2}\left(\frac{C_{2k}^{2l+1}}{C_{2k}^{2l}}+\frac{C_{2k}^{2l-1}}{C_{2k}^{2l}}\right)\right)=\\[1ex]
=C_{2k}^{2l}\left(1+\frac{1}{2}\left(\frac{(2k)!(2l)! (2k-2l)!}{(2l+1)! (2k-2l-1)! (2k)!}\right.\right.+\\[1ex]
+\left.\left.\frac{(2k)!(2l)! (2k-2l)!}{(2l-1)! (2k-2l+1)! (2k)!}\right)\right)=\\[1ex]
=C_{2k}^{2l}\left(1+\frac{1}{2}\left(\frac{2k-2l}{2l+1}+\frac{2l}{2k-2l+1}\right)\right).
\end{gather*}
From the last equality, since $2l+1\geq 3$ and $2k-2l+1\geq 3$, we obtain:
\[
R_{2k}^{2l}\leq C_{2k}^{2l}\left(1+\frac{1}{2}\left(\frac{2k-2l}{3}+\frac{2l}{3}\right)\right)= C_{2k}^{2l}\left(1+\frac{k}{3}\right).
\]
Thus, the following inequality is valid:
\[
\frac{E\left(\xi_1+\xi_2\right)^{2k}}{\left(\psi(2k)\right)^{2k}}\leq \frac{E\left|\xi_1\right|^{2k}}{\left(\psi(2k)\right)^{2k}}+\sum\limits_{l=1}^{k-1}C_{2k}^{2l}\left(1+\frac{k}{3}\right)\frac{\left(\psi(2l)\right)^{2l}\left(\psi(2k-2l)\right)^{2k-2l}}{\left(\psi(2k)\right)^{2k}}\times
\]
\[
\times\left(\frac{E\left|\xi_1\right|^{2l}}{\left(\psi(2l)\right)^{2l}}\right)\left(\frac{E\left|\xi_2\right|^{2k-2l}}{\left(\psi(2k-2l)\right)^{2k-2l}}\right)+\frac{E\left|\xi_2\right|^{2k}}{\left(\psi(2k)\right)^{2k}}\leq \left\|\xi_1\right\|_{2k,\psi,r}^{2k}+
\]
\[
+\sum\limits_{l=1}^{k-1}C_{k}^{l}\left\|\xi_1\right\|_{2k,\psi,r}^{2l}\left\|\xi_2\right\|_{2k,\psi,r}^{2k-2l}+\left\|\xi_2\right\|_{2k,\psi,r}^{2k}=\left(\left\|\xi_1\right\|_{2k,\psi,r}^2+\left\|\xi_2\right\|_{2k,\psi,r}^2\right)^k.
\]
The last inequality implies the inequality \eqref{eq:o1chapter_26}, for $n=2$. Therefore, for any $n$ the inequality \eqref{eq:o1chapter_26} is also valid. Thus, the statement of the Theorem is proved.
\end{proof}

The following consequence follows from the Theorem \ref{th:o1chapter_12} and the inequalities \eqref{eq:o1chapter_17} and \eqref{eq:o1chapter_18}.

\begin{corollary} \label{co:o2chapter-2}
Let $\xi_1,\xi_2,\ldots,\xi_n$ be independent centered random variables from the space $\mathbf{\check{F}}_\psi(\Omega)$. If $\xi_i$ are symmetric random variables and for $k\geq \max (r,2)$ the condition is satisfied
\begin{equation} \label{eq:o1chapter_29}
C_{2k}^{2l}\frac{\left(\psi(2l)\right)^{2l}\left(\psi(2k-2l)\right)^{2k-2l}}{\left(\psi(2k)\right)^{2k}}\leq C_{k}^{l}, \ l=\overline{1,k-1},
\end{equation}
then the inequality holds
\begin{equation} \label{eq:o1chapter_30}
\left\|\sum\limits_{i=1}^n\xi_i\right\|_\psi^2\leq \widehat{\psi}_r^2\sum\limits_{i=1}^n\left\|\xi_i\right\|_\psi^2.
\end{equation}
That is, in this case, for the space $\mathbf{\check{F}}_\psi(\Omega)$, the condition $\mathbf{H}$ with the constant
$C_\psi=\widehat{\psi}_r^2$, where $\widehat{\psi}_r=\max \left(\overline{\psi}_r,\widetilde{C}_r\right)$, $\overline{\psi}_r=\sup\limits_{u\geq r}\frac{\psi(u+2)}{\psi(u)}$, $\widetilde{C}_r=\sup\limits_{1\leq u\leq S_r}\frac{\psi(S_r)}{\psi(u)}$.

If we abandon the symmetry condition, then under the condition \eqref{eq:o1chapter_29} \linebreak the inequality is true
\begin{equation} \label{eq:o1chapter_31}
\left\|\sum\limits_{i=1}^n\xi_i\right\|^2_\psi\leq 4\widehat{\psi}_r^2\sum\limits_{i=1}^n\left\|\xi_i\right\|_\psi^2.
\end{equation}

That is, in this case, for the space $\mathbf{\check{F}}_\psi(\Omega)$, the condition $\mathbf{H}$ with the constant $C_\psi=4\widehat{\psi}_r^2$ is true.

If $\xi_i$ -- are not symmetric and the condition is true
\begin{equation} \label{eq:o1chapter-33}
C_{2k}^{2l}\left(1+\frac{k}{3}\right)\frac{\left(\psi(2l)\right)^{2l}\left(\psi(2k-2l)\right)^{2k-2l}}{\left(\psi(2k)\right)^{2k}}\leq C_{k}^{l}, \ l=\overline{1,k-1},
\end{equation}
then we have the inequality
\begin{equation} \label{eq:o1chapter-34}
\left\|\sum\limits_{i=1}^n\xi_i\right\|^2_\psi\leq \widehat{\psi}_r^2\sum\limits_{i=1}^n\left\|\xi_i\right\|_\psi^2.
\end{equation}

That is, in this case, for the space $\mathbf{\check{F}}_\psi(\Omega)$ the condition $\mathbf{H}$ with the constant $C_\psi=\widehat{\psi}_r^2$ is satisfied.
\end{corollary}

Let us consider examples of spaces $\mathbf{F}_\psi(\Omega)$ for which the conditions of the previous theorems are satisfied. First of all, let us prove the following lemma.

\begin {lemma} \label{le:o1chapter_2}
The inequality holds
\[
C_{2k}^{2l}\leq C_k^l \frac{k^k}{l^l (k-l)^{k-l}}
\]
for $k\geq 2, \ 1\leq l \leq k-1$.
\end{lemma}

\begin{proof}
Consider the equality
\[
C_{2k}^{2l} = C_k^l \frac{C_{2k}^{2l}}{C_k^l}.
\]
By Stirling's formula $n!=\sqrt{2\pi n} n^n e^{-n} e^{\theta_n}$, where $\left|\theta_n\right|<\frac{1}{12n}$ we obtain:\\
$
\frac{C_{2k}^{2l}}{C_k^l}=\frac{(2k)! l! (k-l)!}{(2l)! (2k-2l)! k!}=\frac{k^{2k}l^l(k-l)^{k-l}}{\sqrt{2} l^{2l}(k-l)^{2(k-l)}k^k} \exp\left\{\theta_{2k}+\theta_{2l}+\theta_k+\theta_l+\theta_{k-l}\right\}\leq $ \\
$\leq\frac{k^k}{l^l (k-l)^{k-l}}\frac{1}{\sqrt{2}}\exp \left\{\frac{1}{24k}+\frac{1}{24l}+\frac{1}{24(k-l)}+\frac{1}{12k}+\frac{1}{12l}+\frac{1}{12(k-l)}\right\}\leq $\\
$\leq\frac{k^k}{l^l (k-l)^{k-l}}\frac{1}{\sqrt{2}}\exp \left\{\frac{1}{8}\left(\frac{1}{k}+\frac{1}{k-1}+1\right)\right\}\leq \frac{k^k}{l^l (k-l)^{k-l}}$,\\
which had to be proven.
\end{proof}

\begin{example}
\label{ex:o1chapter_4}
Consider the space $\mathbf{F}_\psi(\Omega)$, where $\psi(u)=u^\theta$, $\theta\geq \frac{1}{2}$.
Let us prove that in this case the condition \eqref{eq:o1chapter_29} holds. The following relations are valid:
\begin{gather*}
C_{2k}^{2l}\frac{(2l)^{2l\theta}(2k-2l)^{(2k-2l)\theta}}{(2k)^{2k\theta}}=C_{2k}^{2l}\left(\frac{l^{2l}(k-l)^{(2k-2l)}}{k^{2k}}\right)^\theta\leq\\[1ex]
\leq C_k^l \frac{k^k}{l^l (k-l)^{k-l}}\left(\frac{l^{2l}(k-l)^{(2k-2l)}}{k^{2k}}\right)^\theta= C_k^l \left(\frac{l^l(k-l)^{k-l}}{k^k}\right)^{2\theta-1},
\end{gather*}
but since
\[
\left(\frac{l^l(k-l)^{k-l}}{k^k}\right)^{2\theta-1}\leq \left(\frac{l^l}{k^l}\frac{(k-l)^{k-l}}{k^{k-l}}\right)^{2\theta-1}\leq 1,
\]
then $\left(\frac{l^l}{k^l}\right)^{2\theta-1}\leq 1$ and $\left(\frac{(k-l)^{k-l}}{k^{k-l}}\right)^{2\theta-1}\leq 1$.

It is obvious that for $\theta\geq \frac{1}{2}$ and $k>2$ the inequality \eqref{eq:o1chapter_29} holds, i.e. for the space $\mathbf{F}_\psi(\Omega)$,
where $\psi(u)=u^\theta$ the condition $\bf{H}$ is satisfied and the following inequality holds:
\[
\left\|\sum\limits_{i=1}^n\xi_i\right\|^2_\psi\leq 4 \cdot 9^\theta \sum\limits_{i=1}^n\left\|\xi_i\right\|_\psi^2.
\]
Note that for $\theta < \frac{1}{2}$ the condition $\bf{H}$ is not satisfied for this space.
\end{example}

\begin{example}
Consider the space $\mathbf{F}_\psi(\Omega)$, where $\psi(u)=u^\theta$. Find such $\theta$ for which the condition \eqref{eq:o1chapter-33} is satisfied.
That is, from the lemma \ref{le:o1chapter_2} for $k\geq 2$ it follows that
\begin{equation} \label{eq:o1chapter-35}
\left(\frac{l^l(k-l)^{k-l}}{k^k}\right)^{2\theta-1}\left(1+\frac{k}{3}\right)\leq 1.
\end{equation}
The following relations hold:
\[
\left(\frac{l^l(k-l)^{k-l}}{k^k}\right)=\left(\frac{l}{k}\right)^l \left(\frac{k-l}{k}\right)^{k-l}=\left(\frac{l}{k}\right)^l \left(1-\frac{l}{k}\right)^{k-l}=A(k,l),
\]
\[
\ln A(k,l)= l(\ln l-\ln k)+(k-l)(\ln (k-l)-\ln k).
\]
Consider $l$, which varies on the interval $[1,k-1]$, then
\[
\left(\ln A(k,l)\right)^{'}=(\ln l-\ln k)+1-(\ln(k-l)-\ln k)-1=\ln \frac{l}{k}-\ln\frac{k-l}{k}.
\]
Moreover, $\left(\ln A(k,l)\right)^{'}=0$ when $l=\frac{k}{2}$. Therefore, $A(k)=\sup\limits_{l>1}A(k, l)=\left(\frac{1}{2}\right)^k$, then
\eqref{eq:o1chapter-35} is satisfied when
\[
\frac{1}{2^{k(2\theta-1)}}\left(1+\frac{k}{3}\right)\leq 1
\]
for $k\geq 2$ and $\theta>\frac{1}{2}$ and $2\theta-1\geq\sup\limits_{k\geq 2}\frac{\ln\left(1+\frac{k}{3}\right)}{k\ln 2}$.
Since when calculating \linebreak $\sup\limits_{k\geq 2}\frac{\ln\left(1+\frac{k}{3}\right)}{k\ln 2}\leq 0.369$, then $\theta\geq 0.6845$.
\end{example}

\begin{theorem} \label{th:o1chapter_13}
Let us consider the space $\mathbf{F}_\psi(\Omega)$, where the function $\psi(u)$, $u\geq1$ is such that $\varphi(u)=\frac{\psi(u)}{u^{1/2}}$ -- monotonically increasing as $u\geq u_0$.
If $u_0=1$, then for $\mathbf{F}_\psi(\Omega)$ the condition $\bf{H}$ with the constant $C_\psi=4\widehat{\psi}_r^2$ is satisfied, where $\widehat{\psi}_r$
is given in the Corollary \ref{co:o2chapter-2}, and if $u_0>1$, then for $\mathbf{F}_\psi(\Omega)$ the condition $\bf{H}$ with the constant
$C_\psi=4\widehat{\psi}_r^2 S_\psi^2 S_{\hat{\psi}}^2$ is satisfied, where

\[
S_\psi=\max \left(1,\max\limits_{1\leq u < u_0}\frac{\psi\left(u_0\right)}{\psi(u)}\left(\frac{u}{u_0}\right)^{1/2}\right),
\]
\[
S_{\hat{\psi}}=\max \left(1,\max\limits_{1\leq u < u_0}\frac{\psi\left(u\right)}{\psi(u_0)}\left(\frac{u_0}{u}\right)^{1/2}\right).
\]
\end{theorem}

\begin{proof}
First, we prove the Theorem when $u_0=1$. We check the condition \eqref{eq:o1chapter_29}. From the lemma \ref{le:o1chapter_2} it follows
\begin{multline*}
C_{2k}^{2l}\frac{\left(\psi(2l)\right)^{2l}\left(\psi(2k-2l)\right)^{2k-2l}}{\left(\psi(2k)\right)^{2k}}\leq \\[1ex]
\leq C_{k}^{l}\frac{\frac{\left(\psi(2l)\right)^{2l}}{(2l)^l} \cdot \frac{\left(\psi(2k-2l)\right)^{2k-2l}} {(2k-2l)^{k-l}}}  {\frac{\left(\psi(2k)\right)^{2k}}{(2k)^k}}=C_{k}^{l}\frac{\left(\varphi(2l)\right)^{2l}\left(\varphi(2k-2l)\right)^{2k-2l}}{\left(\varphi(2k)\right)^{2k}}.
\end{multline*}
Since the function $\varphi(u)$ monotonically increases for $u\geq 1$, then
\[
\frac{\left(\varphi(2l)\right)^{2l}\left(\varphi(2k-2l)\right)^{2k-2l}}{\left(\varphi(2k)\right)^{2k}}=
\left(\frac{\varphi(2l)}{\varphi(2k)}\right)^{2l}\left(\frac{\varphi(2k-2l)}{\varphi(2k)}\right)^{2k-2l}\leq 1.
\]

From the Corollary \ref{co:o2chapter-2} it follows that for the space $\mathbf{F}_\psi(\Omega)$, which is generated by the function $\psi(u)$, the inequality \eqref{eq:o1chapter_31} holds, i.e. the statement of the Theorem is fulfilled.

Let us prove the Theorem in the case $u_0>1$. Consider two spaces $\mathbf{F}_\psi(\Omega)$ and $\mathbf{F}_{\hat{\psi}}(\Omega)$, where
\[
\hat{\psi}(u)=\begin{cases}
\psi(u),&\text{if $u\geq u_0$;}\\
\frac{\psi(u_0)}{u_o^{1/2}}u^{1/2},&\text{if $1\leq u < u_0$.}
\end{cases}
\]
Since
\[
\left\|\xi\right\|_\psi=\max \left(\sup\limits_{1\leq u<u_0 }\frac{\left(E\left|\xi\right|^{u}\right)^{1/u}}{\psi(u)},\sup\limits_{u\geq u_0}\frac{\left(E\left|\xi\right|^{u}\right)^{1/u}}{\psi(u)}\right)
\]
and
\[
\sup\limits_{1\leq u<u_0 }\frac{\left(E\left|\xi\right|^{u}\right)^{1/u}}{\psi(u)}=\sup\limits_{1\leq u<u_0 }\frac{\left(E\left|\xi\right|^{u}\right)^{1/u}}{\frac{\psi(u_0)}{u_o^{1/2}}u^{1/2}}\cdot \frac{\psi(u_0)u^{1/2}}{\psi(u)u_o^{1/2}},
\]
then
\[
\left\|\xi\right\|_\psi\leq
S_\psi\max \left(\sup\limits_{1\leq u<u_0 }\frac{\left(E\left|\xi\right|^{u}\right)^{1/u}}{\frac{\psi(u_0)}{u_o^{1/2}}u^{1/2}},\sup\limits_{u\geq u_0}\frac{\left(E\left|\xi\right|^{u}\right)^{1/u}}{\psi(u)}\right),
\]
that is $\left\|\xi\right\|_\psi\leq S_\psi\left\|\xi\right\|_{\hat{\psi}}$.

Now let's consider
\[
\left\|\xi\right\|_{\hat{\psi}}=\max \left(\sup\limits_{1\leq u<u_0 }\frac{\left(E\left|\xi\right|^{u}\right)^{1/u}}{\frac{\psi(u_0)}{u_o^{1/2}}u^{1/2}},\sup\limits_{u\geq u_0}\frac{\left(E\left|\xi\right|^{u}\right)^{1/u}}{\psi(u)}\right).
\]
Since
\[
\sup\limits_{1\leq u<u_0 }\frac{\left(E\left|\xi\right|^{u}\right)^{1/u}}{\frac{\psi(u_0)}{u_o^{1/2}}u^{1/2}}=\sup\limits_{1\leq u<u_0 }\frac{\left(E\left|\xi\right|^{u}\right)^{1/u}}{\psi(u)}\cdot \frac{\psi(u)u_o^{1/2}}{\psi(u_0)u^{1/2}},
\]
then
\[
\left\|\xi\right\|_{\hat{\psi}}\leq
S_{\hat{\psi}}\max \left(\sup\limits_{1\leq u<u_0 }\frac{\left(E\left|\xi\right|^{u}\right)^{1/u}}{\psi(u)},\sup\limits_{u\geq u_0}\frac{\left(E\left|\xi\right|^{u}\right)^{1/u}}{\psi(u)}\right),
\]
that is $\left\|\xi\right\|_{\hat{\psi}}\leq S_{\hat{\psi}}\left\|\xi\right\|_\psi$.

Therefore, the spaces $\mathbf{F}_\psi(\Omega)$ and $\mathbf{F}_{\hat{\psi}}(\Omega)$ contain the same elements and the norms in these spaces are equivalent.

The function $\hat{\psi}(u)$ such that $\frac{\hat{\psi}(u)}{u^{1/2}}$ is monotone, then the first part of the Theorem implies that
\[
\left\|\sum\limits_{i=1}^n\xi_i\right\|_{\hat{\psi}}^2\leq 4\widehat{\psi}_r^2\sum\limits_{i=1}^n\left\|\xi_i\right\|_{\hat{\psi}}^2.
\]
Therefore
\[
\left\|\sum_{i=1}^{n} \xi_i\right\|_\psi^2\leq  S_\psi^2 \left\|\sum_{i=1}^{n} \xi_i\right\|_{\hat{\psi}}^2 \leq
4\widehat{\psi}_r^2 S_\psi^2 \sum_{i=1}^{n}\left\|\xi_i\right\|_{\hat{\psi}}^2\leq 4\widehat{\psi}_r^2 S_\psi^2 S_{\hat{\psi}}^2 \sum_{i=1}^{n}\left\|\xi_i\right\|_\psi^2.
\]
The last inequality implies the statement of the Theorem.
\end{proof}

\begin{remark}
It is obvious that $S_\psi\leq \frac{\psi\left(u_0\right)}{\psi(u)}$ and $S_{\hat{\psi}}\leq \left(u_0\right)^{1/2}$.
\end{remark}

\begin{corollary} \label{co:o1chapter_3}

If for $u\leq u_0$ the function $\varphi(u)$ decreases monotonically, then $S_\psi=1$ and $S_{\hat{\psi}}=\max \left(1,\frac{\psi(1)}{\psi(u_0)}u_0^{1/2}\right)$.
\end{corollary}

The proof of the Corollary is trivial.

\begin{theorem}\label{th:o2chapter-14}
Let $\mathbf{F}_\psi(\Omega)$ be a space such that the function $\psi(u)=e^{au^\beta}$, where $a>0$, $0<\beta<1$. If $\frac{1}{(2a\beta)^{1/\beta}}=1$, then the condition $\mathbf{H}$ with the constant $C_\psi=4e^{2^\beta a}$ is satisfied for the space $\mathbf{F}_\psi(\Omega)$, and if $\frac{1}{(2a\beta)^{1/\beta}}>1$, then the condition $\mathbf{H}$ with the constant $C_\psi=\frac{4e^{a\left(2^\beta+1\right)-\frac{1}{2\beta}}}{(2a\beta)^{1/{2\beta}}}$ is satisfied for $\mathbf{F}_\psi(\Omega)$.
\end{theorem}

\begin{proof}
Statement of the Theorem follows from Theorem \ref{th:o1chapter_13}. Indeed, according to the Corollary \ref{co:o1chapter_3} in our case $S_\psi=1$, and $S_{\hat{\psi}}=\frac{e^{a-\frac{1}{2\beta}}}{(2a\beta)^{1/{2\beta}}}$ $\left(u_0=\frac{1}{(2a\beta)^{1/\beta}}\right)$, and in this case $\frac{e^{a-\frac{1}{2\beta}}}{(2a\beta)^{1/{2\beta}}}>1$.
\end{proof}

\begin{corollary} \label{co:o1chapter_4}
Let us assume in the Theorem \ref{th:o2chapter-14} $a=1$, then for $\beta=\frac{1}{2}$ the condition $\mathbf{H}$ with the constant $C_\psi\approx 16.453$ is satisfied for the space $\mathbf{F}_\psi(\Omega)$, and for $0<\beta<\frac{1}{2}$ for $\mathbf{F}_\psi(\Omega)$ the condition $\mathbf{H}$ with the constant $C_\psi=\frac{4e^{1+2^\beta-\frac{1}{2\beta}}}{(2\beta)^{1/{2\beta}}}$ is satisfied.
\end{corollary}

\begin{proof}
The Corollary follows from Theorem \ref{th:o2chapter-14}.
\end{proof}

\section*{Conclusions to chapter~\ref{ch:o2series}}

In section \ref{ch:o2series}, the Banach space $\mathbf{F}_\psi(\Omega)$ is introduced. Its main properties are studied and examples of random variables from this space are considered. Large deviation inequalities and a majorizing characteristic for random variables from this space are found. Classes of spaces $\mathbf{F}_\psi(\Omega)$ are introduced, where $\psi(u)=u^\alpha$, $\psi(u)=e^{au^\beta}$, $\psi(u)=\left(\ln (u+1)\right)^\lambda$ and majorizing characteristics of the corresponding spaces are found. The spaces $\mathbf{F}_{S_k,\psi,r}(\Omega)$ and the conditions under which these spaces are equivalent to the space $\mathbf{F}_\psi(\Omega)$ are considered. The conditions under which the condition $\mathbf{H}$ is satisfied for the spaces $\mathbf{F}_\psi(\Omega)$ are also studied.

\chapter{Random processes from spaces $\mathbf{F}_\psi^*(\Omega)$}
\label{ch:o3series}

In this chapter, $\mathbf{F}_\psi^*(\Omega)$ is one of the Banach spaces of random variables $\mathbf{F}_\psi(\Omega)$, $\mathbf{\check{F}}_\psi(\Omega)$ or $\mathbf{F}_{S_k,\psi,r}(\Omega)$.
The norm in this space is denoted by $\left\|\cdot\right\|$.
Estimates of the distribution of suprema of random processes on a compact set from this space are found.
The probabilities of large deviations for sums of independent random processes from spaces $\mathbf{F}^*_\psi(\Omega)$ are considered.
Estimates are found for the distribution of suprema of increments from spaces $\mathbf{F}_\psi(\Omega)$.

\section[Estimates for distribution of suprema of random processes from $\mathbf{F}_\psi^*(\Omega)$]
{Estimates for the distribution of suprema of random processes from spaces $\mathbf{F}_\psi^*(\Omega)$}

\begin{definition}
We say that a random process $X = \{X(t), \ t \in T\}$, where $T$ is a parametric set, belongs to the space $\mathbf{F}_\psi^*(\Omega)$ if for any $t \in T$ the random variable $X(t)$ belongs to the space $\mathbf{F}_\psi^*(\Omega)$.
\end{definition}

\begin{example}
Consider a random process $X = \{X(t),\ t \in T\}$ such that
\begin{equation} \label{eq:o2chapter_1-1}
X(t)=\sum \limits_{k=1}^{\infty}\xi_k L_k(t),
\end{equation}
where $\xi_k$ belong to the space $\mathbf{F}_\psi^*(\Omega)$, $L_k(t)$ are some functions.
Since  $$\left\|X(t)\right\|\leq\sum \limits_{k=1}^{\infty}\left\|\xi_k\right\| \left|L_k(t)\right|$$ and if for all $t \in T$
\begin{equation} \label{eq:o2chapter_1-2}
\sum \limits_{k=1}^{\infty}\left\|\xi_k\right\| \left|L_k(t)\right|<\infty,
\end{equation}
then, it is obvious that $X(t)$ belongs to the space $\mathbf{F}_\psi^*(\Omega)$, i.e. $X \in \mathbf{F}_\psi^*(\Omega)$ and the series \eqref{eq:o2chapter_1-1} converges in the norm of the space $\mathbf{F}_\psi^*(\Omega)$ for each $t \in T$.

Note that the following inequality holds true:

\begin{equation} \label{eq:o2chapter_1-3}
\bigl\|X(t) - X(s) \bigr\|= \bigl\|\sum \limits_{k=1}^{\infty}\xi_k\left(L_k(t)-L_k(s)\right)\bigr\|\leq \sum \limits_{k=1}^{\infty} \bigl\|\xi_k\bigr\| \left|L_k(t)-L_k(s)\right|.
\end{equation}

We can write the series \eqref{eq:o2chapter_1-1} in the form
\[
X(t)=\sum \limits_{k=1}^{\infty}\eta_k \lambda_k L_k(t),
\]
where $\left\|\eta_k\right\|=1$, $\lambda_k>0$ are constants. In this case the condition \eqref{eq:o2chapter_1-2} has the form:
\[
\sum \limits_{k=1}^{\infty} \lambda_k \left|L_k(t)\right|<\infty.
\]
\end{example}

\begin{example}
Let $\eta_k$ be independent centered random variables from the space $\mathbf{F}_\psi^*(\Omega)$ for which the condition $\mathbf{H}$ is satisfied, where $\left\|\eta_k\right\|=1$.
Denote $X(t)=\sum \limits_{k=1}^{\infty}\eta_k \lambda_k L_k(t)$, $\lambda_k>0$.
Let for each $t \in T$ the series $\sum \limits_{k=1}^{\infty}\lambda_k^2 L_k^2(t)$ converge. Then $X(t)$ belongs to the space $\mathbf{F}_\psi^*(\Omega)$, $E X(t)=0$ and the inequality holds
\[
\left\|X(t)\right\|^2\leq C_\psi \sum \limits_{k=1}^{\infty}\lambda_k^2 L_k^2(t),
\]
where $C_\psi$ is a constant from the Definition \ref{de:o1chapter_5}. Since
\begin{equation} \label{eq:o2chapter_1-4}
\left\|X(t)\right\|^2\leq C_\psi \sum \limits_{k=1}^{\infty} \left\|\eta_k \lambda_k L_k(t)\right\|^2=C_\psi \sum \limits_{k=1}^{\infty}\lambda_k^2 L_k^2(t),
\end{equation}
then in this case
\[
\bigl\|X(t) - X(s) \bigr\|^2\leq C_\psi \sum \limits_{k=1}^{\infty}\lambda_k^2 \left(L_k(t)-L_k(s)\right)^2.
\]
\end{example}

\begin{example}
Let $\xi_1(t)$, $\xi_2(t)$, $t \in T$ be Gaussian random processes, $E\xi_1(t)=0$, $E\xi_2(t)=0$, $E\xi_1^2(t)=\sigma_1^2(t)$, $E\xi_2^2(t)=\sigma_2^2(t)$. Consider the process
\[
\eta(t)=b_1(t)\exp\left\{\xi_1(t) c_1(t)\right\}+b_2(t)\exp\left\{\xi_2(t) c_2(t)\right\},
\]
where $b_i(t)$, $c_i(t)$, $i=1,2$ are bounded functions. Let us show that $\eta(t)$ belongs to the space $\mathbf{F}_\psi^*(\Omega)$,
where $\psi(u)=e^{u^{\alpha}}$, $\alpha>1$ and estimate the norm of this process. Thus,
\[
\left\|\eta(t)\right\|=\left\|b_1(t)\exp\left\{\xi_1(t) c_1(t)\right\}+b_2(t)\exp\left\{\xi_2(t) c_2(t)\right\}\right\|\leq
\]
\[
\leq \left|b_1(t)\right|\left\|\exp\left\{\xi_1(t) c_1(t)\right\}\right\|+\left|b_2(t)\right|\left\|\exp\left\{\xi_2(t) c_2(t)\right\}\right\|.
\]
Let $\eta_1(t)=\exp\left\{\xi_1(t) c_1(t)\right\}$, and $\eta_2(t)=\exp\left\{\xi_2(t) c_2(t)\right\}$. Let us consider $\eta_1(t)$ (for $\eta_2(t)$ everything is done similarly). From the inequality \eqref{eq:o1chapter_1} we obtain that
\begin{equation} \label{eq:o2chapter_1-5}
\left\|\eta_1(t)\right\|=\sup\limits_{u\geq1}\frac{\left(E\left|\eta_1(t)\right|^{u}\right)^{1/u}}{e^{u^{\alpha}}}.
\end{equation}
We have $\left(E\left(\eta_1(t)\right)^{u}\right)^{1/u}=\exp\left\{\frac{u\sigma_1^2(t)c^2_1(t)}{2}\right\}$. We substitute this value into the equality \eqref{eq:o2chapter_1-5} and get:
\begin{equation} \label{eq:o2chapter_1-6}
\left\|\eta_1(t)\right\|=\sup\limits_{u\geq1}\frac{\exp\left\{\frac{u\sigma_1^2(t)c^2_1(t)}{2}\right\}}{e^{u^{\alpha}}}.
\end{equation}
From the equality
\[
\left(\frac{u\sigma_1^2(t)c^2_1(t)}{2}-u^{\alpha}\right)^{'}=\frac{\sigma_1^2(t)c^2_1(t)}{2}-\alpha u^{\alpha-1}=0
\]
it follows that the supremum is reached at the point $u=\left(\frac{\sigma_1^2(t)c^2_1(t)}{2\alpha}\right)^{\frac{1}{\alpha-1}}$.
Substituting the value of $u$ into the equality \eqref{eq:o2chapter_1-6}, we obtain:
\[
\left\|\eta_1(t)\right\|= \exp\left\{\frac{\alpha-1}{\alpha^{\frac{\alpha}{\alpha-1}}}\left(\frac{\sigma_1^2(t)c^2_1(t)}{2}\right)^{\frac{\alpha}{\alpha-1}}\right\}.
\]
Therefore,
\[
\left\|\eta(t)\right\|\leq \left|b_1(t)\right|\exp\left\{\frac{\alpha-1}{\alpha^{\frac{\alpha}{\alpha-1}}}\left(\frac{\sigma_1^2(t)c^2_1(t)}{2}\right)^{\frac{\alpha}{\alpha-1}}\right\} +
\]
\[
+\left|b_2(t)\right|\exp\left\{\frac{\alpha-1}{\alpha^{\frac{\alpha}{\alpha-1}}}\left(\frac{\sigma_2^2(t)c^2_2(t)}{2}\right)^{\frac{\alpha}{\alpha-1}}\right\}.
\]
\end{example}

\begin{definition}\label{de:o3chapter-1}
The metric massiveness $N(u)$ of a compact metric space $(T,\rho)$ is the smallest number of closed balls of radius at most $u$ that cover the set $T$.
\end{definition}

\begin{theorem} \label{th:o2chapter_1}
Let $T=(T,\rho)$ be a compact metric space, $N(u)$ be the metric massiveness of the space $\left(T,\rho\right)$, $X = \{X(t), t \in T\}$ be a separable random process from the space $\mathbf{F}_\psi^*(\Omega)$, $\varkappa(n)$ be a majorizing characteristic of the space $\mathbf{F}_\psi^*(\Omega)$, and $\varkappa(u)$, $u\geq 1$ be any monotonically increasing function that coincides with $\varkappa(n)$ for integers $n\geq 1$.
Let there be a function
\[
\sigma=\left\{\sigma(h), \ 0\leq h\leq \sup\limits_{t,s \in T}\rho(t, s)\right\},
\]
such that $\sigma(h)$ is continuous, monotonically increasing and $\sigma(0)=0$ and
\[
\sup\limits_{\rho(t, s) \leq h} \bigl\|X(t) - X(s) \bigr\|\leq \sigma(h).
\]

If for any $z>0$ the condition
\begin{equation}\label{equ1}
\int\limits_{0}^{z}\varkappa\left(N\left(\sigma^{(-1)}(u)\right)\right)du<\infty,
\end{equation}
where $\sigma^{(-1)}(u)$ is the inverse function of $\sigma(u)$,
then with probability one the random variable $\sup\limits_{t \in T}\left|X(t)\right|$ belongs to the space $\mathbf{F}_\psi^*(\Omega)$ and
\begin{equation}\label{equ2}
\left\|\sup\limits_{t \in T}\left|X(t)\right|\right\|\leq B(p),
\end{equation}
where $B(p)=\inf \limits_{t \in T}\left\|X(t)\right\|+\frac{1}{p(1-p)}\int\limits_{0}^{\gamma p}\varkappa\left(N\left(\sigma^{(-1)}(u)\right)\right)du$, $\gamma=\sigma\left(\sup\limits_{t,s \in T}\rho(t, s)\right)$, the number $p$ is such that $0<p<1$.
\end{theorem}

\begin{proof}
For all $u>1$ the following inequalities hold:
\[
P\left\{\left|X(t)-X(s)\right|>\varepsilon\right\}\leq \frac{E\left|X(t)-X(s)\right|^u}{\varepsilon^u}\leq
\]
\[
\leq\frac{\bigl\|X(t) - X(s) \bigr\|^u(\psi(u))^u}{\varepsilon^u}\leq \frac{\sigma(\rho(t,s))(\psi(u))^u}{\varepsilon^u}.
\]

Therefore, for any $\varepsilon>0$ $P\left\{\left|X(t)-X(s)\right|>\varepsilon\right\}\rightarrow 0$ for $\rho(t,s)\rightarrow 0$. Thus, $X(t)$ is a probability-continuous random process on the space $(T,\rho)$. Therefore, any countable everywhere dense set in the space $(T,\rho)$ can be the separability set of the process $X(t)$. Let $\varepsilon_k=\sigma^{(-1)}\left(\gamma p^k \right)$, where $k\geq 0$.

Let $V_{\varepsilon_k}$ be the set of centers of closed balls of radius at most $\varepsilon_k$ covering $T$, and the number of these balls is minimal (minimal $\varepsilon_k$ mesh). Let us denote $V=\bigcup\limits_{k=0}^{\infty} V_{\varepsilon_k}$. It is clear that $V$ is a countable everywhere dense set, therefore $V$ is the separability set of the process $X(t)$. Therefore, with probability one
\[
\sup\limits_{t \in T}\left|X(t)\right|=\sup\limits_{t \in V}\left|X(t)\right|.
\]

Let us introduce the mapping $\alpha_k (t)$ onto $V$ in the following way. If $t \in V_{\varepsilon_m}$ ($m$ is some number), then $\alpha_{m-1} (t)$ is a point in $V_{\varepsilon_{m-1}}$ such that $\rho(t,\alpha_{m-1} (t))\leq \varepsilon_{m-1}$. Such a point exists. If there are several such points, then we fix one of them, if $t \in V_{\varepsilon_{m-1}}$, then $\alpha_{m-1} (t)=t$.

Let $t$ be an arbitrary point in $V$, then it is clear that $t$ belongs to some $V_{\varepsilon_m}$. Let us denote $t_m=t$, $t_{m-1}=\alpha_{m-1}(t_m),t_{m-2}=\alpha_{m-2}(t_{m-1}), \ldots, t_1=\alpha_1(t_2), t_0=\alpha_0(t_1)$. Therefore
\begin{multline*}
X(t)=X(t_m)=X(t_m)-X(t_{m-1})+X(t_{m-1})-X(t_{m-2})+\\[1ex]
+X(t_{m-2})-\ldots-X(t_1)-X(t_0)+X(t_0),
\end{multline*}
then
\begin{multline*}
\left|X(t)\right|\leq \left|X(t_m)-X(t_{m-1})\right|+\left|X(t_{m-1})-X(t_{m-2})\right|+\ldots\\[1ex]
\ldots+\left|X(t_1)-X(t_0)\right|+\left|X(t_0)\right|\leq\max \limits_{t \in V_{\varepsilon_m}} \left|X(t)-X(\alpha_{m-1}(t))\right| + \\[1ex]
+\max \limits_{t \in V_{\varepsilon_{m-1}}} \left|X(t)-X(\alpha_{m-2}(t))\right|+\ldots
\ldots +\max \limits_{t \in V_{\varepsilon_1}} \left|X(t)-X(t_0)\right|+\left|X(\alpha_0(t))\right|.
\end{multline*}

Therefore,
\begin{multline*}
\sup\limits_{t \in T}\left|X(t)\right|=\sup\limits_{t \in V}\left|X(t)\right|\leq \left|X(t_0)\right|+ \\[1ex]
+\sum \limits_{l=1}^{m} \max \limits_{t \in V_{\varepsilon_l}} \left|X(t)-X(\alpha_{l-1}(t))\right|\leq\left|X(t_0)\right|+ \sum \limits_{l=1}^{\infty} \max \limits_{t \in V_{\varepsilon_l}} \left|X(t)-X(\alpha_{l-1}(t))\right|,
\end{multline*}
that is why
\[
\left\|\sup\limits_{t \in T}\left|X(t)\right|\right\|\leq \left\|X(t_0)\right\|+\sum \limits_{l=1}^{\infty} \left\|\max \limits_{t \in V_{\varepsilon_l}} \left|X(t)-X(\alpha_{l-1}(t))\right|\right\|.
\]
From the definition of \ref{eq:o2chapter-3} it follows that
\begin{multline} \label{eq:o2chapter-1}
\left\|\sup\limits_{t \in T}\left|X(t)\right|\right\|\leq \inf \limits_{t \in T}\left\|X(t)\right\|+\sum \limits_{l=1}^{\infty} \varkappa\left(N(\varepsilon_l)\right)\max \limits_{t \in V_{\varepsilon_l}}\left\| X(t)-X(\alpha_{l-1}(t))\right\|\leq \\[1ex]
\leq\inf \limits_{t \in T}\left\|X(t)\right\|+\sum \limits_{l=1}^{\infty} \varkappa\left(N(\varepsilon_l)\right)\sigma(\varepsilon_{l-1})\leq \\[1ex]
\leq\inf \limits_{t \in T}\left\|X(t)\right\|+\sum \limits_{l=1}^{\infty} \varkappa\left(N\left(\sigma^{(-1)}(\gamma p^l)\right)\right)\gamma p^{l-1}.
\end{multline}
Since
\begin{multline*}
\int\limits_{\gamma p^{l+1}}^{\gamma p^l}\varkappa\left(N\left(\sigma^{(-1)}(u)\right)\right)du \geq \varkappa\left(N\left(\sigma^{(-1)}(\gamma p^l)\right)\right)\left(\gamma p^l-\gamma p^{l+1}\right)=\\[1ex]
=\varkappa\left(N\left(\sigma^{(-1)}(\gamma p^l)\right)\right)\gamma p^{l-1}(p-p^2),
\end{multline*}
then we get that $$\gamma p^{l-1}\varkappa\left(N\left(\sigma^{(-1)}(\gamma p^l)\right)\right)\leq\frac{\int\limits_{\gamma p^{l+1}}^{\gamma p^l}\varkappa\left(N\left(\sigma^{(-1)}(u)\right)\right)du}{p(1-p)}.$$ Thus, from the inequality \eqref{eq:o2chapter-1} we have
\begin{multline*}
\left\|\sup\limits_{t \in T}\left|X(t)\right|\right\|\leq \inf \limits_{t \in T}\left\|X(t)\right\|+\sum \limits_{l=1}^{\infty}\frac{1}{p(1-p)}\int\limits_{\gamma p^{l+1}}^{\gamma p^l}\varkappa\left(N\left(\sigma^{(-1)}(u)\right)\right)du=\\[1ex]
=\inf \limits_{t \in T}\left\|X(t)\right\|+\frac{1}{p(1-p)}\int\limits_{0}^{\gamma p}\varkappa\left(N\left(\sigma^{(-1)}(u)\right)\right)du,
\end{multline*}
which was what had to be proved.
\end{proof}

\begin{corollary}\label{co:o2chapter-1}
Let the process $X = \{X(t), \ t \in T\}$, which belongs to the space $\mathbf{F}_\psi(\Omega)$ satisfy the conditions of Theorem \ref{th:o2chapter_1}.
Then for any $\varepsilon>0$ the inequality holds:
\[
P\left\{\sup\limits_{t \in T}\left|X(t)\right|>\varepsilon\right\}\leq \inf\limits_{u\geq1}\frac{B^u(p)(\psi(u))^u}{\varepsilon^u}.
\]
\end{corollary}

\begin{proof}
The Corollary follows from Theorem \ref{th:o1chapter_1}.
\end{proof}

\begin{remark}
The Corollary \ref{co:o2chapter-1} for other spaces $\mathbf{F}_\psi^*(\Omega)$ is obvious. For example, for spaces $\mathbf{F}_{S_k,\psi,r}(\Omega)$ the inequality holds:
\[
P\left\{\sup\limits_{t \in T}\left|X(t)\right|>\varepsilon\right\}\leq \inf\limits_{k\geq r}\frac{B^{S_k}(p)(\psi(S_k))^{S_k}}{\varepsilon^{S_k}}.
\]
\end{remark}

\begin{example}\label{ex:o3chapter-4}
Consider the space $\mathbf{F}_\psi(\Omega)$, where $\psi(u)=u^\alpha$, $\alpha>0$, then from  Corollary \ref{co:o2chapter-1} and Theorem \ref{th:o1chapter_2} for $\varepsilon\geq e^\alpha B(p)$ we obtain:
\[
P\left\{\sup\limits_{t \in T}\left|X(t)\right|>\varepsilon\right\}\leq \exp \left\{-\frac{\alpha}{e}\left(\frac{\varepsilon}{B(p)}\right)^{1/\alpha}\right\}.
\]
\end{example}

\begin{example}\label{ex:o3chapter-5}
Consider the space $\mathbf{F}_\psi(\Omega)$, where $\psi(u)=e^{au^\beta}$, $a>0$, $\beta>0$, then from Corollary \ref{co:o2chapter-1} and Theorem \ref{th:o1chapter_3} for $\varepsilon\geq e^{a(\beta+1)} B(p)$ we obtain:
\[
P\left\{\sup\limits_{t \in T}\left|X(t)\right|>\varepsilon\right\}\leq \exp \left\{-\frac{\beta}{a^{1/\beta}}\left(\frac{\ln \frac{\varepsilon}{B(p)}}{\beta+1}\right)^{\frac{\beta+1}{\beta}}\right\}.
\]
\end{example}

\begin{example}
Consider the space $\mathbf{F}_\psi(\Omega)$, where $\psi(u)=\left(\ln (u+1)\right)^\lambda$, $\lambda>0$, then from  Corollary \ref{co:o2chapter-1} and Theorem \ref{th:o1chapter_4} for $\varepsilon\geq\left(e \ln 2\right)^\lambda B(p)$ we have:
\[
P\left\{\sup\limits_{t \in T}\left|X(t)\right|>\varepsilon\right\}\leq e^\lambda \exp \left\{-\lambda  \exp\left\{\left(\frac{\varepsilon}{B(p)}\right)^{1/\lambda}\frac{1}{e}\right\}\right\}.
\]
\end{example}

\begin{corollary}
\label{co:o2chapter_2}
Let $X = \{X(t), \ t \in \left[c,d\right]\}$, $-\infty<c<d<+\infty$ be a separable random process in the space $\mathbf{F}_\psi(\Omega)$. Let the condition hold
\begin{equation} \label{eq:o2chapter_1}
\sup\limits_{\substack{\left|t-s\right| \leq h\\ t,s \in \left[c,d\right]}}\bigl\|X(t) - X(s) \bigr\|_\psi\leq \sigma(h),
\end{equation}
where $\sigma=\left\{\sigma(h),\ 0\leq h\leq d-c\right\}$ is a continuous, monotonically increasing function and $\sigma(0)=0$. If for any $z>0$ the condition
\[
\int\limits_{0}^{z}\varkappa\left(\frac{d-c}{2\sigma^{(-1)}(u)}+1\right)du<\infty,
\]
then with probability one $\sup\limits_{t \in \left[c,d\right]}\left|X(t)\right| \in \mathbf{F}_\psi(\Omega)$ and for any $0<p<1$ the inequality holds true
\begin{equation} \label{eq:o2chapter_2}
\left\|\sup\limits_{t \in \left[c,d\right]}\left|X(t)\right|\right\|_\psi\leq \widetilde{B}(p),
\end{equation}
where $$\widetilde{B}(p)=\inf \limits_{t \in \left[c,d\right]}\left\|X(t)\right\|_\psi+\frac{1}{p(1-p)}\int\limits_{0}^{\gamma p}\varkappa\left(\frac{d-c}{2\sigma^{(-1)}(u)}+1\right)du,$$ $\gamma=\sigma(d-c)$, $\varkappa(u)$ is the majorizing characteristic of the space $\mathbf{F}_\psi(\Omega)$, $\sigma^{(-1)}(u)$ is the inverse function of $\sigma(u)$.

In addition, for any $\varepsilon>0$ the following inequality holds:
\begin{equation} \label{eq:o2chapter_3}
P\left\{\sup\limits_{t \in \left[c,d\right]}\left|X(t)\right|>\varepsilon\right\}\leq \inf\limits_{u\geq1}\frac{\widetilde{B}^u(p)(\psi(u))^u}{\varepsilon^u}.
\end{equation}
\end{corollary}

\begin{proof}
The Corollary follows from Theorem \ref{th:o2chapter_1}, since the metric massiveness of the interval $\left[c,d\right]$ is estimated as follows:
\[
N(u)\leq \frac{d-c}{2u}+1.
\]
The inequality \eqref{eq:o2chapter_3} follows from Corollary \ref{co:o2chapter-1}.
\end{proof}

\begin{corollary}
\label{co:o2chapter_3}
Let $X = \{X(t), \ t \in \left[c,d\right]\}$, $-\infty<c<d<+\infty$ be a separable random process from the space $\mathbf{F}_\psi(\Omega)$ and for some $0<\mu<1$ the condition
\begin{equation} \label{eq:o2chapter_4}
\sup\limits_{\substack{\left|t-s\right| \leq h\\ t,s \in \left[c,d\right]}} \bigl\|X(t) - X(s) \bigr\|_\psi\leq \frac{C}{\left(\varkappa\left(\frac{d-c}{2h}+1\right)\right)^{1/\mu}},
\end{equation}
where $C>0$ is some constant, $h<d-c$. Then with probability one \linebreak $\sup\limits_{t \in \left[c,d\right]}\left|X(t)\right|$ belongs to the space $\mathbf{F}_\psi(\Omega)$ and the inequality holds
\[
\left\|\sup\limits_{t \in \left[c,d\right]}\left|X(t)\right|\right\|_\psi\leq \inf \limits_{t \in \left[c,d\right]}\left\|X(t)\right\|_\psi+C\left(\varkappa\left(\frac{3}{2}\right)\right)^{(\mu-1)/\mu}\frac{(1+\mu)^{\mu+1}}{\mu^\mu(1-\mu)}=\widetilde{B}.
\]

Furthermore, for any $\varepsilon>0$, the inequality is true
\begin{equation} \label{eq:o2chapter_5}
P\left\{\sup\limits_{t \in \left[c,d\right]}\left|X(t)\right|>\varepsilon\right\}\leq \inf\limits_{u\geq1}\frac{\widetilde{B}^u(\psi(u))^u}{\varepsilon^u}.
\end{equation}
\end{corollary}

\begin{proof}
Our statement follows from the Corollary \ref{co:o2chapter_2}. Indeed,
\[
\sigma(h)=\frac{C}{\left(\varkappa\left(\frac{d-c}{2h}+1\right)\right)^{1/\mu}},
\]
then it is established that \\
\[
\frac{1}{p(1-p)}\int\limits_{0}^{\gamma p}\varkappa\left(\frac{d-c}{2\sigma^{(-1)}(u)}+1\right)du=
\]
\[
=\frac{1}{p(1-p)}\int\limits_{0}^{\gamma p}\frac{C^\mu}{u^\mu}du=\frac{C^\mu}{1-\mu}\gamma^{1-\mu}\frac{1}{(1-p)p^\mu}.
\]
If the expression $\frac{C^\mu}{1-\mu}\gamma^{1-\mu}\frac{1}{(1-p)p^\mu}$ is minimized with respect to $p$, then from the inequality \eqref{eq:o2chapter_2} we obtain our statement, and the inequality \eqref{eq:o2chapter_5} follows from the Corollary \ref{co:o2chapter-1}.
\end{proof}

\begin{example} \label{ex:o2chapter-7}
Consider the space $\mathbf{F}_\psi(\Omega)$, where $\psi(u)=u^\alpha$, $\alpha>0$. According to Theorem \ref{th:o1chapter_7} $\varkappa(n)=\left(\ln n\right)^\alpha\left(\frac{e}{\alpha}\right)^\alpha$, therefore, substituting the value of $\varkappa(n)$ into the inequality \eqref{eq:o2chapter_4}, we obtain:
\[
\sigma(h)=\frac{C}{\left(\frac{e}{\alpha}\ln\left(\frac{d-c}{2h}+1\right)\right)^{\alpha/\mu}}.
\]
Then
\begin{multline*}
\left\|\sup\limits_{t \in \left[c,d\right]}\left|X(t)\right|\right\|_\psi\leq \inf \limits_{t \in \left[c,d\right]}\sup\limits_{u\geq1}\frac{\left(E\left|X(t)\right|^{u}\right)^{1/u}}{u^\alpha}+\\[1ex]
+C\left(\frac{e}{\alpha}\ln\frac{3}{2}\right)^{\alpha(\mu-1)/\mu}\frac{(1+\mu)^{\mu+1}}{\mu^\mu(1-\mu)}=B_1.
\end{multline*}

From the Corollary \ref{co:o2chapter_3} and Theorem \ref{th:o1chapter_2} for $\varepsilon\geq e^\alpha B_1$ we obtain:
\begin{equation} \label{eq:o2chapter-13}
P\left\{\sup\limits_{t \in \left[c,d\right]}\left|X(t)\right|>\varepsilon\right\}\leq \exp \left\{-\frac{\alpha}{e}\left(\frac{\varepsilon}{B_1 }\right)^{1/\alpha}\right\}.
\end{equation}
\end{example}

\begin{example} \label{ex:o2chapter-8}
Consider the space $\mathbf{F}_\psi(\Omega)$, where $\psi(u)=e^{au^\beta}$, $a>0$, $\beta>0$. According to Theorem \ref{th:o1chapter_6}
$$\varkappa(n)=\frac{1}{e^a}\exp \left\{S(a,\beta)(\ln n)^{\frac{\beta}{\beta+1}}\right\},$$ therefore, substituting the value of $\varkappa(n)$ into inequality \eqref{eq:o2chapter_4}, we obtain:

\[
\sigma(h)=\frac{C}{\left(\frac{1}{e^a}\exp \left\{S(a,\beta)(\ln \left(\frac{d-c}{2h}+1\right))^{\frac{\beta}{\beta+1}}\right\}\right)^\frac{1}{\mu}},
\]
where $S(a,\beta)=(\beta a)^{\frac{1}{\beta+1}}(\beta^{-1}+1)$. Then
\begin{multline*}
\left\|\sup\limits_{t \in \left[c,d\right]}\left|X(t)\right|\right\|_\psi\leq \inf \limits_{t \in \left[c,d\right]}\sup\limits_{u\geq1}\frac{\left(E\left|X(t)\right|^{u}\right)^{1/u}}{e^{au^\beta}}+\\[1ex]
+C\left(\frac{1}{e^a}\exp \left\{S(a,\beta) \left(\ln\frac{3}{2}\right)^{\frac{\beta}{\beta+1}}\right\}\right)^{(\mu-1)/\mu}\frac{(1+\mu)^{\mu+1}}{\mu^\mu(1-\mu)}=B_2.
\end{multline*}
From the Corollary \ref{co:o2chapter_3} and  Theorem \ref{th:o1chapter_3} for any $\varepsilon\geq e^{a(\beta+1)}B_2$ we obtain:
\begin{equation} \label{eq:o2chapter-14}
P\left\{\sup\limits_{t \in \left[c,d\right]}\left|X(t)\right|>\varepsilon\right\}\leq \exp \left\{-\frac{\beta}{a^{1/\beta}}\left(\frac{\ln \frac{x}{B_2}}{\beta+1}\right)^{\frac{\beta+1}{\beta}}\right\}.
\end{equation}
\end{example}

\begin{example} \label{ex:o2chapter-9}
Consider the space $\mathbf{F}_\psi(\Omega)$, where $\psi(u)=\left(\ln (u+1)\right)^\lambda$, $\lambda>0$. According to the Theorem \ref{th:o1chapter-8} $\varkappa(n)=e \left(\frac{\ln(\ln n+2)}{\ln 2}\right)^\lambda$, therefore, substituting the value of $\varkappa(n)$ into the inequality \eqref{eq:o2chapter_4}, we obtain:
\[
\sigma(h)=\frac{C}{\left(e\left(\frac{\ln(\ln \left(\frac{d-c}{2h}+1\right)+2)}{\ln 2}\right)^\lambda \right)^\frac{1}{\mu}}.
\]
Then
\begin{multline*}
\left\|\sup\limits_{t \in \left[c,d\right]}\left|X(t)\right|\right\|_\psi\leq \inf \limits_{t \in \left[c,d\right]}\sup\limits_{u\geq1}\frac{\left(E\left|X(t)\right|^{u}\right)^{1/u}}{\left(\ln (u+1)\right)^\lambda}+\\[1ex]
+C\left(e \left(\frac{\ln(\ln \frac{3}{2}+2)}{\ln 2}\right)^\lambda\right)^{(\mu-1)/\mu}\frac{(1+\mu)^{\mu+1}}{\mu^\mu(1-\mu)}=B_3.
\end{multline*}
From  Corollary \ref{co:o2chapter_3} and Theorem \ref{th:o1chapter_4} for $\varepsilon>0$ we obtain:

\begin{equation} \label{eq:o2chapter-15}
P\left\{\sup\limits_{t \in \left[c,d\right]}\left|X(t)\right|>\varepsilon\right\}\leq
e^\lambda \exp \left\{-\lambda  \exp\left\{\left(\frac{x}{B_3}\right)^{1/\lambda}\frac{1}{e}\right\}\right\}.
\end{equation}
\end{example}

In many cases, for the process $X(t)$, the estimates \eqref{eq:o2chapter-13} and \eqref{eq:o2chapter-14} are difficult to find. It is much easier to obtain an estimate of the form $\sigma(h)=\overline{C}\left|h\right|^\delta$, where $\overline{C}$ is some constant, $0<\delta\leq1$. We will show how to obtain the estimates from  Example \ref{ex:o2chapter-7} and Example \ref{ex:o2chapter-8}.

\begin{example}\label{ex:o2chapter-10}
Consider the space $\mathbf{F}_\psi(\Omega)$, where $\psi(u)=u^\alpha$, $\alpha>1$. In this case $\varkappa(n)=\left(\frac{e}{\alpha}\right)^\alpha\left(\ln n\right)^\alpha$. For $x>0$ and $0\leq\tau\leq1$ the inequality holds:
\[
\ln(1+x)=\frac{\tau}{\tau}\ln(1+x)= \frac{1}{\tau}\ln(1+x)^\tau\leq \frac{1}{\tau}\ln(1+x^\tau)\leq \frac{x^\tau}{\tau},
\]
which will be used further. For $\tau<\frac{\delta}{\alpha}$ it follows from Corollary \ref{co:o2chapter_2} that
\begin{multline*}
\varkappa\left(\frac{d-c}{2\sigma^{(-1)}(u)}+1\right)=\left(\frac{e}{\alpha}\right)^\alpha\left(\ln\left(\frac{d-c}{2\sigma^{(-1)}(u)}+1\right)  \right)^\alpha\leq \\
\leq\left(\frac{e}{\alpha}\right)^\alpha\left(\left(\frac{d-c}{2\sigma^{(-1)}(u)}\right)^\tau \frac{1}{\tau} \right)^\alpha = \left(\frac{e}{\alpha\tau}\right)^\alpha\left(\frac{d-c}{2}\right)^{\alpha\tau}\left(\left(\frac{\overline{C}}{u}\right)^{\frac{1}{\delta}}\right)^{\alpha\tau}=\\
=\left(\frac{e}{\alpha\tau}\right)^\alpha\left(\frac{d-c}{2}\right)^{\alpha\tau} \overline{C}^{\frac{\alpha\tau}{\delta}} u^{-\frac{\alpha\tau}{\delta}}=A(\alpha,\tau,\delta)u^{-\frac{\alpha\tau}{\delta}},
\end{multline*}
where $A(\alpha,\tau,\delta)=\left(\frac{e}{\alpha\tau}\right)^\alpha\left(\frac{d-c}{2}\right)^{\alpha\tau} \overline{C}^{\frac{\alpha\tau}{\delta}}$.
Then
\[
\int\limits_{0}^{\gamma p}\varkappa\left(\frac{d-c}{2\sigma^{(-1)}(u)}+1\right)du\leq A(\alpha,\tau,\delta)\int\limits_{0}^{\gamma p}u^{-\frac{\alpha\tau}{\delta}}du=A(\alpha,\tau,\delta)\frac{\left(\gamma p\right)^{1-\frac{\alpha\tau}{\delta}}}{1-\frac{\alpha\tau}{\delta}}.
\]
Therefore, the inequality \eqref{eq:o2chapter_2} has the form
\[
\left\|\sup\limits_{t \in \left[c,d\right]}\left|X(t)\right|\right\|_\psi\leq \widetilde{B}(p,\alpha,\tau,\delta),
\]
where $$\widetilde{B}(p,\alpha,\tau,\delta)=\inf \limits_{t \in \left[c,d\right]}\sup\limits_{u\geq1}\frac{\left(E\left|X(t)\right|^{u}\right)^{1/u}}{u^\alpha}+\frac{A(\alpha,\tau,\delta)\left(\gamma p\right)^{1-\frac{\alpha\tau}{\delta}}}{p(1-p)\left(1-\frac{\alpha\tau}{\delta}\right)},$$
$\gamma=\overline{C}\left|d-c\right|^\delta$.

In addition, from Theorem \ref{th:o1chapter_2} for $\varepsilon\geq e^\alpha \widetilde{B}(p,\alpha,\tau,\delta)$ we obtain:
\[
P\left\{\sup\limits_{t \in \left[c,d\right]}\left|X(t)\right|>\varepsilon\right\}\leq \exp \left\{-\frac{\alpha}{e}\left(\frac{\varepsilon}{\widetilde{B}(p,\alpha,\tau,\delta)}\right)^{1/\alpha}\right\}.
\]
\end{example}

\begin{example}\label{ex:o2chapter-11}
Consider the space $\mathbf{F}_\psi(\Omega)$, where $\psi(u)=e^{au^\beta}$, $a>0$, $\beta>0$. In this case, $$\varkappa(n)=\frac{1}{e^a}\exp \left\{S(a,\beta)(\ln n)^{\frac{\beta}{\beta+1}}\right\},$$ where $S(a,\beta)=(\beta a)^{\frac{1}{\beta+1}}(\beta^{-1}+1)$. From Corollary \ref{co:o2chapter_2} it follows that
\[
\varkappa\left(\frac{d-c}{2\sigma^{(-1)}(u)}+1\right)=\frac{1}{e^a}\exp \left\{S(a,\beta)\left(\ln \left(\frac{D}{u^{\frac{1}{\delta}}}+1\right)\right)^{\frac{\beta}{\beta+1}}\right\},
\]
where $D=\frac{\overline{C}^{\frac{1}{\delta}}(d-c)}{2}$. Then, substituting this value into inequality \eqref{eq:o2chapter_2}, we have:
\begin{multline*}
\int\limits_{0}^{\gamma p}\varkappa\left(\frac{d-c}{2\sigma^{(-1)}(u)}+1\right)du=\int\limits_{0}^{\gamma p} \frac{1}{e^a}\exp \left\{S(a,\beta)\left(\ln \left(\frac{D}{u^{\frac{1}{\delta}}}+1\right)\right)^{\frac{\beta}{\beta+1}}\right\}du\leq \\
\leq \int\limits_{0}^{\gamma p} \frac{1}{e^a}\exp \left\{S(a,\beta)\ln \left(\frac{D}{u^{\frac{1}{\delta}}}+1\right)\left(\ln \left(\frac{D}{u^{\frac{1}{\delta}}}+1\right)\right)^{-\frac{1}{\beta+1}}\right\}du=\\
=\frac{1}{e^a}\int\limits_{0}^{\gamma p}\left(\frac{D}{u^{\frac{1}{\delta}}}+1\right)^{S(a,\beta)\left(\ln \left(\frac{D}{u^{\frac{1}{\delta}}}+1\right)\right)^{-\frac{1}{\beta+1}}}du=I(a,\beta,\delta).
\end{multline*}
The integral $I(a,\beta,\delta)$ converges for any $a>0$, $\beta>0$, and $\delta>0$, since $\left(\ln \left(\frac{D}{u^{\frac{1}{\delta}}}+1\right)\right)^{\frac{1}{\beta+1}}\rightarrow 0$ for $u\rightarrow 0$. Therefore, the inequality \eqref{eq:o2chapter_2} can be rewritten as:
\[
\left\|\sup\limits_{t \in \left[c,d\right]}\left|X(t)\right|\right\|_\psi\leq \widetilde{B}(p,a,\beta,\delta),
\]
where $$\widetilde{B}(p,a,\beta,\delta)=\inf \limits_{t \in \left[c,d\right]}\sup\limits_{u\geq1}\frac{\left(E\left|X(t)\right|^{u}\right)^{1/u}}{e^{au^\beta}}+\frac{I(a,\beta,\delta)}{p(1-p)},$$ $\gamma=\overline{C}\left|d-c\right|^\delta$.

In addition, from Theorem \ref{th:o1chapter_3} for $\varepsilon\geq e^{a(\beta+1)}$ we obtain:
\[
P\left\{\sup\limits_{t \in \left[c,d\right]}\left|X(t)\right|>\varepsilon\right\}\leq \exp \left\{-\frac{\beta}{a^{1/\beta}}\left(\frac{\ln \frac{x}{\widetilde{B}(p,a,\beta,\delta)}}{\beta+1}\right)^{\frac{\beta+1}{\beta}}\right\}.
\]
\end{example}

\begin{example} \label{le:o2chapter_10}
Consider a random process $X = \{X(t), t \in T\}$, where $X(t)=\sum\limits_{k=1}^\infty \xi_k L_k(t)$, $\xi_k \in \mathbf{F}_\psi(\Omega)$, $t \in \left[c,d\right]$. The functions $L_k(t)$ satisfy the Lipschitz condition:
\[
\left|L_k(t)-L_k(s)\right|\leq C_k \left|t-s\right|^\gamma,
\]
where $\gamma$ is a constant, $0<\gamma\leq 1$, and $C_k>0$,  then
\[
\bigl\|X(t) - X(s) \bigr\|_\psi\leq \sum\limits_{k=1}^\infty \left\|\xi_k\right\|_\psi\left|L_k(t)-L_k(s)\right|\leq \left(\sum\limits_{k=1}^\infty \left\|\xi_k\right\|_\psi C_k\right)\left|t-s\right|^\gamma.
\]
Therefore, if the series $\sum\limits_{k=1}^\infty \left\|\xi_k\right\|C_k$ is convergent, then
\[
\sup\limits_{\substack{\left|t-s\right| \leq h\\ t,s \in \left[c,d\right]}} \bigl\|X(t) - X(s) \bigr\|_\psi\leq \widehat{C} h^\gamma,
\]
where $\widehat{C}=\sum\limits_{k=1}^\infty \left\|\xi_k\right\|_\psi C_k$.
\end{example}

We found that estimates from  Example \ref{ex:o2chapter-10} and Example \ref{ex:o2chapter-11} can be applied to the process $X(t)$.
If the series $\sum\limits_{k=1}^\infty \left\|\xi_k\right\|C_k$ is divergent, then we illustrate how such estimates can be obtained.
To this end, consider the following Lemma. Let $\gamma=1$ for simplicity.

\begin{lemma}\label{le:o2chapter-1}
Let the functions $Y_k(t)$, $k=1,2,\ldots$, $t \in \left[c,d\right]$ satisfy the following conditions:
\begin{itemize}
\item [1.] $\sup\limits_{t \in \left[c,d\right]}\left|Y_k(t)\right|\leq B_k$;
\item [2.] $\sup\limits_{\substack{\left|t-s\right| \leq h\\ t,s \in \left[c,d\right]}}\left|Y_k(t) - Y_k(s) \right|\leq \check{C_k} h$.
\end{itemize}
Let, in addition, the function $\varphi(\lambda)$, $\lambda>0$ is continuous and increasing, and the function $\frac{\lambda}{\varphi(\lambda)}$ increases for $\lambda>0$. Then the inequality holds:
\[
\sup\limits_{\substack{\left|t-s\right| \leq h\\ t,s \in \left[c,d\right]}}\left|Y_k(t) - Y_k(s) \right|\leq \max \left(1,2B_k\right)\frac{\varphi(\check{C_k})}{\varphi\left(\frac{1}{\left|h\right|}\right)}.
\]
\end{lemma}

\begin{proof}
From the conditions of the Lemma we have that $\left|Y_k(t)-Y_k(s)\right|\leq \check{C_k} \left|t-s\right|$ and $\left|Y_k(t)\right|< B_k$. Let us consider two cases.
In the first case, where $\check{C_k}>\frac{1}{\left|t-s\right|}$, we have
\[
\frac{\left|Y_k(t)-Y_k(s)\right|}{2B_k}\leq \frac{\left|Y_k(t)\right|+\left|Y_k(s)\right|}{2B_k}\leq 1\leq \frac{\varphi(\check{C_k})}{\varphi\left(\frac{1}{\left|t-s\right|}\right)}.
\]
In the second, where $\check{C_k}\leq\frac{1}{\left|t-s\right|}$, we have
\[
\left|Y_k(t)-Y_k(s)\right|\leq \check{C_k} \left|t-s\right|=\frac{\check{C_k}}{\left|t-s\right|^{-1}}\leq \frac{\varphi(\check{C_k})}{\varphi\left(\frac{1}{\left|t-s\right|}\right)}.
\]
From the inequalities obtained in both cases, the validity of this Lemma follows.
\end{proof}

\begin{example}
If in  Lemma \ref{le:o2chapter-1} we put $Y_k(t)=L_k(t)$ and $\varphi(\lambda)=\lambda^\alpha$, $0<\alpha<1$, then (since in this case $\frac{\lambda}{\varphi(\lambda)}$ increases as $\lambda>0$) the inequality holds:
\[
\sup\limits_{\left|t-s\right| \leq h}\left|L_k(t) - L_k(s) \right|\leq \max \left(1,2B_k\right) \check{C_k}^\alpha h^\alpha,
\]
where $B_k=\sup\limits_{t \in \left[c,d\right]}\left|L_k(t)\right|$.
Therefore,
\[
\bigl\|X(t) - X(s) \bigr\|_\psi\leq \left(\sum\limits_{k=1}^\infty  \left\|\xi_k\right\|_\psi\check{C_k}^\alpha\max \left(1,2B_k\right)\right) h^\alpha.
\]
\end{example}
Denote $\sum\limits_{k=1}^\infty \left\|\xi_k\right\|\check{C_k}^\alpha\max \left(1,2B_k\right)= \widetilde{C}$. In this case we have $\sigma(h)= \widetilde{C} \left|h\right|^\alpha$.
Therefore, to estimate the suprema of $X(t)$, we can use formulas from Example \ref{ex:o2chapter-10}, Example \ref{ex:o2chapter-11} and Example \ref{le:o2chapter_10}.

\begin{example}\label{ex:o2chapter-12}
Let in Example \ref{le:o2chapter_10} the random variables be independent, and for the space $\mathbf{F}_\psi(\Omega)$ the condition $\mathbf{H}$ is satisfied.
Then
\[
\bigl\|X(t) - X(s) \bigr\|_\psi^2\leq C_\psi\sum\limits_{k=1}^\infty \left\|\xi_k\right\|_\psi^2\left(L_k(t)-L_k(s)\right)^2\leq
\]
\[
\leq C_\psi\sum\limits_{k=1}^\infty \left\|\xi_k\right\|_\psi^2C_k^2\left|t-s\right|^{2\gamma}= C_\psi \left|t-s\right|^{2\gamma}\sum\limits_{k=1}^\infty \left\|\xi_k\right\|_\psi^2C_k^2.
\]
If the series $\sum\limits_{k=1}^\infty \left\|\xi_k\right\|_\psi^2C_k^2$ converges, then

$
\sup\limits_{\left|t-s\right| \leq h} \bigl\|X(t) - X(s) \bigr\|_\psi\leq \breve{C} h^\gamma,
$
where $\breve{C}=\sqrt{C_\psi\sum\limits_{k=1}^\infty \left\|\xi_k\right\|_\psi^2C_k^2}$.
\end{example}
Therefore, to estimate the distributions of suprema of the process $X(t)$, we can apply formulas from the Section \ref{ch:o2series} for spaces where the condition $\mathbf{H}$ is satisfied.

\begin{example}\label{ex:o2chapter-15-1}
Let $X(t)$ be a Gaussian centered process, $t \in \left[c,d\right]$,  $EX(t)X(s)=B(t,s)$. Let $\sigma^2(t)=B(t,t)$. Consider $e^{X(t)}-e^{X(s)}$. Let $X(t)>X(s)$. Then
\[
e^{X(t)}-e^{X(s)}=e^{X(t)}\left(1-e^{-(X(t)-X(s))}\right)\leq e^{X(t)}\left(X(t)-X(s)\right).
\]
So, always
\begin{multline*}
\left|e^{X(t)}-e^{X(s)}\right|\leq \\[1ex]
\leq\max \left(e^{X(t)},e^{X(s)}\right)\left|X(t)-X(s)\right|\leq  \left(e^{X(t)}+e^{X(s)}\right)\left|X(t)-X(s)\right|.
\end{multline*}
Consider the expression
\[
I(u)=\frac{\left(E\left|e^{X(t)}-e^{X(s)}\right|^u\right)^{\frac{1}{u}}}{\psi(u)}\leq \frac{\left(E\left(e^{X(t)}+e^{X(s)}\right)^u\left|X(t)-X(s)\right|^u\right)^{\frac{1}{u}}}{\psi(u)}.
\]
Using H\"older's inequality, where $\frac{1}{q}+\frac{1}{p}=1$, for $q>1$, we have
\[
I(u)\leq \frac{1}{\psi(u)}\left(E\left(e^{X(t)}+e^{X(s)}\right)^{pu}\right)^{\frac{1}{pu}} \left(E\left|X(t)-X(s)\right|^{qu}\right)^{\frac{1}{qu}}.
\]
Let $\psi(u)=\exp\left\{u^\beta\right\}$ and let $\omega>\frac{1}{\beta}$, where $\beta>0$. Then
\[
I(u)\leq \frac{1}{(\psi(u))^{1-\omega}}\times
\]
\[ \times\left[\frac{\left(E\exp\left\{puX(t)\right\}\right)^{\frac{1}{pu}}}{(\psi(u))^\omega}+\frac{\left(E\exp\left\{puX(s)\right\}\right)^{\frac{1}{pu}}}{(\psi(u))^\omega}\right]\left(E\left|X(t)-X(s)\right|^{qu}\right)^{\frac{1}{qu}}.
\]
From inequality
$$\left(E\exp\left\{puX(t)\right\}\right)^{\frac{1}{pu}}=\left(\exp\left\{\sigma^2(t)\frac{p^2 u^2}{2}\right\}\right)^{\frac{1}{pu}}\leq \exp\left\{\frac{\sigma^2 pu}{2}\right\},$$ where $\sigma^2=\sup\limits_{t \in \left[c,d\right]}\sigma^2(t)$, we get

\[
\frac{\left(E\exp\left\{puX(t)\right\}\right)^{\frac{1}{pu}}}{\exp\left\{\omega u^\beta\right\}}\leq \exp\left\{\frac{\sigma^2 pu}{2}-\omega u^\beta\right\}.
\]
Let us find the supremum of the right-hand side of the last inequality. The supremum is reached at the point $u=\left(\frac{\sigma^2 p}{2\omega \beta}\right)^{\frac{1}{\beta-1}}$. So
\[
\frac{\left(E\exp\left\{puX(t)\right\}\right)^{\frac{1}{pu}}}{\exp\left\{\omega u^\beta\right\}}\leq \exp\left\{\left(\frac{\sigma^2 p}{2}\right)^{\frac{\beta}{\beta-1}}\left(\frac{1}{\beta \omega}\right)^{\frac{1}{\beta-1}}\left(1-\frac{1}{\beta\omega}\right)\right\}.
\]
The second term is estimated similarly. Therefore,
\[
I(u)\leq 2 \exp\left\{\left(\frac{\sigma^2 p}{2}\right)^{\frac{\beta}{\beta-1}}\left(\frac{1}{\beta \omega}\right)^{\frac{1}{\beta-1}}\left(1-\frac{1}{\beta\omega}\right)\right\} \frac{\left(E\left|X(t)-X(s)\right|^{qu}\right)^{\frac{1}{qu}}}{\exp\left\{u^\beta(1-\omega)\right\}}.
\]
Let $\sigma^2(t,s)=E\left(X(t)-X(s)\right)^2$. Consider the expression
\[
\frac{\left(E\left|X(t)-X(s)\right|^{qu}\right)^{\frac{1}{qu}}}{\exp\left\{ u^\beta(1-\omega)\right\}}=\frac{\left(\int\limits_{-\infty}^{+\infty}\left|x\right|^{qu}\frac{1}{\sqrt{2\pi}\sigma(t,s)}\exp\left\{-\frac{x^2}{2\sigma^2(t,s)}\right\}dx\right)^{\frac{1}{qu}}}{\exp\left\{ u^\beta(1-\omega)\right\}}.
\]
Let's make the substitution $\frac{x}{\sigma}=v$, then we get:
\[
\frac{\left(\int\limits_{-\infty}^{+\infty}\left|x\right|^{qu}\frac{1}{\sqrt{2\pi}\sigma(t,s)}\exp\left\{-\frac{x^2}{2\sigma^2(t,s)}\right\}dx\right)^{\frac{1}{qu}}}{\exp\left\{u^\beta(1-\omega)\right\}}=
\]
\[
=\frac{\left((\sigma(t,s))^{qu}\frac{1}{\sqrt{2\pi}}\int\limits_{-\infty}^{+\infty}\left|v\right|^{qu}\exp\left\{-\frac{v^2}{2}\right\}dv\right)^{\frac{1}{qu}}}{\exp\left\{ u^\beta(1-\omega)\right\}}=
\]
\[
=\frac{\sigma(t,s)\left(\sqrt{\frac{2}{\pi}}\int\limits_{0}^{\infty}v^{qu}\exp\left\{-\frac{v^2}{2}\right\}dv\right)^{\frac{1}{qu}}}{\exp\left\{ u^\beta(1-\omega)\right\}}=\sigma(t,s)\frac{\left(\sqrt{\frac{2}{\pi}}2^{\frac{qu}{2}+\frac{1}{2}}\Gamma\left(\frac{qu}{2}+\frac{1}{2}\right)\right)^{\frac{1}{qu}}}{\exp\left\{ u^\beta(1-\omega)\right\}}.
\]
Let $$C(\beta,\omega)=\sup\limits_{u\geq1}\frac{\left(\sqrt{\frac{2}{\pi}}2^{\frac{qu}{2}+\frac{1}{2}}\Gamma\left(\frac{qu}{2}+\frac{1}{2}\right)\right)^{\frac{1}{qu}}}{\exp\left\{ u^\beta(1-\omega)\right\}}.$$ We have
\[
\left\|e^{X(t)}-e^{X(s)}\right\|=\sup\limits_{u\geq1}I(u)\leq Z\sigma(t,s),
\]
where $$Z=2C(\beta,\omega)\exp\left\{\left(\frac{\sigma^2 p}{2}\right)^{\frac{\beta}{\beta-1}}\left(\frac{1}{\beta \omega}\right)^{\frac{1}{\beta-1}}\left(1-\frac{1}{\beta\omega}\right)\right\},$$ $\sigma(t,s)=\sqrt{E\left(X(t)-X(s)\right)^2}$.
Let
\begin{equation}\label{eq:o2chapter-5-1}
\sigma(t,s)\leq\Delta\left|t-s\right|^\gamma,
\end{equation}
where $\Delta$ is a constant. In this case  we can use results of Example \ref{ex:o2chapter-8} and Example \ref{ex:o2chapter-11} to find the supremum of the process $X(t)$.
\end{example}

\begin{example}
Let in Example \ref{ex:o2chapter-15-1} $X(t)$ be a stationary Gaussian process, $EX(t)=0$, $EX(t)X(s)=B(t-s)$. It is obvious that
\[
E\left|X(t)-X(s)\right|^2=EX^2(t)+EX^2(s)-2EX(t)X(s)=2\left(B(0)-B(t-s)\right).
\]
Let $F(\lambda)$ be the spectral function of the process $X(t)$, i.e. $$B(\tau)= \int\limits_{0}^{\infty}\cos \lambda\tau dF(\lambda).$$ Then
\[
B(0)-B(\tau)=\int\limits_{0}^{\infty}(1-\cos \lambda\tau) dF(\lambda)=\int\limits_{0}^{\infty}2\sin^2 \frac{\lambda\tau}{2} dF(\lambda).
\]
We will use the inequality $\left|\sin u\right|\leq\left|u\right|^\gamma$, where $\gamma$ is any number, $0<\gamma<1$. So,
\[
B(0)-B(\tau)\leq \int\limits_{0}^{\infty}2\left(\frac{\lambda\left|\tau\right|}{2}\right)^{2\gamma}dF(\lambda)=2^{1-2\gamma}\left|\tau\right|^{2\gamma}\int\limits_{0}^{\infty}\lambda^{2\gamma}dF(\lambda).
\]
If $\int\limits_{0}^{\infty}\lambda^{2\gamma}dF(\lambda)<\infty$, then
\[
\left(E\left|X(t)-X(s)\right|^2\right)^{\frac{1}{2}}\leq \acute{C}\left|t-s\right|^\gamma,
\]
where $$\acute{C}=2^{\frac{1-2\gamma}{2}}\left(\int\limits_{0}^{\infty}\lambda^{2\gamma}dF(\lambda)\right)^{\frac{1}{2}}.$$ That is, the condition \eqref{eq:o2chapter-5-1} is satisfied.
\end{example}

\section{Probabilities of large deviations for sums of independent random processes from spaces $\mathbf{F}_\psi^*(\Omega)$}

Let us prove a Lemma that will be used in the following text.

\begin {lemma} \label{le:o2chapter_1}
Let $\xi \in \mathbf{F}_\psi(\Omega)$, $p\geq 1$. Then
\begin{equation} \label{eq:o2chapter_5_1}
\left\|E\left|\xi\right|\right\|_\psi\leq \left\|\xi\right\|_\psi.
\end{equation}
\end {lemma}

\begin{proof}
If $m$ is an arbitrary constant, then
\begin{equation}\label{eq:o2chapter_5_2}
\left\|m\right\|_\psi=\sup\limits_{u\geq1}\frac{\left(E\left|m\right|^{u}\right)^{1/u}}{\psi(u)}= \sup\limits_{u\geq1}\frac{m}{\psi(u)}=\frac{m}{\psi(1)}.
\end{equation}
From the definition of the norm $\left\|\xi\right\|_\psi$ it follows that
\begin{equation} \label{eq:o2chapter_5_3}
E\left|\xi\right|\leq \left\|\xi\right\|_\psi \psi(1).
\end{equation}
Therefore, the validity of the Lemma follows from  equality \eqref{eq:o2chapter_5_2} and  inequality \eqref{eq:o2chapter_5_3}.
\end{proof}

\begin{theorem}
\label{th:o2chapter_2}
Let $(T, \rho)$ be a compact metric space, $Y = \left\{Y(t), \right.$ $\left. t \in T\right\}$ be a random process that belongs to the space $\mathbf{F}_\psi^*(\Omega)$ and for which the condition $\mathbf{H}$ with the constant $C_\psi$ is satisfied. Let $Y$ be a separable process on $(T, \rho)$. Furthermore, let there be a continuous monotonically increasing function
$\sigma(h)$, $\sigma(0)=0$, such that
\[
\sup\limits_{\rho(t, s) \leq h} \bigl\|Y(t) - Y(s) \bigr\|_\psi \leq \sigma(h)
\]
and for any $z > 0$ the condition
\[
\int\limits_{0}^{z} \varkappa\left(N\left(\sigma^{(-1)}(u)\right)\right) du < \infty,
\]
where $\varkappa(n)$ is the majorizing characteristic of the space $\mathbf{F}_\psi(\Omega)$.

Let $X(t) = Y(t) - m(t)$, where $m(t) = EX(t)$ and $X_k (t), k=1,2,\dots$ are independent copies of the process $X(t)$. Let $S_n (t) = \frac{1}{\sqrt{n}} \sum\limits_{k = 1}^{n} X_{k}(t)$.

Then for any $0<p<1$ the inequality holds:
\[
\left\|\sup\limits_{t \in T}\left|S_n (t)\right|\right\|_\psi\leq \widehat{B}(p),
\]
where
$$\widehat{B}(p)=2\sqrt{C_\psi} \inf \limits_{t \in T}\left\|Y(t)\right\|_\psi+\frac{1}{p(1-p)}\int\limits_{0}^{\gamma p}\varkappa\left(N\left(\sigma_1^{(-1)}(u)\right)\right)du,$$ $$\gamma=\sigma_1\left(\sup\limits_{t, s \in T} \rho(t, s)\right)=2\sqrt{C_\psi}\sigma\left(\sup\limits_{t, s \in T} \rho(t, s)\right).$$
\end{theorem}

\begin{proof}
From Definition \ref{de:o1chapter_5} we have that
\[
\left\|S_n (t)-S_n(s)\right\|_\psi^2=\left\|\frac{1}{\sqrt{n}}\sum \limits_{k=1}^n\left(X_k(t)-X_k(s)\right)\right\|_\psi^2\leq
\]
\[
\leq C_\psi \frac{1}{n}\sum \limits_{k=1}^n\left\|X_k(t)-X_k(s)\right\|_\psi^2= C_\psi \left\|X(t)-X(s)\right\|_\psi^2
\]
and
\begin{multline*}
\left\|X(t)-X(s)\right\|_\psi^2 = \left\|Y(t)-Y(s)-(m(t)-m(s))\right\|_\psi^2 \leq \\[1ex]
\leq\left(\left\|Y(t)-Y(s)\right\|_\psi+\left\|m(t)-m(s)\right\|_\psi\right)^2 \leq   \\[1ex]
\leq\left(\left\|Y(t)-Y(s)\right\|_\psi+\frac{\left|m(t)-m(s)\right|}{\psi(1)}\right)^2.
\end{multline*}
From Lemma \ref{le:o2chapter_1} it follows that
\[
\left|m(t)-m(s)\right|\leq E\left\|Y(t)-Y(s)\right\|_\psi\leq \left\|Y(t)-Y(s)\right\|_\psi\psi(1).
\]
Then
\[
\left\|S_n (t)-S_n(s)\right\|_\psi^2\leq 4 C_\psi \left\|Y(t)-Y(s)\right\|_\psi^2.
\]
Therefore,
\[
\sup\limits_{\rho(t, s) \leq h} \bigl\|S_n (t)-S_n(s) \bigr\|_\psi \leq \sigma_1(h),
\]
where $\sigma_1(h)=2\sqrt{C_\psi}\sigma(h)$.

It is obvious that for any $t_0 \in T$
\[
\left\|X\left(t_0\right)\right\|_\psi^2\leq 4 C_\psi \left\|Y\left(t_0\right)\right\|_\psi^2.
\]
Therefore, the validity of the Theorem follows from Theorem \ref{th:o2chapter_1}.

\end{proof}

\begin{corollary}\label{th:o3chapter-5}
Let a process $X = \{X(t), \ t \in T\}$ belongs to the space $\mathbf{F}_\psi(\Omega)$. Let the conditions of Theorem \ref{th:o2chapter_2} be satisfied.
Then for any $\varepsilon>0$ the inequality holds:
\[
P\left\{\sup\limits_{t \in T}\left|S_n (t)\right|>\varepsilon\right\}\leq \inf\limits_{u\geq1}\frac{\widehat{B}^u(p)(\psi(u))^u}{\varepsilon^u}.
\]
\end{corollary}

\begin{proof}
The validity of Corollary \ref{th:o3chapter-5} follows from Corollary \ref{co:o2chapter-1}.
\end{proof}

\section{Distribution of suprema of increments of random processes from spaces $\mathbf{F}_\psi(\Omega)$}

Let $(T,\rho)$ be a compact metric space, $V_{\theta_l}$ be a minimal covering of the set $T$ by balls of radius at most $\theta_l$, where $\theta_l=\sigma^{(-1)}(bp^l)$, $b= \sup\limits_{t, s \in T}\rho(t,s)$, $0<p<1$.

\begin{theorem}\label{th:o2chapter_7}
\label{the1}
Let $X = \{X(t), t \in T\}$ be a separable random process on $(T,\rho)$ from the space $\mathbf{F}_\psi(\Omega)$ and the condition
\begin{equation}\label{eq:o2chapter_6}
\sup\limits_{\rho(t, s) \leq h} \bigl\|X(t) - X(s) \bigr\|_\psi\leq \sigma(h),
\end{equation}
be satisfied, where $\sigma(h)$ is a continuous monotonically increasing function such that $\sigma(0)=0$. If for any $z>0$ the condition holds
\begin{equation}\label{eq:o2chapter_7}
\int\limits_{0}^{z}\varkappa\left(N\left(\sigma^{(-1)}(u)\right)\right)du<\infty,
\end{equation}
where $\varkappa(n)$ is the majorizing characteristic, $\sigma^{(-1)}(u)$ is the inverse of the function to $\sigma(u)$, $0<\theta\leq z$ is some number, then the inequality holds true
\begin{multline*}
\left\|\sup\limits_{\rho(t, s) \leq \theta} \left|X(t) - X(s)\right| \right\|_\psi\leq \varkappa\left(D_k(\theta,p )\right)\sigma(\theta) \frac{3-p}{1-p}+\\[1ex]
+2\frac{1}{p(1-p)}\int\limits_{0}^{\sigma(\theta)p}\varkappa\left(N\left(\sigma^{(-1)}(u)\right)\right)du=S_k(\theta,p),
\end{multline*}
where $D_k(\theta,p )$ is the number of points $u, v$ in $V_{\theta_k}$ such that $$\left\|X(u)-X(v)\right\|\leq\sigma(\theta) \frac{3-p}{1-p}$$ and for any $\varepsilon>0$ the inequality holds:
\begin{equation} \label{eq:o2chapter_8}
P\left\{\sup\limits_{\rho(t, s) \leq \theta} \left|X(t) - X(s)\right| > \varepsilon\right\}\leq \inf\limits_{u\geq1}\frac{\left(S_k(\theta,p)\right)^u(\psi(u))^u}{\varepsilon^u}.
\end{equation}
\end{theorem}

\begin{proof}
Note that the following inequalities hold:
\[
\sigma(\theta_k)<\sigma(\theta)\leq\sigma(\theta_{k-1}),\ bp^k <\sigma(\theta)\leq bp^{k-1},
\]
where $k$ is an integer such that $\theta_k<\theta\leq\theta_{k-1}$.
Let us denote $V_k=\bigcup\limits_{j=k}^{\infty} V_{\theta_j}$. Since under condition \eqref{eq:o2chapter_6} the process is continuous in probability (see the proof of Theorem \ref{th:o2chapter_1}), then $V_k$ is the separant, therefore
\[
\sup\limits_{\substack{\rho(t, s)\leq \theta\\ t,s \in T}} \left|X(t) - X(s)\right| = \sup\limits_{\substack{\rho(t, s)\leq \theta\\ t,s \in V_k}} \left|X(t) - X(s)\right|.
\]

Let $t$ and $s$ be points such that $t,s \in V_k$ and $\rho(t, s) \leq \theta$. Since $t$ and $s$ belong to $V_k$, there exist $m\geq k$ and $r\geq k$ such that $t \in V_{\theta_m}$ and $s \in V_{\theta_r}$. Denote $t_m=t$, $t_{m-1}=\alpha_{m-1}(t_m)$, $t_{m-2}=\alpha_{m-2}(t_{m-1})$, \ldots, $t_k=\alpha_k(t_{k+1})$ and $s_r=s$, $s_{r-1}=\alpha_{r-1}(s_r)$. $s_{r-2}=\alpha_{r-2}(s_{r-1})$, \ldots, $s_k=\alpha_k(s_{k+1})$. It is obvious that the equality holds:
\begin{multline} \label{eq:o2chapter_9}
X(t) - X(s)=\\[1ex]
=\sum \limits_{l=k}^{m-1}\left(X(t_{l+1}) - X(t_l)\right)+\left(X(t_k) - X(s_k)\right)-\sum \limits_{l=k}^{r-1}\left(X(s_{l+1}) - X(s_l)\right).
\end{multline}
From the equality \eqref{eq:o2chapter_9} we have:
\begin{multline*}
\left|X(t_k) - X(s_k)\right|\leq \left|X(t) - X(s)\right|+\sum \limits_{l=k}^{m-1}\left|X(t_{l+1}) - X(t_l)\right|+ \\
+\sum \limits_{l=k}^{r-1}\left|X(s_{l+1}) - X(s_l)\right|.
\end{multline*}
Therefore,
\[
\left\|X(t_k) - X(s_k)\right\|\leq \left\|X(t) - X(s)\right\|+\sum \limits_{l=k}^{m-1}\left\|X(t_{l+1}) - X(t_l)\right\|+
\]
\[
+\sum \limits_{l=k}^{r-1}\left\|X(s_{l+1}) - X(s_l)\right\|\leq \sigma(\theta)+\sum \limits_{l=k}^{m-1}\sigma(\theta_l)+\sum \limits_{l=k}^{r-1}\sigma(\theta_l)\leq \sigma(\theta)+2\sum \limits_{l=k}^{\infty}\sigma(\theta_l)=
\]
\[
=\sigma(\theta)+2\sum \limits_{l=k}^{\infty} bp^l=\sigma(\theta)+2b\frac{p^k}{1-p}\leq \sigma(\theta)+2\sigma(\theta)\frac{1}{1-p}=\sigma(\theta)\frac{3-p}{1-p}.
\]
Then from the equality \eqref{eq:o2chapter_9} we obtain the inequality:
\begin{multline} \label{eq:o2chapter_11}
\left|X(t) - X(s)\right|\leq\sum \limits_{l=k}^{m-1}\left|X(t_{l+1}) - X(t_l)\right|+ \sum \limits_{l=k}^{r-1}\left|X(s_{l+1}) - X(s_l)\right|+ \\ +\left|X(t_k) - X(s_k)\right|\leq \sum \limits_{l=k}^{m-1}\max\limits_{u \in V_{\theta_{l+1}}}\left|X(u) - X(\alpha_l(u))\right|+ \\
+\sum \limits_{l=k}^{r-1}\max\limits_{u \in V_{\theta_{l+1}}}\left|X(u) - X(\alpha_l(u))\right|+\max\limits_{\substack{u,v \in V_{\theta_k}\\\left\|X(u)-X(v)\right\|\leq\sigma(\theta) \frac{3-p}{1-p}}}\left|X(u) - X(v))\right|\leq \\[1ex]
\leq \max\limits_{\substack{u,v \in V_{\theta_k}\\\left\|X(u)-X(v)\right\|\leq\sigma(\theta) \frac{3-p}{1-p}}}\left|X(u) - X(v))\right|+2\sum \limits_{l=k}^{\infty}\max\limits_{u \in V_{\theta_{l+1}}}\left|X(u) - X(\alpha_l(u))\right|.
\end{multline}

Since the last expression in the inequality \eqref{eq:o2chapter_11} does not depend on $t, s \in V_k$, such that $\rho(t, s) \leq \theta$, then, taking into account the separability of the process $X$, with probability one we obtain the inequality:
\begin{multline} \label{eq:o2chapter_12}
\sup\limits_{\substack{\rho(t, s)\leq \theta\\ t,s \in T}} \left|X(t) - X(s)\right| = \sup\limits_{\substack{\rho(t, s)\leq \theta\\ t,s \in V_k}} \left|X(t) - X(s)\right|\leq \\
\leq\max\limits_{\substack{u,v \in V_{\theta_k}\\ \left\|X(u)-X(v)\right\|\leq\sigma(\theta) \frac{3-p}{1-p}}}\left|X(u) - X(v)\right|+ 2\sum \limits_{l=k}^{\infty}\max\limits_{u \in V_{\theta_{l+1}}}\left|X(u) - X(\alpha_l(u))\right|.
\end{multline}
From the inequality \eqref{eq:o2chapter_12} and Definition \ref{eq:o2chapter-3} we obtain the following inequalities:
\[
\left\|\sup\limits_{\rho(t, s) \leq \theta} \left|X(t) - X(s)\right| \right\|\leq
\]
\[
\leq\left\| \max\limits_{\substack{u,v \in V_{\theta_k}\\ \left\|X(u)-X(v)\right\|\leq\sigma(\theta) \frac{3-p}{1-p}}}\left|X(u) - X(v)\right|\right\|+2\sum \limits_{l=k}^{\infty} \left\|\max\limits_{u \in V_{\theta_{l+1}}}\left|X(u) - X(\alpha_l(u))\right|\right\|\leq
\]
\[
\leq\varkappa\left(D_k(\theta,p )\right)\sigma(\theta) \frac{3-p}{1-p}+2\sum \limits_{l=k}^{\infty}\varkappa\left(N(\theta_{l+1})\right)\sigma(\theta_l).
\]

Let's transform the second term on the right side of the last inequality
\[
\sum \limits_{l=k}^{\infty}\varkappa\left(N(\theta_{l+1})\right)\sigma(\theta_l)=\sum \limits_{l=k}^{\infty}\varkappa\left(N\left(\sigma^{(-1)}(bp^{l+1})\right)\right)bp^l,
\]
then we get:
\[
\int\limits_{bp^{l+2}}^{bp^{l+1}}\varkappa\left(N\left(\sigma^{(-1)}(u)\right)\right)du\geq \varkappa\left(N\left(\sigma^{(-1)}(bp^{l+1})\right)\right)(bp^{l+1}-bp^{l+2}),
\]
\[
bp^l\varkappa\left(N\left(\sigma^{(-1)}(bp^{l+1})\right)\right)\leq \frac{1}{p(1-p)}\int\limits_{bp^{l+2}}^{bp^{l+1}}\varkappa\left(N\left(\sigma^{(-1)}(u)\right)\right)du.
\]
Therefore,
\begin{multline*}
\sum \limits_{l=k}^{\infty}\varkappa\left(N(\theta_{l+1})\right)\sigma(\theta_l)\leq \frac{1}{p(1-p)}\int\limits_{0}^{bp^{k+1}}\varkappa\left(N\left(\sigma^{(-1)}(u)\right)\right)du \leq \\[1ex] \leq\frac{1}{p(1-p)}\int\limits_{0}^{p\sigma(\theta)}\varkappa\left(N\left(\sigma^{(-1)}(u)\right)\right)du.
\end{multline*}
To complete the proof of the Theorem, we note that inequality \eqref{eq:o2chapter_8} follows from Theorem \ref{th:o1chapter_1}.
\end{proof}

\begin{corollary}\label{co:o2chapter_6}
Let the conditions of  Theorem~\ref{th:o2chapter_7} be satisfied and the condition $\sigma(\theta)\varkappa\left(D_k(\theta,p )\right)\rightarrow 0$ be true as $\theta\rightarrow 0$. Then
\begin{equation} \label{eq:o2chapter_13}
\left\|\sup\limits_{\rho(t, s) \leq \theta} \left|X(t) - X(s)\right| \right\|\rightarrow 0
\end{equation}
for $\theta\rightarrow 0$. Furthermore, the random process $X(t)$ is uniformly continuous on $(T,\rho)$ with probability one.
\end{corollary}

\begin{proof}
Since from the Corollary condition for $\theta\rightarrow 0$ we have $S_k(\theta,p)\rightarrow 0$, then the relation \eqref{eq:o2chapter_13} is obvious.
From the inequality \eqref{eq:o2chapter_8} it follows that for any $u>1$ the estimate holds true
\[
P\left\{\sup\limits_{\rho(t, s) \leq \theta} \left|X(t) - X(s)\right| > \varepsilon\right\}\leq \frac{\left(S_k(\theta,p)\right)^u(\psi(u))^u}{\varepsilon^u}.
\]
That is why for $\theta\rightarrow 0$ $$\sup\limits_{\rho(t, s) \leq \theta} \left|X(t) - X(s)\right|\rightarrow 0.$$
Therefore, there exists a sequence $\check{\theta}_n$ such that $\check{\theta}_{n+1}<\check{\theta}_n$, $\check{\theta}_n\rightarrow 0$ for $n\rightarrow \infty$, and
$$\sup\limits_{\rho(t, s) \leq \check{\theta}_n} \left|X(t) - X(s)\right|\rightarrow 0$$ as $n\rightarrow \infty$ with probability one.
Since $$\sup\limits_{\rho(t, s) \leq \theta} \left|X(t) - X(s)\right|$$ monotonically decreases in $\theta$, then $$\sup\limits_{\rho(t, s) \leq \theta} \left|X(t) - X(s)\right|\rightarrow 0$$ as $\theta\rightarrow 0$ with probability one.
\end{proof}

\begin{theorem} \label{th:o2chapter_8}
Let $X = \{X(t),\ t \in T\}$ be a separable random process on $(T,\rho)$ from the space $\mathbf{F}_\psi(\Omega)$, where $\psi(u)$ is a weight function such that for the majorizing characteristic of the space $\mathbf{F}_\psi(\Omega)$ the inequality holds:
\begin{equation} \label{eq:o2chapter_14}
\varkappa(n^2)\leq C \varkappa(n),
\end{equation}
where $C>0$ is a constant.

If the conditions \eqref{eq:o2chapter_6} and \eqref{eq:o2chapter_7} are satisfied, then
\begin{equation} \label{eq:o2chapter_15}
\left\|\sup\limits_{\rho(t, s) \leq \theta} \left|X(t) - X(s)\right| \right\|_\psi\leq A(C)\int\limits_{0}^{\sigma(\theta)}\varkappa\left(N\left(\sigma^{(-1)}(u)\right)\right)du=\check{S}_k(\theta),
\end{equation}
where $$A(C)=\frac{4\left(3C^2-12C+4\right)}{3C+2}\cdot\left(\frac{C+2}{C-2}\right)^2.$$
Moreover, for any $\varepsilon>0$ the inequality holds:
\begin{equation} \label{eq:o2chapter_16}
P\left\{\sup\limits_{\rho(t, s) \leq \theta} \left|X(t) - X(s)\right| > \varepsilon\right\}\leq \inf\limits_{u\geq1}\frac{\left(\check{S}_k(\theta)\right)^u(\psi(u))^u}{\varepsilon^u}.
\end{equation}
In addition, the random process $X(t)$ is uniformly continuous on the space $(T,\rho)$ with probability one.
\end{theorem}

\begin{proof}
Theorem~\ref{th:o2chapter_8} follows from Theorem~\ref{th:o2chapter_7}. Indeed,
\[
D_k(\theta,p )\leq \left(N(\theta_k)\right)^2\leq \left(N\left(\sigma^{(-1)}(bp^k)\right)\right)^2,
\]
therefore
\[
\sigma(\theta)\varkappa\left(D_k(\theta,p )\right)\leq \varkappa\left(N^2\left(\sigma^{(-1)}(bp^k)\right)\right)\sigma(\theta).
\]
Since the inequality \eqref{eq:o2chapter_14} holds, then
\begin{equation} \label{eq:o2chapter_17}
\varkappa\left(N^2\left(\sigma^{(-1)}(bp^k)\right)\right)\sigma(\theta)\leq C\varkappa\left(N\left(\sigma^{(-1)}(bp^k)\right)\right)\sigma(\theta).
\end{equation}
It is obvious that
\[
\int\limits_{bp^{k+1}}^{bp^k}\varkappa\left(N\left(\sigma^{(-1)}(u)\right)\right)du\geq \varkappa\left(N\left(\sigma^{(-1)}(bp^k)\right)\right) bp^k (1-p),
\]
therefore
\begin{multline*}
\varkappa\left(N\left(\sigma^{(-1)}(bp^k)\right)\right)\leq \frac{1}{bp^{k-1}p(1-p)}\int\limits_{bp^{k+1}}^{bp^k}\varkappa\left(N\left(\sigma^{(-1)}(u)\right)\right)du\leq \\[1ex] \leq\frac{1}{\sigma(\theta)p(1-p)}\int\limits_{0}^{\sigma(\theta)}\varkappa\left(N\left(\sigma^{(-1)}(u)\right)\right)du.
\end{multline*}
Therefore, from the inequality \eqref{eq:o2chapter_17} it follows:
\[
\varkappa\left(D_k(\theta,p )\right)\sigma(\theta) \frac{3-p}{1-p}+\frac{2}{p(1-p)}\int\limits_{0}^{\sigma(\theta)p}\varkappa\left(N\left(\sigma^{(-1)}(u)\right)\right)du \leq
\]
\[
\leq\varkappa\left(N^2\left(\sigma^{(-1)}(bp^k)\right)\right)\sigma(\theta)\frac{3-p}{1-p}+\frac{2}{p(1-p)}\int\limits_{0}^{\sigma(\theta)}\varkappa\left(N\left(\sigma^{(-1)}(u)\right)\right)du \leq
\]
\[
\leq C\varkappa\left(N\left(\sigma^{(-1)}(bp^k)\right)\right)\sigma(\theta)\frac{3-p}{1-p}+\frac{2}{p(1-p)}\int\limits_{0}^{\sigma(\theta)}\varkappa\left(N\left(\sigma^{(-1)}(u)\right)\right)du \leq
\]
\[
\leq\frac{1}{p(1-p)}\left(C\cdot\frac{3-p}{1-p}+2\right)\int\limits_{0}^{\sigma(\theta)}\varkappa\left(N\left(\sigma^{(-1)}(u)\right)\right)du.
\]
The minimum of the expression $\frac{1}{p(1-p)}\left(C\cdot\frac{3-p}{1-p}+2\right)$ is reached at the point $p=\frac{3C+2}{2C+4}$. Substituting this value into the last inequality, we obtain the inequality \eqref{eq:o2chapter_15}. And the inequality \eqref{eq:o2chapter_16} follows from  Theorem \ref{th:o1chapter_1}. The assertion about sample continuity is proved similarly to the proof of Corollary \ref{co:o2chapter_6}.
\end{proof}

\begin{remark}\label{re:o2chapter-6}
The conditions of the Theorem are satisfied for the functions $\psi(u)=u^\alpha$, $\psi(u)=\left(\ln (u+1)\right)^\lambda$.
If $\psi(u)=u^\alpha$, then $\varkappa(n^2)=2^\alpha \varkappa(n)$, i.e. $C=2^\alpha$, and if $\psi(u)=\left(\ln (u+1)\right)^\lambda$, then $\varkappa(n^2)=\left(\frac{\ln(2(\ln n+1))}{\ln 2}\right)^\lambda e\leq 2^\lambda \varkappa(n)$, i.e. $C=2^\lambda$.
\end{remark}

\begin{remark}
Let $W(h), \ h>0$ be a continuous, monotonically decreasing function, and $N(h)\leq W(h)$. Then the condition \eqref{eq:o2chapter_7} is satisfied if
\begin{equation} \label{eq:o2chapter_18}
\int\limits_{0}^{z}\varkappa\left(W\left(\sigma^{(-1)}(u)\right)\right)du<\infty.
\end{equation}
\end{remark}

\begin{corollary} \label{co:o2chapter_1}
Let in Theorem~\ref{th:o2chapter_8}
\[
\sigma(h)=\frac{R}{\left(\varkappa\left(W(h)\right)\right)^{1/\gamma}},
\]
where $\gamma$ is a number such that $0<\gamma<1$. Then
\[
\left\|\sup\limits_{\rho(t, s) \leq \theta} \left|X(t) - X(s)\right| \right\|\leq
A(C) \frac{R}{1-\gamma}\left(\varkappa\left(W(\theta)\right)\right)^{\frac{\gamma-1}{\gamma}}= Q(\theta,\gamma)
\]
and for any $\varepsilon>0$ the inequality holds:
\begin{equation} \label{eq:o2chapter_19}
P\left\{\sup\limits_{\rho(t, s) \leq \theta} \left|X(t) - X(s)\right| > \varepsilon\right\}\leq \inf\limits_{u\geq1}\frac{\left(Q(\theta,\gamma)\right)^u(\psi(u))^u}{\varepsilon^u}.
\end{equation}
\end{corollary}

\begin{proof}
Indeed, in this case $\sigma^{(-1)}(u)=W^{(-1)}\left(\varkappa^{(-1)}\left(\frac{R^\gamma}{u^\gamma}\right)\right)$. Then
\begin{equation} \label{eq:o2chapter_20}
\int\limits_{0}^{\sigma(\theta)}\varkappa\left(W\left(\sigma^{(-1)}(u)\right)\right)du=\int\limits_{0}^{\sigma(\theta)}\frac{R^\gamma}{u^\gamma}du=\frac{R^\gamma}{1-\gamma}\sigma(\theta)^{1-\gamma}<\infty.
\end{equation}
\end{proof}

\begin{example}
If in Corollary \ref{co:o2chapter_1} the weight function is represented as: $\psi(u)=u^\alpha$, where $\alpha>0$, then for
$\varepsilon\geq e^\alpha Q(\theta,\gamma)$

\[
P\left\{\sup\limits_{\rho(t, s) \leq \theta} \left|X(t) - X(s)\right|>\varepsilon\right\}\leq \exp \left\{-\frac{\alpha}{e}\left(\frac{\varepsilon}{Q(\theta,\gamma)}\right)^{1/\alpha}\right\}
\]
and, using the value of the majorizing characteristic, we obtain:
\[
\sigma(h)=R\left(\frac{e}{\alpha}\ln W(h)\right)^{-\alpha/\gamma}, \ 0<\gamma<1.
\]
\end{example}

\begin{example}
If in  Corollary \ref{co:o2chapter_1} the weight function has the form $\psi(u)=e^{au^\beta}$, where $a>0$, $\beta>0$, then for
$\varepsilon>e^{a(\beta+1)}Q(\theta,\gamma)$
$$
P\left\{\sup\limits_{\rho(t, s) \leq \theta} \left|X(t) - X(s)\right|>\varepsilon\right\}\leq \exp \left\{-\frac{\beta}{a^{1/\beta}}\left(\frac{\ln \frac{\varepsilon}{Q(\theta,\gamma)}}{\beta+1}\right)^{\frac{\beta+1}{\beta}}\right\}
$$
and using the value of the major characteristic we obtain:
$$
\sigma(h)=R\cdot\exp \left\{-\frac{2}{\gamma}\sqrt{a \ln W(h)}+\frac{a}{\gamma}\right\}, \ 0<\gamma<1.
$$
\end{example}

\begin{example}
If in  Corollary \ref{co:o2chapter_1} the weight function is written as: $\psi(u)=\left(\ln (u+1)\right)^\lambda$, where $\lambda>0$,
then for $\varepsilon\geq\left(e \ln 2\right)^\lambda Q(\theta,\gamma)$
\[
P\left\{\sup\limits_{\rho(t, s) \leq \theta} \left|X(t) - X(s)\right|>\varepsilon\right\}\leq
e^\lambda \exp \left\{-\lambda \exp\left\{\left(\frac{\varepsilon}{Q(\theta,\gamma)}\right)^{1/\lambda}\frac{1}{e}\right\}\right\}
\]
and, using the value of the majorizing characteristic, we will get:
\[
\sigma(h)=\frac{R}{\left(\frac{1}{e^a}\exp \left\{S(a,\beta)(\ln W(h))^{\frac{\beta}{\beta+1}}\right\}\right)^\frac{1}{\mu}},
\]
where $$S(a,\beta)=(\beta a)^{\frac{1}{\beta+1}}(\beta^{-1}+1),\, 0<\gamma<1.$$
\end{example}

\begin{theorem} \label{th:o2chapter_10}
Let $X = \{X(t), \ t \in \left[c,d\right]\}$, $c<d$ be a separable random process from the space $\mathbf{F}_\psi(\Omega)$,  the condition \eqref{eq:o2chapter_14} holds true and
\[
\sup\limits_{\substack{\left|t-s\right| \leq h\\ t,s \in \left[c,d\right]}} \bigl\|X(t) - X(s) \bigr\|_\psi\leq \sigma(h),
\]
where $\sigma=\left\{\sigma(h),\ h>0\right\}$ is a continuous, monotonically increasing function such that $\sigma(0)=0$. If for any $z>0$ there is satisfied the condition
\[
\int\limits_{0}^{z}\varkappa\left(\frac{d-c}{2\sigma^{(-1)}(u)}+1\right)du<\infty,
\]
then the inequality holds:
\[
\left\|\sup\limits_{t \in \left[c,d\right]} \left|X(t) - X(s)\right| \right\|_\psi\leq
A(C) \frac{R}{1-\gamma}\left(\varkappa\left(\frac{d-c}{2\theta}+1\right)\right)^{\frac{\gamma-1}{\gamma}}= \tilde{Q}(\theta,\gamma).
\]
Moreover, for any $\varepsilon>0$, the following estimate holds
\begin{equation} \label{eq:o2chapter_21}
P\left\{\sup\limits_{t \in \left[c,d\right]} \left|X(t) - X(s)\right| > \varepsilon\right\}\leq \inf\limits_{u\geq1}\frac{\left(\tilde{Q}(\theta,\gamma)\right)^u(\psi(u))^u}{\varepsilon^u}.
\end{equation}
Furthermore, the random process $X(t)$ is uniformly continuous on the space $(T,\rho)$ with probability one.
\end{theorem}

\begin{proof}
Theorem~\ref{th:o2chapter_10} follows from Theorem~\ref{th:o2chapter_8}, since the metric massiveness of the segment $\left[c,d\right]$ is estimated as follows:
\[
N(u)\leq \frac{d-c}{2u}+1.
\]
\end{proof}

\begin{example}
Consider the space $\mathbf{F}_\psi(\Omega)$, where $\psi(u)=u^\alpha$, $\alpha>0$, then for \\ $\varepsilon\geq e^\alpha A(C) \frac{R}{1-\gamma}\left(\frac{e}{\alpha}\ln\left(\frac{d-c}{2\theta}+1\right)\right)^{\frac{\alpha(\gamma-1)}{\gamma}}$ we have that
\begin{multline*}
P\left\{\sup\limits_{t \in \left[c,d\right]} \left|X(t) - X(s)\right|>\varepsilon\right\}\leq \\[1ex]
\leq\exp \left\{-\frac{\alpha}{e}\left(\frac{\varepsilon}{A(C) \frac{R}{1-\gamma}\left(\frac{e}{\alpha}\ln\left(\frac{d-c}{2\theta}+1\right)\right)^{\frac{\alpha(\gamma-1)}{\gamma}}}\right)^{1/\alpha}\right\}.
\end{multline*}
\end{example}

Similar estimates can be obtained in all examples given above.

\section{Conditions for convergence of wavelet decompositions of random processes from spaces $\mathbf{F}_\psi(\Omega)$}

{Wavelet bases and decompositions of functions in these bases}

Let $\varphi=\left\{\varphi(x), x \in \mathbb{R}\right\} \in L_2(\mathbb{R})$, and let $$\widehat{\varphi}(y)=\int\limits_{\mathbb{R}} e^{-iyx}\varphi(x)dx$$ be the Fourier transform of the function $\varphi:\varphi_{0k}(x)=\varphi(x-k)$.

\begin {definition}[see \cite{daub92,koz04}]
A function $\varphi$ is called an $f$-wavelet if the following conditions are met:
\begin{itemize}
\item [1)] $\sum\limits_{k \in \mathbb{Z}}\left|\widehat{\varphi}(y + 2 \pi k)\right|^2=1$ almost everywhere;
\item[2)] there exists a $2\pi$-periodic function $m_0(x) \in L_2\left(\left[0;2\pi\right]\right)$ such that
$\widehat{\varphi}(y)=m_0\left(\frac{y}{2}\right)\widehat{\varphi}\left(\frac{y}{2}\right)$ almost everywhere;
\item[3)] $\widehat{\varphi}(0)\neq 0$ and $\widehat{\varphi}(y)$ are continuous at zero.
\end{itemize}
\end{definition}

\begin {definition}[see see \cite{daub92,koz04}]
The function $\delta(x)$ is called an $m$-wavelet corresponding to an $f$-wavelet $\varphi$ if its Fourier transform has the form:
\[
\widehat{\delta}(y)=\overline{m_0\left(\frac{y}{2}+\pi\right)} \exp \left\{-i \frac{y}{2}\right\}\widehat{\varphi}\left(\frac{y}{2}\right).
\]
\end{definition}

Let $\varphi_{jk}(x)=2^{\frac{j}{2}}\varphi\left(2^j x-k\right)$, $\delta_{jk}(x)=2^{\frac{j}{2}}\delta\left(2^j x-k\right)$, $j,k \in \mathbb{Z}$.

It is known (see \cite{Hardle98,koz04}) that the system of functions $\left\{\varphi_{0k}, \delta_{jk}, j=0,\ldots, k \in \mathbb{Z}\right\}$
is an orthonormal basis in $L_2(\mathbb{R})$. Any function $f \in L_2(\mathbb{R})$
can be represented as a series that converges in the mean square

\begin{equation} \label{Mlavetse3}
f(x)=\sum\limits_{k \in \mathbb{Z}}\alpha_{0k}\varphi_{0k}(x)+\sum\limits_{j=0}^{\infty}\sum\limits_{k \in \mathbb{Z}} \beta_{jk}\delta_{jk}(x),
\end{equation}
where
\begin{equation} \label{Mlavetse4}
\alpha_{0k}=\int\limits_\mathbb{R} f(x) \overline{\varphi_{0k}(x)}dx, \ \beta_{jk}=\int\limits_\mathbb{R} f(x) \overline{\delta_{jk}(x)}dx
\end{equation}
and
\[
\sum\limits_{k \in \mathbb{Z}}\left|\alpha_{0k}\right|^2+\sum\limits_{j=0}^{\infty}\sum\limits_{k \in \mathbb{Z}} \left|\beta_{jk}\right|^2<\infty.
\]
The representation \eqref{Mlavetse3} is called a wavelet representation.

\begin{remark}
Since the integrals defined in the equalities \eqref{Mlavetse4} for $\alpha_{0k}$ and $\beta_{jk}$ exist not only for functions in $L_2(\mathbb{R})$,
it is possible to obtain wavelet decompositions for a wider class of functions than the space $L_2(\mathbb{R})$, which will coincide in certain norms.
\end{remark}

\begin {definition} [see \cite{daub92,koz04}]
Let $\varphi$ be an $f$-wavelet. The condition $S$ is satisfied for $\varphi$, if there exists an even function $\Phi=\left\{\Phi(x),\ x \in \mathbb{R}\right\}$
such that $\Phi(0)<\infty$, $\Phi(x)$ -- monotonically decreases at
$x\geq 0$, $\int\limits_{\mathbb{R}}\Phi(\left|x\right|)<\infty$ and $\left|\varphi(x)\right|\leq \Phi(\left|x\right|)$ for $x \in \mathbb{R}$.
\end{definition}

\begin{lemma}[see \cite{dari11}]
Let the $f$-wavelet $\varphi$ satisfy the condition $S$ with the function $\Phi$, and $\delta(x)$ is the $m$-wavelet
corresponding to $\varphi$. Then for all $x \in \mathbb{R}$ the inequality holds
\[
\left|\delta(x)\right|\leq B \Phi \left(\left|\frac{2x-1}{4}\right|\right),
\]
where $0<B<\infty$ is some constant.
\end{lemma}

\begin{theorem}[see \cite{dari11}] \label{Mlavetst4}
Let the $f$-wavelet $\varphi$ satisfy the condition $S$ with the function $\Phi$, $c=\left\{c(x), \ x \in \mathbb{R}\right\}$ -- an even function such that
$c(x)>1$, $x \in \mathbb{R}$, $c(x)$ -- monotonically increases for $x>0$ and $\int\limits_{\mathbb{R}} c(x) \Phi(\left|x\right|)<\infty$.
In addition, there exists a function $0<A(u)<\infty$, $u>0$ such that for sufficiently large $x$
\[
c(ax)\leq c(x)\cdot A(a), \ a>0.
\]

If $f=\left\{f(x), x \in \mathbb{R} \right\}$ is a measurable function on $\mathbb{R}$ such that $\left|f(x)\right|\leq c(x)$, $x \in \mathbb{R}$; $f(x)$ is continuous on the interval $(a,b)$, $-\infty<a<b<+\infty$, then
\[
f_m\left(x\right)=\sum\limits_{k \in \mathbb{Z}}\alpha_{0k}\varphi_{0k}(x)+\sum\limits_{j=0}^{m-1}\sum\limits_{k \in \mathbb{Z}} \beta_{jk}\delta_{jk}(x)\xrightarrow [m\rightarrow\infty]  {}f(x),
\]
\[
\alpha_{0k}=\int\limits_\mathbb{R} f(x) \overline{\varphi_{0k}(x)}dx, \ \beta_{jk}=\int\limits_\mathbb{R} f(x) \overline{\delta_{jk}(x)}dx,
\]
uniformly on each segment $\left[\alpha,\beta\right]\subset (a,b)$.
\end{theorem}

{Conditions for uniform convergence of wavelet decompositions of random processes from spaces $\mathbf{F}_\psi(\Omega)$}

The following theorem gives general conditions for uniform convergence of wavelet decompositions for random processes from spaces $\mathbf{F}_\psi(\Omega)$.

\begin{theorem} \label{Mlavetst5}
Let $X = \{X(t), t \in \mathbb{R}\}$ be a separable, measurable, random process from
$\mathbf{F}_\psi(\Omega)$, $T_k=\left[a_k,a_{k+1}\right]$, $-\infty<a_k<a_{k+1}<+\infty$, $k \in \mathbb{Z}$.
For each $T_k$ there exists a continuous strictly monotone increasing function $\sigma_k(h)$,
$0\leq h\leq\left(a_{k+1}-a_k\right)$, $\sigma_k(0)=0$, such that
\[
\sup\limits_{\substack{\left|t-s\right| \leq h\\ t,s \in T_k}}\bigl\|X(t) - X(s) \bigr\|_\psi\leq \sigma_k(h).
\]
Let the following conditions be satisfied:
\begin{enumerate}
\item [1)] $\int\limits_{0}^{\gamma_k}\varkappa\left(\frac{a_{k+1}-a_k}{2\sigma_k^{(-1)}(u)}+1\right)du<\infty$, where $\varkappa(n)$ is the majorizing characteristic, $\gamma_k=\sigma_k\left(\frac{a_{k+1}-a_k}{2}\right)$;
\item [2)] there exists a continuous function $c=\left\{c(t), \ t \in \mathbb{R}\right\}$ such that
\[
c(t)>1, \ r_k=\inf\limits_{t \in T_k}c(t);
\]
\item [3)] for any $\varepsilon>0$ the series $\sum\limits_{k \in \mathbb{Z}} \inf\limits_{u\geq1}\frac{B_k^u(\psi(u))^u}{\left(\varepsilon r_k\right)^u}$ are convergent,
 where $$B_k=\inf \limits_{t \in T_k}\left\|X(t)\right\|_\psi + \frac{1}{p_k(1-p_k)}\int\limits_{0}^{\gamma_k p_k}\varkappa\left(\frac{a_{k+1}-a_k}{2\sigma_k^{(-1)}(u)}+1\right)du,$$
 $p_k$, $0<p_k<1$, are some numbers;
\item [4)] the function $\psi(u)$ is such that for the majorizing characteristic of the space $\mathbf{F}_\psi(\Omega)$ the condition is satisfied
\begin{equation} \label{Mlavetse5}
\varkappa(n^2)\leq C_{\varkappa} \varkappa(n),
\end{equation}
where $C_{\varkappa}>0$ is a constant;
\item [5)] for a process $X$ from an interval $(a,b)$ the condition is satisfied:
\[
\sup\limits_{\rho(t, s) \leq h} \bigl\|X(t) - X(s) \bigr\|_\psi \leq \sigma(h),
\]
where $\sigma(h)$, $0\leq h\leq b-a$, is a continuous, monotonically increasing function, $\sigma(0)=0$ and for any $z>0$ the condition is satisfied
$$\int\limits_{0}^{z}\varkappa\left(\frac{b-a}{2\sigma^{(-1)}(u)}+1\right)du<\infty;$$
\item [6)] let $\varphi$ be a $f$-wavelet, and $\delta$ be the corresponding $m$-wavelet, and for $\varphi$ the condition $S$ with the function $\Phi$ holds;
\item [7)] for the function $c(x)$ there exists a function $0<A(u)<\infty$, $u>0$, such that for sufficiently large $x$, $a>0$, $c(ax)\leq c(x) A(a)$ and $\int\limits_R c(x) \Phi(\left|x\right|)<\infty$.
\end{enumerate}
Then, for any interval $\left[\alpha,\beta\right]\subset (a,b)$, $X_n(t)\rightarrow X(t)$ as $n\rightarrow \infty$ uniformly over $t \in \left[\alpha,\beta\right]$ with probability one, where
\[
X_n\left(t\right)=\sum\limits_{k \in \mathbb{Z}}\xi_{0k}\varphi_{0k}(x)+\sum\limits_{j=0}^{n-1}\sum\limits_{k \in \mathbb{Z}} \eta_{jk}\delta_{jk}(x),
\]
\[
\xi_{0k}=\int\limits_{\mathbb{R}} X(t) \overline{\varphi_{0k}(t)}dt, \ \eta_{jk}=\int\limits_{\mathbb{R}} X(t) \overline{\delta_{jk}(t)}dt.
\]
\end{theorem}

\begin{proof}
The proof of the theorem follows from Theorem {\rm\ref{Mlavetst4}}.
Since, by Theorem 5.5.19 \cite{mlav143}
there exists a function $c(t)$ and a random variable $\xi$ such that with probability one $\left|X(t)\right|<c(t)\xi$.
In addition, according to Theorem \ref{th:o2chapter_8} the random process $X(t)$ with probability one is uniformly continuous on the interval $[a,b]$.
\end{proof}

\begin{remark}
The condition \eqref{Mlavetse5} of Theorem {\rm\ref{Mlavetst5}} is satisfied for the functions $\psi(u)=u^\alpha$, $\psi(u)=\left(\ln (u+1)\right)^\lambda$.
If $\psi(u)=u^\alpha$, then $\varkappa(n^2)=2^\alpha \varkappa(n)$, i.e. $C_{\varkappa}=2^\alpha$, and if $\psi(u)=\left(\ln (u+1)\right)^\lambda$, then
$\varkappa(n^2)\leq 2^\lambda \varkappa(n)$, i.e. $C_{\varkappa}=2^\lambda$.
\end{remark}

\section*{Conclusions to chapter~\ref{ch:o3series}}

In section~\ref{ch:o3series}, estimates of the distribution of suprema of random processes from the spaces $\mathbf{F}^*_\psi(\Omega)$ are found.
The probabilities of large deviations for sums of independent random processes from the spaces $\mathbf{F}^*_\psi(\Omega)$ are considered and examples are given.
Estimates and conditions of sample continuity with probability one from these spaces are found.


\chapter{Stochastic processes from the space $SF_\psi(\Omega)$ of random variables}
\label{ch:o4series}

The space of strictly  $\mathbf{F}_\psi(\Omega)$ stochastic
processes is introduced in the same manner as strictly Orlicz
stochastic processes  defined in  \cite{Barrasa:1995}, as well as
strictly  $Sub_\varphi(\Omega)$ stochastic processes in
\cite{vasyl2008}.

\section{Strictly $F_\psi(\Omega)$ space of random variables. Basic properties}

\begin{definition}\label{def2}
A family $\Delta$ of random variables $\xi$ from the space
$\mathbf{F}_\psi(\Omega)$ such that $E\xi=0$ is called strictly
$\mathbf{F}_\psi(\Omega)$ (denoted as $SF_\psi(\Omega)$), if there exists a constant
$C_\Delta$ such that for at most countable set $I$ of random
variables $\xi_i \in \Delta$, $i \in I$ and for all $\lambda_i \in
R^1$ the following inequality holds true
\begin{equation*}
\left\|\sum_{i \in I}\lambda_i\xi_i\right\|_\psi\leq
C_\Delta\left(E\left(\sum_{i \in
I}\lambda_i\xi_i\right)^2\right)^{1/2}.
\end{equation*}
\end{definition}

\begin{remark}
The constant $C_\Delta$ is called determining constant of the family
$\Delta$.
\end{remark}

\begin{remark}\label{rem1}
It follows from Theorem \ref{the1} that there exists a constant
$B=\psi(2)$ such that the following inequality holds true
\[
\left(E\left(\sum_{i \in I}\lambda_i\xi_i\right)^2\right)^{1/2}\leq
B \left\|\sum_{i \in I}\lambda_i\xi_i\right\|_\psi.
\]
\end{remark}

\begin{theorem}\label{the6}
Let for the space  $\mathbf{F}_\psi(\Omega)$ the condition $\mathbf{H}$ hold
true. Let $\Delta$ be a family of independent random variables
$\xi_i$, $i \in I$, from $\mathbf{F}_\psi(\Omega)$  and there exists a
constant $R$, $R>0$, such that $\left\|\xi_i\right\|_\psi\leq
R\left(E\xi_i^2\right)^{1/2}$. Then $\Delta$ is the strictly
$\mathbf{F}_\psi(\Omega)$ family of random variables with the determining
constant $C_\Delta=\sqrt{C_\psi}R$, where $C_\psi$ is a scale
constant from Definition \ref{de:o1chapter_5}.
\end{theorem}

\begin{proof}
It follows from Definition  \ref{de:o1chapter_5} and Definition \ref{def2} that
\[
\left\|\sum_{i \in I}\lambda_i\xi_i\right\|_\psi^2\leq C_\psi
\sum_{i \in I}\left\|\lambda_i\xi_i\right\|_\psi^2\leq C_\psi
\sum_{i \in I} \lambda_i^2 R^2 E\xi_i^2= C_\psi R^2 E \left(\sum_{i
\in I}\lambda_i\xi_i\right)^2.
\]
\end{proof}

\begin{remark}\label{rem4}
It follows from Theorem \ref{the6} and Example \ref{ex1} that a
family of i.i.d random variables with normal distribution
$N(0,\sigma^2)$ is  strictly $\mathbf{F}_\psi(\Omega)$,
$\psi(u)=\sqrt{u}$, family of random variables with the determining
constant $C_\Delta=\sqrt{C_\psi}R=2\sqrt{2}e^{-5/6}$.

\end{remark}

\begin{remark}\label{rem5}
It follows from Theorem \ref{the6} and Example \ref{ex2} that a
family of i.i.d random variables with symmetric exponential
distribution  is  strictly $\mathbf{F}_\psi(\Omega)$, $\psi(u)=u$, family
of random variables with the determining constant
$C_\Delta=\sqrt{C_\psi}R=32\pi e^{-2}.$

\end{remark}

\begin{theorem}\label{the7}
Let $\Delta$ be a strictly $\mathbf{F}_\psi(\Omega)$ family of random
variables with the determining constant $C_\Delta$. The linear
closure of the family $\Delta$ in the space $\mathbf{F}_\psi(\Omega)$ is
strictly $\mathbf{F}_\psi(\Omega)$ family with the same determining constant.
\end{theorem}

\begin{proof}
Let $\xi_1,\xi_2,\ldots,\xi_n$ be random variables from $\Delta$.
Let $\eta_i=\sum_{j=1}^n a_{ij}\xi_j$, $i=1,2,\ldots,m$ be an
element of the linear closure of the set $\Delta$. Then
\[
\left\|\sum_{i=1}^m\lambda_i\eta_i\right\|_\psi=\left\|\sum_{i=1}^m\lambda_i\sum_{j=1}^n
a_{ij}\xi_j\right\|_\psi=\left\|\sum_{j=1}^n\left(\sum_{i=1}^m\lambda_i
a_{ij}\right)\xi_j\right\|_\psi\leq
\]
\[
C_\Delta\left(E\left(\sum_{j=1}^n\left(\sum_{i=1}^m\lambda_i
a_{ij}\right)\xi_j\right)^2\right)^{1/2}=C_\Delta\left(E\left(\sum_{i=1}^m\lambda_i
\eta_i\right)^2\right)^{1/2}.
\]
\end{proof}

For different elements of the linear closure of the set $\Delta$ we
verify the conclusion of the theorem by passing to the limit (see
Remark  \ref{rem1}).

\begin{definition}
A stochastic process $X=\left\{X(t), t \in T\right\}$ from the space
$\mathbf{F}_\psi(\Omega)$ is called strictly
$\mathbf{F}_\psi(\Omega)$ if the family of random variables
$\left\{X(t), t \in T\right\}$ is strictly
$\mathbf{F}_\psi(\Omega)$. The stochastic processes $X=\left\{X(t),
t \in T\right\}$ and $Y=\left\{Y(t), t \in T\right\}$ are called
jointly strictly $\mathbf{F}_\psi(\Omega)$ if the family of random
variables $\left\{X(t), Y(t), t \in T\right\}$ is strictly
$\mathbf{F}_\psi(\Omega)$.
\end{definition}

\begin{example}
Let $\xi_k, k=1,\ldots,\infty$ be a sequence of independent random
variables $\xi_k \in F_\psi(\Omega)$, such that
$\left\|\xi_k\right\|_\psi\leq R\left(E\xi_k^2\right)^{1/2}$ and let
$\mathbf{F}_\psi(\Omega)$ satisfy condition $H$. If a process $X(t)$ admits
the expansion
\begin{equation}\label{equ56}
X(t)=\sum_{k=1}^\infty \varphi_k(t)\xi_k, t \in T,
\end{equation}
and the series in the expansion \eqref{equ56} converges in the mean
square then the process $X(t)$ is strictly $\mathbf{F}_\psi(\Omega)$ with
determining constant $\sqrt{C_\psi}R$ (see Theorem \ref{the6}).
\end{example}

\begin{example}
Let stochastic processes $X(t)$ and $Y(t)$ admit the following expansions
\[
X(t)=\sum_{k=1}^\infty \varphi_k(t)\xi_k, Y(t)=\sum_{k=1}^\infty
\psi_k(t)\xi_k, t \in T,
\]
where $\left\{\xi_k, k=1,\ldots,\infty\right\}$ is strictly
$\mathbf{F}_\psi(\Omega)$ family of random variables and the series in the expansions above converge in
 mean square. It follows from Theorem \ref{the7} that $X(t)$ is
strictly $\mathbf{F}_\psi(\Omega)$ process, $Y(t)$ is strictly
$\mathbf{F}_\psi(\Omega)$ process and $X(t)$, $Y(t)$ are the jointly strictly
$\mathbf{F}_\psi(\Omega)$ processes.
\end{example}

\begin{theorem}\label{the8}
Let  $X=\left\{X(t), t \in T\right\}$ be a strictly $\mathbf{F}_\psi(\Omega)$
stochastic process and let $\varphi_k(t), k=1,\ldots,\infty$ be a
denumerable set of measurable functions on a measurable space
$\left\{T,B,\mu\right\}$. Let $\Xi=\xi_k, k=1,\ldots,\infty$, where
\[
\xi_k=\int_T\varphi_k(t)X(t)d\mu(t)
\]
and the integral are interpreted in the mean square sense. Then
$\Xi$ is strictly $\mathbf{F}_\psi(\Omega)$ family. The determining constant
of the family is equal to the determining constant of the process
$X(t)$.
\end{theorem}

\begin{proof}
Obviously, Theorem  \ref{the8} follows from Theorem \ref{the7}.
\end{proof}

\begin{corollary}
Note, that in the case when $X=\left\{X(t), t \in T\right\}$ and
$Y=\left\{Y(t), t \in T\right\}$ are  jointly $SF_\psi(\Omega)$
processes with the same determining constant the families $\xi_k,
k=1,\ldots,\infty$ and $\eta_k, k=1,\ldots,\infty$, where
\[
\xi_k=\int_T\varphi_k(t)X(t)d\mu(t),
\eta_k=\int_T\varphi_k(t)Y(t)d\mu(t),
\]
are $SF_\psi(\Omega)$. The determining constants of the families are
equal to the determining constant of the processes $X(t)$ and
$Y(t)$.
\end{corollary}

\section{Approximation of stochastic processes from the space  $SF_\psi(\Omega)$}\label{sec4}

\begin{definition}\label{def-approx}
Let  $\left\{T, \mathcal{F}, v\right\}$ be a measurable space,
$X=\left\{X(t), t \in T\right\}$ be a stochastic process. We say
that stochastic process $\tilde{X}=\left\{\tilde{X}(t), t \in
T\right\}$ approximates the process $X$  with respect to the norm of
some Banach space  $B(T)$ (usually the spaces $L_p(T)$, $C(T)$ are
considered, $\left\|\cdot\right\|_B$ stands for the norm in this
space) with given reliability  $1-\delta$, $\delta>0$ and accuracy
$\varepsilon>0$, if the following condition
\[
P\left\{\left\|\tilde{X}(t)-X(t)\right\|_B>\varepsilon\right\}\leq
\delta
\]
holds true.
\end{definition}

Let $X=\left\{X(t), t \in T\right\}$ be a stochastic process from
the space  $SF_\psi(\Omega)$ with the determining constant $C_\psi$.
We also suppose that  condition  $H$ with the constant  $C_B$ is
fulfilled.

\begin{theorem}\cite{Kozachenko:2011-1}\label{the9}
Let $X(t)$, $t \in T$, be a second-order stochastic process,
$EX(t)=0$, with correlation function $R(t,s)=EX(t)\overline{X(s)}$.
Let $\left(\Lambda, \mathcal{B}, \mu\right)$ be a measurable space
with $\sigma$-finite measure $\mu$, functions $f_i(t,\lambda)$, $t
\in T$, $i=1,\ldots,n$, belong to $L_2\left(\Lambda,\mu\right)$,
$g_k(\lambda), k \in Z$ be an orthonormal basis (ONB) in
$L_2\left(\Lambda,\mu\right)$.

Suppose that the correlation function $R(t,s)$ can be represented as
\[
R(t,s)=\sum_{i=1}^n\int_\Lambda f_i(t,\lambda)
\overline{f_i(s,\lambda)}d\mu(\lambda),
\]
then the process $X(t)$ can be represented as a convergent in mean
square series
\[
X(t)=\sum_{i=1}^n\sum_{k \in Z}a_{ik}(t)\xi_{ik},
\]
where
\[
a_{ik}(t)=\int_\Lambda f_i(t,\lambda)
\overline{g_k(\lambda)}d\mu(\lambda),
\]
\[
E\xi_{ik}=0,
E\xi_{ik}\overline{\xi_{jl}}=\delta_{ij}\cdot\delta_{kl},
\]
and $\delta_{ij}=\begin{cases}
                   0,& i\neq j;\\
                   1,& $i=j$.\\
                   \end{cases}$
\end{theorem}

\begin{remark}\cite{Kozachenko:2011-1}\label{rem3}
Let $H_X$ be a closure in mean square of the sums $\sum_{k=1}^m
c_kX\left(t_k\right)$, $t_k \in T$, $c_k \in C$ and $H_\xi$ be a
closure in mean square of the sums $\sum_{i=1}^m
\sum_{k=1}^lc_{ik}\xi_{ik}$, $c_{ik} \in C$.

If the system of functions $\left\{f_i(t,\lambda), i=1,\ldots,n, t
\in T\right\}$ is complete in $L_2\left(\Lambda,\mu\right)$ then
$H_X=H_\xi$.

Hence, if in Theorem 9 the process $X(t)$ is Gaussian then
$\xi_{ik}$, $i=\overline{1,n}$, $k \in Z$, are independent Gaussian
random variables.
\end{remark}

\begin{theorem} \cite{Kozachenko:2011-1}\label{the10}
Let $X=\left\{X(t), t \in T\right\}$ be a stochastic process which
satisfies conditions of Theorem \ref{the9} for $n=1$, i.e., its
covariance function $R(t,s)$ can be  represented as
\[
R(t,s)=EX(t)\overline{X(s)}=\int_\Lambda f(t,\lambda)
\overline{f(s,\lambda)}d\mu(\lambda),
\]
$f(t,\cdot) \in L_2\left(\Lambda,\mu\right)$ and $X(t)$ can be
expanded according to Theorem \ref{the9} as a convergent in mean square
series
\begin{equation}\label{equ71}
X(t)=\sum_{k \in Z}a_k(t)\xi_k,
\end{equation}
where
\[
a_k(t)=\int_\Lambda f(t,\lambda)
\overline{g_k(\lambda)}d\mu(\lambda),
\]
$\xi_k$ are uncorrelated random variables,
$\left\{g_k(\cdot)\right\}$, $k \in Z$, is ONB in
$L_2\left(\Lambda,\mu\right)$. Let $\left\{\Lambda, \mathcal{U},
\mu\right\}$ and $\left\{T, \mathcal{F}, v\right\}$ be measurable
spaces with $\sigma$-finite measures $\mu(\cdot)$ and $v(\cdot)$
correspondingly. Suppose that $f \in L_2 (T\times\Lambda)$, there
exist functions $b_l(t)$, $l \in Z$, such that
\[
\int_T \vert b_l(t)\vert^2 dv(t)<\infty
\]
and $b_l(t)$ is a solution of integral equation
\[
g_l(\lambda)=\int_T b_l(t) f(t,\lambda)dv(t).
\]
Then the following relationship holds for random variables in
\eqref{equ71}:
\[
\xi_l=\int_T X(t)b_l(t)dv(t).
\]
\end{theorem}

\begin{corollary}
Let  $X(t)$, $t \in R$, be a stationary process with spectral
density $h(\lambda)>0$, that is
\[
R(t,s)=R(t-s)=\int_R h(\lambda)\exp^{i\lambda(t-s)}d\lambda=\int_R
\sqrt{h(\lambda)}\exp^{i\lambda
t}\overline{\sqrt{h(\lambda)}\exp^{i\lambda s}}d\lambda,
\]
in this case  $f(t,\lambda)=\sqrt{h(\lambda)}\exp^{i\lambda t}$. Let
$-\infty<t<\infty$ and  $\nu(\cdot)$ be a Lebesgue measure,
$d\nu(t)=dt$. Then
\begin{equation}\label{equ81}
g_l(\lambda)=\int_{-\infty}^\infty
b_l(t)\sqrt{h(\lambda)}\exp^{i\lambda t}
dt=\sqrt{h(\lambda)}\int_{-\infty}^\infty b_l(t)\exp^{i\lambda t}
dt=\sqrt{h(\lambda)} \hat{b}_l(\lambda),
\end{equation}
where $\hat{b}_l(\lambda)$ is a Fourier transform of the function
$b_l(\lambda)$, it means
\[
\hat{b}_l(\lambda)=\int_{-\infty}^\infty b_l(t)\exp^{i\lambda t} dt.
\]
Then $\hat{b}_l(\lambda)=\frac{g_l(\lambda)}{\sqrt{h(\lambda)}}$. It
$\frac{g_l(\lambda)}{\sqrt{h(\lambda)}} \in L_1(R)$, then
$b_l(t)=\frac{1}{2\pi}\int_{-\infty}^\infty\frac{g_l(\lambda)}{\sqrt{h(\lambda)}}
e^{-i\lambda t}d\lambda$ is the inverse Fourier transform.
\end{corollary}

\begin{theorem}\label{the-11}
Let $X=\left\{X(t), t \in R\right\}$ be a real-valued centered
stationary process, $EX(t)=0$, from the space $SF_\psi(\Omega)$ such
that
\[
R(t,s)=R(t-s)=EX(t)X(s)=\int_{-\infty}^\infty h(\lambda) \cos
\lambda(t-s)d\lambda,
\]
where $ h(\lambda) \in L_1(R)$ is an even function, $ h(\lambda)>0$ and
$ h(\lambda)=0$ as $\vert \lambda\vert>A>0$ ( it could be also
that  $A=\infty$). Then the covariance function can be written as

\[
R(t,s)=\int_{-A}^A h(\lambda) \exp
\left\{i\lambda(t-s)\right\}d\lambda=\int_{-A}^A f(t,\lambda)
\overline{f(s,\lambda)}d\lambda,
\]
where $f(t,\lambda)= exp \left\{i\lambda
t\right\}\sqrt{h(\lambda)}$. The process $X(t)$ can be represented
as a convergent in mean square series
\[
X(t)=\sum_{l \in Z}a_l(t)\xi_l,
\]
where $a_l(t)=\int_{-A}^A \sqrt{h(\lambda)}\exp^{it\lambda}
\overline{g_l(\lambda)}d\mu(\lambda)$, $g_l(\lambda), l \in Z$ is an
ONB in $L_2\left(-A,A\right)$, $\xi_l$ are uncorrelated random
variables, $E\xi_l\xi_k=\delta_{l,k}$.

Let $T=\left[a,b\right]$, $-\infty<a<b<\infty$, $\xi_l \in
SF_\psi(\Omega)$,
\[
b_l(t)=\frac{1}{2\pi}\int_{-A}^A\frac{g_l(\lambda)}{\sqrt{h(\lambda)}}
e^{-i\lambda t}d\lambda.
\]
Then
\begin{equation}\label{equ57}
\xi_l=\int_a^b X(t)b_l(t)dt.
\end{equation}
\end{theorem}

\begin{proof}
The theorem follows from Theorem \ref{the10}, taking into account
 $\Lambda=\left[-A,A\right]$, $T=\left[a,b\right]$, and
$b_l(t)$ is a solution of equation \eqref{equ81}.
\end{proof}

\begin{remark}
Equality \eqref{equ57} gives a possibility to simulate random
variable $\xi_l$ from the space $SF_\psi(\Omega)$ under known
information on  sample path  $X(t)$.
\end{remark}

 The following theorem states the conditions for expansion of a
stochastic process from the space $SF_\psi(\Omega)$ in the form of
series giving a possibility to approximate the process by  the
model. Moreover, the rate of convergence in uniform metrics is also
investigated.

\begin{theorem}\label{the11}
Let $X=\left\{X(t), t \in T\right\}$ be a stochastic process from the
space $SF_\psi(\Omega)$, where $(T,\rho)$ is a compact metric
separable space, and suppose
\[ R(t,s)=EX(t)X(s)=\int_\Lambda
f(t,\lambda) f(s,\lambda)d\mu(\lambda),
\]
where $(\Lambda, \mathcal{A}, \mu)$ is some measurable space,
$f(t,\lambda) \in L_2\left(\Lambda,\mu\right)$ for any $t \in T$,
the system of the functions $f(t,\lambda)$, $t \in T$, is complete in
the space $L_2\left(\Lambda,\mu\right)$.

Assume that for the process  $X$ the conditions of Theorem
\ref{the10} are fulfilled, hence $X(t)=\sum_{k=1}^\infty \xi_k
a_k(t)$, where the series converges in mean squares, therefore, it is
also convergent in the norm of the space $SF_\psi(\Omega)$,
\[
a_k(t)=\int_\Lambda f(t,\lambda) \overline{g_k(\lambda)}
d\mu(\lambda),
\]
where $g_k(\lambda)$ is ONB in $L_2\left(\Lambda,\mu\right)$,
$\xi_k$ are random variables such that $E\xi_k=0$, $E\xi_k
\xi_l=\delta_{kl}$ and it follows from Remark \ref{rem3} that random
variables  $\xi_k$ belong to  the space $SF_\psi(\Omega)$ too.

Let $X_N(t)=\sum_{k=1}^N \xi_k a_k(t)$, $N\geq1$ and
$Y_N(t)=X(t)-X_N(t)$. If the condition
\begin{equation}\label{equ7}
\sup\limits_{\rho(t, s) \leq h}
\left(E\vert Y_N(t)-Y_N(s)\vert^2\right)^{1/2}\leq \sigma_N(h)
\end{equation}
is satisfied, where for any $N$, $\sigma_N(h)$, $h>0$, is a continuous
monotonically increasing function with initial state $\sigma_N(0)=0$,
and for some  $z>0$ (hence, for any value $z>0$)
\begin{equation}\label{equ8}
\int\limits_{0}^{z}\varkappa\left(\mathbf{N}\left(\sigma_N^{(-1)}(\frac{u}{C_\Delta})\right)\right)du<\infty,
\end{equation}
where $\sigma_N^{(-1)}(u)$ is an inverse function of $\sigma_N(u)$,
$\varkappa(s)$ is majorizing characteristics of the space
$SF_\psi(\Omega)$,
then with probability one $$\sup\limits_{t \in T}Y_N(t) \in
SF_\psi(\Omega)$$
 and  $$\left\|\sup\limits_{t \in
T}\vert Y_N(t)\vert\right\|\leq B_N(p),\quad 0<p<1,$$
 where
\begin{equation}\label{eq-Bn}
B_N(p)=\inf \limits_{t \in T}C_\Delta\left(E\vert
Y_N(t)\vert^2\right)^{1/2}+\frac{1}{p(1-p)}\int\limits_{0}^{\gamma
p}\varkappa\left(\mathbf{N}\left(\sigma_N^{(-1)}(\frac{u}{C_\Delta})\right)\right)du,
\end{equation}
$\gamma=C_\Delta\sigma_N\left(\sup\limits_{t,s \in T}\rho(t,
s)\right)$. Moreover, for  $ \varepsilon>0$ the following inequality
holds:
\begin{equation}
P\left\{\sup\limits_{t \in
T}\vert Y_N(t)\vert>\varepsilon\right\}\leq
\inf\limits_{u\geq1}\frac{B_N^u(p)(\psi(u))^u}{\varepsilon^u}.
\end{equation}
\end{theorem}

\begin{proof}
We use Theorem \ref{th:o2chapter_1} and Corollary \ref{co:o2chapter-1} to prove the
statement of given theorem. Really, inequality \eqref{equ7} implies that
 $$\sup\limits_{\rho(t, s) \leq h}
\left\|Y_N(t)-Y_N(s)\right\|\leq C_\Delta\sigma(h),$$ so
$s(h)=C_\Delta\sigma(h)$, and the inverse function then equals
$s^{(-1)}(h)=\sigma^{(-1)}\left(\frac{h}{C_\Delta}\right)$. Hence,
condition  \eqref{equ1} follows from  \eqref{equ8}.  Inequality
\eqref{equ2} is also carried out, since in our case  $B(p)$ equals
 $B_N(p)$.
\end{proof}

\begin{corollary}\label{cor-diff}
Let in Theorem \ref{the11} we study such $a_k(t)$, that
\[
\vert a_k(t)-a_k(s)\vert<c_k z(\rho(t-s)),
\]
where $z(u)>0$, $u>0$, is a monotonically increasing function,
$z(0)=0$. Then
\[
E\vert Y_N(t)-Y_N(s)\vert^2\leq \sum_{k=N+1}^\infty c^2_k
z^2(\rho(t-s)).
\]
If the series  $\sum_{k=N+1}^\infty c^2_k$ is convergent, then the
conditions of Theorem \ref{the11} are fulfilled with
$\sigma_N(u)=\left(\sum_{k=N+1}^\infty c^2_k z^2(u)\right)^{1/2}$.
\end{corollary}
Let clarify the statement of Theorem \ref{the11} in the case when
 parametric set $T$ is a real-valued domain $[a,b]$,  $\psi$ is a power  function $\psi(u)=u^\alpha,\, u\geq 1, \alpha>0$
and the standard deviation of increments are majorized by
\begin{equation}\label{eq9}
\sup\limits_{\vert t- s\vert \leq h}
\left(E\vert Y_N(t)-Y_N(s)\vert^2\right)^{1/2}\leq C_N h^\beta,
\quad, \beta>\alpha>0.
\end{equation}
\begin{remark}\label{rem-Cn}
It's obvious that $C_N\to 0$ as $N\to \infty$.
\end{remark}

\begin{theorem}\label{the-simulation}
Let $X=\left\{X(t), t \in [a,b] \right\}$ be a centered stochastic
process from the space $SF_\psi(\Omega)$, $\psi(u)=u^\alpha, \,
u\geq 1, \alpha>0$, with covariance function
\[ R(t,s)=EX(t)X(s)=\int_\Lambda
f(t,\lambda) f(s,\lambda)d\mu(\lambda),
\]
 and suppose the
conditions of Theorem \ref{the11} are fulfilled.

Hence
\begin{equation}\label{eq-process}
X(t)=\sum_{k=1}^\infty \xi_k a_k(t),
\end{equation}
 where the series
converges in mean squares, therefore, it is also convergent in the
norm of the space $SF_\psi(\Omega)$,
\begin{equation}\label{eq-ak}
a_k(t)=\int_\Lambda f(t,\lambda) g_k(\lambda)d\mu(\lambda),
\end{equation}
 $g_k(\lambda)$ is ONB in $L_2\left(\Lambda,\mu\right)$,
$\xi_k$ are random variables such that $E\xi_k=0$, $E\xi_k
\xi_l=\delta_{kl}$ and it follows from Remark \ref{rem3} that random
variables  $\xi_k$ belong to  the space $SF_\psi(\Omega)$ too.

We consider an approximated model of the process $X$  and the
accuracy process as
 $$X_N(t)=\sum_{k=1}^N
\xi_k a_k(t),\quad Y_N(t)=X(t)-X_N(t).$$

If the condition \eqref{eq9} is satisfied, then with probability one
$$\sup\limits_{t \in T}Y_N(t) \in SF_\psi(\Omega)$$
 and  $$\left\|\sup\limits_{t \in
T}\vert Y_N(t)\vert\right\|\leq B_N,$$
 where
\begin{equation}\label{eq-Bn-the}
B_N= C_\Delta \inf_{t\in [a,b]}\left(\sum_{k=N+1}^\infty
a_k^2(t)\right)^\frac12+ \frac{ C_\Delta C_N\left(b-a\right)^\beta
(\alpha+\beta)^{1+\alpha/\beta}}{\beta-\alpha} \left(\frac{e
 }{2\alpha^{1+1/\beta}} \right)^\alpha,
\end{equation}
$\beta>\alpha>0$. Moreover, for  $ \varepsilon>0$ the following
inequality on the rate of convergence of approximated model in the
space of continuous functions holds true:
\begin{equation}\label{eq-ineq}
P\left\{\sup\limits_{t \in
T}\vert Y_N(t)\vert>\varepsilon\right\}\leq
\inf\limits_{u\geq1}\frac{B_N^u u^{\alpha u}}{\varepsilon^u}.
\end{equation}
\end{theorem}

\begin{remark}
Since the value $C_N$ tends to 0 as $N\to \infty$  (see Remark
\ref{rem-Cn}) and the series $\sum_{k=1}^\infty a_k^2(t)$ is
absolutely convergent then $B_N$ as a function of $N$ tends to 0 as
$N\to \infty$.
\end{remark}

\begin{proof}
The proof employs the results of Theorem \ref{the11}. First of all,
find expectation and variance for the accuracy process. Since
$E\xi_k=0$, $E\xi_k \xi_l=\delta_{kl}$ then $E Y_N(t)=0$ and

$$
Var(X(t)-X_N(t))=Var Y_N(t)=E Y_N^2(t)=$$
\begin{equation}\label{eq-var}
=E\left(\sum_{k=N+1}^\infty a_k(t)\xi_k \right)^2=
\sum_{k=N+1}^\infty a_k^2(t).
\end{equation}

It  follows from Theorem \ref{th:o1chapter_7} that in the case of
$\psi(u)=u^\alpha,\, u\geq 1, \alpha>0$ the majorizing
characteristics is
\begin{equation}\label{eq-kappa2}
 \varkappa(n)=
e^\alpha \left(\frac{\ln n}{\alpha}\right)^\alpha, n\geq 2.
\end{equation}
In monograph \cite{byld98}, it is shown that the metric
masiveness $\mathbf{N}_T(u)$ on the interval $T=[a,b]$ is estimated
in the following way
\begin{equation}\label{eq-mass}
\mathbf{N}_{[a,b]}(u)\leq \frac{b-a}{2u}+1.
\end{equation}
In the notation of Theorem \ref{the11} from \eqref{eq9} we have that
the function $\sigma_N(t)=C_N\cdot t^\beta$, then the inverse
function equals

\begin{equation}\label{eq-inv sigma}
\sigma^{(-1)}(t)=\left(\frac{t}{C_N}\right)^{1/\beta}.
\end{equation}

Substituting \eqref{eq-kappa2}, \eqref{eq-mass} and \eqref{eq-inv
sigma} in the integral from  \eqref{eq-Bn} and applying trivial
inequality $\ln(1+x)\leq x,\, x>0,$ we obtain
\begin{eqnarray}\label{eq-int}
&\int\limits_{0}^{\gamma
p}\varkappa\left(\mathbf{N}\left(\sigma_N^{(-1)}(\frac{u}{C_\Delta})\right)\right)du\leq\nonumber\\
&\leq
\int\limits_{0}^{\gamma
p}\frac{e^\alpha}{\alpha^\alpha}\ln^\alpha\left(\frac{b-a}{2\sigma_N^{(-1)}(\frac{u}{C_\Delta})}+1\right)du\leq\nonumber\\
&\leq
\int\limits_{0}^{\gamma
p}\frac{e^\alpha}{\alpha^\alpha}\left(\frac{(b-a)(C_N\cdot C_\Delta)^{1/\beta} }{2u^{1/\beta}}\right)^\alpha du=\nonumber\\
&=
\left(\frac{e (b-a) }{2\alpha} \right)^\alpha (C_N\cdot
C_\Delta)^{\alpha /\beta} \int\limits_{0}^{\gamma p}
u^{-\alpha/\beta} du=\nonumber\\
&= \frac{\beta}{\beta-\alpha} \left(\frac{e (b-a) }{2\alpha}
\right)^\alpha (C_N\cdot C_\Delta)^{\alpha /\beta} (\gamma
p)^\frac{\beta-\alpha}{\beta}.
\end{eqnarray}

Taking into account \eqref{eq-var} and \eqref{eq-int}, the value
$B_N(p)$ from \eqref{eq-Bn} can be estimated as
$$
B_N(p)\leq C_\Delta \inf_{t\in [a,b]}\left(\sum_{k=N+1}^\infty
a_k^2(t)\right)^\frac12+ \frac{\beta \gamma}{\beta-\alpha}
\left(\frac{e (b-a) }{2\alpha}  \left(\frac{C_N\cdot
C_\Delta}{\gamma}\right)^{1 /\beta} \right)^\alpha \frac{
p^{-\frac{\alpha}{\beta}}}{1-p}.
$$
The value $\gamma=C_\Delta\sigma_N\left(\sup\limits_{t,s \in
T}\rho(t, s)\right)$ under the condition of Theorem
\ref{the-simulation} equals
$$
\gamma=C_\Delta C_N\left(b-a\right)^\beta.
$$
Then

$$
B_N(p)\leq C_\Delta \inf_{t\in [a,b]}\left(\sum_{k=N+1}^\infty
a_k^2(t)\right)^\frac12+ \frac{\beta C_\Delta
C_N\left(b-a\right)^\beta}{\beta-\alpha} \left(\frac{e  }{2\alpha}
\right)^\alpha \frac{ p^{-\frac{\alpha}{\beta}}}{1-p}.
$$

 Let us minimize $B_N(p), p\in (0,1)$ over $p$. It is easy to check
that the argument $p$ that minimizes this expression equals
$$
p_{\min}= \text{argmin}_{p\in (0,1)} B_N(p)=\frac{\alpha}{\alpha+\beta}.
$$
Therefore,
\begin{eqnarray*}
B_N=\min_{p\in (0,1)} B_N(p)&= &C_\Delta \inf_{t\in
[a,b]}\left(\sum_{k=N+1}^\infty a_k^2(t)\right)^\frac12+\\
&+& \frac{ C_\Delta C_N\left(b-a\right)^\beta
(\alpha+\beta)^{1+\alpha/\beta}}{\beta-\alpha} \left(\frac{e
 }{2\alpha^{1+1/\beta}} \right)^\alpha,
\end{eqnarray*}
that completely coincides with \eqref{eq-Bn-the}. So, the theorem is
 proved.
\end{proof}

\begin{corollary}\label{coroll-1}
Let $X=\left\{X(t), t \in [a,b] \right\}$ be a centered stochastic
process from the space $SF_\psi(\Omega)$, $\psi(u)=u^\alpha, \,
u\geq 1, \alpha>0$, with covariance function
\[ R(t,s)=EX(t)X(s)=\int_\Lambda
f(t,\lambda) f(s,\lambda)d\mu(\lambda),
\]
and suppose that all conditions of Theorem \ref{the-simulation} are met. Then
the process
 $$X_N(t)=\sum_{k=1}^N
\xi_k a_k(t),$$
where $a_k(t),\, k\geq1,\,$ are from \eqref{eq-ak}, approximates
stochastic process $X(t)$ with given accuracy $\varepsilon, \,
\varepsilon>0,$ and reliability $1-\delta, \, \delta\in (0,1),$ in
uniform norm if the cutting-off level $N$ satisfies inequality:
\begin{equation}\label{eq-approx-Bn}
B_N\leq \frac{\varepsilon
\delta^{-\frac1{\ln\delta}}}{\ln^\alpha\frac1\delta},
\end{equation}
where $B_N$ is defined in \eqref{eq-Bn-the}.
\end{corollary}

\begin{proof}
By Definition \ref{def-approx} the process (model) $X_N(t), \, t\in
T$, approximates a stochastic process $X(t), t\in T,$  in uniform
norm  with given reliability  $1-\delta$, $\delta>0$, and accuracy
$\varepsilon>0$, if
\[
P\left\{\sup_{t\in
T}\vert X_N(t)-X(t)\vert>\varepsilon\right\}=P\left\{\sup_{t\in
T}\vert Y_N(t)\vert>\varepsilon\right\}\leq \delta.
\]
 By \eqref{eq-ineq} the left-hand side of the inequality above is
bounded by $\frac{B_N^u u^{\alpha u}}{\varepsilon^u},\, $ $u\geq 1.$
So, to approximate the process $X(t)$ by $X_N(t)$ it is enough to find such $N$ that
$$
\frac{B_N^u u^{\alpha u}}{\varepsilon^u}<\delta \quad
\Leftrightarrow \quad B_N  \leq  \frac{\varepsilon
\delta^{\frac1{u}}}{u^{\alpha }}.
$$
Taking the maximum of $ \frac{\varepsilon
\delta^{\frac1{u}}}{u^{\alpha }}$ over $u$ we get the point of
maximum $ u_{\max}=-\ln \delta=\ln\frac1\delta$ and obtain
\eqref{eq-approx-Bn}.
\end{proof}

\section{Case study}\label{sec5}
Consider a centered real-valued stationary stochastic process
$X=\left\{X(t), t \in [a,b] \right\}$  from the space
$SF_\psi(\Omega)$ with covariance function

\begin{equation}\label{eq-case-cov}
R(\tau)=\int_{-A}^A h(\lambda) \exp
\left\{i\lambda\tau\right\}d\lambda,\quad\quad
h(\lambda)=e^{-\vert \lambda\vert},  \vert \lambda\vert\leq A
\end{equation}

and consider $g_k(\lambda)=\frac{1}{A}e^{\frac{ik\lambda\pi}{A}}$
that are ONB on  $\left[-A,A\right]$.

Then by Theorem \ref{the11} $$X(t)=\sum_{k=1}^\infty \xi_k
a_k(t),\quad Y_N(t)=\sum_{N+1}^\infty \xi_k a_k(t),\quad t\in
[a,b],$$
 where in our case
\begin{eqnarray}\label{eq-case-ak}
a_k(t)& =& \int_\Lambda f(t,\lambda)
\overline{g_k(\lambda)}d\mu(\lambda)=\frac1{A}\int_{-A}^A e^{-
\frac{\vert \lambda \vert }2} e^{i\lambda t} \overline{e^{\frac{ik\lambda
\pi}{A}}} d\lambda\nonumber\\
&=&\frac1{A}\int_{-A}^A e^{- \frac{\vert \lambda \vert }2} \cos(\lambda
t-\frac{k\lambda \pi}{A}) d\lambda.
\end{eqnarray}

The following Theorem provides the method of construction and
generation of random variable $\xi_k$ from the space
$SF_\psi(\Omega)$ if the information on trajectory of stochastic
process $X(t)$ is known.

\begin{theorem} Assume that a centered real-valued stationary stochastic process $X=\left\{X(t), t
\in [a,b] \right\}$  from the space $SF_\psi(\Omega)$ has the
covariance function \eqref{eq-case-cov} and all conditions of
Theorem \ref{the-11} are satisfied. Then
\begin{equation}\label{eq-case-57}
\xi_k=\int_a^b X(t)b_k(t)dt,
\end{equation}
where
\[
b_k(t)=\frac{2e^{\frac{A}{2}}(\cos( At-k\pi)+2(t-\frac{k\pi}{A})\sin
( At-k\pi))-2}{\pi A(1+4(t-\frac{k\pi}{A})^2)}.
\]

\end{theorem}
\begin{proof}
It follows from Theorem \ref{the-11} that
\[
b_k(t)=\frac{1}{2\pi}\int_{-A}^A\frac{g_k(\lambda)}{\sqrt{h(\lambda)}}
e^{-i\lambda t}d\lambda.
\]

Then substituting the functions $g_k(\lambda)$ and $h(\lambda)$ we
obtain
\[
b_k(t)=\frac{1}{2\pi A} \int_{-A}^A e^{\frac{ik\lambda\pi}{A}}
e^{\frac{\vert \lambda\vert}{2}} e^{-i\lambda t}d\lambda=
\frac{1}{2\pi A} \int_{-A}^A
e^{-i\lambda\left(t-\frac{k\pi}{A}\right)+\frac{\vert \lambda\vert}{2}}d\lambda=
\]
\[
\frac{1}{2\pi A} \int_{-A}^A \cos
\lambda\left(t-\frac{k\pi}{A}\right)
e^{\frac{\vert \lambda\vert}{2}}d\lambda= \frac{1}{\pi A}
\int_{0}^A \cos \lambda\left(t-\frac{k\pi}{A}\right)
e^{\frac{\lambda}{2}}d\lambda
\]

Let's take the integration above two times by parts.  Then finally
\[
b_k(t)=\frac{2e^{\frac{A}{2}}(\cos( At-k\pi)+2(t-\frac{k\pi}{A})\sin
( At-k\pi))-2}{\pi A(1+4(t-\frac{k\pi}{A})^2)},
\]
which completely proves Theorem.

\end{proof}
\begin{lemma}\label{lem-1}
The following inequalities hold true:
\begin{equation}\label{eq-case-ak-est}
\vert a_k(t)\vert\leq
\frac{\sqrt{\tilde{c}_1}}{k}+\frac{\sqrt{\tilde{c}_3}}{k^2},\quad
k>\frac{bA}{\pi}.
\end{equation}

\begin{equation}\label{eq-case-VarY-est}
Var(Y_N(t))= \sum_{k=N+1}^\infty a_k^2(t)\leq\frac{\tilde{c}_1}{N}+
\frac{\tilde{c}_2}{2N^2}+ \frac{\tilde{c}_3}{3N^3},
\end{equation}
where
\begin{equation}\label{eq-case-c}
\tilde{c}_1=\frac{2^8}{e^{A}\pi^2},\quad \tilde{c}_2= \frac{2^7A(1+
e^{-A/2})}{e^{A/2}\pi^3},\quad \tilde{c}_3= \frac{2^4A^2(1+
e^{-A/2})^2}{\pi^4}
\end{equation}

\end{lemma}
\begin{proof}
Taking the integration by parts two times  in \eqref{eq-case-ak} and
deriving integral itself we have
\begin{eqnarray}\label{eq-case-ak2}
a_k(t)& =& \frac2{A}\int_{0}^A e^{-\frac{\lambda}2} \cos(\lambda
t-\frac{k\lambda \pi}{A}) d\lambda\nonumber \\
&=&\frac2{A e^{A/2}} \frac{4 (t-\frac{k\pi}{A})\sin(At-k\pi)
-2\cos(At-k\pi) +2e^{A/2}}{1+ 4(t-\frac{k\pi}{A})^2}.
\end{eqnarray}
Since $t\in [a,b]$ and  $k>\frac{bA}{\pi}$ then $a_k(t)$ can be
estimated as follows:
\begin{eqnarray*}
\vert a_k(t)\vert &  \leq & \frac2{A e^{A/2}} \frac{4 \vert t-\frac{k\pi}{A}\vert +2
+2e^{A/2}}{1+ 4(t-\frac{k\pi}{A})^2}\\
&\leq &\frac4{A e^{A/2}}\frac{4 \frac{k\pi}{A} +1
+e^{A/2}}{(\frac{k\pi}{A})^2}\\
&=& \frac{16}{e^{A/2}\pi}\cdot \frac1{k}+ \frac{4A(1+
e^{-A/2})}{\pi^2}\cdot
\frac1{k^2}=\frac{\sqrt{\tilde{c}_1}}{k}+\frac{\sqrt{\tilde{c}_3}}{k^2}.
\end{eqnarray*}
The variance of accuracy process $Y_N(t)$ by \eqref{eq-var} and
\eqref{eq-case-ak-est} is bounded as
$$
Var(Y_N(t))= \sum_{k=N+1}^\infty a_k^2(t)\leq
\tilde{c}_1\sum_{k=N+1}^\infty\frac1{k^2}+
\tilde{c}_2\sum_{k=N+1}^\infty\frac1{k^3}+\tilde{c}_3\sum_{k=N+1}^\infty\frac1{k^4}.
$$
Since
\begin{eqnarray}\label{eq-case-sum}
\sum_{k=N+1}^\infty \frac1{k^s}&=& \sum_{k=N+1}^\infty\int_{k-1}^k
\frac1{k^s} dx \leq \sum_{k=N+1}^\infty\int_{k-1}^k
\frac1{x^s} dx \nonumber  \\
& = &\int_N^\infty \frac1{x^s} dx=\frac{1}{(s-1)N^{s-1}},\quad
\text{as} \quad s>1,
\end{eqnarray}
then
$$
Var(Y_N(t))\leq \frac{\tilde{c}_1}{N}+  \frac{\tilde{c}_2}{2N^2}+
\frac{\tilde{c}_3}{3N^3}.
$$
\end{proof}
We show now that in the case study framework Corollary \ref{cor-diff}
holds true.

\begin{lemma}\label{lem-2}

The functions $a_k(t)$ satisfy the H\"{o}lder condition:

\begin{equation}\label{eq-case-ak-est_sec5}
\vert a_k(t)- a_k(s) \vert\leq\left(
\frac{\hat{c}_1}{k}+\frac{\hat{c}_2}{k^2}+
\frac{\hat{c}_4}{k^4}+\frac{\hat{c}_3}{k^3}\right)\cdot
z(\vert t-s\vert),\quad k>\frac{bA}{\pi},
\end{equation}
where

\begin{equation}\label{eq-case-z}
z(h)=\left\{
       \begin{array}{ll}
         h^\beta, & h\in [0,1],\quad \beta\in (0,1] ; \\
         h, & h>1.
       \end{array}
     \right. ,
\end{equation}

\begin{eqnarray}\label{eq-case-c-hat}
\hat{c}_1&=& \frac{2^{2-\beta}e^{-A/2}A^\beta}{\pi},\nonumber\\
\hat{c}_2&=& \frac{2 e^{-A/2}A ( 2b(A/2)^\beta+2(A/2)^\beta +3)}{\pi^2}, \nonumber\\
\hat{c}_3&=& \frac{2e^{-A/2}A^{2}(1+4b+e^{A/2})}{\pi^3},\\
\hat{c}_4&=& \frac{2 e^{-A/2}A^{3}b(1+2b+e^{A/2})}{\pi^4}.\nonumber
\end{eqnarray}

\end{lemma}

\begin{proof}
Really, by \eqref{eq-case-ak2}
\begin{equation}\label{eq-case-a-begin}
\vert a_k(t)-a_k(s)\vert\leq I_1+I_2+I_3,
\end{equation}
where
\begin{eqnarray}
I_1& = &\frac8{A e^{A/2}}\vert  \frac{
(t-\frac{k\pi}{A})\sin(At-k\pi) }{1+ 4(t-\frac{k\pi}{A})^2}-
\frac{(s-\frac{k\pi}{A})\sin(As-k\pi)
}{1+ 4(s-\frac{k\pi}{A})^2} \vert, \label{I-1}\\
I_2& = & \frac4{A e^{A/2}} \vert \frac{\cos(At-k\pi) }{1+ 4(t-\frac{k\pi}{A})^2} -\frac{\cos(As-k\pi) }{1+ 4(s-\frac{k\pi}{A})^2} \vert, \label{I-2}\\
I_3& = & \frac4{A } \vert \frac{1}{1+ 4(t-\frac{k\pi}{A})^2}-
\frac{1}{1+ 4(s-\frac{k\pi}{A})^2}\vert. \label{I-3}
\end{eqnarray}

Estimate now each summand $I_j, \, j=1,2,3,$ separately. We will use
the following transformation approach:
\begin{equation}\label{eq-case-1}
\vert \frac{f_i(t)}{g(t)}-\frac{f_i(s)}{g(s)} \vert\leq
\frac{\vert f_i(t)\vert}{\vert g(t)g(s)\vert}\vert g(t)-g(s)\vert+\frac{\vert f_i(t)-f_i(s)\vert}{g(s)}
\end{equation}
for some functions $f_i(t), i=1,2,3,\, g(s), \, t,s\in [a,b]$.

For $I_1$ the functions in \eqref{eq-case-1} are $$f_1(t)=\frac8{A
e^{A/2}} (t-\frac{k\pi}{A})\sin(At-k\pi)\quad \text{and}\quad g(t)=
1 + 4(t-\frac{k\pi}{A})^2. $$
Then the difference $\vert g(t)-g(s)\vert$ can be estimated as
\begin{equation}\label{eq-case-2}
\vert g(t)-g(s)\vert=8 \vert t-s\vert \cdot\vert t+s-\frac{2 k\pi}{A}\vert \leq 4(b+\frac{
k\pi}{A})\vert t-s\vert .
\end{equation}

 To evaluate the increment $\vert f_1(t)-f_1(s)\vert $, we use
such inequality: $\vert \sin(h)\vert \leq \vert h\vert ^\beta, \, \beta\in(0,1].$ Then

\begin{eqnarray}\label{eq-case-3}
\vert f_1(t)-f_1(s)\vert &=& \frac8{A e^{A/2}}
\vert (t-\frac{k\pi}{A})\sin(At-k\pi)-
(s-\frac{k\pi}{A})\sin(As-k\pi)\vert\nonumber\\
&\leq &   \frac8{A e^{A/2}} \left(\vert t-\frac{k\pi}{A}\vert \vert \sin(At-k\pi)-
\sin(As-k\pi)\vert  \right.\nonumber \\
&+&\left. \,\,\,\,\,\,\,\quad\quad
\vert \sin(As-k\pi)\vert \vert (t-\frac{k\pi}{A})-(s-\frac{k\pi}{A})\vert
\right)\nonumber\\
&\leq & \frac8{A e^{A/2}}
\left((b+\frac{k\pi}{A})2\vert \sin\frac{At-As}{2}\vert +\vert t-s\vert \right)\nonumber\\
&\leq & \frac8{A e^{A/2}}
\left(2^{1-\beta}A^\beta(b+\frac{k\pi}{A})\vert t-s\vert ^\beta+\vert t-s\vert \right)\nonumber\\
&\leq & \frac{8z(\vert t-s\vert )}{A e^{A/2}}
\left(2^{1-\beta}A^\beta(b+\frac{k\pi}{A})+1\right).
\end{eqnarray}

It follows from \eqref{I-1}, \eqref{eq-case-1}-\eqref{eq-case-3}
that
\begin{eqnarray}\label{eq-case-I1}
I_1 & \leq&  \frac8{A e^{A/2}}\left( \frac{8(b+\frac{
k\pi}{A})^2\vert t-s\vert }{16 \left(k\pi/A\right)^4 }+
\frac{(2^{1-\beta}A^\beta(b+\frac{k\pi}{A})+1) z(\vert t-s\vert )}{4
\left(k\pi/A\right)^2}\right)\nonumber\\
&\leq &  \left(\frac{4A^3b^2}{e^{A/2}\pi^4k^4}+
\frac{8A^2b}{e^{A/2}\pi^3k^3}+ \frac{2A(3+A^\beta
2^{1-\beta})}{e^{A/2}\pi^2k^2}+\frac{2^{2-\beta}A^\beta}{e^{A/2}\pi
k}\right)\cdot z(\vert t-s\vert ).
\end{eqnarray}

Let us estimate now $I_2$ from \eqref{I-2}. For $I_2$ the functions
in \eqref{eq-case-1} are $$f_2(t)=\frac4{A e^{A/2}}
\cos(At-k\pi)\quad \text{and}\quad g(t)= 1 + 4(t-\frac{k\pi}{A})^2.
$$

 Similarly to $I_1$, we  use
 inequality $\vert \sin(h)\vert \leq \vert h\vert ^\beta, \, \beta\in(0,1],$  to estimate the increment $\vert f_2(t)-f_2(s)\vert $. Then

\begin{eqnarray}\label{eq-case-4}
\vert f_2(t)-f_2(s)\vert &=& \frac4{A e^{A/2}} \vert \cos(At-k\pi)-
\cos(As-k\pi)\vert\nonumber\\
&= &   \frac4{A e^{A/2}}
\left(2\vert \sin\frac{At-As}{2}\vert \vert \sin (At+As-2k\pi)\vert \right)\nonumber\\
&\leq & \frac4{A e^{A/2}}
\left(2^{1-\beta}A^\beta\vert t-s\vert ^\beta\right).
\end{eqnarray}

It follows from \eqref{I-2}, \eqref{eq-case-1},\eqref{eq-case-2} and
\eqref{eq-case-4} that
\begin{eqnarray}\label{eq-case-I2}
I_2 & \leq&  \frac4{A e^{A/2}}\left( \frac{8(b+\frac{
k\pi}{A})\vert t-s\vert }{16 \left(k\pi/A\right)^4 }+
\frac{(2^{1-\beta}A^\beta) \vert t-s\vert ^\beta}{4
\left(k\pi/A\right)^2}\right)\nonumber\\
&\leq &  \left(\frac{2A^3b}{ e^{A/2}\pi^4k^4}+ \frac{2A^2}{
e^{A/2}\pi^3k^3}+ \frac{A^{1+\beta} 2^{1-\beta}}{
e^{A/2}\pi^2k^2}\right)\cdot z(\vert t-s\vert ).
\end{eqnarray}

And for the third summand $I_3$ applying \eqref{eq-case-2} we have
\begin{eqnarray}\label{eq-case-I3}
I_3& = & \frac4{A } \vert \frac{1}{1+ 4(t-\frac{k\pi}{A})^2}-
\frac{1}{1+ 4(s-\frac{k\pi}{A})^2}\vert\nonumber\\
&\leq &  \frac4{A }   \frac{8(b+\frac{k\pi}{A})\vert t-s\vert }{16
(k\pi/A)^4}=\left(\frac{2A^3 b}{\pi^4 k^4}+\frac{2A^2}{\pi^3
k^3}\right) \vert t-s\vert .
\end{eqnarray}
Substituting \eqref{eq-case-I1}, \eqref{eq-case-I2},
\eqref{eq-case-I3} in \eqref{eq-case-a-begin} Lemma will be
completed.
\end{proof}

\begin{theorem}\label{the-case}
 Let $X=\left\{X(t), t \in [a,b] \right\}$ be a centered stochastic
process from the space $SF_\psi(\Omega)$, $\psi(u)=u^\alpha, \,
u\geq 1, \alpha>0$, with covariance function \eqref{eq-case-cov}.
Then the process
 $$X_N(t)=\sum_{k=1}^N
\xi_k a_k(t),$$
where $a_k(t),\, k\geq1,\,$ are from  \eqref{eq-case-ak2},
approximates stochastic process $X(t)$ with given accuracy
$\varepsilon, \, \varepsilon>0,$ and reliability $1-\delta, \,
\delta\in (0,1),$ in uniform norm if the cutting-off level $N$
satisfies inequality:
\begin{equation}\label{eq-approx-Bn_sec5}
\hat{B}_N\leq \frac{\varepsilon
\delta^{-\frac1{\ln\delta}}}{\ln^\alpha\frac1\delta},
\end{equation}
where

\begin{equation}\label{eq-case-Bn}
\hat{B}_N= C_\Delta \left(\frac{\tilde{c}_1}{N}+
\frac{\tilde{c}_2}{2N^2}+ \frac{\tilde{c}_3}{3N^3}\right)^\frac12+
\hat{C}_N \cdot K,
\end{equation}

$$
K=\frac{ C_\Delta \left(b-a\right) (\alpha+1)^{1+\alpha}}{1-\alpha}
\left(\frac{e
 }{2\alpha^{2}} \right)^\alpha,
$$
\begin{equation}\label{eq-case-c_n} \hat{C}_N=\left(
\frac{\hat{c}_1^2}{N}+\frac{\hat{c}_1\hat{c}_2}{N^2}+
\frac{\hat{c}_1^2+2\hat{c}_1\hat{c}_3}{3N^3}+\frac{\hat{c}_1\hat{c}_4+\hat{c}_2\hat{c}_3}{2N^4}+
\frac{\hat{c}_3^2+2\hat{c}_2\hat{c}_4}{5N^5}+
\frac{\hat{c}_3\hat{c}_4}{3N^6}+\frac{\hat{c}_4^2}{7N^7}
 \right)^\frac12
\end{equation}
and the fixed values $\tilde{c}_i,\, i=\overline{1,3}$,
$\hat{c}_j,\, j=\overline{1,4},$ are from \eqref{eq-case-c} and
\eqref{eq-case-c-hat} respectively.

Moreover for $ \varepsilon>0$ the following inequality on the rate
of convergence of approximated model in the space of continuous
functions holds true:
\begin{equation}\label{eq-ineq_sec5}
P\left\{\sup\limits_{t \in
T}\vert X(t)-X_N(t)\vert>\varepsilon\right\}=P\left\{\sup\limits_{t
\in T}\vert Y_N(t)\vert>\varepsilon\right\}\leq
\inf\limits_{u\geq1}\frac{\hat{B}_N^u u^{\alpha u}}{\varepsilon^u}.
\end{equation}

\end{theorem}

\begin{proof}

 The proof follows from Corollaries \ref{cor-diff} and
\ref{coroll-1}, Lemmas \ref{lem-1} and \ref{lem-2} for $\beta=1$.
Indeed, the function $\sigma_N(h)$ that majorizes $ \sup\limits_{\vert t- s\vert  \leq
h} \left(E\vert Y_N(t)-Y_N(s)\vert^2\right)^{1/2}\leq \sigma_N(h)$
equals
$$
\sigma_N(h)=\hat{C}_N h=\sqrt{\sum_{k=N+1}^\infty \left(
\frac{\hat{c}_1}{k}+\frac{\hat{c}_2}{k^2}+
\frac{\hat{c}_4}{k^4}+\frac{\hat{c}_3}{k^3}\right)^2}\cdot h,
$$
therefore, applying \eqref{eq-case-sum} we obtain \eqref{eq-case-c_n}.

\end{proof}

\begin{example}
Theorem \ref{the-case} is fulfilled for centered Gaussian stochastic processes.
In this case according to Example \ref{ex1} and Remark \ref{rem4} we have $\psi(u)=u^\frac12,$ $\alpha=\frac12$ and
$C_\Delta=2\sqrt{2}e^{-5/6}.$

If we take accuracy
$\varepsilon=0.5$ and reliability $1-\delta=0.95$, parameter $\beta=1$, then in the case of $T=[a,b]$ with $a=0$, $b=1$ and  $\Lambda=[-A,A]$ with
$A= \pi/(2b)$, the cutting-off level which
satisfies inequality (\ref{eq-approx-Bn_sec5}) is $N=518956$.

If we consider accuracy
$\varepsilon=0.05$ and reliability $1-\delta=0.9$,  parameter $\beta=0.9$, then in the case of $T=[a,b]$ with $a=0$, $b=2$ and
$A= \pi/(2b)$, the cutting-off level which
satisfies inequality (\ref{eq-approx-Bn_sec5}) is $N=985676$.

In the case of accuracy
$\varepsilon=0.1$ and reliability $1-\delta=0.9$,  parameters $\beta=0.9$,   $a=0$, $b=1$ and
$A= \pi/(1.1b)$, the cutting-off level which
satisfies inequality (\ref{eq-approx-Bn_sec5}) is $N=59867$.

If we take the accuracy
$\varepsilon=0.1$ and reliability $1-\delta=0.9$, the parameter $\beta=0.9$, then for $a=0$, $b=0.5$ and
$A= \pi/(1.1b)$, the cutting-off level which
satisfies inequality (\ref{eq-approx-Bn_sec5}) is $N=3112$.

All calculation were performed using the \emph{R} environment.
\end{example}

\section*{Conclusion}\label{sec7}

In this chapter we investigate approximation problems for stochastic processes from the space ${F}_\psi(\Omega)$.
 As approximating model the cutting-off series expansion is considered.
 We investigate the rate of convergence in uniform metrics that
gives us possibility to construct the approximating model with given
accuracy and reliability in the space of continuous functions.
We have considered in details the case of power function
$\psi(u)=u^\alpha,\, \alpha>0,$ and presented some results of calculations in the case of Gaussian processes.

\chapter{Random processes from Orlicz spaces of random variables}
\label{ch:o5series}

In this Chapter \ref{ch:o5series} the main properties of Orlicz spaces, in particular Orlicz spaces of exponential type, are described.
The connection of Orlicz spaces with the spaces $\mathbf{F}_\psi(\Omega)$ is investigated.
Large deviation inequalities for random processes from Orlicz spaces are found.

\section[Basic properties of Orlicz spaces and spaces $\mathbf{F}_\psi(\Omega)$]
{Basic properties of Orlicz spaces and connection with $\mathbf{F}_\psi(\Omega)$ spaces}

\begin{definition}\cite{dari11}
We say that a $C$-function $U$ satisfies the $\mathtt{g}$-condition if there exist constants $z_0\geq 0,$ $ K>0$ and $A>0$ such that for $x\geq z_0$, $y\geq z_0$ the inequality holds:
\[
U(x)U(y)\leq AU(Kxy).
\]
\end{definition}

\begin{definition} \cite{dari11}
A $C$-function $U=\left\{U(x), \ x \in \mathbb{R}\right\}$ is called an $N$-function of Orlicz ($N$-function) if the following conditions are satisfied:
\[
\lim\limits_{x \rightarrow \infty}\frac{U(x)}{x}=0; \ \lim\limits_{x \rightarrow 0}\frac{U(x)}{x}=\infty.
\]
\end{definition}

\begin{definition}\cite{giul03}
Let $\varphi=\left\{\varphi(x), \ x \in \mathbb{R} \right\}$ be an $N$-function. The function $\varphi^*$ is defined by the condition
\[
\varphi^*(x)=\sup\limits_{y \in \mathbb{R}}(xy-\varphi(y))
\]
and is called the Young-Fenchel transform with respect to $\varphi$.
\end{definition}

\begin{example}
The function $U(x)=a\left|x\right|^\alpha, x \in \R, a>0, \alpha\geq 1$ satisfies the $\mathtt{g}$-condition for $K=1$, $A=a$ and $z_0=0$.
\end{example}
The
$C$-function $U(x)=\exp\left\{\varphi(x)\right\}-1, x \in \R $, where $\varphi=\left(\varphi(x),x \in \R\right)$ is an arbitrary $C$-function that satisfies the $\mathtt{g}$-condition for $K=1$, $A=1$, $z_0=2$
$( \text{ if } \varphi(x)=\left|x\right|^\alpha,\alpha\geq 1, \text{ then } z_0=2^{1/\alpha})$.

\begin{lemma}
\label{le:o5chapter-2}
Let $m$ be a constant. Then for any Orlicz space $m \in L_U(\Omega)$ and $\left\|m\right\|_U=\frac{\left|m\right|}{U^{(-1)}(1)}$.
\end{lemma}

Lemma~\ref{le:o5chapter-2} is obvious.

\begin{lemma}
\label{le:o5chapter-3}
Let $\xi$ belong to the space $L_U(\Omega)$. Then there exists a constant $d_U$ such that $E\left|\xi\right|\leq d_U \left\|\xi\right\|_U$.
\end{lemma}

This Lemma is a consequence of Theorem 2.3.2 from \cite{byld98}.

\begin{example}
For the function $U(x)=\left|x\right|^p, p\geq2$ we have $ d_U=1$. For $U(x)=\exp\left\{\varphi(x)\right\}-1$, where $\varphi(x)$ is an $N$-function we have $d_U=\frac{2}{\varphi^{*(-1)}(1)}$. Here $\varphi^{*(-1)}(x)$ is the inverse function of $\varphi^*(x)$, $x \in \R $ which is the Young-Fenhel transform of the function $\varphi(x)$ (see Lemma 2.3.3 from \cite{byld98}).
\end{example}

For example: if $\varphi(x)=\frac{\left|x\right|^\alpha}{\alpha}, \alpha > 1$, then $\varphi^*(x)=\frac{\left|x\right|^\beta}{\beta}$, where $\frac{1}{\beta}+\frac{1}{\alpha}=1$, $\varphi^{*(-1)}(x)=(x\beta)^{1/\beta}$ and $\varphi^{*(-1)}(1)=(\beta)^{1/\beta}$. If $\alpha=2$, then $\beta=2$ and $d_U=\sqrt{2}$. If $\alpha=4$, then $\beta=\frac{4}{3}$ and $d_U=\frac{3^{3/4}}{4^{1/4}}$.

\begin{corollary} \rm{\cite{byld98}}
Let a $C$-function satisfy the $\mathtt{g}$-condition with constants $A$, $K$ and $z_0$. Then the sequence $(\varkappa(n), \ n\geq 1)$ is a majorizing characteristic of the space $L_U(\Omega)$ if
\[
\varkappa(n) =\begin{cases}
        n,&\text{if \; $n \leq U(z_0)$,}\\
        K(1+U(z_0))\max\left\{1,A\right\} U^{(-1)} (n),& \text {if \; $n > U(z_0)$.}
        \end{cases}
\]
\end{corollary}

\begin{definition} \label{de:o3chapter_6}
For the Orlicz space $L_U(\Omega)$, the condition $\bf{H}$ is satisfied if for any centered independent random variables  $\xi_1,\xi_2,\ldots,\xi_n$ from the space $L_U(\Omega)$ the following inequality holds:
\[
\left\|\sum_{k=1}^{n} \xi_k\right\|_U ^2\leq C_U \sum_{k=1}^{n}\left\|\xi_k\right\|_U ^2,
\]
where $C_U$ is an absolute constant.
\end{definition}

Examples of Orlicz spaces for which the condition $\mathbf{H}$ is satisfied:
\begin{itemize}
\item spaces $L_p(\Omega), p\geq2$, where $C_U=C_p=\sqrt{2}\left(\Gamma(p+1)/2\sqrt{\pi}\right)^{1/p}$ \cite{mats88};
\item spaces $L_U(\Omega)$, where $U(x)$ are such $C$-functions that there exist $p>q\geq2$, for which $U\left(\sqrt[q]{x}\right)$ is convex, and $U\left(\sqrt[p]{x}\right)$ is concave, and $C_U=2B_p$ \cite{bagr92}, where $B_p=2 k^{\frac{1}{2}}$, and $2k$ is the smallest even number not less than $p$ (see \cite{zygm65} p. 341);
\item Orlicz spaces are generated by the $C$-function $U(x)=exp\left\{\left|x\right|^\alpha\right\}-1$, where $1\leq\alpha\leq2$ (see \cite{giul03} and \cite{byld98}). In the work of Yu.~V. Kozachenko \cite{koz85} it is shown that for $\alpha\geq2$ the condition $\mathbf{H}$ is not satisfied for these spaces.
\end{itemize}

Furthermore, we will show that the condition $\bf{H}$ is satisfied for some Orlicz spaces $L_U(\Omega)$ such that the function $U(x)$ grows exponentially as $x\rightarrow\infty$.

Consider the Orlicz $C$-function
\begin{equation} \label{eq:o3chapter_1}
U(x)=\begin{cases}
\left(\frac{e\alpha}{2}\right)^{2/\alpha}x^2,&\text{if $\left|x\right|\leq x_\alpha$;}\\
\exp \left\{\left|x\right|^\alpha\right\},&\text{if $\left|x\right|> x_\alpha$,}
\end{cases}
\end{equation}
where $x_\alpha=\left(\frac{2}{\alpha}\right)^{1/\alpha}$, $0<\alpha<1$. $L_U(\Omega)$ is the Orlicz space generated by the function $U(x)$.

Consider the function $U_1(x)=\exp \left\{\left|x\right|^\alpha\right\}$, $0<\alpha\leq1$. Denote by the symbol $\mathcal{S}_{U_1}(\Omega)$ the family $\xi$ for which there exists $r$ such that $EU_1\left(\frac{\xi}{r}\right)<\infty$. Introduce on $\mathcal{S}_{U_1}(\Omega)$ the functional
\[
\left\langle \left\langle \xi\right\rangle\right\rangle_{U_1}=\inf\left\{r>0;EU_1\left(\frac{\xi}{r}\right)\leq2\right\}.
\]

\begin{lemma} \label{le:o1chapter_3}
In order for $\xi\in L_U(\Omega)$ to be true, it is necessary and sufficient that $\xi\in\mathcal{S}_{U_1}(\Omega)$ and the inequalities are true:
\begin{equation} \label{eq:o5chapter-2}
\left\|\xi\right\|_U\leq \left(e^{2/\alpha+2}\right)\left\langle \left\langle \xi\right\rangle\right\rangle_{U_1},
\end{equation}
\begin{equation} \label{eq:o5chapter-3}
\left\langle \left\langle \xi\right\rangle\right\rangle_{U_1}\leq \left\|\xi\right\|_U\left(e^{2/\alpha}+1\right)^{1/\alpha}.
\end{equation}
\end{lemma}

\begin{proof}
First, we prove the inequality \eqref{eq:o5chapter-2}. Let $r>0$, then
\begin{multline*}
EU\left(\frac{\xi}{r}\right)=EU\left(\frac{\xi}{r}\right)\mathbb{I}\left\{\frac{\left|\xi\right|}{r}\leq x_{\alpha}\right\}+EU\left(\frac{\xi}{r}\right)\mathbb{I}\left\{\frac{\left|\xi\right|}{r}> x_{\alpha}\right\}\leq\\[1ex]
\leq U\left(x_\alpha\right)+ E\exp \left\{\left|\frac{\xi}{r}\right|^\alpha\right\}= e^{2/\alpha}+E\exp \left\{\left|\frac{\xi}{r}\right|^\alpha\right\}.
\end{multline*}

Now let $r=\left\langle \left\langle \xi\right\rangle\right\rangle_{U_1}$, then $EU\left(\frac{\xi}{\left\langle \left\langle \xi\right\rangle\right\rangle_{U_1}}\right)\leq e^{2/\alpha}+2$. Since for $0<\alpha<1$ the inequality $U(\alpha x)\leq \alpha \ U(x)$ holds (Lemma 2.2.2 from the book \cite{byld98}), then
\[
EU\left(\frac{\xi}{\left\langle \left\langle \xi\right\rangle\right\rangle_{U_1}\left(e^{2\alpha}+2\right)}\right)\leq \frac{1}{e^{2\alpha}+2}EU\left(\frac{\xi}{\left\langle \left\langle \xi\right\rangle\right\rangle_{U_1}}\right)\leq 1.
\]
Therefore,
\begin{equation}
\left\|\xi\right\|_U\leq \left(e^{2/\alpha+2}\right)\left\langle \left\langle \xi\right\rangle\right\rangle_{U_1}.
\end{equation}

To verify the validity of the inequality \eqref{eq:o5chapter-3}, we note that

\[
E\exp \left\{\left|\frac{\xi}{r}\right|^\alpha\right\}=E\exp \left\{\left|\frac{\xi}{r}\right|^\alpha\right\}\mathbb{I}\left\{\frac{\left|\xi\right|}{r}< x_{\alpha}\right\}+
\]
\[
+E\exp \left\{\left|\frac{\xi}{r}\right|^\alpha\right\}\mathbb{I}\left\{\frac{\left|\xi\right|}{r}\geq x_{\alpha}\right\}\leq \exp \left\{\left(x_\alpha\right)^\alpha\right\}+ EU\left(\frac{\xi}{r}\right).
\]

Let $r=\left\|\xi\right\|_U$, then we have:
\begin{equation} \label{eq:o3chapter_22}
E\exp \left\{\left|\frac{\xi}{\left\|\xi\right\|_U}\right|^\alpha\right\}\leq e^{2/\alpha}+1.
\end{equation}
Using the inequality $\exp \left\{\left|ax\right|\right\}-1\leq a\left(\exp\left\{\left|x\right|\right\}-1\right)$, when $0<a\leq1$, we have:
\[
E\exp \left\{\left|\frac{\xi}{\left\|\xi\right\|_U}\right|^\alpha\frac{1}{e^{2/\alpha}+1}\right\}-1\leq \frac{1}{e^{2/\alpha}+1}\left(E\exp \left\{\left|\frac{\xi}{\left\|\xi\right\|_U}\right|^\alpha\right\}-1\right).
\]
Therefore, from the inequality \eqref{eq:o3chapter_22} it follows that
\begin{multline*}
E\exp \left\{\left|\frac{\xi}{\left\|\xi\right\|_U}\right|^\alpha\frac{1}{e^{2/\alpha}+1}\right\}-1\leq \frac{1}{e^{2/\alpha}+1}\left(E\exp \left\{\left|\frac{\xi}{\left\|\xi\right\|_U}\right|^\alpha\right\}-1\right)\leq\\[1ex]
\leq \frac{1}{e^{2/\alpha}+1}\left(e^{2/\alpha}+1-1\right)=\frac{e^{2/\alpha}}{e^{2/\alpha}+1}.
\end{multline*}

We received
\[
E\exp \left\{\left|\frac{\xi}{\left\|\xi\right\|_U\left(e^{2/\alpha}+1\right)^{1/\alpha}}\right|^\alpha\right\}\leq \frac{e^{2/\alpha}}{e^{2/\alpha}+1} +1\leq 2.
\]
It follows that $\left\langle \left\langle \xi\right\rangle\right\rangle_{U_1}\leq \left\|\xi\right\|_U\left(e^{2/\alpha}+1\right)^{1/\alpha}.$
\end{proof}

\begin {lemma} \label{le:o1chapter_4}
The inequality is valid
\[
\left\langle \left\langle \xi\right\rangle\right\rangle_{U_1}\geq \alpha^{1/\alpha}e^{1/\alpha}\left(\sup\limits_{n\geq1}\frac{\left(E\left|\xi\right|^n\right)^{1/n}}{n^{1/\alpha}}\right)
\]
for $0<\alpha<1$.
\end{lemma}

\begin{proof}
From the inequalities
\[
x^n\exp\left\{-x^\alpha\right\}\leq\left(\frac{n}{\alpha}\right)^{n/\alpha}\exp\left\{-\frac{n}{\alpha}\right\},
\]
\[
x^n\leq\exp\left\{x^\alpha\right\}\left(\frac{n}{\alpha}\right)^{n/\alpha}\exp\left\{-\frac{n}{\alpha}\right\}
\]
it follows that
\[
\frac{E\left|\xi\right|^n}{r^n}\leq E\exp\left\{\left(\frac{\left|\xi\right|}{r}\right)^\alpha\right\}\left(\frac{n}{\alpha}\right)^{n/\alpha}\exp\left\{-\frac{n}{\alpha}\right\},
\]
\[
E\left|\xi\right|^n\leq \left\langle \left\langle \xi\right\rangle\right\rangle_{U_1}^n 2 \left(\frac{n}{\alpha}\right)^{n/\alpha}\exp\left\{-\frac{n}{\alpha}\right\},
\]
\[
\left(E\left|\xi\right|^n\right)^{1/n}\leq \left\langle \left\langle \xi\right\rangle\right\rangle_{U_1} 2^{1/n} \left(\frac{n}{\alpha}\right)^{1/\alpha}\exp\left\{-\frac{1}{\alpha}\right\}.
\]
Since $\left\langle \left\langle \xi\right\rangle\right\rangle_{U_1}=\inf\left\{r>0;E\exp\left\{\left(\frac{\xi}{r}\right)^\alpha\right\}\right\}\leq 2$, then
\[
\left\langle \left\langle \xi\right\rangle\right\rangle_{U_1}\geq \left(E\left|\xi\right|^n\right)^{1/n}\frac{1}{2^{1/n} \left(\frac{n}{\alpha}\right)^{1/\alpha}\exp\left\{-\frac{1}{\alpha}\right\}}\geq \frac{\left(E\left|\xi\right|^n\right)^{1/n}}{n^{1/\alpha}}\alpha^{1/\alpha}e^{1/\alpha},
\]
which was to be proved.
\end{proof}

\begin {lemma} \label{le:o1chapter_5}
For $0<\alpha<1$ the inequality holds
\[
\left\langle \left\langle \xi\right\rangle\right\rangle_{U_1}\leq \left(1+\frac{e^{1/12}}{\sqrt{2\pi}}\right)^{1/\alpha}e^{1/\alpha}\left(\sup\limits_{n\geq1}\frac{\left(E\left|\xi\right|^n\right)^{1/n}}{n^{1/\alpha}}\right).
\]
\end{lemma}

\begin{proof}
From the Lyapunov inequality it follows that for $0<\alpha<1$ the inequality holds
\[
E\left|\xi\right|^{n\alpha}\leq\left(E\left|\xi\right|^n\right)^\alpha.
\]

We denote by $J_\alpha=E\exp\left\{\frac{\left|\xi\right|^\alpha}{r^\alpha}\right\}-1$. Therefore,
\[
J_\alpha=\sum\limits_{n=1}^\infty\frac{E\left|\xi\right|^{n\alpha}}{n!r^{\alpha n}}\leq\sum\limits_{n=1}^\infty\frac{\left(E\left|\xi\right|^n\right)^\alpha}{n!r^{\alpha n}}.
\]
Let $\hat{z}=\sup\limits_{n\geq1}\frac{\left(E\left|\xi\right|^n\right)^{1/n}}{n^{1/\alpha}}$. Since $E\left|\xi\right|^{n}\leq \hat{z}^n n^{n/\alpha}$, then $J_\alpha\leq\sum\limits_{n=1}^\infty\frac{\hat{z}^{n\alpha} n^n}{n!r^{\alpha n}}$.
From Stirling's formula we have that

\[
J_\alpha\leq\sum\limits_{n=1}^\infty\left(\frac{\hat{z}^\alpha e}{r^\alpha}\right)^n\frac{e^{1/{12n}}}{\sqrt{2\pi n}}\leq \frac{e^{1/{12}}}{\sqrt{2\pi}}\sum\limits_{n=1}^\infty\left(\frac{\hat{z}^\alpha e}{r^\alpha}\right)^{n}.
\]
Let $r=\frac{1}{s^{1/\alpha}}\hat{z} e^{1/\alpha}$, where $0<s<1$, then
\[
J_\alpha\leq\frac{e^{1/{12}}}{\sqrt{2\pi}}\sum\limits_{n=1}^\infty s^n=\frac{e^{1/{12}}}{\sqrt{2\pi}}\frac{s}{1-s}.
\]
If we put $s=\frac{1}{\left(1+\frac{e^{1/12}}{\sqrt{2\pi}}\right)}$, then we get:
\[
E\exp\left\{\frac{\left|\xi\right|^\alpha}{\left(\frac{1}{s^{1/\alpha}}\hat{z} e^{1/\alpha}\right)^\alpha}\right\}\leq 2.
\]
From these considerations it follows that the statement of the lemma is valid.
\end{proof}

\begin{theorem}
\label{th:o3chapter_1}
The Orlicz spaces $L_U(\Omega)$, where the function $U(x)$ is given by \eqref{eq:o3chapter_1}, contain the same elements as the spaces $\mathbf{F}_\psi(\Omega)$, where $\psi(u)=u^{1/\alpha}$, $\alpha>0$, and the norms in these spaces are equivalent and the following inequalities hold:
\begin{equation} \label{eq:o3chapter_4}
\left\|\xi\right\|_U\leq C_{\psi U} \left\|\xi\right\|_\psi,
\end{equation}
\begin{equation}\label{eq:o3chapter_5}
\left\|\xi\right\|_U\geq C_{U\psi} \left\|\xi\right\|_\psi,
\end{equation}
where
$$C_{\psi U}=e^{2/\alpha+2} \left(1+\frac{e^{1/12}}{\sqrt{2\pi}}\right)^{1/\alpha}e^{1/\alpha},$$
$$C_{U\psi}=\frac{1}{2^{1/\alpha}}\left(e^{2/\alpha}+1\right)^{-1/\alpha}\alpha^{1/\alpha}e^{1/\alpha}.$$
\end{theorem}

\begin{proof}
The Theorem follows from Lemma \ref{le:o1chapter_3}, Lemma \ref{le:o1chapter_4} and Lemma  \ref{le:o1chapter_5}. From Lemma \ref{le:o1chapter_3}, Lemma \ref{le:o1chapter_5} we have that
\[
\left\|\xi\right\|_U\leq e^{2/\alpha+2} \left(1+\frac{e^{1/12}}{\sqrt{2\pi}}\right)^{1/\alpha}e^{1/\alpha}\left(\sup\limits_{n\geq1}\frac{\left(E\left|\xi\right|^n\right)^{1/n}}{n^{1/\alpha}}\right),
\]
and from Lemma \ref{le:o1chapter_3}, Lemma \ref{le:o1chapter_4} we have inequality
\[
\left\|\xi\right\|_U \geq \left(e^{2/\alpha}+1\right)^{-1/\alpha}\alpha^{1/\alpha}e^{1/\alpha}\left(\sup\limits_{n\geq1}\frac{\left(E\left|\xi\right|^n\right)^{1/n}}{n^{1/\alpha}}\right).
\]
It's easy to see that
$$\sup\limits_{n\geq1}\frac{\left(E\left|\xi\right|^n\right)^{1/n}}{n^{1/\alpha}}\leq \sup\limits_{u\geq1}\frac{\left(E\left|\xi\right|^u\right)^{1/u}}{u^{1/\alpha}}=\left\|\xi\right\|_\psi.$$
 Then
\[
\sup\limits_{u\geq1}\frac{\left(E\left|\xi\right|^u\right)^{1/u}}{u^{1/\alpha}}\leq \sup\limits_{n\geq 2} \sup\limits_{n-1\leq u\leq n}\frac{\left(E\left|\xi\right|^u\right)^{1/u}}{u^{1/\alpha}}\leq
\]
\[
\leq\sup\limits_{n\geq 2} \frac{\left(E\left|\xi\right|^n\right)^{1/n}}{(n-1)^{1/\alpha}}\leq 2^{1/\alpha} \sup\limits_{n\geq1}\frac{\left(E\left|\xi\right|^n\right)^{1/n}}{n^{1/\alpha}},
\]
therefore the inequalities \eqref{eq:o3chapter_4} and \eqref{eq:o3chapter_5} hold true.
\end{proof}

\begin{theorem}
For the Orlicz space $L_U(\Omega)$, where $U(x)$ is given by \eqref{eq:o3chapter_1}, the condition $\bf{H}$ holds with the constant
\[
C_U=4\cdot 9^{\frac{1}{\alpha}}\left(\frac{C_{\psi U}}{C_{U\psi}}\right)^2,
\]
where $C_{\psi U}$ and $C_{U\psi}$ are defined in Theorem \ref{th:o3chapter_1}.
\end{theorem}

\begin{proof}
Let $\xi_k$, $k=1,2,\ldots,n$ be independent random variables from the Orlicz space $L_U(\Omega)$. Then from Theorem \ref{th:o3chapter_1} it follows that
\begin{equation}\label{eq:o3chapter_6}
\left\|\sum\limits_{k=1}^n\xi_k\right\|^2_U\leq C_{\psi U}^2 \sum\limits_{k=1}^n\left\|\xi_k\right\|_\psi^2,
\end{equation}
where $\psi(u)=u^{1/\alpha}$. From the calculations in Example \ref{ex:o1chapter_4} we have:
\begin{equation}\label{eq:o3chapter_7}
\left\|\sum\limits_{i=1}^n\xi_i\right\|^2_\psi\leq 4 \cdot 9^{\frac{1}{\alpha}} \sum\limits_{i=1}^n\left\|\xi_i\right\|_\psi^2.
\end{equation}
From Theorem \ref{th:o3chapter_1} and inequalities \eqref{eq:o3chapter_6} and \eqref{eq:o3chapter_7} it is obvious that the inequality holds:
\[
\left\|\sum\limits_{k=1}^n\xi_k\right\|^2_U\leq C_{\psi U}^2 \cdot 4 \cdot 9^{\frac{1}{\alpha}} \frac{1}{C_{U\psi}^2}\sum\limits_{k=1}^n\left\|\xi_k\right\|_U^2,
\]
which had to be proved.
\end{proof}

\section{Orlicz spaces of exponential type}

\begin{definition}\cite{byld98, giul03}
Let $\psi$ be an arbitrary $N$-function. The Orlicz space generated by the $N$-function
\[
U(x)=\exp \left\{\psi(x)\right\}-1, \ x \in \mathbb{R}
\]
is called the Orlicz space of exponential type.
\end{definition}

Let us denote this space by $Exp_\psi\left(\Omega\right)$, and the norm by $\left\|\cdot\right\|_{U(\psi)}$.

\begin{definition}\cite{giul03}
For an $N$-function $\psi$, the condition $Q$ is satisfied if
\[
\liminf \limits_{x \rightarrow 0} \frac{\psi(x)}{x^2}=C>0,
\]
where $C$ can be equal to $+\infty$.
\end{definition}

Examples of $N$-functions for which the condition $Q$ is satisfied are the following functions:
\begin{enumerate}
\item[1)] $\psi(x)=C\left|x\right|^{\alpha}, \ C>0, \ 1<\alpha\leq2$;
\item[2)] $\psi(x)=\begin{cases}
C\left|x\right|^2, \ \left|x\right|\leq1, \\
C\left|x\right|^\alpha, \ \left|x\right|>1, \ \alpha>2.
\end{cases}$
\end{enumerate}

\begin{definition}\cite{giul03}
Let $\varphi$ be an $N$-function for which the condition $Q$ holds. The random variable $\xi$ belongs to the space $Sub_\varphi(\Omega)$  ($\varphi$-sub-Gaussian), if $E\xi=0$, $E\exp\left\{\lambda \xi\right\}$ exists for all $\lambda \in \mathbb{R}$ and there exists a constant $a>0$ such that the following inequality holds for all $\lambda \in \mathbb{R}$
\[
E\exp\left\{\lambda \xi\right\}\leq \exp\left\{\varphi(\lambda a)\right\}.
\]
\end{definition}

In the works \cite{giul03, koz85} it is shown that the space $Sub_\varphi(\Omega)$ is a Banach space with respect to the norm
\[
\tau_\varphi(\xi)=\inf\left(a\geq0 : E\exp\lambda \xi\leq \exp\left\{\varphi(\lambda a)\right\}, \ \lambda \in \mathbb{R} \right).
\]

Let $L_U(\Omega)$ be an Orlicz space of exponential type generated by the $N$ - function $U(x)=\exp \left\{\psi(x)\right\}-1$, where $\psi(x)$ is an $N$-function such that for the $N$-function $\psi^*(x)$ the condition $Q$ is satisfied.

\begin{example}
If $\psi(x)=C \left|x\right|^\alpha$, $\alpha\geq 2$, then $\psi^*(x)=C_\beta \left|x\right|^\beta$, where $C_\beta=\frac{\alpha-1}{\alpha}\left(\frac{1}{C\alpha}\right)^{\frac{1}{\alpha-1}}$, $\beta>1$ is a number such that $\frac{1}{\alpha}+\frac{1}{\beta}=1$. Moreover, if $C=\frac{1}{\alpha}$, then $C_\beta=\frac{1}{\beta}$.
\end{example}

Consider the Orlicz space of exponential type $L_U(\Omega)$, where $U(x)=\exp \left\{\psi(x)\right\}-1$, such that for $\psi^*(x)$ the condition $Q$ is satisfied, as well as the space $Sub_{\psi^*}(\Omega)$.

\begin{theorem}\rm{\cite{giul03}}
\label{th:o3chapter_3}
In order for the random variable $\xi$, $E\xi=0$, to belong to the space $Exp_\psi\left(\Omega\right)$, it is necessary and sufficient that $\xi$ belongs to the space $Sub_{\psi^*}(\Omega)$, and the norms $\left\|\xi\right\|_{U(\psi)}$ and $\tau_{\psi^*}(\xi)$ are equivalent, i.e. the inequalities hold
\[
\left\|\xi\right\|_{U(\psi)}\leq 3 \tau_{\psi^*}(\xi),
\]
\[
\tau_{\psi^*}(\xi)\leq R_\psi \left\|\xi\right\|_{U(\psi)},
\]
where $R_\psi=S_{\psi^*} e^{\frac{49}{48}}$, $S_{\psi^*}=\max\limits_{i=\overline{1,3}}\gamma^{-1}_i$, and $\gamma_i=\gamma_i(\lambda_0)$ are defined as follows:
$\gamma_1$ is a solution of the equation $\gamma=\lambda_0 \sqrt{c_0(1-\gamma)}$, where $\lambda_0>0$, $c_0=\inf\limits_{0<\left|\lambda\right|\leq \lambda_0}\frac{\psi^*(\lambda)}{\lambda^2}$, $\gamma_2$ is a solution of the equation $\gamma^3-2(1-\gamma)=0$, $\gamma_3$ is a solution of the equation $\gamma=\psi^{*(-1)}(2)\sqrt{c_0(1-\gamma)}$.
\end{theorem}

\begin{example}
If the function $\psi^*(x)=\frac{\left|x\right|^\beta}{\beta}$, where $1<\beta\leq 2$, then from Theorem \ref{th:o3chapter_3} we have that
$c_0=\frac{1}{\beta}\left|\lambda_0\right|^{\beta-2}$ and $\gamma_1=\lambda_0^{\frac{\beta}{2}}\frac{1}{\sqrt{\beta}}\sqrt{1-\gamma}$, $\gamma_2=0.770917$, $\gamma_3=2^{\frac{1}{\beta}}\beta^{\frac{1}{\beta}-\frac{1}{2}}\lambda_0^{\frac{\beta}{2}-1}\sqrt{1-\gamma}$.
If we choose $\lambda_0$ such that $\gamma_1>\gamma_2$ and $\gamma_3>\gamma_2$, then $S_{\psi^*}= \frac{1}{\gamma_2}=1.2972$.
\end{example}

\begin{theorem}\rm{\cite{giul03}}
\label{th:o3chapter_4}
Let the space $Sub_{\psi^*}(\Omega)$ be such that the function $\psi^*\left(\sqrt{x}\right)$, $x>0$ is convex and $\xi_1,\xi_2,\ldots,\xi_n$ are independent centered random variables from the space $Sub_{\psi^*}(\Omega)$, then the inequality holds:
\[
\tau_{\psi^*}^2(\sum \limits_{k=1}^n\xi_k)\leq \sum \limits_{k=1}^n \tau_{\psi^*}^2(\xi_k).
\]
\end{theorem}

\begin{theorem}
\label{th:o3chapter_5}
Let $Exp_\psi\left(\Omega\right)$ be a space such that the function $\psi^*\left(\sqrt{x}\right)$, $x>0$, is convex and let $\xi_1,\xi_2,\ldots,\xi_n$ be independent centered random variables from this space. Then the inequality holds
\[
\left\|\sum \limits_{k=1}^n\xi_k\right\|^2_{U(\psi)}\leq 9 R^2_\psi \sum \limits_{k=1}^n\left\|\xi_k\right\|^2_{U(\psi)},
\]
that is, for this space the condition $\mathbf{H}$ with the constant $9 R^2_\psi$ is satisfied, where $R_\psi$ is determined in Theorem \ref{th:o3chapter_3}.
\end{theorem}

\begin{proof}
From Theorem \ref{th:o3chapter_3} and  Theorem\ref{th:o3chapter_4} we obtain that
\[
\left\|\sum \limits_{k=1}^n\xi_k\right\|^2_{U(\psi)}\leq 9 \tau_{\psi^*}^2(\sum \limits_{k=1}^n\xi_k)\leq 9 \sum \limits_{k=1}^n \tau_{\psi^*}^2(\xi_k)\leq 9 R^2_\psi \sum \limits_{k=1}^n\left\|\xi_k\right\|^2_{U(\psi)}.
\]
From these inequalities follows the statement of the Theorem.
\end{proof}

Theorem \ref{th:o3chapter_5} holds for Orlicz spaces $Exp_\psi\left(\Omega\right)$ when
\begin{equation}\label{eq:o3chapter_8}
\psi^*(x)=\begin{cases}
\frac{\left|x\right|^2}{\alpha}, \ \left|x\right|\leq1; \\
\frac{\left|x\right|^\alpha}{\alpha}, \ \left|x\right|>1, \ \alpha>2,
\end{cases}
\end{equation}
since
\begin{equation}\label{eq:o3chapter_9}
\psi(x)=\begin{cases}
\frac{\beta\left|x\right|^2}{4(\beta-1)}, \ x \leq 2\left(1-\frac{1}{\beta}\right); \\
\left|x\right|-\left(1-\frac{1}{\beta}\right), \ 2\left(1-\frac{1}{\beta}\right)<x<1; \\
\frac{\left|x\right|^\beta}{\beta}, \ x\geq 1,
\end{cases}
\end{equation}
where $\beta>1$ is a number such that $\frac{1}{\alpha}+\frac{1}{\beta}=1$.

Note that for $\alpha\geq 2$ we have $1<\beta\leq 2$.

\begin{remark}
It is easy to see that the condition $\mathbf{H}$ is also satisfied for the spaces $Exp_{\widetilde{\psi}}\left(\Omega\right)$, where $\widetilde{\psi}(x)=C \left|x\right|^\beta$, $1<\beta\leq 2$.
\end{remark}

It is obvious that the function $\widetilde{\psi}(x)=C \left|x\right|^\beta$ is equivalent to the function \eqref{eq:o3chapter_9}. Therefore, the spaces $Exp_\psi\left(\Omega\right)$ and $Exp_{\widetilde{\psi}}\left(\Omega\right)$ contain the same elements and their norms are equivalent. Indeed, let $\xi$ belong to the space $Exp_{\widetilde{\psi}}\left(\Omega\right)$. For simplicity, let $C=\frac{1}{\beta}$, i.e. ${\widetilde{\psi}}(x)=\frac{\left|x\right|^\beta}{\beta}$. Therefore, ${\widetilde{\psi}}(x)=\psi(x)$ for $\left|x\right|\geq 1$. Let $r>0$, then
\begin{multline*}
E\exp \left\{\psi\left(\frac{\left|\xi\right|}{r}\right)\right\}-1=E \ \mathbb{I}\left\{\frac{\left|\xi\right|}{r}\leq 1\right\}
\exp \left\{\psi\left(\frac{\left|\xi\right|}{r}\right)\right\}+\\[1ex]
+E \ \mathbb{I}\left\{\frac{\left|\xi\right|}{r}> 1\right\}
\exp \left\{\psi\left(\frac{\left|\xi\right|}{r}\right)\right\}-1 \leq \exp \left\{\psi(1)\right\}+\\[1ex]
+E \ \mathbb{I}\left\{\frac{\left|\xi\right|}{r}> 1\right\}
\exp \left\{\widetilde{\psi}\left(\frac{\left|\xi\right|}{r}\right)\right\}-1 \leq e^{\frac{1}{\beta}}+E\exp \left\{\widetilde{\psi}\left(\frac{\left|\xi\right|}{r}\right)\right\}-1.
\end{multline*}
If we put $r=\left\|\xi\right\|_{U(\widetilde{\psi})}$, then we get that
\[
E\exp \left\{\psi\left(\frac{\left|\xi\right|}{\left\|\xi\right\|_{U(\widetilde{\psi})}}\right)\right\}-1\leq e^{\frac{1}{\beta}}+1.
\]
Therefore,
\[
\left\|\xi\right\|_{U(\psi)}\leq \left\|\xi\right\|_{U(\widetilde{\psi})} \left(1+e^{\frac{1}{\beta}}\right).
\]
Similarly, we obtain that
\[
\left\|\xi\right\|_{U(\widetilde{\psi})} \leq  \left\|\xi\right\|_{U(\psi)}\left(1+e^{\frac{1}{\beta}}\right).
\]

Let $\xi_1,\xi_2,\ldots,\xi_n$ be independent centered random variables from the space $Exp_{\widetilde{\psi}}\left(\Omega\right)$, then the following inequality follows from Theorem \ref{th:o3chapter_5}:
\[
\left\|\sum \limits_{k=1}^n\xi_k\right\|^2_{U(\widetilde{\psi})}\leq \left(1+e^{\frac{1}{\beta}}\right)^2 \left\|\sum \limits_{k=1}^n\xi_k\right\|^2_{U(\psi)}\leq
\]

\[
\leq9 R^2_\psi \left(1+e^{\frac{1}{\beta}}\right)^2 \sum \limits_{k=1}^n\left\|\xi_k\right\|^2_{U(\psi)} \leq 9 R^2_\psi \left(1+e^{\frac{1}{\beta}}\right)^4 \sum \limits_{k=1}^n\left\|\xi_k\right\|^2_{U(\widetilde{\psi})}.
\]

In other words, for the space $Exp_{\widetilde{\psi}}\left(\Omega\right)$ condition $\mathbf{H}$ holds with respect to the constant $C_{\widetilde{\psi}}= R^2_\psi \left(1+e^{\frac{1}{\beta}}\right)^4$.

\section{Random processes from Orlicz spaces}

\begin{definition}
Let $X = \{X(t),\ t \in T\}$ be a random process, where $T$ is a parametric set. We say that the process $X$ belongs to the Orlicz space $L_U(\Omega)$ if for any $t \in T$ the random variable $X(t)$ belongs to the space $L_U(\Omega)$.
\end{definition}

Let $\rho(t,s)=\left\|X(t)-X(s)\right\|_U$ be the pseudometric generated in $T$ by the process $X = \{X(t),\ t \in T\}$, which belongs to the Orlicz space $L_U(\Omega)$.
Consider the pseudometric space $\left(T,\rho\right)$. Let $N(u)$ be the metric massiveness of the space $\left(T,\rho\right)$ (see Definition in \ref{de:o3chapter-1}).

\begin{theorem} \label{th:o3chapter_6}
Let $X = \{X(t),\ t \in T\}$ be a random process from the space $L_U(\Omega)$, where $U(x)$ is a $C$-function that satisfies the $\mathtt{g}$-condition and the process $X$ is separable on $\left(T,\rho\right)$. If
\begin{enumerate}
\item $\sup\limits_{t \in T} \bigl\|X(t) \bigr\|_U < \infty$;
\item condition is satisfied
\begin{equation}
\int\limits_{0}^{\varepsilon_0} U^{(-1)} (N(\varepsilon)) d\varepsilon < \infty,
\end{equation}
\end{enumerate}
where $\varepsilon_0=\sup\limits_{t,s \in T} \rho (t,s)$, then:
\begin{enumerate}
\item $\sup\limits_{t \in T} \left|X(t) \right| \in L_U(\Omega)$;
\item $\left\|\sup\limits_{t \in T} \left|X(t) \right|\right\|_U \leq B$, where $B=\inf\limits_{t \in T} \left\|X(t)\right\|_U +\inf\limits_{0<\theta<1}\frac{1}{\theta(1-\theta)} \int\limits_{0}^{\theta\varepsilon_0} \varkappa(N(\varepsilon)) d\varepsilon$, where
$\varkappa(n)$ is the majorizing characteristic of the space $L_U(\Omega)$;
\item for all $x>0$
   \[
   P\left\{\sup\limits_{t \in T} \left|X(t)\right|\geq x \right\}\leq \frac{1}{U\left(\frac{x}{B}\right)}.
   \]
\end{enumerate}
\end{theorem}

\begin{proof}
Theorem~\ref{th:o3chapter_6} is a modified Corollary 3.3.1 from the book \cite{byld98}.
\end{proof}

\begin{theorem} \label{th:o3chapter_7}
Let $(T,w)$ be a compact metric space, $N_w(u)$ be the metric massiveness of the space $(T,w)$, $X = \{X(t), t \in T\}$ be a separable random process from the Orlicz space $L_U(\Omega)$, $\varkappa(n)$ be the majorizing characteristic of the space $L_U(\Omega)$. Let the function $U(x)$ satisfies the $\mathtt{g}$-condition. Let there be a function
\[
\sigma = \{\sigma(h),\,\,\,0 \leq h \leq \sup\limits_{t, s \,\in\, T} w(t, s)\},
\]
such that $\sigma(h)$ is monotonically increasing, continuous, and $\sigma(h)\rightarrow 0$ as $h\rightarrow 0$, and
\[
\sup\limits_{w(t, s) \leq h} \bigl\|X(t) - X(s) \bigr\|_U \leq \sigma(h).
\]
If
\[
\int\limits_{0}^{\delta_0} U^{(-1)}(N_w(\sigma^{(-1)} (u))) du < \infty,
\]
where $\delta_0=\sigma\left(\sup\limits_{t, s \,\in\, T} w(t,s)\right)$, $\sigma^{(-1)} (u)$ is an inverse to $\sigma(u)$ function, then
\begin{enumerate}
	\item $\sup\limits_{t \in T} \left|X(t) \right| \in L_U(\Omega)$;
	\item $\left\|\sup\limits_{t \in T} \left|X(t) \right|\right\|_U \leq \widetilde{B}$,\\
where $\widetilde{B}=\inf\limits_{t \in T} \left\|X(t)\right\|_U +\inf\limits_{0<\theta<1} \frac{1}{\theta(1-\theta)} \int\limits_{0}^{\theta\delta_0} \varkappa(N_w(\sigma^{(-1)} (u))) du$.
\end{enumerate}
\end{theorem}

\begin{proof}
This Theorem follows from Theorem~\ref{th:o3chapter_6}, since the process $X$ is separable on $(T,w)$, then the process $X(t)$ is separable on $(T,\rho)$ and the inequality holds
\[
N(u)\leq N_w \left(\sigma^{-1}(u)\right).
\]
\end{proof}

Consider the space $\mathbb R^{d}$ with the metric $m\left(\vec{x},\vec{y}\right)=\max\limits_{1\leq i \leq d} \left|x_i - y_i\right|$.

\begin{corollary}
Let $T$ be the cube $\left\{0\leq x_i<T, i=\overline{1,d}\right\},T>0$. Then from Theorem~\ref{th:o3chapter_7}
\[
N_w (u)\leq \left(\frac{T}{2u}+1\right)^d
\]
and the following inequality holds:
\[
\widetilde{B}\leq\widehat{B}=\inf\limits_{t \in T} \left\|X(t)\right\|_U +\inf\limits_{0<\theta<1} \frac{1}{\theta(1-\theta)} \int\limits_{0}^{\theta\delta_0} \varkappa\left(\left(\frac{T}{2\sigma^{(-1)}(u)}+1\right)^d\right) du.
\]
\end{corollary}

\begin{theorem}\label{th:o3chapter_8}
Let $Y = \left\{Y(t), t \in T\right\}$ be a random process belonging to the Orlicz space $L_U(\Omega)$.
Let for the function $U$ the $\mathtt{g}$-condition holds true and for the space $L_U(\Omega)$ the condition $\mathbf{H}$ with the constant $C_U$ holds true.
Let $(T, w)$ be a compact metric space and $Y$ be a separable process on $(T, w)$, $N_w(u)$ be the metric massiveness of the space $(T, w)$. Let there be a function
\[
\sigma = \{\sigma(h),\,\,\,0 \leq h \leq \sup\limits_{t, s \,\in\, T} w(t, s)\},
\]
such that $\sigma(h)$ is monotonically increasing continuous, $\sigma(h)\rightarrow 0$, when $h\rightarrow 0$ and
\[
\sup\limits_{w(t, s) \leq h} \bigl\|Y(t) - Y(s) \bigr\|_U \leq \sigma(h)
\]
and
\begin{equation}
\int\limits_{0}^{\delta_0} U^{(-1)}(N_w(\sigma^{(-1)} (u))) du < \infty.
\end{equation}

Let $X(t) = Y(t) - m(t)$, where $m(t) = EX(t)$ and $X_k (t)$ are independent copies of $X(t)$,
$$S_n (t) = \frac{1}{\sqrt{n}} \sum\limits_{k = 1}^{n} X_{k}(t).$$
Then for all $\varepsilon > 0$ the inequality holds true
\begin{equation}\label{eq:o5chapter-12}
P\{\sup_{t \in T} |S_n (t) | > \varepsilon\} \leq \frac{1}{U \left(\frac{\varepsilon}{B(t_0, \theta)}\right)},
\end{equation}
where $t$ is an arbitrary point from $[0,T]$, $0 < \theta < 1$,
\begin{equation}
B(t_0, \theta) = \|X(t_0)\|_U + \frac{1}{\theta(1 - \theta)}\int\limits_{0}^{\delta_{0}\theta} \varkappa (N_w (\sigma^{(-1)}_1 (u)) du,
\end{equation}
where $\sigma_1 (h)=C_U\left(1+\frac{d_U}{U^{(-1)}(1)}\right)\sigma (h)$, $d_U$ is a constant defined in Lemma \ref{le:o5chapter-3}, $\delta_0 = \sigma_1 \left(\sup\limits_{t, s \,\in\,T } w(t,s)\right)$, $\varkappa(n)$ is the major characteristic of the space $L_U(\Omega)$.
\end{theorem}

\begin{proof}
From Definition~\ref{de:o3chapter_6} we have that
\begin{multline} \label{de:o3chapter_14}
\left\|S_n (t)-S_n(s)\right\|_U^2=\left\|\frac{1}{\sqrt{n}}\sum \limits_{k=1}^n\left(X_k(t)-X_k(s)\right)\right\|_U^2\leq \\[1ex]
\leq C_U \frac{1}{n}\sum \limits_{k=1}^n\left\|X_k(t)-X_k(s)\right\|_U^2= C_U \left\|X(t)-X(s)\right\|_U^2.
\end{multline}

From Lemma~\ref{le:o5chapter-2} and Lemma~\ref{le:o5chapter-3} we obtain:
\begin{multline} \label{de:o3chapter_15}
\left\|X(t)-X(s)\right\|_U = \left\|Y(t)-Y(s)-(m(t)-m(s))\right\|_U\leq \left\|Y(t)-Y(s)\right\|_U+ \\[1ex]
+\frac{1}{U^{(-1)}(1)}\left|m(t)-m(s)\right|\leq  \left\|Y(t)-Y(s)\right\|_U+ \\[1ex]
+\frac{1}{U^{(-1)}(1)}E\left|Y(t)-Y(s)\right|\leq  \left\|Y(t)-Y(s)\right\|_U + \\[1ex]
+\frac{d_U}{U^{(-1)}(1)}\left\|Y(t)-Y(s)\right\|_U = \left\|Y(t)-Y(s)\right\|_U \left(1+\frac{d_U}{U^{(-1)}(1)}\right).
\end{multline}

Therefore, from \eqref{de:o3chapter_14} it follows that
\[
\sup\limits_{w(t, s)\,\leq\, h} \|X(t) - X(s)\|_U \leq \sigma_1(h).
\]
Now the statement of Theorem~\ref{th:o3chapter_8} is established on the basis of Theorem~\ref{th:o3chapter_7}.
\end{proof}

\begin{remark}\label{re:o5chapter-2}
It is easy to prove (as in~\eqref{de:o3chapter_15}) that
\begin{equation}
\left\|X(t_0)\right\|_U\leq\left(1+\frac{d_U}{U^{(-1)}(1)}\right)\left\|Y(t_0)\right\|.
\end{equation}
\end{remark}

\section*{Conclusions to chapter \ref{ch:o5series}}

In this chapter we consider the main properties of Orlicz spaces and exponential-type Orlicz spaces. The connection of Orlicz spaces with the spaces $\mathbf{F}_\psi(\Omega)$ is investigated.
Large deviation inequalities are found for stochastic processes from Orlicz spaces.



\chapter{Spaces $D_{V,W}(\Omega)$ of random variables}
\label{omch:dvw1}

In this chapter we consider the pre-Banach $K_\sigma$-spaces $D_{V,W}(\Omega)$. Basic properties of these spaces are studied, as well as the conditions for the convergence of series in $D_{V,W}(\Omega)$..

\section{Definition and basic properties of the space $D_{V,W}(\Omega).$}

\begin{definition}
\label{omomd2.1}
Let $W=\{W(x),x\in R\}$ and $V=\{V(x),x\in R\}$, $W(x)>0$, $V(x)>0$, $x\neq 0$ be some even monotonically increasing continuous functions at $x>0$, which have the following properties:

\begin{enumerate}
\item There is a constant $C>0$ for which  $$W^{(-1)}(x+y)\leq C(W^{(-1)}(x)+ W^{(-1)}(y)),\, x>0,\, y>0,$$ where $W^{(-1)}(x)$ is the function inverse at $x>0$ to $W(x)$, $V^{(-1)}(x)$ is a function inverse to $V(x)$ for $x>0$;
\item There is a continuous function $Z=\{Z(x),x>0\}$ such that $0<Z(x)<\infty$, $|x|<\infty$, and for any constant $a>0$ the inequality $V(ax)\leq Z(a)V(x)$ for $x>0$ holds;
\item $W(0)=0$.
\end{enumerate}
A random
variable $\xi$ belongs to the space $D_{V,W}(\Omega)$ if
\begin{equation}
\label{om2.1}
\sup\limits_{x\geq0} V(x)W^{(-1)}(P\{|\xi|>x\})<\infty.
\end{equation}
\end{definition}

Examples of functions $W$ and $V$ are the following functions:
$$W(x)=|x|^a,\, V(x)=|x|^b\,\,\, (Z(x)=|x|^b),$$
$$W(x)=\exp\{|x|^a\}-1,\, a>0,\, b>0.$$

We will say that for functions $V$ and $W$ condition $B1$ is satisfied if:

\bigskip
{\bf B1)} The function $W$ is a $C$-Orlicz function, and $V$ is a function inverse to
a $C$-Orlicz function
\bigskip
\begin{theorem}
\label{omt2.1} The space $D_{V,W}(\Omega)$ is a space with prenorm
$$||\xi||_{V,W} = \left(\sup\limits_{x>0}V(x)W^{-1}(P\{|\xi|>x\})\right)^{1/2}.$$
If $\xi_n$ is a sequence of random variables from $D_{V,W}(\Omega)$ such that $||\xi_n-\xi_m||_{V,W}\to\infty$ when $n,m\to\infty$ and
$\sup_{n}||\xi_n||_{V,W}<\infty$, then there exists a random variable $\xi\in
D_{V,W}(\Omega)$ such that $||\xi_n-\xi||_{V,W}\to 0$, $n\to\infty$.
In addition, the prenorm $||\cdot||_{V,W}$ is subject to the function
$J(\lambda)=(Z(\lambda))^{1/2}$.
\\
If condition $B1$ is satisfied for $W$ and $V$, then the functional $||\cdot||$ is a quasi-norm,
and the space is complete with respect to this quasi-norm.
\\
Finely,
\begin{equation}
\label{om2.2} P\{|\xi|>x\}\leq
W\left(\frac{||\xi||^2_{V,W}}{V(x)}\right).
\end{equation}
for all $x>0$.
\end{theorem}

\begin{dov}
It is clear that for the space $D_{V,W}(\Omega)$ the conditions $a1, a2, a3$ of the Definition \ref{om-prebanach} are satisfied,
i.e. $D_{V,W}(\Omega)$ is a pre-$K_\sigma$-space. Also,
$||\xi||_{V,W}=0$ if and only if $\xi=0$ with probability one. The inequality \ref{om2.2} is obvious and follows from the definition of the prenorm.

Let us now prove that the space $D_{V,W}(\Omega)$ is linear. For this
it is sufficient to prove that, if $\xi$ and $\eta$ belong to $D_{V,W}$, then
$a\cdot\xi\in D_{V,W}(\Omega)$ and $\xi+\eta\in D_{V,W}(\Omega)$. If $a\neq0$, then

$$||a\xi||_{V,W} = \left(\sup_{x>0}V(x)W^{(-1)}(P\{|a\xi|>x\})\right)^{1/2} =$$

$$=  \left(\sup_{x>0}V(x)W^{(-1)}(P\{|\xi|>\frac{x}{|a|}\})\right)^{1/2} = $$

$$ = \left(\sup_{x>0}V(x\frac{|a|}{|a|})W^{(-1)}(P\{|\xi|>\frac{x}{|a|}\})\right)^{1/2}.$$

Putting $y:=x/|a|$, we get:

$$||a\xi||_{V,W} = \sup_{y>0}V(|a|y)W^{(-1)}(P\{|\xi|>y\}))^{1/2}\leq $$
$$\leq (Z(|a|)\cdot\sup_{y>0}V(y)W^{(-1)}(P\{|\xi|>y\}))^{1/2}=Z(|a|)^{1/2}||\xi||_{V,W}<\infty.$$
Here we also showed that the prenorm $||\cdot||_{V,W}$
is subject to the function $J(\lambda)=(Z(\lambda))^{1/2}$.

Let us now prove that the sum of the elements of the space $D_{V,W}$ also belongs to this space:

$$||\xi+\eta||_{V,W} = \left(\sup_{x>0}V(x)W^{(-1)}(P\{|\xi+\eta|>x\})\right)^{1/2}\leq$$

$$\leq (\sup_{x>0}V(x)W^{(-1)}(P\{|\xi|>x/2\}+P\{|\eta|>x/2\}))^{1/2}\leq$$

$$\leq (\sup_{x>0}V(x)\cdot C\cdot(W^{(-1)}(P\{|\xi|>x/2\})+W^{(-1)}(P\{|\eta|>x/2\})))^{1/2}\leq$$

$$\leq(C\sup_{x>0}V(x)\cdot(W^{(-1)}(P\{|\xi|>x/2\}))+$$
$$+\sup_{x>0}V(x)\cdot(W^{(-1)}(P\{|\eta|>x/2\})))^{1/2} \leq $$

\begin{equation}
\label{om2.4}
\leq(CZ(2)(||\xi||_{V,W}+||\eta||_{V,W}))<\infty.
\end{equation}

Let $\xi_n \in D_{V,W}(\Omega)$ be a sequence such that
$||\xi_n-\xi_m||_{V,W}\to0$ as $n,m\to\infty$ and let $\sup_n||\xi_n||_{V,W}<\infty$. 
From inequality (\ref{om2.2}) it follows
$$P\{|\xi_n-\xi_m|>\epsilon\}\leq W\left(
\frac{||\xi_n-\xi_m||^2_{V,W}}{V(\epsilon)} \right)\to0.$$
Therefore, there exists a random variable $\xi$ such that $\xi_n\to\xi$ in
probability, hence, $P\{|\xi_n|>x\}\to P\{|\xi|>x\}$ at points where
the function $P\{|\xi|>x\}$ is continuous. Then at points of continuity of the function $P\{|\xi|>x\}$
$$P\{|\xi|>x\}\leq \sup_{n\geq1}P\{|\xi_n|>x\}\leq$$
$$\leq \sup_{n\geq 1} W\left(\frac{||\xi_n||_{V,W}^2}{V(x)}\right)= W\left(\frac{\sup_{n\geq1}||\xi_n||_{V,W}^2}{V(x)}\right).$$
Since the function $W(\frac{a}{V(x)})$ is continuous in $x$ (where $a>0$ is
any constant), and $P\{|\xi|>x\}$ does not monotonically increase,
the latter inequality is also true for all $x>0$. Therefore, $\xi\in D_{V,W}$
and $||\xi||_{V,W}\leq\sup_{n\geq 1}||\xi_n||_{V,W}$.

Now, instead of the sequence $\xi_n$, consider the sequence
$\xi_n-\xi_m$, $n>m$. As in the previous case, we have
$$||\xi-\xi_m||_{V,W}\leq\sup_{n\geq m}||\xi_n-\xi_m||_{V,W},$$
that is $||\xi_m-\xi||_{V,W}\to0$ as $m\to\infty$.

Let $0<\alpha<1$ and condition $B1$ be satisfied. From the properties of the $C$-functions: $V(ax)\leq aV(x)$, $a>1$, and $W(x+y)\leq W(x)+W(y)$, it follows that
$$||\xi+\eta||^2_{V,W} = \sup\limits_{x>0}V(x)W^{-1}(P\{|\xi+\eta|>x\})\leq$$
$$\leq\sup\limits_{x>0}(V(x)W^{-1}(P\{|\xi|+|\eta|>x\}))\leq$$ $$\leq\sup\limits_{x>0}(V(x)W^{-1}(P\{|\xi|>\alpha x\}+P\{|\eta|>(1-\alpha)x\}))\leq$$
$$\leq \sup\limits_{x>0}V\left(\frac{\alpha x}{\alpha}\right) W^{-1}(P\{|\xi|>\alpha x\}) + \sup\limits_{x>0}V\left(\frac{(1-\alpha) x}{(1-\alpha)}\right)\times$$
 $$\times W^{-1}(P\{|\eta|>(1-\alpha) x\})\leq\frac{1}{\alpha}||\xi||^2_{V,W} + \frac{1}{1-\alpha} ||\eta||^2_{V,W}.$$
The triangle inequality follows from the last inequality if we put
 $$\alpha =\frac{
||\xi||_{V,W}}{(||\xi||_{V,W}+||\eta||_{V,W})}.$$

Let us prove that under condition $B1$ the space $D_{V,W}(\Omega)$ is complete.

Since $|$ $||\xi_n||_{V,W}-||\xi_m||_{V,W}$ $|$ $\leq ||\xi_n-\xi_m||_{V,W}\to0$ for $n,m\to\infty$, then there exists a limit $\sigma=\lim\limits_{n\to\infty}||\xi_n||_{V,W}$.

If in the inequality $$P\{|\xi_n|>x\}\leq
W\left(\frac{||\xi_n||^2_{V,W}}{V(x)}\right).$$
we go to the limit as $n\to\infty$ (it is enough to consider the points $x$,
which are the points of continuity of the distribution function of the random variable
$\xi$), then we obtain the inequality $$P\{|\xi|>x\}\leq
W\left(\frac{\sigma^2}{V(x)}\right).$$
Therefore, the random variable $\xi$ belongs to the space $D_{V,W}(\Omega)$ and $||\xi||_{V,W}\leq\sigma$.

Now, instead of the sequence $\xi_n$, consider the sequence
$\xi_m-\xi_n$, $m\to\infty$. As in the previous case,
we obtain $$||\xi-\xi_n||_{V,W}
\leq\lim\limits_{m\to\infty}||\xi_m-\xi_n||_{V,W}.$$
 Therefore,
$||\xi-\xi_n||_{V,W}\to 0$ as $n\to\infty$, i.e., the space
$D_{V,W}$ is complete.
\end{dov}

Now, let us find the majorizing characteristic of the space $D_{V,W}(\Omega)$.

\begin{theorem}
\label{omt2.2} The sequence $$\kappa(n) =
\sup\limits_{0<t<1/n}\left( \frac{W^{(-1)}(tn)}{W^{(-1)}(t)}
\right)^{1/2}$$ is a majorizing characteristic of the space
$D_{V,W}(\Omega)$.
\end{theorem}

\begin{dov} Let $\xi_i$, $i=1,\ldots,n$ be random variables
from the space $D_{V,W}(\Omega)$, $a=\max\limits_{1\leq i\leq n}
||\xi_i||_{V,W}$. The following relations hold true
$$||\max\limits_{1\leq i\leq n} |\xi_i|||^2_{V,W} =
\sup\limits_{x>0}(V(x)W^{(-1)}(P\{\max\limits_{1\leq i\leq
n}|\xi_i|>x\}))\leq$$
$$\leq \sup\limits_{x>0} (V(x)W^{(-1)}(\min\{1,\sum_{i=1}^nP\{|\xi_i|>x\}\}))\leq$$
$$\leq\sup_{x>0}\left(V(x)W^{(-1)}\left(\min\left\{1,\sum_{i=1}^n W\left(\frac{||\xi_i||^2_{V,W}}{V(x)}\right)\right\}\right)\right)\leq$$
\begin{equation}
\label{om2.3}
\leq\sup_{x>0}\left(V(x)W^{(-1)}\left(\min\left\{1,nW\left(\frac{a^2}{V(x)}\right)\right\}\right)\right).
\end{equation}

Let $t=W(\frac{a^2}{V(x)})$.
Then $V(x) =\frac{a^2}{W^{(-1)}(t)}$. So, from (\ref{om2.3}) we get that

$$||\max_{1\leq i\leq n}|\xi_i|||^2_{V,W}\leq \sup_{t>0}\left( \frac{a^2}{W^{(-1)}(t)} W^{(-1)}(\min\{1,nt\}) \right)=$$

$$=a^2\sup_{0<t\leq \frac{1}{n}}\left( \frac{W^{(-1)}(nt)}{W^{(-1)}(t)} \right).$$
\end{dov}

In the future, in this section, the prenorm $||\cdot||_{V,W}$ will be denoted by $||\cdot||$.

\section{Convergence of series in spaces $D_{V,W}(\Omega)$}

\begin{theorem}
\label{omt3.1} Let $\xi_k$ be random variables from the space $D_{V,W}(\Omega)$, let $||\cdot||$ be the prenorm, $||\xi_k||>0$, $f(x)=
xV(W(x))$, $x>0$, $f^{(-1)}(x)$ be the inverse function of $f(x)$. In order for the series
\begin{equation}
\label{om3.1}
\sum_{k=1}^\infty \xi_k
\end{equation}
to converge in probability, it is sufficient that the series
\begin{equation}
\label{om3.2}
\sum_{k=1}^\infty \alpha^*_k,
\end{equation}
converges. Here $\alpha_k^*= V^{(-1)}\left( \frac{||\xi_k||^2}{f^{(-1)}(||\xi_k||^2)} \right)$.
In this case, for $$x\geq\mu = \sum_{k=1}^\infty V^{(-1)}\left( \frac{||\xi_k||^2}{f^{(-1)}(||\xi_k||^2)} \right)$$
the inequality holds
\begin{equation}
\label{om3.3}
P\left\{ \left| \sum_{k=1}^\infty \xi_k \right| \geq x \right\} \leq \sum_{k=1}^\infty W\left( \frac{||\xi_k||^2}{V(\frac{\alpha^*_kx}{\mu})} \right)<\infty.
\end{equation}

\end{theorem}
\begin{remark}
\label{omr3.1} The function $x/f^{(-1)}(x)$ is monotonically increasing.
This follows from the fact that the function $f(x)/x =
V(W(x))$ is monotonically increasing.
\end{remark}

\begin{dov} Let us prove that the convergence of the series (\ref{om3.2}) implies the convergence of the series on the right-hand side (\ref{om3.3}). It is obvious that for $x\geq\mu$ the series (\ref{om3.3}) converges when the series

\begin{equation}
\label{om3.4}
\sum_{k=1}^\infty W\left( \frac{||\xi_k||^2}{V(\alpha^*_k)} \right),
\end{equation}
converges, since the functions $W$ and $V$ are monotonically increasing.

Let us now prove the following equality. Given $||\xi_k||^2\neq 0$

\begin{equation}
\label{om3.5} W\left( \frac{||\xi_k||^2}{V(\alpha^*_k)} \right) =
V^{(-1)}\left( \frac{||\xi_k||^2}{f^{(-1)}(||\xi_k||^2)} \right).
\end{equation}

\noindent Really,

$$\frac{||\xi_k||^2}{f^{(-1)}(||\xi_k||^2)} = \frac{f(f^{(-1)}(||\xi_k||^2))}{f^{(-1)}(||\xi_k||^2)}=$$

$$=\frac{f^{(-1)}(||\xi_k||^2) V(W(f^{(-1)}(||\xi_k||^2))) }{f^{(-1)}(||\xi_k||^2)} = $$

 $$=V\left( W \left( \frac{||\xi_k||^2f^{(-1)}(||\xi_k||^2)}{||\xi_k||^2} \right) \right) = $$

$$ =  V\left( W \left( \frac{||\xi_k||^2}{V(V^{(-1)}(||\xi_k||^2/f^{(-1)}(||\xi_k||^2)))} \right) \right) =$$

$$ = V\left( W \left( \frac{||\xi_k||^2}{V(\alpha^*_k)} \right) \right).$$

\noindent From the last equality we obtain \eqref{om3.5}.

Now consider the probability $P\{|\sum_{k=l}^m\xi_k|>x\}$. Let us fix $x>0$. Then

$$P\left\{\left| \sum_{k=l}^m \xi_k \right|>x\right\}\leq \sum_{k=l}^m P\left\{ |\xi_k|>\alpha_kx \right\},$$
where $\sum_{k=l}^m \alpha_k = 1$, $\alpha_k>0$.
Since, according to Theorem \ref{omt2.1}, for all $\xi\in D_{V,W}(\Omega)$ the inequality
$$P\{|\xi|>x\}\leq W\left(\frac{||\xi||^2}{V(x)}\right)$$ holds, then

\begin{equation}
\label{omdvw-conv-t0}
\sum_{k=l}^m P\{ |\xi_k|>\alpha_kx \}\leq \sum_{k=l}^m W\left( \frac{||\xi_k||^2}{V(\alpha_kx)} \right).
\end{equation}

\noindent Putting $\alpha_k = \alpha^*_k/\mu_{lm}$, where
$\mu_{lm} = \sum_{k=l}^m V^{(-1)}(||\xi_k||^2/f^{(-1)}(||\xi_k||^2))$,
we will get

\begin{equation}
\label{omdvw-conv-t1}
\sum_{k=l}^m W\left( \frac{||\xi_k||^2}{V(\alpha_k x)} \right) = \sum_{k=l}^m W\left( \frac{||\xi_k||^2}{V(\alpha^*_k x / \mu_{lm})} \right).
\end{equation}

Since the series (\ref{om3.4}) converges, it follows from (\ref{om3.5}) that
$\mu_{lm}<x$ for sufficiently large $l,m$. Then
$$P\{ |\sum_{k=l}^m\xi_k|>x \}\leq
\sum_{k=l}^m W(||\xi_k||^2/V(\alpha^*_kx)).$$
From the convergence of the series (\ref{om3.4}) it follows that $ P\{ |\sum_{k=l}^m \xi_k|>x \}\to 0$ as $l,m\to\infty$,
i.e. the series (\ref{om3.1}) converges in probability.

The inequality (\ref{om3.3}) follows from (\ref{omdvw-conv-t0}) and (\ref{omdvw-conv-t1}),
when we set $l:=1$ and direct m to infinity.
\end{dov}

Let us find the conditions for convergence of series of random variables from the spaces $D_{V,W}(\Omega)$, when $V(x)=|x|^b$, $W(x)=|x|^a$, $0<b<1$, $a>1$.
The following statement holds.

\begin{theorem}
\label{omt3.2}
Let $W(x) = |x|^a$ and $V(x)=|x|^b$, $a>0,b>0$.
Then the series (\ref{om3.1}) converges in probability when the series converges

$$\mu = \sum_{k=1}^\infty ||\xi_k||^{\frac{2a}{ab+1}},$$
and for $x\geq\mu$ the inequality holds

$$P\{|\sum_{k=1}^\infty \xi_k|>x\}\leq \frac{1}{x^{ab}}\left( \sum_{k=1}^\infty ||\xi_k||^{\frac{2a}{ab+1}} \right)^{ab+1},$$
i.e. $\sum_{k=1}^\infty \xi_k$ belongs to the space $D_{V,W}$ and

$$||\sum_{k=1}^\infty \xi_k ||\leq \left( \sum_{k=1}^\infty ||\xi_k||^{\frac{2a}{ab+1}} \right)^{\frac{ab+1}{2a}}.$$

\end{theorem}

\begin{dov}$$P\{|\sum_{k=1}^\infty\xi_k|>x\}\leq \sum_{k=1}^\infty P\{|\xi_k|>\alpha_k x\},$$
where $\sum_{k=1}^\infty \alpha_k=1$.

Since, according to Theorem \ref{omt3.1}, for all $\xi\in D_{V,W}(\Omega)$ the inequality holds
$P\{|\xi|>x\}\leq W(\frac{||\xi||^2}{V(x)})$, then

$$\sum_{k=1}^\infty P\{ |\xi_k|>\alpha_kx \}\leq \sum_{k=1}^\infty W\left( \frac{||\xi_k||^2}{V(\alpha_kx)} \right).$$

\noindent Substituting the explicit form of $W$ and $V$, we obtain
$$\sum_{k=1}^\infty W\left( \frac{||\xi_k||^2}{V(\alpha_kx)} \right) = \frac{1}{x^{ab}}\sum_{k=1}^\infty\frac{||\xi_k||^{2a}}{\alpha_k^{ab}}.$$

\noindent The minimum of the last expression is achieved at the following $\alpha_k$:

$$\alpha_k = \frac{||\xi_k||^{\frac{2a}{ab+1}}}{\sum_{i=1}^\infty ||\xi_i||^{\frac{2a}{ab+1}}}.$$

\noindent For such values of $\alpha_k$

$$\frac{1}{x^{ab}}\sum_{k=1}^\infty\frac{||\xi_k||^{2a}}{\alpha_k^{ab}} = \frac{1}{x^{ab}}\sum_{k=1}^\infty\frac{||\xi_k||^{2a}}{\frac{||\xi_k||^{\frac{2a^2b}{ab+1}}}{( \sum_{i=1}^\infty ||\xi_i||^{\frac{2a}{ab+1}})^{ab}}} = $$

$$ = \frac{1}{x^{ab}}\left(\sum_{k=1}^\infty ||\xi_k||^{\frac{2a}{ab+1}} \right)^{ab}\sum_{i=1}^\infty ||\xi_i||^{\frac{2a}{ab+1}} = \frac{1}{x^{ab}}\left(\sum_{k=1}^\infty ||\xi_k||^{\frac{2a}{ab+1}} \right)^{ab+1},$$ in other words, $$P\left\{ \left| \sum_{k=1}^\infty \xi_k \right|>x \right\}\leq \frac{1}{x^{ab}}\left( \sum_{k=1}^\infty ||\xi_k||^\frac{2a}{ab+1} \right)^{ab+1}.$$
The statement of the Theorem follows from the last inequality.

\end{dov}

Now substitute the expressions for the functions $W$ and $V$ into (\ref{om3.2}):

$$\sum_{k=1}^\infty V^{(-1)}\left( \frac{||\xi_k||^2}{f^{(-1)}(||\xi_k||^2)} \right) = \sum_{k=1}^\infty \frac{||\xi_k||^{2/b}}{||\xi_k||^{2/(b(ab+1))}}  =\sum_{k=1}^\infty ||\xi_k||^{\frac{2a}{ab+1}}<\infty.$$

As can be seen, the results coincide, i.e. the conditions of Theorem \ref{omt3.1} cannot be improved.

In the case when the prenorm $||\cdot||$ is a quasi-norm, the following obvious statement follows from Theorem \ref{omt2.1}.

\begin{theorem}
\label{omt3.3} Let the random variables $\xi_k\in D_{V,W}$,
$W$ and $V$ satisfy condition $B1$, and
the series
$$\sum_{k=1}^\infty ||\xi_k||$$ converges.
Then the series
$$\sum_{k=1}^\infty \xi_k$$
converses in probability, and its sum belongs to the space
$D_{V,W}$. In this case, the inequality holds

$$P\left\{ \left|\sum_{k=1}^\infty\xi_k\right|\geq x \right\}\leq W\left(\frac{(\sum_{k=1}^\infty||\xi_k||)^2}{V(x)}\right).$$
\end{theorem}

\begin{remark}
\label{omr3.2} In the case where $W(x)=|x|^a$, $a\geq1$, $V(x)=|x|^b$, $0<b\leq 1$, condition B1 is satisfied for the functions $V$ and $W$,
then $\frac{2a}{ab+1}\geq1$, that is, if the condition
of Theorem \ref{omt3.2} is satisfied, the condition of Theorem \ref{omt3.3} may not be satisfied.
\end{remark}

\section{Conditions for convergence of infinite series of random variables with given distributions in spaces $D_{V,W}(\Omega)$}

Let for the space $D_{V,W}(\Omega)$ the functions $V(x)=|x|^b$, $W(x)=|x|^a$, $a>0$, $b>0$. Let $\xi_k$ be random variables from the space $D_{V,W}(\Omega) $.
Then, as it is known from Theorem \ref{omt3.2},  the convergence condition of the series has the form
\begin{equation}
\label{om5.3}
\sum_{k=1}^\infty ||\xi_k||^{\frac{2a}{ab+1}}<\infty.
\end{equation}
Let $\xi_k=a_k\hat{\xi}_k$, where $\hat{\xi}_k$ are independent identically distributed symmetric random variables.

\begin{example}
Let $P\{|\hat{\xi}_k|>x\}=\frac{1}{x^c+1}$, $x>0$, $c>0$. Then

$$||\xi_k|| = \left( \sup_{x>0}\frac{x^b}{(x^c/a_k^c+1)^{1/a}} \right)^{1/2}$$

\noindent and (\ref{om5.3}) can be rewritten as
\begin{equation}
\label{om5.4}
\sum_{k=1}^\infty ||\xi_k||^{\frac{2a}{ab+1}} = \sum_{k=1}^\infty \left(\sup_{x>0}\frac{x^b}{(x^c/a_k^c+1)^{1/a}}\right)^{\frac{a}{ab+1}}<\infty
\end{equation}
For $b> c/a$ the supremum of the function $\frac{x^b}{(x^c+a_k^c)^{1/a}}$ is infinite. For $b< c/a$ this function has an extremum at the point
$$x = a_k\left( \frac{ab}{c-ab} \right)^{1/c}.$$
\noindent Substituting this value into (\ref{om5.4}), for $b<c/a$ we have
$$\sum_{k=1}^\infty \left( a_k^{c/a}\sup_{x>0}\frac{x^b}{(x^c+a_k^c)^{1/a}}\right)^{\frac{a}{ab+1}} = \sum_{k=1}^\infty \left( a_k^b\frac{(ab)^{b/c}}{c^{1/a}}(c-ab)^\frac{c-ab}{ac} \right)^{\frac{a}{ab+1}}<\infty$$
\noindent If $b=c/a$, then the extremum will be reached at $+\infty$ and will be equal to $a_k^b$.
Therefore, the condition for convergence of the series at $b<c/a$ will be
$$\sum_{k=1}^\infty ||\xi_k||^\frac{2a}{ab+1} =  \left( \frac{(ab)^{b/c}}{c^{1/a}}(c-ab)^\frac{c-ab}{ac} \right)^{\frac{a}{ab+1}}\sum_{k=1}^\infty a_k^{\frac{ab}{ab+1}} <\infty,$$
\noindent that is
$$\sum_{k=1}^\infty a_k^{\frac{ab}{ab+1}} <\infty,$$
\noindent The estimate for the distribution function of the sum of the series will be of the form
$$P\{|\sum_{k=1}^\infty \xi_k|>x\}\leq \frac{1}{x^{ab}}{R}_c\left( \sum_{k=1}^\infty  a_k^\frac{ab}{ab+1} \right)^{ab+1},$$
\noindent where ${R}_c = \left( \frac{(ab)^{ab/c}}{c}(c-ab)^\frac{c-ab}{c}\right)$ for $b< c/a$.
It is easy to see that ${R}_c\to 1$, $c\to ab$, since $\frac{(ab)^{ab/c}}{c}\to 1$, $c\to ab$, and
$(c-ab)^\frac{c-ab}{c}\to 1$, $c\to ab$.
Therefore, for $b=c/a$ the estimate of the distribution function of the sum of the series will be
$$P\{|\sum_{k=1}^\infty \xi_k|>x\}\leq \frac{1}{x^c}\left( \sum_{k=1}^\infty  a_k^\frac{c}{c+1} \right)^{c+1}.$$
\end{example}

\begin{example}
Let $\hat{\xi}_k$ have a normal distribution, $\xi_k=a_k\hat{\xi}_k$.
Then
$$P\{|\xi_k|>x\}=2-\frac{\sqrt{2}}{\sigma\sqrt{\pi}}\int_{-\infty}^x \exp\left(-\frac{u^2}{2\sigma^2a_k^2}\right)du.$$

\noindent Substituting into the expression for the prenorm, we get

$$||\xi_k|| = \left( \sup_{x>0} x^b \left(2-\frac{\sqrt{2}}{\sigma\sqrt{\pi}}\int_{-\infty}^x \exp\left(-\frac{u^2}{2\sigma^2a_k^2}\right)du\right)^{1/a} \right)^{1/2}\leq$$

$$\leq \sup_{x>0}\left(x^b\left( \frac{\sqrt{2}a_k\sigma}{x\sqrt{\pi}}e^{x^2/a_k^2\sigma^2} \right)^{1/a}\right)^{1/2}$$

\noindent The supremum is reached at the point $x=a_k\sigma\sqrt{(ab-1)/2}$, and the desired prenorm takes the form

$$||\xi_k|| = (\sigma^2(ab-1))^{b/2}\left( \frac{e^{ab-1}2^{2-ab}}{\pi(ab-1)} \right)a_k^b.$$
Therefore, the convergence condition can be written as
$$\sum_{k=1}^\infty ||\xi_k||^\frac{2a}{ab+1} =  \left((\sigma^2(ab-1))^{ab}\left( \frac{e^{ab-1}2^{2-ab}}{\pi (ab-1)} \right)\right)^{1/2ab+2}\sum_{k=1}^\infty a_k^{\frac{ab}{ab+1}} <\infty,$$

\noindent that is,
$$\sum_{k=1}^\infty a_k^{\frac{ab}{ab+1}} <\infty.$$

The estimate for the distribution function of the sum of the series will be of the form

$$P\{|\sum_{k=1}^\infty \xi_k|>x\}\leq \frac{1}{x^{ab}}{N}_c\left( \sum_{k=1}^\infty  a_k^\frac{ab}{ab+1} \right)^{ab+1},$$

\noindent where $${N}_c=\left((\sigma^2(ab-1))^{ab}\left(
\frac{e^{ab-1}2^{2-ab}}{\pi (ab-1)} \right)\right)^{1/2}.$$
\end{example}

\begin{example}
Let the quantity $\hat{\xi}_k$ have a Cauchy distribution, $\xi_k=a_k\hat{\xi}_k$. Then

$$P\{|\xi_k|>x\}=1-\frac{2}{\pi}\arctan\left(\frac{x}{\gamma a_k}\right)=\frac{2}{\pi}\arctan\left(\frac{\gamma a_k}{x}\right)$$

\noindent for $x>0$. Substituting $t=\frac{x}{\gamma a_k}$, we have

$$\frac{2}{\pi}\arctan\left(\frac{\gamma a_k}{x}\right)= \frac{2}{\pi}\arctan\frac{1}{t} = \frac{2}{\pi}\int_t^\infty \frac{du}{1+u^2}\leq$$

$$\leq \frac{2}{\pi}\int_{t}^\infty \frac{du}{u^2}= \frac{2}{\pi t}=\frac{2\gamma a_k}{\pi x},$$

\noindent and, therefore,

$$P\{|\xi_k|>x\}\leq \frac{2\gamma a_k}{\pi x}.$$

Since $P\{|\xi_k|>x\}\leq 1$ $\forall x\in R$, we need to consider two cases: when $2\gamma a_k/\pi x>1$ and $2\gamma a_k/\pi x\leq1$.
In the first case, we can put $P\{|\xi_k|>x\}=1$ provided that $x<2\gamma a_k/\pi x$. Then
$$||\xi_k||=\left(\sup_{x>0}V(x)W^{(-1)}(P\{|\xi_k|>x\})\right)^{1/2} = \left(\sup_{x>0}x^b1^{1/a}\right)^{1/2} = \sup_{x>0}x^{b/2}.$$

\noindent Since $x<2\gamma a_k/\pi x$, the supremum will be reached precisely at this point, and then
$$||\xi_k||=\left(\frac{2\gamma a_k}{\pi }\right)^{b/2}.$$

Let us now consider the second case, when $2\gamma a_k/\pi x\leq1$.
In this case, we will have that under the condition $x\geq 2\gamma a_k/\pi x$
$$||\xi_k||=\left(\sup_{x>0}x^{b}\left(\frac{2\gamma a_k}{\pi x}\right)^{1/a}\right)^{1/2}.$$

\noindent Since $b/a>0$ is always the case, such a function will be increasing, and its supremum will be infinite, and, therefore, the prenorm $\xi_k$ will be undefined.

Finally, we will have

$$||\xi_k||=\left(\frac{2\gamma a_k}{\pi }\right)^{b/2}.$$

\noindent The corresponding convergence condition is written as

$$\sum_{k=1}^\infty ||\xi_k||^\frac{2a}{ab+1} = \left(\frac{2\gamma}{\pi}\right)^{ab/(ab+1)}\sum_{k=1}^\infty a_k^{ab/(ab+1)}<\infty,$$

\noindent that is

$$\sum_{k=1}^\infty a_k^{ab/(ab+1)}<\infty.$$
The estimate for the distribution function of the sum of the series will be of the form

$$P\{|\sum_{k=1}^\infty \xi_k|>x\}\leq \left(\frac{2\gamma}{\pi x}\right)^\frac{ab}{ab+1}\left( \sum_{k=1}^\infty  a_k^\frac{ab}{ab+1} \right)^{ab+1}$$
provided that $x<2\gamma a_k/\pi$.
\end{example}

\section*{Conclusions to chapter \ref{omch:dvw1}}

In section \ref{omch:dvw1}, the spaces $D_{V,W}$ are introduced: pre-Banach spaces of random variables with a given prenorm.
Their basic properties are studied, and a majorizing characteristic of these spaces is found.
A Theorem on the conditions for convergence of infinite sums of random variables from the vspace $D_{V,W}$ is proved in general case, as well as for the cases
where $V$ and $W$ are power functions, and where $W$ is a $C$-function, and in $V$ there is a function inverse to the $C$-function.
Several examples of random variables from $D_{V,W}$ spaces are considered.



\chapter{Random processes from spaces $D_{V,W}(\Omega)$ and their modeling}
\label{omch:dvw2}

In this chapter we introduce the concept of a random process from $D_{V,W}(\Omega)$ spaces. The following properties of these processes are studied: distribution of supremum, sample continuity, uniform convergence. Based on the studied properties, the accuracy and reliability of building models of processes from $D_{V,W}(\Omega)$ in this space are estimated.

\section{Random processes from spaces $D_{V,W}(\Omega)$}

\begin{definition}
\label{omd3.1} We will say that a random process $X(t) = \left\{X(t), t\in T\right\}$ belongs to the space $D_{V,W}$ if  $X(t)\in D_{V,W}$ for each $t\in T$.
\end{definition}
Examples of random processes from the space $D_{V,W}$ are random processes that can be represented as series
\begin{equation}
\label{om3.6}
X(t) = \sum_{k=1}^\infty \xi_k \phi_k(t),\mbox{   }t \in T,
\end{equation}
where $\xi_k\in D_{V,W}$ and this series converges in the space $D_{V,W}$.

The conditions for convergence of the series (\ref{om3.6}) are given by the following theorem.
\begin{theorem}
\label{omt3.5} Let $W(x)=|x|^a$, $V(x)=|x|^b$, $a\geq1$,
$0<b\leq 1$. Then the series (\ref{om3.6}) converges in probability if
 converges the series
$$\mu = \sum_{k=1}^\infty \phi_k^{\frac{ab}{ab+1}}(t)||\xi_k||_{V,W}^{\frac{2a}{ab+1}}.$$
Furthermore, for all $x\geq\mu$,
$$P\{|\sum_{k=1}^\infty \phi_k(t)\xi_k|\geq x\}\leq \frac{1}{x^{ab}}\left( \sum_{k=1}^\infty \phi^{\frac{ab}{ab+1}}_k(t)||\xi_k||_{V,W}^{\frac{2a}{ab+1}} \right)^{ab+1},$$
that is $X(t)\in D_{V,W}$.
\end{theorem}

\begin{dov}
The statement of the Theorem follows from Theorem \ref{omt3.2} and Theorem \ref{omt2.1}.
\end{dov}

In the following, in this section, we will denote the prenorm $||\cdot||_{V,W}$ of the space $D_{V,W}$ by $||\cdot||$.

Let $X=\{X(t),t\in T\}$ be a random process from the space
$D_{V,W}$, and let $\rho_X(t,s)=||X(s)-X(t)||$ be a premetric,
generated by the process $X$.

We will say that for the process $X(t)$ condition A1 is satisfied if

{\bf A1)} $\sup_{t\in T}||X(t)||<\infty$.

We will say that for a process $X(t)$ condition A2 is satisfied if

{\bf A2)} The space $(T,\rho_X)$ is separable and the process $X$ is separable
on $(T,\rho_X)$.

Denote by $\varepsilon_0=\sup_{t,s\in T}\rho_X(t,s)$. From condition (A1)
and Theorem \ref{omt2.1} it follows that $\varepsilon_0<\infty$.
Denote by
$\varepsilon_k=\varepsilon_0 \theta^k$, $\theta\in (0,1)$. Denote by $N(\varepsilon)$ the metric massiveness of the space
$(T,\rho_X)$, i.e. the minimum number of closed balls of radius $\varepsilon$ that cover
$(T,\rho_X)$.

Denote by $S_n$ the smallest $\varepsilon_n$-grid of the set $T$ with respect to the pseudometric $\rho_X$, and set $S=\bigcup_{n=0}^{\infty} S_n$.
The set $S_0$ consists of only one point (any); denote it by $t_0$. The set $S$ is countable and everywhere dense in $T$ with respect to the pseudometric $\rho_X$.
Since the process $X$ is continuous in probability and separable, any countable everywhere dense set in $(T,\rho_X)$ is a separability set for the process $X$.
Then
$$\sup_{t\in T}|X(t)| = \sup_{t\in S}|X(t)|$$
almost surely.

\begin{definition}
\label{omalpha-proc}
A family of mappings $\alpha_k(t)$, $k=0,1,\ldots$ is called an $\alpha$-procedure if each point in $S$ is associated with one point $\alpha_k$ in $S_k$ such that
$\rho_X(t,\alpha_k(t))\leq \varepsilon_k$.
\end{definition}

The following Theorem gives us the conditions under which $\sup_{t\in T}X(t)$ is finite with probability 1, and estimates for the distribution of this supremum.

\begin{theorem}
\label{omt4.1} Let a process $X$ satisfy conditions A1 and A2. If the series
\begin{equation}
\label{omum4.1}
\sum_{n=1}^\infty V^{(-1)} \left( \frac{\kappa(N(\varepsilon_n))^2\varepsilon_{n-1}^2}{f^{(-1)}(\kappa(N(\varepsilon_n))^2\varepsilon_{n-1}^2)} \right),
\end{equation}
where $f=xV(W(x))$, converges, then the following inequality holds true
$$P\{\sup_{t\in T}|X(t)| \geq x\}\leq$$
\begin{equation}
\label{om4.0}  \leq
W\left(\frac{\inf_{t\in T}||X(t)||^2}{V(\psi_0
x)}\right)+\sum_{k=1}^\infty W\left(
\frac{\kappa(N(\varepsilon_k))^2\varepsilon_{k-1}^2}{V(\psi_k x)}
\right)
\end{equation}

\noindent for $x>\Psi$. Here
$$\psi_0 = \frac{1}{\Psi}V^{(-1)}\left(\frac{\inf_{t\in
T}||X(t)||^2}{f^{(-1)}(\inf_{t\in T}||X(t)||^2)}\right),$$
$$\psi_k=\frac{1}{\Psi}V^{(-1)}\left(\frac{\kappa(N(\varepsilon_k))^2\varepsilon_{k-1}^2}{f^{(-1)}(\kappa(N(\varepsilon_k))^2\varepsilon_{k-1}^2)}\right),$$
$$\Psi = V^{(-1)}\left(\frac{\inf_{t\in T}||X(t)||^2}{f^{(-1)}(\inf_{t\in T}||X(t)||^2)}\right)+\sum_{k=1}^\infty V^{(-1)}\left(\frac{\kappa(N(\varepsilon_k))^2\varepsilon_{k-1}^2}{f^{(-1)}(\kappa(N(\varepsilon_k))^2\varepsilon_{k-1}^2)}\right).$$
\end{theorem}
\begin{dov}

Using the $\alpha$-procedure for selecting points of the sets $S_n$,
we obtain a set of points
$t=t_m,t_{m-1}=\alpha_{m-1}(t_m),...,t_1=\alpha_1(t_2),t_0=\alpha_0(t_1)$,
such that $t_n\in S_n,n=0,1,\ldots,m$, and $S_0=t_0$. Since

$$X(t)=X(t_0)+\sum_{n=1}^m (X(t_n)-X(t_{n-1})),$$

\noindent we have
\begin{equation}
\label{om4.1} \sup_{t\in T} |X(t)|\leq |X(t_0)|+\sum_{k=1}^\infty
\max_{s\in S_k}|X(s)-X(\alpha_{k-1}(s))|,
\end{equation}
therefore,
$$P\{\sup_{t\in T}|X(t)| \geq x\}\leq P\{|X(t_0)|\geq \psi_0x\}+$$
$$+\sum_{k=1}^\infty P\{ \max_{s\in S_k} |X(s)-X(\alpha_{k-1}(s))|\geq \psi_kx \},$$
where $\psi_k$ are numbers such that $\psi_k>0$,
$\sum_{k=0}^\infty\psi_k=1$. Using the inequality from Theorem
\ref{omt2.1}, we have
$$P\{|X(t_0)|\geq \psi_0x\}+\sum_{k=1}^\infty P\{ \max_{s\in S_k} |X(s)-X(\alpha_{k-1}(s))|\geq \psi_kx \}\leq $$
\begin{equation}
\label{om4.2} \leq W\left(\frac{||X(t_0)||^2}{V(\psi_0
x)}\right)+\sum_{k=1}^\infty W\left(\frac{||\max_{s\in S_k}
|X(s)-X(\alpha_{k-1}(s))|||^2}{V(\psi_kx)}\right).
\end{equation}

Let us consider separately $||\max_{s\in S_k} |X(s)-X(\alpha_{k-1}(s))|||$. We have

$$||\max_{s\in S_k} |X(s)-X(\alpha_{k-1}(s))||| \leq \kappa(N(\varepsilon_k))\max_{s\in S_k} ||X(s)-X(\alpha_{k-1}(s))||,$$
and since
$\rho_X(s,\alpha_{k-1}(s)) =
||X(s)-X(\alpha_{k-1}(s))||\leq\varepsilon_{k-1}$, we have
$$\kappa(N(\varepsilon_k))\max_{s\in S_k} ||X(s)-X(\alpha_{k-1}(s))||\leq \kappa(N(\varepsilon_k))\varepsilon_{k-1}.$$
Substituting the obtained estimate into inequality (\ref{om4.2}), we obtain
$$P\{\sup_{t\in T}|X(t)| \geq x\}\leq W\left(\frac{||X(t_0)||^2}{V(\psi_0 x)}\right)+\sum_{k=1}^\infty W\left(\frac{\kappa(N(\varepsilon_k))^2\varepsilon_{k-1}^2}{V(\psi_kx)}\right).$$
Since $t_0$ can be chosen arbitrarily, then
\begin{equation}
\label{om4.x} P\{\sup_{t\in T}|X(t)| \geq x\}\leq
W\left(\frac{\inf_{t\in T}||X(t)||^2}{V(\psi_0
x)}\right)+$$ $$+\sum_{k=1}^\infty
W\left(\frac{\kappa(N(\varepsilon_k))^2\varepsilon_{k-1}^2}{V(\psi_kx)}\right).
\end{equation}

Just as in the proof of Theorem \ref{omt3.1}, it is easy to show that the series (\ref{om4.x}) converges if we choose $\psi_0$ and $\psi_k$ as follows:
$$\psi_0 = \frac{1}{\Psi}V^{(-1)}\left(\frac{\inf_{t\in T}||X(t)||^2}{f^{(-1)}(\inf_{t\in T}||X(t)||^2)}\right),$$
$$\psi_k = \frac{1}{\Psi}V^{(-1)}\left(\frac{\kappa(N(\varepsilon_k))^2\varepsilon_{k-1}^2}{f^{(-1)}(\mu(N(\varepsilon_k))^2\varepsilon_{k-1}^2)}\right),$$
\noindent where
$$\Psi = V^{(-1)}\left(\frac{\inf_{t\in T}||X(t)||^2}{f^{(-1)}(\inf_{t\in T}||X(t)||^2)}\right)+\sum_{k=1}^\infty V^{(-1)}\left(\frac{\kappa(N(\varepsilon_k))^2\varepsilon_{k-1}^2}{f^{(-1)}(\kappa(N(\varepsilon_k))^2\varepsilon_{k-1}^2)}\right)$$
and $f=xV(W(x))$.

\noindent Since in (\ref{om4.x}) the majorizing series is convergent and $P\{\sup_{t\in T}|X(t)|\geq x\}\to 0$ as $x\to\infty$, then $\sup_{t\in T}|X(t)|$
is bounded almost surely.
\end{dov}

\begin{theorem}
\label{omt4.2} Let the process $X=\{X(t),t\in T\}$ be such that $X\in D_{V,W}$, $W=|x|^a$, $V=|x|^b$, $a>0,b>0$; moreover, $X$ satisfies conditions A1 and A2.
If the condition
\begin{equation}
\label{om4.3}
\int_0^{\Delta_0p}(N(u^{\frac{ab+1}{2a}}))^{\frac{1}{ab+1}}du<\infty,
\end{equation}
\noindent holds true, where $p:=\theta^{\frac{2a}{ab+1}}$,
$\theta$ is any number, $0\leq\theta\leq 1$,
$\Delta_0:=\varepsilon_0^{\frac{2a}{ab+1}}$, \\ $\varepsilon_0=\sup_{t,s\in T}\rho_X(t,s)$, then $\sup_{t\in T}|X(t)|\in D_{V,W}$, and, furthermore,
$$P\{\sup_{t\in T}|X(t)|\geq x\}\leq \frac{1}{x^{ab}}\left( \inf_{t\in T}||X(t)||^{\frac{2a}{ab+1}}+\right.$$
\begin{equation}
\label{om4.3b}
+\left.\frac{1}{p(1-p)}\int_0^{\Delta_0p}(N(u^{\frac{ab+1}{2a}}))^{\frac{1}{ab+1}}du \right).
\end{equation}
\end{theorem}

\begin{dov} Given $W$ and $V$, the condition \eqref{omum4.1} will have the form
$$\sum_{k=1}^\infty \kappa(N(\varepsilon_k))^{\frac{2a}{ab+1}}\varepsilon_{k-1}^{\frac{2a}{ab+1}} = $$
$$ = \varepsilon_0^{\frac{2a}{ab+1}}\sum_{k=1}^\infty \kappa(N(\varepsilon_0\theta^k))^{\frac{2a}{ab+1}}(\theta^{k-1})^{\frac{2a}{ab+1}}<\infty,$$

\noindent since $\varepsilon_k=\theta^k\varepsilon_0$. If we put
$\Delta_0=\varepsilon_0^{\frac{2a}{ab+1}}$,
$p=\theta^{\frac{2a}{ab+1}}$, then we get
$$ \varepsilon_0^{\frac{2a}{ab+1}}\sum_{k=1}^\infty \kappa(N(\varepsilon_0\theta^k))^{\frac{2a}{ab+1}}(\theta^{k-1})^{\frac{2a}{ab+1}} =$$
$$ = \sum_{k=1}^\infty (\kappa(N((\Delta_0p^k)^{\frac{ab+1}{2a}})))^\frac{2a}{ab+1}\Delta_0p^{k-1}.$$

\noindent Therefore
$$\sum_{k=1}^\infty (\kappa(N((\Delta_0p^k)^{\frac{ab+1}{2a}})))^\frac{2a}{ab+1}\Delta_0p^{k-1} \leq$$
$$\leq \sum_{k=1}^\infty \int_{\Delta_0p^{k+1}}^{\Delta_0p^k}(\kappa(N(u^{\frac{ab+1}{2a}})))^\frac{2a}{ab+1}du \frac{1}{p(1-p)} = $$
$$ = \frac{1}{p(1-p)} \int_0^{\Delta_0p}(\kappa(N(u^{\frac{ab+1}{2a}})))^\frac{2a}{ab+1}du.$$

\noindent If this integral is convergent, then $\sup_{t\in T}|X(t)|\in D_{V,W}$. This follows from (\ref{om4.3b}).

The statement of the Theorem follows from the fact that for the data $V,W$

$$\kappa(n) = \sup_{0<t<1/n} \left( \frac{W^{(-1)}(tn)}{W^{(-1)}(t)} \right)^{1/2} = n^{1/2a}.$$
\end{dov}

\begin{theorem}
\label{omt4.3}
Let the process $X=\{X(t),t\in [0,T]\}$ be such that $X\in D_{V,W}$, $W(x)=|x|^a$, $V(x)=|x|^b$, $a>0$, $b>0$, the process $X$ satisfies condition (A1), and $X$ is separable on $[0,T]$.
Let
$$\sup_{|t-s|\leq h}||X(t)-X(s)||\leq Dh^\zeta=\delta(h),$$
$D>0$, $\zeta>\frac{1}{2a}$.
Then $\sup_{t\in [0,T]}|X(t)|\in
D_{V,W}$, and, moreover, for any $x>0$
$$P\{\sup_{t\in [0,T]}|X(t)|>x\}\leq$$
$$\leq\frac{1}{x^{ab}}\left( \inf_{t\in [0,T]}||X(t)||^{\frac{2a}{ab+1}}+\int_0^{\Delta_0p} \left(\frac{TD^{1/\zeta}}{2u^{\frac{ab+1}{2a\zeta}}}+1\right)^{\frac{1}{ab+1}}du\right),$$
$\Delta_0$ and $p$ are given in Theorem \ref{omt4.2}.

\end{theorem}

\begin{dov}
From the previous Theorem $\sup_{t\in T}|X(t)|\in D_{V,W}$, if the integral (\ref{om4.3}) is finite. Under the conditions of this Theorem,

$$N(\varepsilon)\leq \frac{T}{2\delta^{(-1)}(\varepsilon)}+1.$$

\noindent Then the condition (\ref{om4.3}) can be rewritten as

$$\int_0^{\Delta_0p} \left( \frac{T}{2\delta^{(-1)}(u^{(ab+1)/2a})} +1\right)^{1/(ab+1)}du<\infty,$$
\noindent or, substituting $\delta(h)$,

$$\int_0^{\Delta_0p} \left( \frac{D^{1/\zeta}T}{2u^{(ab+1)/2a\zeta}} +1\right)^{1/(ab+1)}du<\infty.$$
For this integral to be finite, it is sufficient that the integral be finite

$$\int_0^{\Delta_0p} \frac{1}{u^{1/2a\zeta}}du.$$

\noindent This is done when $\zeta>\frac{1}{2a}$.
\end{dov}

\section{Continuity of processes from the space $D_{V,W}(\Omega)$}

Let $X = \{X(t), t\in T\}$ be a random process from the space $D_{V,W}$ such that $\sup_{t\in T} ||X(t)||<\infty$.
Let $\rho_X$ be a quasimetric generated by the process $X$, let the space $(T,\rho_X)$ be separable, and let the process $X(t)$ be separable on $(T,\rho_X)$.
Let $\theta\in(0,1)$ and $\varepsilon_l=\varepsilon_0\theta^l, l\geq1$, $\varepsilon_0=\sup_{t,s\in T}||X(t)-X(s)||$.
Denote by $V_{\varepsilon_k}$ the set of centers of closed balls of radius $\varepsilon_k$ that form a minimal cover of the space $(T,\rho_X)$.
The number of points of the set $V_{\varepsilon_k}$ is equal to $N(\varepsilon_k)$.
Assume that $N(\varepsilon)<\infty$ $\forall \varepsilon>0$.
Let $t,s\in T$ be points such that $\rho_X(t,s)<\varepsilon$, $0<\varepsilon\leq\varepsilon_0$.

Choose $k$ such that $\varepsilon_k<\varepsilon\leq\varepsilon_{k-1}$.
The set $V_k=\cup_{j=k}^{\infty}V_{\varepsilon_j}$ is the separability set of the process $X(t)$, because $X(t)$ is continuous in probability in the space $(T,\rho_X)$.

Let $S_n$ denote the smallest $\varepsilon_n$-grid of the set T with respect to the pseudometric $\rho_X$, and let $S=\bigcup_{n=0}^{\infty} S_n$.

\begin{theorem}
\label{omt1}
Let the process $X$ satisfy the above conditions and let the series
$$\sum_{l=k}^\infty V^{(-1)}\left( \frac{(\kappa((N(\varepsilon_l))^2))^2\varepsilon_{l-1}^2}{f^{(-1)}((\kappa(N(\varepsilon_l))^2)^2\varepsilon_{l-1}^2)} \right),$$
\noindent  be convergent, and let $x\geq\Psi$, where
$$\Psi = V^{(-1)}\left(\frac{(\kappa((N(\varepsilon_k))^2))^2\hat{\varepsilon}^2}{f^{(-1)}((\kappa(N^2(\varepsilon_k))^2)^2\hat{\varepsilon}^2)}\right) + $$
$$+2\sum_{l=k}^\infty  V^{(-1)}\left(\frac{(\kappa((N(\varepsilon_l))^2))^2\varepsilon_{l-1}^2}{f^{(-1)}((\kappa(N(\varepsilon_l))^2)^2\varepsilon_{l-1}^2)}\right).$$
\noindent Then
$$P\{\sup_{\rho_X(t,s)<\varepsilon}|X(t)-X(s)|\geq x\}\leq $$
$$\leq W\left(\frac{(\kappa((N(\varepsilon_k))^2))^2\varepsilon_{k-1}^2}{V(\psi_0 x)}\right) + 2\sum_{l=k}^\infty W\left(\frac{(\kappa((N(\varepsilon_l))^2))^2\hat{\varepsilon}^2}{V(\psi_lx)}\right),$$

\noindent where

$$\psi_0 = \frac{1}{\Psi}V^{(-1)}\left(\frac{(\kappa((N(\varepsilon_k))^2))^2\hat{\varepsilon}^2}{f^{(-1)}((\kappa(N^2(\varepsilon_k))^2)^2\hat{\varepsilon}^2)}\right),$$
$$\psi_l = \frac{1}{\Psi}V^{(-1)}\left(\frac{(\kappa((N(\varepsilon_l))^2))^2\varepsilon_{l-1}^2}{f^{(-1)}((\kappa(N(\varepsilon_l))^2)^2\varepsilon_{l-1}^2)}\right),$$
$$\hat{\varepsilon}=\varepsilon_k\frac{5-3\theta}{1-\theta}.$$
\noindent In this case, the process $X(t)$ is sample continuous in the space $(T,\rho_X)$.
\end{theorem}

\begin{dov}
Let $t$ and $s$ be arbitrary points in $V_k$ such that $\rho(t,s)\leq\varepsilon$.
It is clear that there exist points $m$ and $m_1$ such that $t\in V_{\varepsilon_m}$, $s\in V_{\varepsilon_{m_1}}$, $m>k$, $m_1>k$.
Let
$$t_m=\alpha_m(t),t_{m-1}=\alpha_{m-1}(t_m),\ldots,t_k=\alpha_k(t_{k+1}),$$
$$s_{m_1}=\alpha_{m_1}(t),s_{{m_1}-1}=\alpha_{{m_1}-1}(s_{m_1}),\ldots,s_k=\alpha_k(s_{k+1}),$$
\noindent where $\alpha_k(t)$ is the alpha procedure. Then
$$X(t)-X(s)=\sum_{l=k}^{m-1}(X(t_{l+1})-X(t_l))+$$
\begin{equation}
\label{omq1}
+\sum_{l=k}^{m-1}(X(s_{l+1})-X(s_l))+(X(t_k)-X(s_k)).
\end{equation}

From the last equality it follows that
$$||X(t_k)-X(s_k)||\leq ||X(t)-X(\alpha_m(t))||+||X(s)-X(\alpha_m(s))||+$$ $$+\sum_{l=k}^{m-1}||X(t_{l+1})-X(t_l)||+\sum_{l=k}^{m-1}||X(s_{l+1})-X(s_l)||+||X(t)-X(s)||\leq$$
$$\leq 2\sum_{l=k}^{m-1}\max_{u\in V_l}||X(u)-X(\alpha_l(u))||+||X(t)-X(\alpha_k(t))||+$$
$$+||X(s)-X(\alpha_k(s))||+||X(t)-X(s)||\leq2\sum_{l=k}^{m}\varepsilon_l+2\varepsilon_k+\varepsilon\leq$$
$$\leq \varepsilon_k\frac{5-3\theta}{1-\theta}:=\hat{\varepsilon}.$$

From (\ref{omq1}), directing $m$ to infinity, we obtain
$$\sup_{\rho_X(t,s)\leq \varepsilon}|X(t)-X(s)| = \sup_{|t-s|\leq \varepsilon, t,s,\in V_k}|X(t)-X(s)|\leq$$
$$\leq \max_{u,v\in V_k}|X(u)-X(v)| + 2 \sum_{l=k}^\infty \max_{u\in V_{l+1}}|X(u)-X(\alpha_l(u))|.$$

\noindent Therefore, for $\psi_0>0$, $\psi_l>0$, such that $\psi_0+2\sum_{l=k}^\infty \psi_l\leq 1$, the inequality holds
$$P\{\sup_{\rho_X(t,s)\leq \varepsilon}|X(t)-X(s)|\geq x\}\leq$$
$$\leq P\{\max_{t_k,s_k\in V_k}|X(t_k)-X(s_k)|\geq \psi_0 x\} + $$
$$+2\sum_{l=k}^\infty P\{\max_{u\in V_{l+1}}|X(u)-X(\alpha_l(u))|\geq\psi_lx\}.$$
Similarly to the proof of Theorem \ref{omt4.1}, we obtain
$$P\{\sup_{\rho_X(t,s)\leq \varepsilon}|X(t)-X(s)|\geq x\}\leq$$
$$\leq W\left(\frac{(\kappa((N(\varepsilon_k))^2))^2\hat{\varepsilon}^2}{V(\psi_0 x)}\right) + 2\sum_{l=k}^\infty W\left(\frac{(\kappa((N(\varepsilon_l))^2))^2\varepsilon_{l-1}^2}{V(\psi_lx)}\right),$$
\noindent where
$$\psi_0 = \frac{1}{\Psi}V^{(-1)}\left(\frac{(\kappa((N(\varepsilon_k))^2))^2\hat{\varepsilon}^2}{f^{(-1)}((\kappa(N(\varepsilon_k))^2))^2\hat{\varepsilon}^2)}\right),$$
$$\psi_l = \frac{1}{\Psi}V^{(-1)}\left(\frac{\kappa(N(\varepsilon_l))^2\varepsilon_{l-1}^2}{f^{(-1)}((\kappa(N(\varepsilon_l))^2)^2\varepsilon_{l-1}^2}\right),$$
$$\Psi = V^{(-1)}\left(\frac{(\kappa((N(\varepsilon_k))^2))^2\hat{\varepsilon}^2}{f^{(-1)}((\kappa(N(\varepsilon_k))^2))^2\hat{\varepsilon}^2)}\right) + $$ $$+ \sum_{l=k}^\infty  V^{(-1)}\left(\frac{(\kappa(N(\varepsilon_l))^2)^2\varepsilon_{l-1}^2}{f^{(-1)}((\kappa(N(\varepsilon_l))^2)^2\varepsilon_{l-1}^2}\right).$$

The first statement of the Theorem is proved.

Since $W(x)$ monotonically increases for all $x>0$, then for a fixed $x$
$$W\left(\frac{(\kappa((N(\varepsilon_k))^2))^2\hat{\varepsilon}^2}{V(\psi_0 x)}\right)\to0,$$
\noindent and, since the series
$$\sum_{l=k}^\infty W\left(\frac{\kappa(N(\varepsilon_l))^2\varepsilon_{l-1}^2}{V(\psi_lx)}\right)$$
\noindent converges, then if we direct $k\to\infty$, we get
$$W\left(\frac{(\kappa((N(\varepsilon_k))^2))^2\hat{\varepsilon}^2}{V(\psi_0 x)}\right) + \sum_{l=k}^\infty W\left(\frac{\kappa(N(\varepsilon_l))^2\varepsilon_{l-1}^2}{V(\psi_lx)}\right)\to 0,$$
\noindent which automatically means
$$P\{\sup_{\rho_X(t,s)<\varepsilon}|X(t)-X(s)|\geq x\}\to 0, k\to\infty.$$
\noindent From this it is clear that the process is sample continuous on $(T,\rho_X)$.
\end{dov}

\begin{theorem}
Let $W(x)=x^a$, $V(x)=x^b$, $a>1$, $0<b<1$. Then, if
\begin{equation}
\label{omtd4.3}
\int_0^{\Delta_0p^{k}}N(u^{\frac{ab+1}{2a}})^{2/(ab+1)}du<\infty,
\end{equation}
\noindent then
$$P\{\sup_{\rho_X(t,s) <\varepsilon}|X(s)-X(t)|\geq x\}\leq$$
\begin{equation}
\label{omtd4.3-2}
\leq \frac{1}{x^{ab}p(1-p)}\left(2\int_0^{\Delta_0p^{k+1}}N(u^{\frac{ab+1}{2a}})^{2/(ab+1)}du+\right.$$
$$\left.+\left(\frac{5-3\theta}{(1-\theta)\theta}\right)^\frac{1}{ab+1}\int_{\Delta_0p^{k+1}}^{\Delta_0p^{k}}N(u^{\frac{ab+1}{2a}})^{2/(ab+1)}du\right),
\end{equation}
\noindent $\Delta_0$ and $p$ are given in Theorem \ref{omt4.1}, $\theta\in(0,1)$. In this case, the process $X(t)$ is sample continuous on $(T,\rho_X)$.
\end{theorem}

\begin{dov}

$$P\{\sup_{\rho_X(t,s) <\varepsilon}|X(s)-X(t)|\geq x\}\leq$$

$$\leq \frac{1}{x^{ab}p(1-p)}\left( (\kappa((N(\varepsilon_k))^2)\hat{\varepsilon})^{\frac{2a}{ab+1}} + 2\sum_{l=k+1}^\infty(\kappa((N(\varepsilon_l))^2)\varepsilon_{l-1})^{\frac{2a}{ab+1}}\right)\leq$$

$$\leq \frac{1}{x^{ab}p(1-p)}\left( \int_{\Delta_0p^{k+1}}^{\Delta_0p^k}\left(\kappa\left(\left(N(u^{\frac{ab+1}{2a}})\right)^2\frac{5-3\theta}{(1-\theta)\theta}\right)\right)^{\frac{2a}{ab+1}}du +\right.$$

$$\left.+ 2\int_{0}^{\Delta_0p^{k+1}}(\kappa((N(u^{\frac{ab+1}{2a}}))^2))^{\frac{2a}{ab+1}}du\right).$$

The last result was obtained using the proof of Theorem \ref{omt4.1}.

From Theorem \ref{omt2.2} we have that

$$\kappa(n) = \sup\limits_{0<t<1/n}\left( \frac{W^{(-1)}(tn)}{W^{(-1)}(t)}
\right)^{1/2} = n^{1/2a}.$$

\noindent Taking this expression into account, we get the following

$$P\{\sup_{\rho_X(t,s) <\varepsilon}|X(s)-X(t)|\geq x\}\leq$$

$$\leq \frac{1}{x^{ab}p(1-p)}\left( \left(\frac{5-3\theta}{(1-\theta)\theta}\right)^\frac{1}{ab+1}\int_{\Delta_0p^{k+1}}^{\Delta_0p^k}(N(u^{\frac{ab+1}{2a}}))^{2/(ab+1)}du + \right.$$

$$+\left.2\int_{0}^{\Delta_0p^{k+1}}N(u^{\frac{ab+1}{2a}})^{2/(ab+1)}du\right)\leq \frac{1}{x^{ab}p(1-p)}\times $$ $$\times\left(2\int_0^{\Delta_0p^{k+1}}N(u^{\frac{ab+1}{2a}})^{2/(ab+1)}du+\right.$$

$$\left.+\left(\frac{5-3\theta}{(1-\theta)\theta}\right)^\frac{1}{ab+1}\int_{\Delta_0p^{k+1}}^{\Delta_0p^{k}}N(u^{\frac{ab+1}{2a}})^{2/(ab+1)}du\right).$$

\end{dov}

\begin{theorem}
Let the process $X=\{X(t),t\in [0,T]\}$ be such that $X\in D_{V,W}$, $X$-separable on
$[0,T]$, $W(x)=|x|^a$, $V(x)=|x|^b$, $a\geq1$, $0<b\leq1$.
Let
$$\sup_{|t-s|\leq h}||X(t)-X(s)||\leq Dh^\zeta,$$ $D>0$, $\zeta>\frac{1}{2a}$.
Then $\sup_{t\in [0,T]}|X(t)|\in D_{V,W}$, and, furthermore, for any $x>0$

$$P\left\{\sup_{\rho_X(t,s) \leq \varepsilon}|X(s)-X(t)|\geq x\right\}\leq$$
$$\leq \frac{1}{x^{ab}p(1-p)}\left(\left(\frac{5-3\theta}{(1-\theta)\theta}\right)^\frac{1}{ab+1}\int_{\Delta_0p^{k+1}}^{\Delta_0p^{k}}\left( \frac{D^{1/\zeta}T}{2u^\frac{ab+1}{2a\zeta}} +1 \right)^{2/(ab+1)}du+\right.$$
$$+\left.2\int_0^{\Delta_0p^{k+1}}\left( \frac{D^{1/\zeta}T}{2u^\frac{ab+1}{2a\zeta}} +1 \right)^{2/(ab+1)}du\right),$$

\noindent $\Delta_0$ and $p$ are given in Theorem \ref{omt4.1}, $\theta\in(0,1)$.
 In this case, the process $X(t)$ is sample continuous on $(T,\rho_X)$.
\end{theorem}
\begin{dov}
Under the conditions of the Theorem,
$$N(\varepsilon)\leq \frac{D^{1/\zeta}T}{2\varepsilon^{1/\zeta}}+1,$$

\noindent Then the condition (\ref{omtd4.3-2}) can be rewritten as

$$P\{\sup_{\rho_X(t,s)<\varepsilon}|X(s)-X(t)|\geq x\}\leq$$
$$\leq\frac{1}{x^{ab}p(1-p)}\left(2\int_0^{\Delta_0p^{k+1}}N(u^{\frac{ab+1}{2a}})^{2/(ab+1)}du+\right.$$
$$\left.+\left(\frac{5-3\theta}{(1-\theta)\theta}\right)^\frac{1}{ab+1}\int_{\Delta_0p^{k+1}}^{\Delta_0p^{k}}N(u^{\frac{ab+1}{2a}})^{2/(ab+1)}du\right)\leq$$
$$\leq \frac{1}{x^{ab}p(1-p)}\left(\left(\frac{5-3\theta}{(1-\theta)\theta}\right)^\frac{1}{ab+1}\int_{\Delta_0p^{k+1}}^{\Delta_0p^{k}}\left( \frac{D^{1/\zeta}T}{2u^\frac{ab+1}{2a\zeta}} +1 \right)^{2/(ab+1)}du+\right.$$
$$+\left.2\int_0^{\Delta_0p^{k+1}}\left( \frac{D^{1/\zeta}T}{2u^\frac{ab+1}{2a\zeta}} +1 \right)^{2/(ab+1)}du\right).$$
\noindent The integrals in the last expression are convergent when converges the integral
$$\int_0^{\Delta_0p^{k+1}} \frac{1}{u^\frac{1}{a\zeta}}du,$$
\noindent which is achieved when $\zeta>\frac{1}{a}$.
The values of these integrals can be estimated using the hypergeometric function.
\end{dov}

\section{Uniform convergence of functional series in $D_{V,W}(\Omega)$ }

Let the process $X(t)$ be represented in the form \eqref{om3.6}, i.e.
$$X(t) = \sum_{k=1}^\infty \xi_k \phi_k(t),\mbox{ }t \in T,$$ where $\xi_k\in D_{V,W}$ and this series converges in the space $D_{V,W}$.
Denote
$$X_{N,M}(t):=\sum_{k=N+1}^M \xi_k\phi_k(t).$$

Let us consider the case when $V(x)=|x|^b$, $W(x)=|x|^a$, $a>0$, $b>0$, $T=[0,T]$ is a segment. Under these conditions, the metric massiveness will take the form:
$$N(\varepsilon)=\frac{D^{1/\zeta}T}{2\varepsilon^{1/\zeta}}+1,$$
\noindent when for all $N,M$ the following condition is satisfied
\begin{equation}
\label{omrivnep}
\sup_{|t-s|<h}(||X_{N,M}(s)-X_{N,M}(t)||)\leq D_{N,M}|h|^\zeta,
\end{equation}
\noindent $D_{N,M}<D$.

\begin{theorem}\label{omt1.5.2}
Let $V(x)=|x|^b$, $W(x)=|x|^a$, $a>0$, $b>0$, $T=[0,T]$ be a segment, let $||X_{N,M}(t)||\to0$  uniformly in $t\in[0,T]$ as $N,M\to\infty$ and let the condition (\ref{omrivnep}) be satisfied.
Then
$$P\{\sup_{t\in T}|X_{N,M}(t)|\geq x\}\to 0, N,M\to \infty,$$ i.e. $X_{N,M}$ is uniformly convergent in probability.
\end{theorem}
\begin{dov} From Theorem \ref{omt4.3} we have that
$$P\{\sup_{t\in [0,T]}|X_{N,M}(t)|>x\}\leq$$
$$\leq\frac{1}{x^{ab}}\left( \inf_{t\in [0,T]}||X_{N,M}(t)||^{\frac{2a}{ab+1}}+\int_0^{\Delta_{N,M}p} \left(\frac{D^{1/\zeta}T}{2u^{\frac{ab+1}{2a\zeta}}}+1\right)^{\frac{1}{(ab+1)}}du\right),$$
\noindent where $\Delta_{N,M}=\varepsilon_{N,M}^\frac{2a}{ab+1}$, $p=\theta^\frac{2a}{ab+1}$.

Let us denote the right-hand side of the inequality as $Z(N,M,x)$. Under the conditions of the Theorem, $||X_{N,M}(t)||\to0$, $N,M\to\infty$, whence $\inf_{t\in [0,T]}$ $||X_{N,M}(t)||\to0$, $N,M\to\infty$, and also $\Delta_{N,M}\to0$, $N,M\to\infty$. Hence we obtain that $Z(N,M,x)\to0$, $N,M\to\infty$.
\end{dov}

\begin{remark} By directing $M\to\infty$, we obtain that $$\sum\limits_{k=N+1}^\infty \xi_k\phi_k(t)\to 0,\quad N\to\infty,$$ in probability.
\end{remark}

\section{Models of stochastic processes from spaces $D_{V,W}(\Omega)$}

Let a stochastic process $X(t)$ be represented in the form
\begin{equation}
\label{omprocess}
X(t)=\sum_{k=1}^\infty \xi_k \phi_k(t),
\end{equation}
\noindent where $t \in [0,T]$. Let
$$X_{N}(t)=\sum_{k=1}^N\xi_k\phi_k(t).$$
\noindent The expression $X_N(t)$ will be called the model of the process $X$.

Let us denote \begin{equation}
\label{om5.5}
\tilde{X}_N(t):= \sum_{k=N+1}^\infty\xi_k\phi_k(t) = X(t)-X_N(t).
\end{equation}

\begin{theorem} \label{ommainmodel}
Let the process $X=\{X(t),t\in [0,T]\}$ (see (\ref{omprocess})) be such that $\xi_k\in D_{V,W}$, $X$ satisfies conditions A1 and A2; moreover, $W(x)=|x|^a$, $V(x)=|x|^b$, $a\geq1$, $0<b\leq1$.
Let the following conditions be satisfied
$$\sup_{|t-s|<h}|\phi_k(s)-\phi_k(t)|\leq C_k|h|^\zeta,$$
$$\sum_{k=1}^\infty ||\xi_k||^\frac{2a}{ab+1}C_k^\frac{ab}{ab+1}<\infty$$
\noindent and let
$$\zeta>\frac{1}{ab}.$$
\noindent Then $\sup_{t\in T}|\tilde{X}_{N}(t)|$ $\in D_{V,W}$, and, besides,
$$P\left\{\sup_{t\in [0,T]}|\tilde{X}_{N}(t)|>x\right\}\leq \frac{1}{x^{ab}}\left( \inf_{t\in [0,T]}||\tilde{X}_{N}(t)||^{\frac{2a}{ab+1}} +\right.$$
$$\left.+\frac{T^{1/(ab+1)}(\sum_{k=N+1}^\infty C_k^\frac{ab}{ab+1}||\xi_k||^\frac{2a}{ab+1})^\frac{1}{ab\zeta}}{p(1-p)2^{ab/(ab+1)}} \frac{ab\zeta(\Delta_Np)^{1-\frac{1}{ab\zeta}}}{ab\zeta-1}     +\frac{\Delta_N}{1-p}\right),$$
$\Delta_N = \Delta_{N,\infty}$, where $\Delta_{N,M}$ is specified in Theorem \ref{omt1.5.2}, $p=\theta^\frac{2a}{ab+1}$, $\theta\in(0,1)$.
\end{theorem}
\begin{dov}

Since $$\sup_{|t-s|<h}|\phi_k(t)-\phi_k(s)|\leq C_k  |h|^\zeta.$$
From Theorem \ref{omt2.1} and Theorem \ref{omt3.2} it follows that
$$\sup_{|t-s|<h}||\tilde{X}_{N}(t)-\tilde{X}_{N}(s)|| = \sup_{|t-s|<h}||\sum_{k=N+1}^\infty \xi_k(\phi_k(t)-\phi_k(s))||\leq $$
$$\leq \left(\sum_{k=N+1}^\infty||\xi_k||^\frac{2a}{ab+1}\sup_{|t-s|<h}J^\frac{2a}{ab+1}\left(\phi_k(t)-\phi_k(s)\right)\right)^\frac{ab+1}{2a}\leq$$
$$\leq \left(\sum_{k=1}^\infty ||\xi_k||^\frac{2a}{ab+1}\left(C_k^{b/2}h^{b\zeta/2}\right)^\frac{2a}{ab+1}\right)^\frac{ab+1}{2a} = h^{b\zeta/2}\left(\sum_{k=1}^\infty ||\xi_k||^\frac{2a}{ab+1}C_k^\frac{ab}{ab+1}\right)^\frac{ab+1}{2a},$$
\noindent since $J(z) = z^{b/2}$ and
$$||\sum_{k=N+1}^\infty\xi_k||\leq(\sum_{k=N+1}^\infty||\xi_k||^\frac{2a}{ab+1})^\frac{ab+1}{2a}.$$

Since by the condition of the Theorem $t\in [0,T]$, then $N(\varepsilon)\leq\frac{T}{2\delta^{(-1)}(h)}+1$, where $\delta(h)$ can be put
$$\delta(h)=h^{b\zeta/2}\left(\sum_{k=1}^\infty ||\xi_k||^\frac{2a}{ab+1}C_k^\frac{ab}{ab+1}\right)^\frac{ab+1}{2a}$$
\noindent provided that the series $$\sum_{k=N+1}^\infty ||\xi_k||^\frac{2a}{ab+1}C_k^\frac{ab}{ab+1}$$ converges.

By Theorem \ref{omt4.2}
$$P\{\sup_{t\in [0,T]}|\tilde{X}_{N}(t)|>x\} \leq \frac{1}{x^{ab}}\left( \inf_{t\in [0,T]}||\tilde{X}_{N}(t)||^{\frac{2a}{ab+1}}+\right.$$
$$+\left.\frac{1}{p(1-p)}\int_0^{\Delta_Np}(N(u^{\frac{ab+1}{2a}}))^{1/(ab+1)}du \right)\leq \frac{1}{x^{ab}}\left( \inf_{t\in [0,T]}||\tilde{X}_{N}(t)||^{\frac{2a}{ab+1}}+\right.$$
$$\left.+\frac{1}{p(1-p)}\int_0^{\Delta_Np} \left(\frac{T}{2\delta^{(-1)}(u^\frac{ab+1}{2a})}+1\right)^{1/(ab+1)}\right)\leq$$ $$\leq\frac{1}{x^{ab}}\left( \inf_{t\in [0,T]}||\tilde{X}_{N}(t)||^{\frac{2a}{ab+1}}+\right.$$
$$\left.+\frac{1}{p(1-p)}\int_0^{\Delta_Np} \left(\frac{T(\sum_{k=N+1}^\infty C_k^\frac{ab}{ab+1}||\xi_k||^\frac{2a}{ab+1})^\frac{ab+1}{ab\zeta}}{2u^\frac{ab+1}{ab\zeta}}+1\right)^{1/(ab+1)}\right)\leq$$
$$\leq \frac{1}{x^{ab}}\left( \inf_{t\in [0,T]}||\tilde{X}_{N}(t)||^{\frac{2a}{ab+1}}+\frac{T^{1/(ab+1)}(\sum_{k=N+1}^\infty C_k^\frac{ab}{ab+1}||\xi_k||^\frac{2a}{ab+1})^\frac{1}{ab\zeta}}{2^{1/(ab+1)}p(1-p)}\right.$$
$$\left.\int_0^{\Delta_Np} \frac{du}{u^\frac{1}{ab\zeta}}+\frac{\Delta_N}{1-p}\right)\leq \frac{1}{x^{ab}}\left( \inf_{t\in [0,T]}||\tilde{X}_{N}(t)||^{\frac{2a}{ab+1}} +\right.$$
$$\left.+\frac{T^{1/(ab+1)}(\sum_{k=N+1}^\infty C_k^\frac{ab}{ab+1}||\xi_k||^\frac{2a}{ab+1})^\frac{1}{ab\zeta}}{2^{ab/(ab+1)}p(1-p)} \frac{ab\zeta(\Delta_Np)^{1-\frac{1}{ab\zeta}}}{ab\zeta-1}     +\frac{\Delta_N}{1-p}\right),$$
\noindent provided that the integral
$$\int_0^{\Delta_Np} \frac{1}{u^\frac{1}{ab\zeta}}du$$
\noindent is finite, i.e. when $\zeta>\frac{1}{ab}$.

\end{dov}

\begin{corollary}
Let the process $X=\{X(t),t\in [0,T]\}$, which can be represented in the form (\ref{omprocess}), be such that
$\xi_k\in D_{V,W}$, $W(x)=|x|^a$, $V(x)=|x|^b$, $a\geq1$, $0<b\leq1$. Let, in addition, $X$ satisfy conditions A1 and A2.

The model $X_N(t)$ approximates the process $X(t)$ at $t\in [0,T]$ with given reliability $1-\nu$, $0<\nu<1$ and accuracy $\ae>0$ in the space $D_{V,W}(\Omega)$, i.e.
$$P\{\sup_{t\in [0,T]} |\widetilde{X}_N(t)|>\ae\}\leq\nu,$$
provided that $$\sup_{|t-s|<h}|\phi_k(s)-\phi_k(t)|\leq C_k|h|^\zeta,$$ if the following conditions are satisfied:
$$ \frac{1}{\ae^{ab}}\left( \inf_{t\in [0,T]}||\tilde{X}_{N}(t)||^{\frac{2a}{ab+1}} +\frac{\Delta_N}{1-p}+\right.$$
\begin{equation}
\label{ommodel_cond}
\left.+\frac{T^{1/(ab+1)}(\sum_{k=N+1}^\infty C_k^\frac{ab}{ab+1}||\xi_k||^\frac{2a}{ab+1})^\frac{1}{ab\zeta}}{2^{ab/(ab+1)}p(1-p)} \frac{ab\zeta(\Delta_Np)^{1-\frac{1}{ab\zeta}}}{ab\zeta-1}     \right)\leq\nu,
\end{equation}
$$\sum_{k=1}^\infty |C_k|^\frac{ab}{ab+1}||\xi_k||^\frac{2a}{ab+1}<\infty,$$
$$\zeta>\frac{1}{ab},$$
where $p=\theta^\frac{2a}{ab+1}$, $\theta$ is any number such that $0<\theta<1$,

$$\Delta_N=\left(\sup_{t,s\in T}||X(s)- X(t)||\right)^\frac{2a}{ab+1}.$$
\end{corollary}

\begin{remark} The left-hand side (\ref{ommodel_cond}) reaches its minimum when $\theta$ is a solution to the equation
$$\frac{ab\zeta T^{1/(ab+1)}(\sum_{k=N+1}^\infty C_k^\frac{ab}{ab+1}||\xi_k||^\frac{2a}{ab+1})^\frac{1}{ab\zeta}}{2^{ab/(ab+1)}(ab\zeta-1)} \Delta_N((ab+1)\theta^\frac{2a}{ab+1}-1)-$$
$$-\Delta_N^\frac{ab+1}{ab}\theta^\frac{2}{b}((\Delta_N+1)\theta^\frac{1}{ab+1}-\Delta_N-ab-1)=0.$$
\end{remark}

\section{Construction of models of stochastic processes from spaces $D_{V,W}$}

In this section we will consider processes of the form
$$X(t) = \sum_{k=1}^\infty \xi_k \varphi_k(t),$$
\noindent where $\xi_k\in D_{V,W}$, on the interval $[0,T]$.

\begin{example} {\bf Generalized Brownian motion.} Let $V(x)=|x|^b$, $W(x)=|x|^a$, $a>1$, $0<b<1$, and let the process $X(t)$ be represented in the form
$$X(t) = \sum_{k=1}^\infty \xi_k \sqrt{\frac{2}{\pi k}}\sin(\pi k t),$$
\noindent $\xi_k\in D_{V,W}$. In this case $\varphi_k(t)=\sqrt{\frac{2}{\pi k}}\sin(\pi k t).$ For this process
$$\sup_{|t-s|<h}|\varphi_k(t)-\varphi_k(s)| = \sup_{|t-s|<h}\left|\frac{\sqrt{2}(\sin(\pi k t)-\sin(\pi k s))}{\sqrt{\pi k}}\right|\leq$$
$$\leq 2/\pi k\left|\sin\left(\frac{\pi k h}{2} \right)\right|\leq (2/\pi k)^{1/2-\alpha} h^\alpha,$$
\noindent because $|\sin t| \leq t^\alpha$, $0<\alpha<1$. Therefore, $C_k = (2/\pi k)^{1/2-\alpha}$, $\zeta=\alpha$, and $\frac{1}{ab}<\alpha$. This is achieved
for $a\in (1/b\alpha,+\infty)$. For the same process
$$\inf_{t\in [0,T]}||X(t)||^\frac{ab}{ab+1}=0,$$
$$\Delta_N=\left(\sup_{t,s\in [0,T]} ||\tilde{X}_{N}(t)-\tilde{X}_{N}(s)||\right)^\frac{2a}{ab+1} \leq$$
 $$\leq\left(\sup_{t,s\in [0,T]}\left(\sum_{k=N+1}^\infty\left|\left|\sqrt{\frac{2}{\pi k}}(\xi_k (\sin (\pi k t)-\sin (\pi k s)))\right|\right|^\frac{2a}{ab+1}\right)^\frac{ab+1}{2a}\right)^\frac{2a}{ab+1}\leq$$
$$\leq (8/\pi k)^{3ab/4ab+4}\sum_{k=N+1}^\infty ||\xi_k||^\frac{2a}{ab+1}.$$

Having chosen the necessary values of accuracy $\alpha$, reliability $1-\nu$ and calculating the value of the constant $\theta$ at which (\ref{ommodel_cond}) is minimized, from the inequality

$$\nu\geq\frac{1}{\ae^{ab}}\left(\frac{ab\alpha T^{1/(ab+1)} \left(\theta^\frac{2a}{ab+1}(8/\pi k)^{3ab/4ab+4}\sum_{k=N+1}^\infty ||\xi_k||^\frac{2a}{ab+1}\right)^{1-\frac{1}{ab\alpha}}}{\theta^\frac{2a}{ab+1}(1-\theta^\frac{2a}{ab+1})(ab\alpha-1)2^{ab/(ab+1)}}\times\right.$$
$$\left.\times\left(\sum_{k=N+1}^\infty ||\xi_k||^\frac{2a}{ab+1}(2/\pi k)^\frac{ab(1/2-\alpha)}{ab+1}\right)^\frac{1}{ab\alpha} +\right.$$ $$\left.+ \frac{(8/\pi k)^{3ab/4ab+4}\sum_{k=N+1}^\infty ||\xi_k||^\frac{2a}{ab+1}}{1-\theta^\frac{2a}{ab+1}}\right)$$
\noindent under condition
$$\sum_{k=1}^\infty (2/\pi k)^\frac{ab(1/2-\alpha)}{ab+1}||\xi_k||^\frac{2a}{ab+1}<\infty$$
\noindent we find the value of $N$ we need.
\end{example}

\begin{example} Let $V(x)=|x|^b, W(x)=|x|^a$, $a>1$, $0<b<1$, and let the process $X(t)$ be represented in
the form
$$X(t) = \sum_{k=1}^\infty \xi_k A_k(\sin(B_k t)+\cos(B_kt)),$$
\noindent $\xi_k\in D_{V,W}$, $A_k>0, B_k>0$. In this case
$$\sup_{|t-s|<h}|A_k(\sin( B_k t) + \cos( B_k t))-A_k(\sin( B_k s)+\cos( B_k s))|=$$
$$=\sup_{|t-s|<h}\left|A_k\left(2\sin\left(\frac{B_k}{2}\left(t-s\right)\right)\cos\left(\frac{B_k}{2}\left(t+s\right)\right)+\right.\right.$$
$$\left.\left.+2\sin\left(\frac{B_k}{2}(t+s)\right)\sin\left(\frac{B_k}{2}(t-s)\right)\right)\right|\leq \sup_{|t-s|<h}4A_k\left|\sin\left(\frac{B_k}{2}(t-s)\right)\right|\leq$$
$$\leq 2^{2-\alpha}A_kB_k^\alpha h^\alpha,$$
\noindent because $\sin t \leq t^\alpha$, $0<\alpha<1$.
The condition for convergence of the integral is satisfied if $\frac{1}{ab}<\alpha$, i.e. when $a\in (\frac{1}{b\alpha},+\infty)$.
For the same process
$$\inf_{t\in [0,T]}|A_k(\sin(B_kt)+\cos(B_kt))| \leq |A_k|,$$
$$\Delta_N = \sup_{t,s\in [0,T]}\left(\sum_{k=N+1}^\infty ||\xi_k(A_k(\sin(B_kt)+\cos(B_kt))-A_k(\sin(B_ks)+\right.$$
$$+\left.\cos(B_ks)))\right)\leq\sum_{k=N+1}^\infty ||\xi_k 2\sqrt{2}A_k||^\frac{2a}{ab+1}=2^\frac{3ab}{4ab+4}A_k^\frac{ab}{ab+1}\sum_{k=N+1}^\infty ||\xi_k||^\frac{2a}{ab+1}.$$
Choosing reliability $1-\nu$, accuracy $\ae$ and calculating the constant $\theta$, we have
$$\nu\geq\frac{1}{\ae^{ab}}\left(\sum_{k=1}^\infty ||\xi_k||^\frac{2a}{ab+1}|A_k|^\frac{ab}{ab+1}+\left(\sum_{k=N+1}^\infty ||\xi_k||^\frac{2a}{ab+1}(2^{2-\alpha}A_kB_k^\alpha)^\frac{ab}{ab+1}\right)^\frac{1}{ab\alpha}\times\right.$$
$$\left.   \times\frac{ab\alpha T^{1/(ab+1)} \left(\theta^\frac{2a}{ab+1}2^\frac{3ab}{4ab+4}A_k^\frac{ab}{ab+1}\sum_{k=N+1}^\infty ||\xi_k||^\frac{2a}{ab+1}\right)^{1-\frac{1}{ab\alpha}}}{\theta^\frac{2a}{ab+1}(1-\theta^\frac{2a}{ab+1})(ab\alpha-1)2^{ab/(ab+1)}} +\right.$$
 $$\left.+\frac{2^\frac{3ab}{4ab+4}A_k^\frac{ab}{ab+1}\sum_{k=N+1}^\infty ||\xi_k||^\frac{2a}{ab+1}}{1-\theta^\frac{2a}{ab+1}}\right)$$
\noindent under condition
$$\sum_{k=1}^\infty ||\xi_k||^\frac{2a}{ab+1}(2^{2-\alpha}A_kB_k^\alpha)^\frac{ab}{ab+1}<\infty$$
\noindent we find the corresponding value of $N$.
\end{example}

\section*{Conclusions to chapter \ref{omch:dvw2}}

In section \ref{omch:dvw2}, the concept of a random process in the spaces $D_{V,W}$ is introduced, the distribution of the supremum of such a process is studied, as well as the conditions for its sample continuity.
Also, the conditions for uniform convergence of functional series containing elements from $D_{V,W}$ are studied.
In addition, the estimates of the reliability and accuracy of process models from the spaces $D_{V,W}$ is investigated, and examples of such random processes are also considered.

\chapter{Estimates of norms in $L_p(T)$ of random processes from $\mathbf{F}_\psi(\Omega)$ spaces}
\label{ch:o6series}

\section{Estimates of norms in $L_p(T)$ of random processes from $\mathbf{F}_\psi(\Omega)$ spaces}
\label{ch:o6series}

The section \ref{ch:o6series} contains estimates for distributions of norms in $L_p(T)$ of random processes from spaces $\mathbf{F}_\psi(\Omega)$.

\begin{lemma}
\label{le:o6chapter-1}
Let $\xi$ be a random variable belonging to the space  $\mathbf{F}_\psi(\Omega)$. Then for $p\geq 1$
\[
\left\|\xi\right\|_\psi\leq \frac{\psi(p)}{\psi(1)}\sup\limits_{u\geq p}\frac{\left(E\left|\xi\right|^{u}\right)^{1/u}}{\psi(u)}.
\]
\end{lemma}

\begin{proof}
Obviously, from Lyapunov's inequality we have:
\[
\sup\limits_{1\leq u\leq p}\frac{\left(E\left|\xi\right|^{u}\right)^{1/u}}{\psi(u)}\leq \frac{\psi(p)}{\psi(1)} \frac{\left(E\left|\xi\right|^{p}\right)^{1/p}}{\psi(p)}.
\]
Therefore,
\[
\left\|\xi\right\|_\psi= \max \left(\sup\limits_{1\leq u\leq p}\frac{\left(E\left|\xi\right|^{u}\right)^{1/u}}{\psi(u)}, \sup\limits_{u\geq p}\frac{\left(E\left|\xi\right|^{u}\right)^{1/u}}{\psi(u)}\right)\leq
\]
\[
\leq \max \left(\frac{\psi(p)}{\psi(1)} \frac{\left(E\left|\xi\right|^{p}\right)^{1/p}}{\psi(p)}, \sup\limits_{u\geq p}\frac{\left(E\left|\xi\right|^{u}\right)^{1/u}}{\psi(u)}\right)\leq \frac{\psi(p)}{\psi(1)}\sup\limits_{u\geq p}\frac{\left(E\left|\xi\right|^{u}\right)^{1/u}}{\psi(u)}.
\]
The last inequality implies what had to be proved.
\end{proof}

\begin {theorem}\label{th:o6chapter-1}
Let $\nu$ be a $\sigma$-finite measure in a compact metric space $(T,\rho)$, $X = \{X(t), \ t \in T\}$ be a random, measurable process from the space $\mathbf{F}_\psi(\Omega)$.
Let for some $p\geq 1$, the following condition holds true
\begin{equation}\label{eq:o6chapter-1}
\int\limits_{T}\left\|X(t)\right\|_\psi^p d\nu(t)<\infty.
\end{equation}
Then:
\begin{enumerate}
\item [1.] with probability one there exists an integral $\int\limits_{T}\left|X(t)\right|^p d\nu(t)$ and the inequality holds:
\begin{equation} \label{eq:o6chapter-2}
\left\|\left(\int\limits_{T}\left|X(t)\right|^p d\nu(t)\right)^{1/p}\right\|_\psi\leq \frac{\psi(p)}{\psi(1)}\left(\int\limits_{T}\left\|X(t)\right\|_\psi^p d\nu(t)\right)^{1/p};
\end{equation}
\item [2.] for any $\varepsilon>0$ the inequality holds:
\begin{multline} \label{eq:o6chapter-3}
P\left\{\left(\int\limits_{T}\left|X(t)\right|^p d\nu(t)\right)^{1/p}>\varepsilon\right\}\leq \\[1ex] \leq\inf\limits_{u\geq1}\frac{\left(\frac{\psi(p)}{\psi(1)}\right)^u\left(\int\limits_{T}\left\|X(t)\right\|_\psi^p d\nu(t)\right)^{u/p}(\psi(u))^u}{\varepsilon^u}.
\end{multline}
\end{enumerate}
\end {theorem}

\begin{proof}
Since
\[
E\int\limits_{T}\left|X(t)\right|^p d\nu(t)=\int\limits_{T}E\left|X(t)\right|^p d\nu(t)\leq \int\limits_{T}(\psi(p))^p\left\|X(t)\right\|_\psi^p d\nu(t)  <\infty,
\]
then $\int\limits_{T}\left|X(t)\right|^p d\nu(t)$ exists with probability one. From the generalized Minkowski inequality it follows that for $u\geq p$
\begin{multline} \label{eq:o6chapter-4}
E\left(\int\limits_{T}\left|X(t)\right|^p d\nu(t)\right)^{u/p}=\left(\left(E\left(\int\limits_{T}\left|X(t)\right|^p d\nu(t)\right)^{u/p}\right)^{p/u}\right)^{u/p}\leq \\[1ex]
\leq\left(\int\limits_{T}\left(E\left|X(t)\right|^u\right)^{p/u}d\nu(t)\right)^{u/p}\leq
\left(\int\limits_{T}\left\|X(t)\right\|_\psi^p (\psi(u))^p d\nu(t)\right)^{u/p}\leq \\[1ex]
\leq(\psi(u))^u \left(\int\limits_{T}\left\|X(t)\right\|_\psi^p d\nu(t)\right)^{u/p}.
\end{multline}

From Lemma \ref{le:o6chapter-1} and inequality \eqref{eq:o6chapter-4} we obtain that
\[
\left\|\left(\int\limits_{T}\left|X(t)\right|^p d\nu(t)\right)^{1/p}\right\|_\psi \leq \frac{\psi(p)}{\psi(1)} \sup\limits_{u\geq p}\frac{\left(E\left|\int\limits_{T}\left|X(t)\right|^p d\nu(t)\right|^{u/p}\right)^{1/u}}{\psi(u)} \leq
\]
\[
\leq\frac{\psi(p)}{\psi(1)}\sup\limits_{u\geq p}\frac{\psi(u)\left(\int\limits_{T}\left\|X(t)\right\|_\psi^p d\nu(t)\right)^{1/p}}{\psi(u)}=\frac{\psi(p)}{\psi(1)}\left(\int\limits_{T}\left\|X(t)\right\|_\psi^p d\nu(t)\right)^{1/p}.
\]
The inequality \eqref{eq:o6chapter-2} is proved. The inequality  \eqref{eq:o6chapter-3} follows from Theorem~\ref{th:o1chapter_1}.
\end{proof}

\begin{example}
Consider the space $\mathbf{F}_\psi(\Omega)$, where $\psi(u)=u^\alpha$, $\alpha>0$. From Theorem \ref{th:o6chapter-1} and Theorem \ref{th:o1chapter_2} for
$$\varepsilon\geq e^\alpha\frac{\psi(p)}{\psi(1)}\left(\int\limits_{T}\left\|X(t)\right\|_\psi^pd\nu(t)\right)^{1/p}$$ we get:
$$
P\left\{\left(\int\limits_{T}\left|X(t)\right|^p d\nu(t)\right)^{1/p}>\varepsilon\right\}
\leq\exp \left\{-\frac{\alpha}{e}\left(\frac{\varepsilon}{\frac{\psi(p)}{\psi(1)}\left(\int\limits_{T}\left\|X(t)\right\|_\psi^pd\nu(t)\right)^{1/p}}\right)^{1/\alpha}\right\}.
$$
\end{example}

\begin{example}
Consider the space $\mathbf{F}_\psi(\Omega)$, where $\psi(u)=e^{au^\beta}$, $a>0$, $\beta>0$. From Theorem \ref{th:o6chapter-1} and Theorem \ref{th:o1chapter_3} for $$\varepsilon\geq e^{a(\beta+1)} \frac{\psi(p)}{\psi(1)}\left(\int\limits_{T}\left\|X(t)\right\|_\psi^pd\nu(t)\right)^{1/p}$$ it follows that
\begin{multline*}
P\left\{\left(\int\limits_{T}\left|X(t)\right|^p d\nu(t)\right)^{1/p}>\varepsilon\right\}\leq \\[1ex]
\leq\exp \left\{-\frac{\beta}{a^{1/\beta}}\left(\frac{\ln \frac{\varepsilon}{\frac{\psi(p)}{\psi(1)}\left(\int\limits_{T}\left\|X(t)\right\|_\psi^pd\nu(t)\right)^{1/p}}}{\beta+1}\right)^{\frac{\beta+1}{\beta}}\right\}.
\end{multline*}
\end{example}

\begin{example}
Consider the space $\mathbf{F}_\psi(\Omega)$, where $\psi(u)=\left(\ln (u+1)\right)^\lambda$, $\lambda>0$. According to Theorem \ref{th:o6chapter-1} and Theorem\ref{th:o1chapter_4} for $\varepsilon>0$ we can state that
\[
P\left\{\left(\int\limits_{T}\left|X(t)\right|^p d\nu(t)\right)^{1/p}>\varepsilon\right\}\leq
\]
\[
\leq e^\lambda \exp \left\{-\lambda  \exp\left\{\left(\frac{\varepsilon}{\frac{\psi(p)}{\psi(1)}\left(\int\limits_{T}\left\|X(t)\right\|_\psi^pd\nu(t)\right)^{1/p}}\right)^{1/\lambda}\frac{1}{e}\right\}\right\}.
\]
\end{example}

\begin {theorem}\label{th:o6chapter-2}
Let $\nu$ be a $\sigma$-finite measure in a compact metric space $(T,\rho)$, $Y = \{Y(t), \ t \in T\}$ be a random process from the space $\mathbf{F}_\psi(\Omega)$ and for this space the condition $\mathbf{H}$ with the constant $C_\psi$ holds. Let $EY(t)=m(t)$, $Z_n(t)=\frac{1}{n}\sum\limits_{k=1}^{n}Y_k(t)-m(t)=\frac{1}{n}\sum\limits_{k=1}^{n}(Y_k(t)-m(t))$, where $Y_k(t)$ are independent copies of $Y(t)$. Then for all $p\geq 1$ the inequality holds true
\begin{equation}\label{eq:o6chapter-5}
\left\|\left(\int\limits_{T}\left|Z_n(t)\right|^p d\nu(t)\right)^{1/p}\right\|_\psi\leq \frac{2 \sqrt{C_\psi}}{\sqrt{n}}\cdot \frac{\psi(p)}{\psi(1)}\left(\int\limits_{T}\left\|Y(t)\right\|_\psi^p d\nu(t)\right)^{1/p}
\end{equation}
and for any $\varepsilon>0$ the following estimate is valid
\begin{multline} \label{eq:o6chapter-6}
P\left\{\left(\int\limits_{T}\left|Z_n(t)\right|^p d\nu(t)\right)^{1/p}>\varepsilon\right\}\leq \\[1ex]
\leq\inf\limits_{u\geq1}\frac{\left(\frac{2 \sqrt{C_\psi}}{\sqrt{n}}\cdot\frac{\psi(p)}{\psi(1)}\right)^u\left(\int\limits_{T}\left\|Y(t)\right\|_\psi^pd\nu(t)\right)^{u/p}(\psi(u))^u}{\varepsilon^u}.
\end{multline}
\end {theorem}

\begin{proof}
From the Definition \ref{de:o1chapter_5} and Lemma \ref{le:o2chapter_1} it follows that
\[
\left\|Z_n(t)\right\|_\psi^2\leq \frac{1}{n^2}C_\psi \sum\limits_{k=1}^{n}\left\|Y_k(t)-m(t)\right\|_\psi^2= \frac{1}{n}C_\psi \left\|Y(t)-m(t)\right\|_\psi^2\leq
\]
\[
\leq\frac{1}{n}C_\psi \left(\left\|Y(t)\right\|_\psi+\left\|m(t)\right\|_\psi\right)^2\leq\frac{4}{n}C_\psi  \left\|Y(t)\right\|_\psi^2.
\]

Since from Theorem \ref{th:o6chapter-1} the relation is established
\[
\left\|\left(\int\limits_{T}\left|Z_n(t)\right|^p d\nu(t)\right)^{1/p}\right\|_\psi \leq\frac{\psi(p)}{\psi(1)}\left(\int\limits_{T}\left\|Z_n(t)\right\|_\psi^p d\nu(t)\right)^{1/p}\leq
\]
\[
\leq\frac{\psi(p)}{\psi(1)}\left(\int\limits_{T}\left(\frac{2 \sqrt{C_\psi}}{\sqrt{n}}\left\|Y(t)\right\|_\psi\right)^p d\nu(t)\right)^{1/p}=
\]
\[
=\frac{2 \sqrt{C_\psi}}{\sqrt{n}}\cdot \frac{\psi(p)}{\psi(1)}\left(\int\limits_{T}\left\|Y(t)\right\|_\psi^p d\nu(t)\right)^{1/p},
\]
then inequality \eqref{eq:o6chapter-5} holds, and inequality \eqref{eq:o6chapter-6} follows from Theorem \ref{th:o1chapter_1}.
\end{proof}

\begin{example} \label{ex:o6chapter-4}
Consider the space $\mathbf{F}_\psi(\Omega)$, where $\psi(u)=u^\alpha$, $\alpha>0$, then from Theorem \ref{th:o6chapter-2} and Theorem \ref{th:o1chapter_2} for $$\varepsilon \geq e^\alpha \frac{2 \sqrt{C_\psi}}{\sqrt{n}}\cdot \frac{\psi(p)}{\psi(1)}\left(\int\limits_{T}\left\|Y(t)\right\|_\psi^p d\nu(t)\right)^{1/p}$$ we obtain:
\begin{multline*}
P\left\{\left(\int\limits_{T}\left|Z_n(t)\right|^p d\nu(t)\right)^{1/p}>\varepsilon\right\}\leq \\[1ex]
\leq\exp \left\{-\frac{\alpha}{e}\left(\frac{\varepsilon}{\frac{2 \sqrt{C_\psi}}{\sqrt{n}}\cdot \frac{\psi(p)}{\psi(1)}\left(\int\limits_{T}\left\|Y(t)\right\|_\psi^p d\nu(t)\right)^{1/p}}\right)^{1/\alpha}\right\}.
\end{multline*}
\end{example}

\begin{example} \label{ex:o6chapter-5}
Consider the space $\mathbf{F}_\psi(\Omega)$, where $\psi(u)=e^{au^\beta}$, $a>0$, $\beta>0$, then from Theorem \ref{th:o6chapter-2} and Theorem \ref{th:o1chapter_3} for $$\varepsilon\geq e^{a(\beta+1)} \frac{2 \sqrt{C_\psi}}{\sqrt{n}}\cdot \frac{\psi(p)}{\psi(1)}\left(\int\limits_{T}\left\|Y(t)\right\|_\psi^p d\nu(t)\right)^{1/p}$$ we can conclude that
\begin{multline*}
P\left\{\left(\int\limits_{T}\left|Z_n(t)\right|^p d\nu(t)\right)^{1/p}>\varepsilon\right\}\leq \\[1ex]
\leq\exp \left\{-\frac{\beta}{a^{1/\beta}}\left(\frac{\ln \frac{\varepsilon}{\frac{2 \sqrt{C_\psi}}{\sqrt{n}}\cdot \frac{\psi(p)}{\psi(1)}\left(\int\limits_{T}\left\|Y(t)\right\|_\psi^p d\nu(t)\right)^{1/p}}}{\beta+1}\right)^{\frac{\beta+1}{\beta}}\right\}.
\end{multline*}
\end{example}

\begin{example}
Consider the space $\mathbf{F}_\psi(\Omega)$, where $\psi(u)=\left(\ln (u+1)\right)^\lambda$, $\lambda>0$, then from Theorem \ref{th:o6chapter-2} and Theorem \ref{th:o1chapter_4} for $\varepsilon \geq 0$ we have the following estimate:
\begin{multline*}
P\left\{\left(\int\limits_{T}\left|Z_n(t)\right|^p d\nu(t)\right)^{1/p}>\varepsilon\right\}\leq \\[1ex]
\leq e^\lambda \exp \left\{-\lambda  \exp\left\{\left(\frac{\varepsilon}{\frac{2 \sqrt{C_\psi}}{\sqrt{n}}\cdot \frac{\psi(p)}{\psi(1)}\left(\int\limits_{T}\left\|Y(t)\right\|^p d\nu(t)\right)^{1/p}}\right)^{1/\lambda}\frac{1}{e}\right\}\right\}.
\end{multline*}
\end{example}

\section*{Conclusions to chapter \ref{ch:o6series}}
Estimates for the distributions of norms in $L_p(T)$ of random processes from the spaces $\mathbf{F}_\psi(\Omega)$, as well as estimates of the norms of sums of independent copies of a random process from this space, have been found.

\chapter{Accuracy and reliability of Monte Carlo integral calculations}
\label{ch:o7series}

In this chapter we deal with the Monte Carlo method of calculating multiple integrals given on $\mathbb R^n$ with specified reliability and accuracy. Two approaches to solving this problem are proposed. The first approach is based on the theory of Orlicz spaces of random variables, and the second is based on the theory of $\mathbf{F}_\psi(\Omega)$ spaces. Integrals that depend on a parameter are also considered.

\section{Calculation of integrals by the Monte Carlo method} \label{se:o7series-1}

Let $\left\{\mathcal{S},\mathcal{A},\mu\right\}$ be a measurable space, let $\mu$ be a $\sigma$-finite measure, let $p(s)\geq 0, s\in \mathcal{S}$, be a function such that $\int_{\mathcal{S}}p(s)d\mu(s)=1$, and let $P(A)$, $A \in \mathcal{A}$, be a measure defined as: $P(A)=\int_{A}p(s)d\mu(s)$. Since $P(A)$ is a probability measure, then the space $\left\{\mathcal{S},\mathcal{A},P\right\}$ is a probability space.

Let $f(s)$ be a measurable function on $\left\{\mathcal{S},\mathcal{A},\mu\right\}$. Let us consider the integral (this integral is assumed to exist) $$I=\int\limits_{\mathcal{S}}f(s)p(s)d\mu(s),$$

\begin{remark}
We can consider an integral of the form $\int\limits_{\mathcal{S}}\varphi(s)d\mu(s)$.  If $p(s)> 0$ is the probability density in the space $\left\{\mathcal{S},\mathcal{A},\mu\right\}$, then
\[
\int\limits_{\mathcal{S}}\varphi(s)d\mu(s)=\int\limits_{\mathcal{S}}\frac{\varphi(s)}{p(s)}p(s)d\mu(s)=\int\limits_{\mathcal{S}}f(s)p(s)d\mu(s),
\]
where $f(s)=\frac{\varphi(s)}{p(s)}$.
\end{remark}

We assume that the functions $f(s)=\xi$ are random variables with $\left\{\mathcal{S},\mathcal{A},P\right\}$ and $$\int\limits_{\mathcal{S}}f(s)p(s)d\mu(s)=\int\limits_{\mathcal{S}}f(s)dP(s)=E\xi.$$

Let $\xi_{i}$, $i=1,\ldots,n$ be independent copies of the random variable $\xi$, $$Z_{n}=\frac{1}{n}\sum \limits_{i=1}^{n}\xi_{i},$$
 then by the strong law of large numbers $Z_{n}\rightarrow E\xi_{1}=I$ with probability one. Consider $Z_{n}$ as an estimate for $I$.

\begin{definition}\label{de:o7chapter-1}
We say that $Z_{n}$ approximates $I$ with reliability $1-\delta$ $(0<\delta<1)$ and accuracy $\varepsilon>0$ if the following inequality holds:
\begin{equation}
\label{eq:o7chapter-3}
P\left\{\left|Z_{n}-I\right|>\varepsilon\right\}\leq \delta.
\end{equation}
\end{definition}

\section{Application of the theory of Orlicz spaces}

\subsection{Accuracy and reliability of the calculation of integrals}

\begin{theorem}\label{th:o7chapter-1}
Let $\xi_1, \xi_2, \ldots , \xi_n$ be independent, identically distributed random variables belonging to the Orlicz space $L_U(\Omega)$.
 Let the condition $\mathbf{H}$ be satisfied for the Orlicz space $L_U(\Omega)$.
Let $$Y_n=\frac{1}{\sqrt{n}}\sum \limits_{i=1}^n\left(\xi_i-I\right),$$ where $I=E\xi_1$.

Then for any $\varepsilon>0$ the inequality holds:
\begin{equation}
\label{eq:o7chapter-1}
P\left\{\left|Y_n\right|>\varepsilon\right\}\leq \frac{1}{U\left(\frac{\varepsilon}{L}\right)},
\end{equation}
where $L=\left\|\xi_1-I\right\|_U\sqrt{C_U}$, $C_U$ is a constant from Definition~\ref{de:o3chapter_6}.
\end{theorem}

\begin{proof}
From Definition~\ref{de:o3chapter_6} it follows that
\begin{multline*}
\left\|Y_n\right\|_U^2=\left\|\frac{1}{\sqrt{n}}\sum \limits_{i=1}^n\left(\xi_i-I\right)\right\|_U^2=\frac{1}{n}\left\|\sum \limits_{i=1}^n\left(\xi_i-I\right)\right\|_U^2\leq\\
\leq\frac{1}{n} C_U \sum \limits_{i=1}^n \left\|\xi_i-I\right\|^2_U= C_U \left\|\xi_1-I\right\|^2_U.
\end{multline*}
The inequality~\eqref{eq:o7chapter-1} is established based on Lemma~\ref{le:o4chapter-1}.
\end{proof}

\begin{corollary}
\label{co:o7chapter-1}
Let the conditions of Theorem~\ref{th:o7chapter-1} be satisfied. Then for any $\varepsilon>0$ the inequality holds:
\begin{equation}
\label{eq:o7chapter-2}
P\left\{\left|\frac{1}{n}\sum \limits_{i=1}^n \xi_i-I\right|>\varepsilon\right\}\leq\frac{1}{U\left(\frac{\sqrt{n}\varepsilon}{L}\right)}.
\end{equation}
\end{corollary}

\begin{proof}
The obvious equality
\[
\frac{1}{n}\sum \limits_{i=1}^n \xi_i-I=\frac{1}{n}\sum \limits_{i=1}^n \left(\xi_i-I\right)=\frac{1}{\sqrt{n}}Y_n,
\]
therefore we obtain:
\[
P\left\{\left|\frac{1}{n}\sum \limits_{i=1}^n \xi_i-I\right|>\varepsilon\right\}=P\left\{\frac{1}{\sqrt{n}}\left|Y_n\right|>\varepsilon\right\}=P\left\{\left|Y_n\right|>\sqrt{n}\varepsilon\right\}\leq\frac{1}{U\left(\frac{\sqrt{n}\varepsilon}{L}\right)},
\]
which had to be proved.
\end{proof}

\begin {theorem}\label{th:o7chapter-2}
Let $$I=\int\limits_{\mathcal{S}}f(s)p(s)d\mu(s),$$ let $\xi(s)$ be a random variable, $s \in \left\{\mathcal{S},\mathcal{A},P\right\}$, let $p(s)$ be the density of the distribution of $\xi$, let $\xi_{i}$, $i=1,2,\ldots,n$, be independent copies of the random variable $\xi$, and let
$$Z_{n}=\frac{1}{n}\sum \limits_{i=1}^{n}\xi_{i}.$$
If the random variable $\xi$ belongs to the space $L_U(\Omega)$, for which the condition $\mathbf{H}$ with constant $C_{U}$ is satisfied, then $Z_{n}$ approximates $I$ with reliability $1-\delta$ and accuracy $\varepsilon$ (see definition \ref{de:o7chapter-1}) when the inequality is satisfied:
\begin{equation}
\label{eq:o7chapter-4}
n\geq\left(\frac{LU^{(-1)}\left(\frac{1}{\delta}\right)}{\varepsilon}\right)^{2},
\end{equation}
where $L= \left\|\xi-I\right\|_U\sqrt{C_U}$.
\end{theorem}

\begin{proof}
From Corollary~\ref{co:o7chapter-1} it follows that
\[
P\left\{\left|Z_n - I\right|>\varepsilon\right\}\leq \left(U\left(\frac{\sqrt{n}\varepsilon}{L}\right)\right)^{-1}.
\]
The statement of this Theorem is true if $\left(U\left(\frac{\sqrt{n}\varepsilon}{L}\right)\right)^{-1}\leq \delta$, that is, if the inequality~\eqref{eq:o7chapter-4} holds.
\end{proof}

\begin{remark}\label{re:o7chapter-1}
It is easy to show that
\begin{equation}
\label{eq:o7chapter-5}
\left\|\xi-I\right\|_U\leq\left(1+ \frac{d_U}{U^{(-1)}(1)}\right)\left\|\xi\right\|_U,
\end{equation}
where $d_U$ is defined in Lemma~\ref{le:o5chapter-3}.
\end{remark}

Since
\[
\left\|\xi-I\right\|_U\leq \left\|\xi\right\|_U + \left\|I\right\|_U,
\]
then from Lemma~\ref{le:o5chapter-2} we have that $$\left\|I\right\|_U<\frac{\left|I\right|}{U^{(-1)}(1)},$$ and from Lemma~\ref{le:o5chapter-3} it follows that $\left|I\right|\leq d_U \left\|f\right\|_U$. Thus, inequality~\eqref{eq:o7chapter-5} is true. Therefore, in the inequality ~\eqref{eq:o7chapter-4} instead of $L$ we can substitute $$\widehat{L}=\left\|\xi\right\|_U\sqrt{C_U}\left(1+ \frac{d_U}{U^{(-1)}(1)}\right).$$

\begin{example}\label{ex:o7chapter-1}
The Monte Carlo method is often used to calculate multiple integrals. But for simplicity, let's consider the integral of one variable. Let $$\int\limits_{-\infty}^{+\infty}f(x)\exp\{\frac{-x^2}{2a^2}-bx\}dx=I, a>0,$$ $\left|f(x)\right|<1$, then $$I=\sqrt{2\pi}a\frac{1}{\sqrt{2\pi}a}\int\limits_{-\infty}^{+\infty}f(x)\exp\{\frac{-x^2}{2a^2}-bx\}dx.$$
Let $J=E\exp\{-\xi b\}$, where $\xi=N(0,a^2)$, $\eta_i=f(\xi_i)\exp\{-\xi_i b\}$, $\xi_i$ are independent copies of the random variable $\xi$. Let $$J_n=\frac{1}{n}\sum\limits_{i=1}^n\eta_i=\frac{1}{n}\sum\limits_{i=1}^n f(\xi_i)\exp\{-\xi_i b\}.$$
The estimate for $I$ is $I_n=\sqrt{2\pi}aJ_n$. Let $U(x)=\left|x\right|^p,p\geq2$, then according to Remark~\ref{re:o7chapter-1} to Theorem~\ref{th:o7chapter-2} and the value of $C_U$ for the space $L_U(\Omega)$ we obtain (in this case $\left\|\eta\right\|_U=\left\|\eta\right\|_p=\left(E\left|\eta\right|^p\right)^{1/p}$):
\[
L=\left\|\eta\right\|_p 2 \sqrt{\sqrt{2}\left(\Gamma(p+1)/2\sqrt{\pi}\right)^{1/p}},
\]
where $$\left\|\eta\right\|_p^p=E\left|f(\xi)\right|^p\left(\exp\{-\xi b\}\right)^p\leq E\exp\{-\xi bp\}=\exp\{\frac{p^2b^2a^2}{2}\}.$$

According to inequality~\eqref{eq:o7chapter-4} the integral $I$ will be calculated with an accuracy $\varepsilon$ and a reliability of $1-\delta$ when the inequality is satisfied:
\[
n\geq \frac{a^22\pi L^2}{\varepsilon^2\delta^{2/p}}.
\]

The last inequality must be true for $p\geq2$, that is, it is necessary to find the smallest $n$, that is find the minimum of the right-hand side in $p$, more precisely, find the approximate value of the minimum. According to the Stirling formula $$\Gamma(p)\cong \exp\{-p\}p^{p-1/2}(2\pi)^{1/2}$$ we have that
\[
\frac{L^2}{\delta^{2/p}}\cong \frac{4\sqrt{2}\exp\{pb^2a^2\}p(p/2)^{1/2p}}{\delta^{2/p}}.
\]
It is easy to see that the expression $\frac{L^2}{\delta^{2/p}}$ takes on an approximate minimum value at the point
\[
p\cong\frac{2(-ln\delta)}{1+\sqrt{1+4a^2b^2(-ln\delta)}}.
\]
\end{example}

\begin{theorem}
Let the conditions of Theorem~\ref{th:o7chapter-1} be satisfied. Then with probability one for sufficiently large $n$ the following inequality is valid
\[
\left|Z_n - I\right|\leq \frac{L}{\sqrt{n}}U^{(-1)}\left(\frac{1}{\delta_n}\right),
\]
where $\delta_n>0$ is any sequence such that $\sum_{n=1}^\infty \delta_n < \infty$.
\end{theorem}

\begin{proof}
This Theorem follows from the Borel-Cantelli lemma. Indeed, from the Corollary~\ref{co:o7chapter-1} it follows that

\[
P\left\{\left|Z_n - I\right|> \frac{L}{\sqrt{n}}U^{(-1)}\left(\frac{1}{\delta_n}\right)\right\}\leq \left(U \left(\frac{\sqrt{n}}{L}\frac{L}{\sqrt{n}}U^{(-1)}\left(\frac{1}{\delta_n}\right)\right)\right)^{-1}=\delta_n.
\]
\end{proof}

\begin{example}
If $U(x)=\left|x\right|^p$, $p\geq 2$ and $\delta_n=\frac{1}{n^{1+\kappa}}$, $\kappa>0$, where $\kappa$ is such that $\frac{1+\kappa}{p}<\frac{1}{2}$, then for sufficiently large $n$ we have:
\[
\left|Z_n - I\right|\leq\frac{L}{n^{1/2-1/p-\kappa/p}}.
\]

If $U(x)=exp\left\{\left|x\right|^{\alpha} \right\}-1$, $1\leq\alpha\leq2$ and $\delta_n=\frac{1}{n^{1+\kappa}}$, $\kappa>0$, then also for sufficiently large $n$ we obtain the following estimate:
\[
\left|Z_n - I\right|\leq \check{L}\frac{1}{n^{1/2}}\left(\ln n \right)^{1/\alpha},
\]
where $\check{L}$ is a constant.
\end{example}

\subsection{Reliability and accuracy in the space $C(T)$ of calculation of integrals depending on the parameter}\label{su:o7series-2}\indent

Consider the integral $$\int\limits_{\mathcal{S}}f(s,t)p(s)d\mu(s)=I(t),\, t \in T.$$
We assume that this integral exists. Let all the assumptions of section \ref{se:o7series-1} be fulfilled, while the function $f(s,t)$ depends on $t \in T$, where $(T,\rho)$ is a compact metric space and this function $f(s,t)$ is continuous with respect to $t$.

Consider $f(s,t)$ as a random process on $\left\{\mathcal{S},\mathcal{A},P\right\}$ and denote it by $\xi(s,t)=\xi(t)$ and $$I(t)=\int\limits_{\mathcal{S}}f(s,t)p(s)d\mu(s)=\int\limits_{\mathcal{S}}f(s,t)dP(s)=E\xi(t).$$

Let $\xi_{i}(t)$, $i=1,2,\ldots,n$ be independent copies of the random process $\xi(t)$, $Z_{n}(t)=\frac{1}{n}\sum \limits_{i=1}^{n}\xi_{i}(t)$. Then by the strengthened law of large numbers $Z_{n}(t)\rightarrow E\xi(t)=I(t)$ with probability one for any $t \in T$.

\begin{definition}\label{de:o7chapter-2}
We say that $Z_{n}(t)$ approaches $I(t)$ in the space $C(T)$ with reliability $1-\delta>0$ and accuracy $\varepsilon>0$ if the following inequality holds:
\[
P\left\{\sup\limits_{t \in T}\left|Z_{n}(t)-I(t)\right|>\varepsilon\right\}\leq \delta.
\]
\end{definition}

\begin{theorem}\label{th:o7chapter-4}
Let:
\begin{itemize}
\item [1.] the random process $\xi(t)$ belongs to the space $L_U(\Omega)$, where for the space $L_U(\Omega)$ the condition $\mathbf{H}$ with the constant $C_{U}$ is satisfied and the function $U$ satisfies the condition $\mathtt{g}$;
\item [2.] there is a continuous, increasing function $\sigma = \left(\sigma(h), 0\leq h \leq \delta_{0}\right)$,
$\delta_{0}= \sigma_1\left(\sup\limits_{t, s \,\in\, T} \rho(t, s)\right)$, such that
\begin{equation}\label{eq:o7chapter-6}
\sup\limits_{\rho(t, s) \leq h} \bigl\|\xi(t) - \xi(s) \bigr\|_U \leq \sigma(h)
\end{equation}
and
\begin{equation}
\int\limits_{0}^{\delta_0} U^{(-1)}(N(\sigma^{(-1)} (u))) du < \infty;
\end{equation}

\item [3.] the following inequality holds
\[
\left(U \left(\frac{\varepsilon\sqrt{n}}{\check{B}(\theta)}\right)\right)^{-1}\leq\delta,
\]
i.e
\begin{equation}\label{eq:o7chapter-8}
n\geq \frac{\check{B}^2(\theta) U^{(-1)}\left(\frac{1}{\delta}\right)}{\varepsilon^2},
\end{equation}
where $\check{B}(\theta) = C_U\left(1+\frac{d_U}{U^{(-1)}(1)}\right)\inf\limits_{t \in T}\left\|\xi(t)\right\|_U + \frac{1}{\theta(1 - \theta)}
\int\limits_{0}^{\delta_{0}\theta} \varkappa (N (\sigma^{(-1)}_1 (u)) du$, \\
$\sigma_1 (h)=C_U\left(1+\frac{d_U}{U^{-1}(1)}\right)\sigma (h)$, $0<\theta<1$, $\varkappa(u)$ is the majorizing characteristic of the space $L_U(\Omega)$, $N(u)$ is the metric massiveness of this space.
\end{itemize}
Then $Z_{n}(t)$ approximates $I(t)$ with reliability $1-\delta$ and accuracy $\varepsilon$ in the space $C(T)$ (see definition \ref{de:o7chapter-2}).
\end{theorem}

\begin{proof}
The Theorem follows from Theorem \ref{th:o3chapter_8} and Remark \ref{re:o5chapter-2}.
The function $f(t,s)$ is continuous. That is, the process $\xi(t)$ is separable. Thus, from inequality \eqref{eq:o5chapter-12} and Remark \ref{re:o5chapter-2} we have that
\[
P\{\sup_{t \in T} \sqrt{n}\left|Z_n (t) -m(t) \right| > \varepsilon\} \leq \frac{1}{U \left(\frac{\varepsilon}{\check{B}(\theta)}\right)}.
\]
That is,
\[
P\{\sup_{t \in T} \left|Z_n (t) -m(t) \right| > \varepsilon\}=
\]
\[
=P\{\sup_{t \in T} \sqrt{n}\left|Z_n (t) -m(t) \right| > \sqrt{n}\varepsilon\}\leq \frac{1}{U \left(\frac{\varepsilon\sqrt{n}}{\check{B}( \theta)}\right)}.
\]
From the last inequality we obtain the inequality \eqref{eq:o7chapter-8}.
\end{proof}

\begin{example}
Let $$I(t)=\sqrt{2\pi}a\frac{1}{\sqrt{2\pi}a}\int\limits_{-\infty}^{+\infty}f(x)\exp\{\frac{-x^2}{2a^2}-tx\}dx,\,a>0,\,\left|f(x)\right|<1,\, 0\leq t\leq T.$$ We use the same notation as in  Example~\ref{ex:o7chapter-1}. The estimate for $I(t)$ is $$I_n(t)=\sqrt{2\pi}aJ_n(t),$$ where $$J_n(t)=\frac{1}{n}\sum\limits_{i=1}^n\eta_i,\, \eta_i=f(\xi_i)\exp\left\{-\xi_it\right\}.$$
Let $U(x)=\left|x\right|^p, \ p\geq2$, then, according to Theorem~\ref{th:o7chapter-4}, we have that
\[
\check{B}(\theta) = 2\bigl\|\eta(t)\bigr\|_p + \frac{1}{\theta(1 - \theta)}
\int\limits_{0}^{\delta_{0}\theta} \varkappa (N_w (\sigma^{(-1)}_1 (u)) du,
\]
where $\bigl\|\eta(t)\bigr\|_p\leq\exp\left\{\frac{pt^2a^2}{2}\right\}$, and $\inf\limits_{0 \leq t \leq T}\bigl\|\eta(t)\bigr\|_p=1$.

Since for $u<\delta_0$ the conditions $\frac{T}{2\sigma^{(-1)}_1 (u)}>\frac{1}{2}$ and $\sigma_1(h)=2\sigma(h)$ are valid, then Theorem \ref{th:o7chapter-4} yields the following estimate:
\begin{multline*}
\int\limits_{0}^{\delta_{0}\theta} \varkappa (N_w (\sigma^{(-1)}_1 (u)) du\leq \int\limits_{0}^{\delta_{0}\theta} \left(\frac{T}{2\sigma^{(-1)}_1 (u)}+1\right)^{1/p}du\leq \\[1ex]
\leq\int\limits_{0}^{\delta_{0}\theta} \left(\frac{3}{2}T\right)^{1/p}\left(\frac{1}{\sigma^{(-1)}_1 (u)}\right)^{1/p}du=2\left(\frac{3}{2}T\right)^{1/p}\int\limits_{0}^{{\delta_{0}\theta}/2} \left(\frac{1}{\sigma^{(-1)} (v)}\right)^{1/p}dv.
\end{multline*}

From the inequality ~\eqref{eq:o7chapter-6} we find $\sigma(v)$, $0\leq v\leq T$,
\begin{multline*}
\left\|\exp\{-\xi t\}-\exp\{-\xi s\}\right\|_p^p\leq E\left|\exp\{-\xi t\}-\exp\{-\xi s\}\right|^p=\\[1ex]
=EI\{\xi\geq0\}\left|\exp\{-\xi t\}-\exp\{-\xi s\}\right|^p+\\[1ex]
+EI\{\xi<0\}\left|\exp\{-\xi t\}-\exp\{-\xi s\}\right|^p=\Delta_{+}+\Delta_{-}
\end{multline*}

Let $\frac{1}{\alpha}+\frac{1}{\beta}=1$, $\beta>1$ and $s>t$. Then
\[
\Delta_{-}=EI\{\xi<0\}\left|\exp\{-\xi t\}-\exp\{-\xi s\}\right|^p=
\]
\[
=EI\{\xi<0\}\left|\exp\{-\xi t\}(1-\exp\{-\xi (s-t)\})\right|^p\leq
\]
\[
\leq EI\{\xi<0\}\left|\exp\{-\xi t\}\left|\xi\right|(s-t)\right|^p\leq
\]
\[
\leq\left(EI\{\xi<0\}\exp\{-\xi pt\beta\}\right)^{1/\beta}\left(EI\{\xi<0\}\left|\xi\right|^{p\alpha}\right)^{1/\alpha} (s-t)^p.
\]
So $$\Delta_{-}\leq\left|t-s\right|^p \exp\{\frac{a^2p^2t^2\beta}{2}\}\left(E\left|\xi\right|^{p\alpha}\right)^{1/\alpha}.$$ Finding the estimate for $\Delta_{+}$ is similar.
Therefore, $\sigma(v)=\left|t-s\right|C_p$, where $$C_p=2^{1/p} \exp\left\{\frac{a^2p\beta T^2}{2}\right\}\left(E\left|\xi\right|^{p\alpha}\right)^{1/p\alpha}.$$ Let us estimate the integral
\[
E\left|\xi\right|^{p\alpha}=\frac{1}{\sqrt{2\pi }a}\int\limits_{-\infty}^{+\infty}\left|x\right|^{p\alpha}\exp\left\{\frac{-x^2}{2a^2}\right\}dx=\frac{a^{p\alpha}}{\sqrt{2\pi }}\int\limits_{-\infty}^{+\infty}\left|t\right|^{p\alpha}\exp\left\{\frac{-t^2}{2}\right\}dt.
\]
From the inequality $x^s\leq\left(\frac{s}{e}\right)\exp\{x\}$ it follows that
\[
E\left|\xi\right|^{p\alpha}\leq\frac{a^{p\alpha}}{\sqrt{2\pi}}\left(\frac{p\alpha}{e}\right)^{p\alpha}\int\limits_{-\infty}^{+\infty}\exp\left\{\left|t\right|\right\}\exp\left\{\frac{-t^2}{2}\right\}dt\leq 2 \exp\left\{1/2\right\}a^{p\alpha}\left(\frac{p\alpha}{e}\right)^{p\alpha}.
\]
We obtain that $$C_{p}=\left(2\exp\left\{1/2\right\}\right)^{1/p\alpha}a \frac{p\alpha}{e} 2^{1/p} \exp\left\{\frac{a^2p\beta T^2}{2}\right\}.$$ Substituting the obtained values, we have that
\[
\check{B}(\theta)=1+\frac{1}{\theta(1 - \theta)}\int\limits_{0}^{{\delta_{0}\theta}/2} \left(\frac{C_p}{v}\right)^{1/p}dv
=1+\frac{C_p^{1/p}p}{\theta(1 - \theta)(p-1)}\left(\frac{\delta_o\theta}{2}\right)^{1-1/p}.
\]
Since $\delta_0=\sigma_1 \left(\sup\limits_{t, s \,\in\, T} \rho(t,s)\right)$ and $\sup\limits_{t, s \,\in\, T} \rho(t,s)=T$, then $\delta_0=\sigma_1(T)=2\sigma (T)=2TC_p$ and $$\check{B}(\theta)=2\exp\left\{\frac{a^2pT^2}{2}\right\}+\frac{2^{1/p-1}C_ppT}{(1 - \theta)(p-1)}\left(T\theta\right)^{-1/p}.$$
According to the inequality~\eqref{eq:o7chapter-8} the given integral will be calculated with an accuracy of $\varepsilon$ and a reliability of $1-\delta$ if the inequality holds:
\[
n\geq \inf\limits_{p\geq2, 0<\theta<1}\left(\frac{a^22\pi \check{B}^2(\theta)}{\varepsilon^2\delta^{2/p}}\right).
\]
\end{example}

\section{Application of the theory of spaces $\mathbf{F}_\psi(\Omega)$}

\subsection{Accuracy and reliability of the calculation of integrals}\indent

This subsection retains all the notation, conditions, and definitions of section \ref{se:o7series-1}.
\begin{theorem}\label{th:o7chapter-5}

Let $\xi_1, \xi_2, \ldots , \xi_n$, be independent, identically distributed random variables belonging to the space $\mathbf{F}_\psi(\Omega)$. Let the condition $\mathbf{H}$ is satisfied for the space $\mathbf{F}_\psi(\Omega)$. Let $$Y_n=\frac{1}{\sqrt{n}}\sum \limits_{i=1}^n\left(\xi_i-I\right),\,I=E\xi_1.$$
Then for any $\varepsilon>0$ the inequality holds
\begin{equation}
\label{eq:o7chapter-9}
P\left\{\left|Y_n\right|>\varepsilon\right\}\leq \inf\limits_{u\geq1}\frac{L^u(\psi(u))^u}{\varepsilon^u},
\end{equation}
where $L=\left\|\xi_1-I\right\|_{\psi}\sqrt{C_\psi}$, $C_\psi$ is a constant from Definition~\ref{de:o1chapter_5}.
\end{theorem}

\begin{proof}
From  Definition~\ref{de:o1chapter_5} it follows that
\begin{multline*}
\left\|Y_n\right\|_\psi^2=\left\|\frac{1}{\sqrt{n}}\sum \limits_{i=1}^n\left(\xi_i-I\right)\right\|_\psi^2=\frac{1}{n}\left\|\sum \limits_{i=1}^n\left(\xi_i-I\right)\right\|_\psi^2\leq\\
\leq\frac{1}{n} C_\psi \sum \limits_{i=1}^n \left\|\xi_i-I\right\|^2_\psi= C_\psi \left\|\xi_1-I\right\|^2_\psi.
\end{multline*}
The inequality ~\eqref{eq:o7chapter-9} follows from Theorem ~\ref{th:o1chapter_1}.
\end{proof}

\begin{corollary}
\label{co:o7chapter-2}
Let the conditions of Theorem ~\ref{th:o7chapter-5} be satisfied. Then for any $\varepsilon>0$ the inequality holds
\begin{equation}
\label{eq:o7chapter-21}
P\left\{\left|\frac{1}{n}\sum \limits_{i=1}^n \xi_i-I\right|>
\varepsilon\right\}\leq \inf\limits_{u\geq1}\frac{L^u(\psi(u))^u}{\left(\sqrt{n}\varepsilon\right)^u}.
\end{equation}
\end{corollary}

\begin{proof}
The proof of  Corollary \ref{co:o7chapter-2} is similar to the proof of  Corollary~\ref{co:o7chapter-1}.
\end{proof}

\begin{remark}\label{re:o7chapter-3}
It is obvious that
\begin{equation}
\label{eq:o7chapter-11}
\left\|\xi_1-I\right\|_\psi\leq 2\left\|\xi_1\right\|_\psi.
\end{equation}
\end{remark}

Indeed, from Lemma~\ref{le:o2chapter_1} we have that
\[
\left\|\xi_1-E\xi_1\right\|_\psi\leq \left\|\xi_1\right\|_\psi+\left\|\left|E\xi_1\right|\right\|_\psi\leq 2\left\|\xi_1\right\|_\psi.
\]

\begin{corollary}
\label{co:o7chapter-3}
Let the conditions of Theorem~\ref{th:o7chapter-5} be true, then for any $\varepsilon>0$ the inequality holds
\begin{equation}
\label{eq:o7chapter-12}
P\left\{\left|\frac{1}{n}\sum \limits_{i=1}^n \xi_i-I\right|>
\varepsilon\right\}\leq \inf\limits_{u\geq1}\frac{2^u \widetilde{L}^u(\psi(u))^u}{\left(\sqrt{n}\varepsilon\right)^u},
\end{equation}
where $\widetilde{L}=\left\|\xi_1\right\|_{\psi}\sqrt{C_\psi}$.
\end{corollary}

\begin{proof}
Corollary \ref{co:o7chapter-3} follows from Corollary~\ref{co:o7chapter-2} and Remark~\ref{re:o7chapter-3}.
\end{proof}

\begin{example}\label{ex:o7chapter-4}
Consider the space $\mathbf{F}_\psi(\Omega)$, where $\psi(u)=u^\alpha$, $\alpha>\frac{1}{2}$. Then, taking into account results from Example \ref{ex:o1chapter_4}, we have that for this space the condition $\mathbf{H}$ with the constant $C_\psi=4 \cdot9^\alpha$ is satisfied. Then from  Corollary \ref{co:o7chapter-3} and Theorem \ref{th:o1chapter_2} for $\varepsilon\geq \frac{4(3e)^\alpha\left\|\xi_1\right\|_\psi}{\sqrt{n}}$ it follows that
\[
P\left\{\left|\frac{1}{n}\sum \limits_{i=1}^n \xi_i-I\right|>\varepsilon\right\}\leq \exp \left\{-\frac{\alpha}{3e}\left(\frac{\sqrt{n}\varepsilon}{4\left\|\xi_1\right\|_\psi}\right)^{1/\alpha}\right\}.
\]
\end{example}

\begin{example} \label{ex:o7chapter-5}
Consider the space $\mathbf{F}_\psi(\Omega)$, where $\psi(u)=e^{au^\beta}$, $a>0$, $0<\beta<1$. Based on Theorem \ref{th:o2chapter-14}, we have two cases.
In the first case, where $\frac{1}{(2a\beta)^{1/\beta}}=1$, the condition $\mathbf{H}$ with the constant $C_\psi=4e^{2^\beta a}$ is satisfied for the space $\mathbf{F}_\psi(\Omega)$.
Then, by the Corollary  \ref{co:o7chapter-3} and Theorem \ref{th:o1chapter_3}, we have:
\[
P\left\{\left|\frac{1}{n}\sum \limits_{i=1}^n \xi_i-I\right|>\varepsilon\right\}\leq \exp \left\{-\frac{\beta}{a^{1/\beta}}\left(\frac{\ln \frac{\sqrt{n}\varepsilon}{4e^{2^{\beta-1} a}\left\|\xi_1\right\|_\psi}}{\beta+1}\right)^{\frac{\beta+1}{\beta}}\right\},
\]
in the case where $$\varepsilon\geq\frac{4e^{a\left(2^{\beta-1}+\beta+1\right)}\left\|\xi_1\right\|_\psi}{\sqrt{n}}.$$
In the second case, where $\frac{1}{(2a\beta)^{1/\beta}}>1$, for the space $\mathbf{F}_\psi(\Omega)$ the condition $\mathbf{H}$ with the constant $C_\psi=\frac{4e^{a\left(2^\beta+1\right)-\frac{1}{2\beta}}}{(2a\beta)^{1/{2\beta}}}$ is satisfied. Then from  Corollary \ref{co:o7chapter-3} and Theorem \ref{th:o1chapter_3} it follows that
\[
P\left\{\left|\frac{1}{n}\sum \limits_{i=1}^n \xi_i-I\right|>\varepsilon\right\}\leq \exp \left\{-\frac{\beta}{a^{1/\beta}}\left(\frac{\ln \frac{\sqrt{n}\varepsilon (2a\beta)^{1/{4\beta}}}{4 e^{\frac{a}{2}\left(2^\beta+1\right)-\frac{1}{4\beta}}\left\|\xi_1\right\|_\psi}}{\beta+1}\right)^{\frac{\beta+1}{\beta}}\right\},
\]
if $$\varepsilon \geq \frac{4 e^{a\left(2^{\beta-1}+\beta+\frac{3}{2}\right)-\frac{1}{4\beta}}\left\|\xi_1\right\|_\psi}{\sqrt{n}(2a\beta)^{1/{4\beta}}}.$$
\end{example}

\begin{theorem}\label{th:o7chapter-6}
Let $I=\int\limits_{\mathcal{S}}f(s)p(s)d\mu(s)$, $\xi(s)$ be a random variable, $s \in \left\{\mathcal{S},\mathcal{A},P\right\}$, $p(s)$ be the density of the distribution of $\xi$, $\xi_{i}$, $i=1,2,\ldots,n$ be independent copies of the random variable $\xi$, $Z_{n}=\frac{1}{n}\sum \limits_{i=1}^{n}\xi_{i}$. If the random variable $\xi$ belongs to the space $\mathbf{F}_\psi(\Omega)$, for which the condition $\mathbf{H}$ with constant $C_\psi$ is satisfied and $n$ is such that
\begin{equation}\label{eq:o7chapter-13}
\inf\limits_{u\geq1}\frac{2^u \widetilde{L}^u(\psi(u))^u}{\left(\sqrt{n}\varepsilon\right)^u}\leq \delta,
\end{equation}
then $Z_{n}$ approximates $I$ with reliability $1-\delta$ and accuracy $\varepsilon$ (see definition \ref{de:o7chapter-1}).
In the estimate \ref{eq:o7chapter-13} $\widetilde{L}$ is defined in Corollary \ref{co:o7chapter-3}.
\end{theorem}

\begin{proof}
The theorem follows from  Corollary \ref{co:o7chapter-3} and inequality~\eqref{eq:o7chapter-3}.
\end{proof}

\begin{example}\label{ex:o7chapter-6}
Consider the space $\mathbf{F}_\psi(\Omega)$, where $\psi(u)=u^\alpha$, $\alpha\geq\frac{1}{2}$. Based on the considered Example \ref{ex:o7chapter-4} and Theorem \ref{th:o7chapter-6}
for $$\varepsilon\geq \frac{4(3e)^\alpha\left\|\xi\right\|_\psi}{\sqrt{n}}$$ we have that
\[
\inf\limits_{u\geq1}\frac{2^u \widetilde{L}^u(\psi(u))^u}{\left(\sqrt{n}\varepsilon\right)^u}=\exp \left\{-\frac{\alpha}{3e}\left(\frac{\sqrt{n}\varepsilon}{4\left\|\xi\right\|_\psi}\right)^{1/\alpha}\right\}.
\]
Therefore, the inequality \eqref{eq:o7chapter-13} is satisfied when the following inequality is true
\[
\exp \left\{-\frac{\alpha}{3e}\left(\frac{\sqrt{n}\varepsilon}{4\left\|\xi\right\|_\psi}\right)^{1/\alpha}\right\}\leq \delta.
\]
In this case
\[
n\geq\left(\frac{4 \left\|\xi\right\|_\psi}{\varepsilon}\right)^2\left((-\ln \delta)\frac{3e}{\alpha}\right)^{2\alpha}
\]
and
\[
n\geq\left(\frac{4 (3e)^\alpha\left\|\xi\right\|_\psi}{\varepsilon}\right)^2 \max\left(1,\left(\frac{-\ln \delta}{\alpha}\right)^{2\alpha}\right).
\]
\end{example}

\begin{remark}
If we evaluate the accuracy and reliability using Chebyshev's inequality, we get that
\[
P\left\{\left|\frac{1}{n}\sum \limits_{i=1}^n \xi_i-I\right|>\varepsilon\right\}\leq \frac{D\xi_1}{n\varepsilon^2},
\]
then $Z_{n}$ approximates $I$ with reliability $1-\delta$ and accuracy $\varepsilon$ when
\[
n\geq\frac{D\xi_1}{\delta\varepsilon^2}.
\]
So, if $\delta=0.01$, then $n\geq C_1\frac{100}{\varepsilon^2}$, and from Example \ref{ex:o7chapter-6} we have that
$n\geq C_2\frac{\ln100}{\varepsilon^2}\approx C_2\frac{4.61}{\varepsilon^2}$, where $C_1$, $C_2$ are some constants.
\end{remark}

\begin{example}
Let us consider the space $\mathbf{F}_\psi(\Omega)$, where $\psi(u)=e^{au^\beta}$, $a>0$, $0<\beta<1$. Based on results of Example \ref{ex:o7chapter-5} and Theorem \ref{th:o7chapter-6}, we consider two cases.
In the first case, for $\frac{1}{(2a\beta)^{1/\beta}}=1$, the inequality \eqref{eq:o7chapter-13} holds true when the following inequality is true
\[
\exp \left\{-\frac{\beta}{a^{1/\beta}}\left(\frac{\ln \frac{\sqrt{n}\varepsilon}{4e^{2^{\beta-1} a}\left\|\xi\right\|_\psi}}{\beta+1}\right)^{\frac{\beta+1}{\beta}}\right\}\leq \delta.
\]
In this case
\[
n\geq\left(\frac{4e^{2^{\beta-1} a}\left\|\xi\right\|_\psi}{\varepsilon}\right)^2 \exp\left\{2(\beta+1)\left((-\ln\delta)\frac{a^{1/\beta}}{\beta}\right)^{\frac{\beta}{\beta+1}}\right\}
\]
and
\[
n\geq\left(\frac{4e^{2^{\beta-1} a}\left\|\xi\right\|_\psi}{\varepsilon}\right)^2 \max \left(e^{2(\beta+1)}, \exp\left\{2(\beta+1)\left((-\ln\delta)\frac{a^{1/\beta}}{\beta}\right)^{\frac{\beta}{\beta+1}}\right\}\right).
\]
In the second case, where $\frac{1}{(2a\beta)^{1/\beta}}>1$, the inequality \eqref{eq:o7chapter-13} is true when the inequality holds
\[
\exp \left\{-\frac{\beta}{a^{1/\beta}}\left(\frac{\ln \frac{\sqrt{n}\varepsilon (2a\beta)^{1/{4\beta}}}{4 e^{\frac{a}{2}\left(2^\beta+1\right)-\frac{1}{4\beta}}\left\|\xi\right\|_\psi}}{\beta+1}\right)^{\frac{\beta+1}{\beta}}\right\}\leq \delta.
\]
In this case
\[
n\geq\left(\frac{4 e^{\frac{a}{2}\left(2^\beta+1\right)-\frac{1}{4\beta}}\left\|\xi\right\|_\psi}{\varepsilon (2a\beta)^{1/{4\beta}}}\right)^2 \exp\left\{2(\beta+1)\left((-\ln\delta)\frac{a^{1/\beta}}{\beta}\right)^{\frac{\beta}{\beta+1}}\right\}
\]
and
\[
n\geq\left(\frac{4 e^{\frac{a}{2}\left(2^\beta+1\right)-\frac{1}{4\beta}}\left\|\xi\right\|_\psi}{\varepsilon (2a\beta)^{1/{4\beta}}}\right)^2\times
\]
\[
\times\max \left(e^{2a(\beta+1)}, \exp\left\{2(\beta+1)\left((-\ln\delta)\frac{a^{1/\beta}}{\beta}\right)^{\frac{\beta}{\beta+1}}\right\}\right).
\]
\end{example}

\begin{example}\label{ex:o7chapter-8}
Consider an integral of the form
\[
\int\limits_{0}^{+\infty}\int\limits_{0}^{+\infty}c(x,y)(x+y)^{v-1}e^{-px}e^{-qy}dxdy,
\]
where $\left|c(x,y)\right|\leq 1$, $v>\frac{3}{2}$. Denote
\[
I=\frac{1}{pq}\int\limits_{0}^{+\infty}\int\limits_{0}^{+\infty}c(x,y)(x+y)^{v-1}pe^{-px}qe^{-qy}dxdy.
\]
Let $\xi$ and $\eta$ be independent random variables that are distributed according to the exponential distribution
\[
P\left\{\xi<x\right\}=\begin{cases}
1-e^{-px},& $x>0$;\\
0,& $x<0$,
\end{cases}
\]
\[
P\left\{\eta<y\right\}=\begin{cases}
1-e^{-qx},& $y>0$;\\
0,& $y<0$,
\end{cases}
\]
where $p(x)=pe^{-px}$, $p(y)=qe^{-qy}$.

So,
\[
I=\frac{1}{pq}\int\limits_{0}^{+\infty}\int\limits_{0}^{+\infty}c(x,y)(x+y)^{v-1}pe^{-px}qe^{-qy}dxdy=\frac{1}{pq} E c(\xi,\eta)(\xi+\eta)^{v-1}.
\]
Let the function $\psi(u)=u^{v-1}$, then since $v>\frac{3}{2}$, we get:
\[
\sup\limits_{u\geq1}\frac{\left(E\left(c(x,y)(\xi+\eta)^{v-1}\right)^{u}\right)^{1/u}}{u^{v-1}}\leq
\]
\[
\leq\sup\limits_{u\geq1}\frac{\left(E\left(\left|c(\xi,\eta)\right|(\xi+\eta)^{u(v-1)}\right)^{\frac{1}{u(v-1)}}\right)^{v-1}}{u^{v-1}}\leq
\]
\[
\leq\sup\limits_{u\geq1}\frac{\left(E\left((\xi+\eta)^{u(v-1)}\right)^{\frac{1}{u(v-1)}}\right)^{v-1}}{u^{v-1}}\leq
\]
\[ \leq\sup\limits_{u\geq1}\frac{\left(\left(E\xi^{u(v-1)}\right)^{\frac{1}{u(v-1)}}+\left(E\eta^{u(v-1)}\right)^{\frac{1}{u(v-1)}}\right)^{v-1}}{u^{v-1}}.
\]
Consider $E\xi^{u(v-1)}=\int\limits_{0}^{+\infty}x^{u(v-1)}pe^{-px}dx$. We make a change of variables in this integral, i.e. $px=t$. Then
\[
\int\limits_{0}^{+\infty}x^{u(v-1)}pe^{-px}dx=\frac{1}{p^{u(v-1)}}\int\limits_{0}^{+\infty}e^{-t}t^{u(v-1)}dt=\frac{1}{p^{u(v-1)}}\Gamma(u(v-1)+1).
\]
Therefore, $$\left(E\xi^{u(v-1)}\right)^{\frac{1}{u(v-1)}}=\frac{1}{p}\left(\Gamma(u(v-1)+1)\right)^{\frac{1}{u(v-1)}}.$$ Similarly, we find $$\left(E\eta^{u(v-1)}\right)^{\frac{1}{u(v-1)}}=\frac{1}{q}\left(\Gamma(u(v-1)+1)\right)^{\frac{1}{u(v-1)}}.$$
Thus,
\[
\sup\limits_{u\geq1}\frac{\left(E\left((\xi+\eta)^{v-1}\right)^{u}\right)^{1/u}}{u^{v-1}}\leq \sup\limits_{u\geq1}\frac{\left(\Gamma(u(v-1)+1)\right)^{\frac{1}{u}}\left(\frac{1}{p}+\frac{1}{q}\right)^{v-1}}{u^{v-1}}.
\]
Since $\Gamma(z)\leq e^{-z}z^{z-\frac{1}{2}}C_z$, where $C_z=\sqrt{2\pi}e^{\frac{1}{12z}}$, then
\[
\left(\Gamma(u(v-1)+1)\right)^{\frac{1}{u}}\leq e^{-\left(v-1+\frac{1}{u}\right)}(u(v-1)+1)^{v-1+\frac{1}{2u}}\left(C_z\right)^{\frac{1}{u}},
\]
where $z=u(v-1)+1$. Note that $C_z\leq S=\sqrt{2\pi}e^{\frac{1}{18}}$. Hence, we have that
\begin{multline*}
\sup\limits_{u\geq1}\frac{\left(E\left((\xi+\eta)^{v-1}\right)^{u}\right)^{1/u}}{u^{v-1}}\leq\\ \leq\sup\limits_{u\geq1}\frac{e^{-(v-1)}e^{-\frac{1}{u}}(u(v-1)+1)^{v-1}(u(v-1)+1)^{\frac{1}{2u}}\left(S\right)^{\frac{1}{u}}\left(\frac{1}{p}+\frac{1}{q}\right)^{v-1}}{u^{v-1}}\leq\\
\leq e^{-(v-1)} \left(\frac{1}{p}+\frac{1}{q}\right)^{v-1}
\sup\limits_{u\geq1} \left(\frac{S}{e}\right)^{\frac{1}{u}}\frac{u^{v-1}(v-1+\frac{1}{u})^{v-1} u^{\frac{1}{2u}}(v-1+\frac{1}{u})^{\frac{1}{2u}}}{u^{v-1}} \leq\\
\leq e^{-(v-1)} \left(\frac{1}{p}+\frac{1}{q}\right)^{v-1}
\sup\limits_{u\geq1} \left(\frac{S}{e}\right)^{\frac{1}{u}} u^{\frac{1}{2u}}(v-1+\frac{1}{u})^{v-1+\frac{1}{2u}} \leq\\
\leq e^{-(v-1)+\frac{1}{2e}} v^{v-\frac{1}{2}}\left(\frac{1}{p}+\frac{1}{q}\right)^{v-1}.
\end{multline*}
That is $$\left\|c(\xi,\eta)\frac{1}{pq}(\xi+\eta)^{v-1}\right\|_\psi\leq\frac{1}{pq}e^{-(v-1)+\frac{1}{2e}} v^{v-\frac{1}{2}}\left(\frac{1}{p}+\frac{1}{q}\right)^{v-1}.$$
According to results of Example \ref{ex:o7chapter-6} we can obtain the following inequality:
\[
n\geq\left(\frac{4 (3)^{v-1} e^{\frac{1}{2e}} v^{v-\frac{1}{2}}\left(\frac{1}{p}+\frac{1}{q}\right)^{v-1}}{pq\varepsilon}\right)^2\max\left(1,\left(-\frac{\ln \delta}{v-1}\right)^{2(v-1)}\right).
\]
\end{example}

\subsection{Reliability and accuracy in the space $C(T)$ of calculation of integrals depending on the parameter}\indent

Consider the integral $$\int\limits_{\mathcal{S}}f(s,t)p(s)d\mu(s)=I(t)$$ under the condition that it exists. We use all the conditions and notation of subsection \ref{su:o7series-2}.

\begin{theorem}\label{th:o7chapter-7}
Let a random process $\xi(t)$ belong to the space $\mathbf{F}_\psi(\Omega)$ and for which the condition $\bf{H}$ with constant $C_\psi$ is satisfied, let $$\widetilde{Z}_{n}(t)=\frac{1}{n}\sum \limits_{i=1}^{n}\left(\xi_{i}(t)-I(t)\right),$$ where $\xi_i(t)$ are independent copies of the random process $\xi(t)$.
Moreover, let there exist a continuous monotonically increasing function $\sigma(h)$, $\sigma(0)=0$ such that
\[
\sup\limits_{\rho(t, s) \leq h} \bigl\|\xi(t) - \xi(s) \bigr\|_\psi \leq \sigma(h).
\]
Suppose that for any $z > 0$ the condition is satisfied
\[
\int\limits_{0}^{z} \varkappa\left(N\left(\sigma^{(-1)}(u)\right)\right) du < \infty,
\]
where $\varkappa(u)$ is the majorizing characteristic, $N(u)$ is the metric massiveness of the space $\mathbf{F}_\psi(\Omega)$, then for any $0<p<1$ the inequality holds
\[
\left\|\sup\limits_{t \in T}\left|\widetilde{Z}_n (t)\right|\right\|_\psi\leq \widehat{B}(p),
\]
where
$$\widehat{B}(p)=2\sqrt{C_\psi} \inf \limits_{t \in T}\left\|\xi(t)\right\|_\psi+\frac{1}{p(1-p)}\int\limits_{0}^{\gamma p}\varkappa\left(N\left(\sigma_1^{(-1)}(u)\right)\right)du,$$
 $\sigma_1^{(-1)}(u)$  is the inverse function ot $\sigma_1(u)$, $$\sigma_1(h)=2\sqrt{C_\psi}\sigma(h),\,\, \gamma=2\sqrt{C_\psi}\sigma\left(\sup\limits_{t, s \in T} \rho(t, s)\right).$$
In this case, $Z_{n}(t)$ approximates $I(t)$ with reliability $1-\delta$ and accuracy $\varepsilon$ in the space $C(T)$ (see Definition \ref{de:o7chapter-2}), if the number $n$ is such that the condition is true
\begin{equation}\label{eq:o7chapter-14}
\inf\limits_{u\geq1}\frac{\widehat{B}^u(p)(\psi(u))^u}{(\varepsilon\sqrt{n})^u}\leq\delta.
\end{equation}
\end{theorem}

\begin{proof}
The theorem follows from Theorem \ref{th:o2chapter_2} and Theorem~\ref{th:o7chapter-4}, and Corollary~\ref{th:o3chapter-5}.
\end{proof}

\begin{example}\label{ex:o7chapter-9}
Consider the space $\mathbf{F}_\psi(\Omega)$, where $\psi(u)=u^\alpha$, $\alpha>\frac{1}{2}$. From Example \ref{ex:o3chapter-4} and Theorem \ref{th:o1chapter_7} for $\varepsilon\geq \frac{e^\alpha\widehat{B}(p)}{\sqrt{n}}$ it follows that
\[
\inf\limits_{u\geq1}\frac{\widehat{B}^u(p)(\psi(u))^u}{(\varepsilon\sqrt{n})^u}=\exp \left\{-\frac{\alpha}{e}\left(\frac{\sqrt{n}\varepsilon}{\widehat{B}(p)}\right)^{1/\alpha}\right\},
\]
where
$$\widehat{B}(p)=4 \cdot 3^{\alpha}\inf \limits_{t \in T}\left\|\xi(t)\right\|_\psi+\frac{1}{p(1-p)}\left(\frac{e}{\alpha}\right)^\alpha\int\limits_{0}^{\gamma p}\left(\ln \left(N\left(\sigma_1^{(-1)}(u)\right)\right)\right)^\alpha du.$$
Therefore, the inequality \eqref{eq:o7chapter-14} holds true when the inequality is satisfied
\[
\exp \left\{-\frac{\alpha}{e}\left(\frac{\sqrt{n}\varepsilon}{\widehat{B}(p)}\right)^{1/\alpha}\right\}\leq \delta,
\]
under conditions
\[
n\geq\left(\frac{\widehat{B}(p)}{\varepsilon}\right)^2\left((-\ln \delta)\frac{e}{\alpha}\right)^{2\alpha}
\]
and
\[
n\geq\left(\frac{ e^{2\alpha}\widehat{B}(p)}{\varepsilon}\right)^2 \max\left(1,\left(\frac{-\ln \delta}{\alpha}\right)^{2\alpha}\right).
\]
\end{example}

\begin{example} \label{ex:o7chapter-10}
Consider an integral of the form
\[
\int\limits_{0}^{+\infty}\int\limits_{0}^{+\infty}c(t,x,y)(x+y)^{v-1}e^{-px}e^{-qy}dxdy,
\]
where $0\leq t\leq 1$, $\left|c(t,x,y)\right|\leq 1$. For $t_1, t_2 \in [0,1]$, we have
 $\left|c(t_1,x,y)-c(t_2,x,y)\right|\leq \widetilde{\gamma} \left(t_1-t_2\right)R$,
where $R>0$, $\widetilde{\gamma}(h)$, $0\leq h\leq 1$ is a monotonically increasing, continuous function such that
$\widetilde{\gamma}(0)=0$.
From Example \ref{ex:o7chapter-8} it follows that $$\left\|I(t)\right\|_\psi\leq\frac{1}{pq}e^{-v} v^{v-\frac{1}{2}}S\left(\frac{1}{p}+\frac{1}{q}\right)^{v-1},$$
where $S=\sqrt{2\pi}e^{\frac{1}{18}}$, then
\begin{multline*}
\left\|I(t_1)-I(t_2)\right\|_\psi=\left\|\left(c(t_1,\xi,\eta)-c(t_2,\xi,\eta)\right)\frac{1}{pq}(\xi+\eta)^{v-1}\right\|_\psi\leq \\[1ex] \leq\widetilde{\gamma}\left(t_1-t_2\right)R \frac{1}{pq}e^{-v} v^{v-\frac{1}{2}}S\left(\frac{1}{p}+\frac{1}{q}\right)^{v-1}.
\end{multline*}
That is, in terms of Theorem \ref{th:o7chapter-7} and Example \ref{ex:o7chapter-9}, we have that
\[
\sigma(h)=\widetilde{\gamma}(h)R \frac{1}{pq}e^{-v} v^{v-\frac{1}{2}}S\left(\frac{1}{p}+\frac{1}{q}\right)^{v-1},
\]
and
\[
\sigma_1(h)=4 \cdot 3^{v-1}R \frac{1}{pq}e^{-v} v^{v-\frac{1}{2}}S\left(\frac{1}{p}+\frac{1}{q}\right)^{v-1}\widetilde{\gamma}(h)=D(v)\widetilde{\gamma}(h),
\]
so $\sigma^{(-1)}_1(h)=\widetilde{\gamma}^{(-1)}\frac{h}{D(v)}$.
If, in addition, $N(u)\leq \frac{1}{2u}+1$, then from Example \ref{ex:o7chapter-9} it follows that $Z_{n}(t)=\frac{1}{n}\sum \limits_{i=1}^{n}I_n(t)$
approximates $I(t)$ in the space $C([0,1])$ with reliability $1-\delta$ and accuracy $\varepsilon$ if the inequality holds
\[
n\geq\left(\frac{\widehat{B}(r)}{\varepsilon}\right)^2\left((-\ln \delta)\frac{e}{v-1}\right)^{2(v-1)},
\]
where \begin{multline*}
\widehat{B}(r)=4 \cdot 3^{v-1}\frac{1}{pq}e^{-v} v^{v-\frac{1}{2}}S\left(\frac{1}{p}+
\frac{1}{q}\right)^{v-1}+\\
+\frac{1}{r(1-r)}\left(\frac{e}{v-1}\right)^{v-1}\int\limits_{0}^{\gamma r}\left(\ln \left(\frac{1}{2\widetilde{\gamma}^{(-1)}\frac{x}{D(v)}}+1\right)\right)^{v-1} dx,
\end{multline*}
$\gamma=D(v)\widetilde{\gamma}(1).$
\end{example}

An example of the function $c(t,x,y)$ can be the function
\[
c(t,x,y)=2^\theta\frac{\sin(t(x+1)(y+1))}{((x+1)(y+1))^\theta},
\]
where $0<\theta\leq 1$. Indeed,
$$
\left|c(t_1,x,y)-c(t_2,x,y)\right|=\left|\frac{2^\theta\cdot 2 \sin\left(\frac{(t_1-t_2)(x+1)(y+1)}{2}\right)\cos\left(\frac{(t_1+t_2)(x+1)(y+1)}{2}\right)}{(x+1)^\theta(y+1)^\theta}\right|.
$$
Since $\left|\sin x\right|\leq \left|x\right|^\theta$, where $0<\theta\leq 1$, we have
\[
\left|c(t_1,x,y)-c(t_2,x,y)\right|\leq 2 \left|t_1-t_2\right|^\theta.
\]

\begin{example}
Let in Example \ref{ex:o7chapter-10} $\widetilde{\gamma}(h)=L h^\beta$, then $$\widetilde{\gamma}^{(-1)}\left(\frac{x}{D(v)}\right)=\left(\frac{x}{LD(v)}\right)^{1/\beta}.$$
Thus,
\[
\int\limits_{0}^{\gamma r}\left(\ln \left(\frac{1}{2\widetilde{\gamma}^{(-1)}\frac{x}{D(v)}}+1\right)\right)^{v-1} dx=\int\limits_{0}^{\gamma r}\left(\ln \left(\frac{(L\cdot D(v))^{1/\beta}}{2x^{1/\beta}}+1\right)\right)^{v-1} dx.
\]
From the inequality $\ln(1+x)\leq \frac{x^\tau}{\tau}$, where $x>0$, $0\leq \tau \leq 1$, it follows that
\[
\left(\ln \left(\frac{(L\cdot D(v))^{1/\beta}}{2x^{1/\beta}}+1\right)\right)^{v-1} \leq \frac{1}{2^{\tau(v-1)}}\left(\frac{L\cdot D(v)}{x}\right)^{\frac{\tau(v-1)}{\beta}}\frac{1}{\tau^{v-1}}.
\]
Then the integral $$\int\limits_{0}^{\gamma r}\frac{1}{2^{\tau(v-1)}}\left(\frac{L\cdot D(v)}{x}\right)^{\frac{\tau(v-1)}{\beta}}\frac{1}{\tau^{v-1}}dx$$ will converge when $\frac{\tau(v-1)}{\beta}<1$.
\end{example}

\begin{example}
Consider the space $\mathbf{F}_\psi(\Omega)$, where $\psi(u)=e^{au^\beta}$, where $a>0$, $0<\beta<1$. From Example \ref{ex:o3chapter-5} and Theorem \ref{th:o1chapter_6} for $\varepsilon\geq \frac{e^{a(\beta+1)}\widehat{B}(p)}{\sqrt{n}}$ it follows that
\[
\inf\limits_{u\geq1}\frac{\widehat{B}^u(p)(\psi(u))^u}{(\varepsilon\sqrt{n})^u}=\exp \left\{-\frac{\beta}{a^{1/\beta}}\left(\frac{\ln \frac{\sqrt{n}\varepsilon}{\widehat{B}(p)}}{\beta+1}\right)^{\frac{\beta+1}{\beta}}\right\}.
\]
If $\frac{1}{(2a\beta)^{1/\beta}}=1$, then we have:
\begin{multline*}
\widehat{B}(p)=4e^{2^{\beta-1} a}\inf \limits_{t \in T}\left\|\xi(t)\right\|_\psi+\\[1ex]
+\frac{1}{p(1-p)}\int\limits_{0}^{\gamma p}\frac{1}{e^a}\exp \left\{S(a,\beta)\left(\ln \left(N\left(\sigma_1^{(-1)}(u)\right)\right)\right)^{\frac{\beta}{\beta+1}}\right\},
\end{multline*}
where $S(a,\beta)=(\beta a)^{\frac{1}{\beta+1}}(\beta^{-1}+1)$.
In the case $\frac{1}{(2a\beta)^{1/\beta}}>1$, we have similarly
\begin{multline*}
\widehat{B}(p)=\frac{4e^{\frac{a}{2}\left(2^\beta+1\right)-\frac{1}{4\beta}}}{(2a\beta)^{1/{4\beta}}} \inf \limits_{t \in T}\left\|\xi(t)\right\|_\psi+\\[1ex]
+\frac{1}{p(1-p)}\int\limits_{0}^{\gamma p}\frac{1}{e^a}\exp \left\{S(a,\beta)\left(\ln \left(N\left(\sigma_1^{(-1)}(u)\right)\right)\right)^{\frac{\beta}{\beta+1}}\right\}.
\end{multline*}
Therefore, the inequality \eqref{eq:o7chapter-14} is satisfied when the inequality is true
\[
\exp \left\{-\frac{\beta}{a^{1/\beta}}\left(\frac{\ln \frac{\sqrt{n}\varepsilon}{\widehat{B}(p)}}{\beta+1}\right)^{\frac{\beta+1}{\beta}}\right\}\leq \delta
\]
under conditions
\[
n\geq\left(\frac{\widehat{B}(p)}{\varepsilon}\right)^2 \exp\left\{2(\beta+1)\left((-\ln\delta)\frac{a^{1/\beta}}{\beta}\right)^{\frac{\beta}{\beta+1}}\right\}
\]
and
\[
n\geq\left(\frac{\widehat{B}(p)}{\varepsilon}\right)^2\max\left(e^{a(\beta+1)}, \exp\left\{2(\beta+1)\left((-\ln\delta)\frac{a^{1/\beta}}{\beta}\right)^{\frac{\beta}{\beta+1}}\right\}\right).
\]
\end{example}

\subsection{Reliability and accuracy in the space $L_p(T)$ of calculation of integrals depending on the parameter}\indent

This subsection retains the notation of the previous subsection.
\begin{definition}\label{de:o7chapter-3}
Let us say that $Z_{n}(t)$ approaches $I(t)$ in the space $L_p(T)$ with reliability $1-\delta>0$ and accuracy $\varepsilon>0$ if the following inequality holds:
\[
P\left\{\left(\int\limits_{T}\left|Z_n(t)-I(t)\right|^p d\mu(t)\right)^{1/p}>\varepsilon\right\}\leq \delta.
\]
\end{definition}

\begin{theorem}\label{th:o7chapter-8}
Let $I(t)=E\xi(t)=\int\limits_{\mathcal{S}}f(s,t)p(s)d\mu(s)$, let $\xi(t)$ be a random process belonging to the space $\mathbf{F}_\psi(\Omega)$ for which the condition $\bf{H}$ with the constant $C_\psi$ is satisfied, let $\widetilde{Z}_{n}(t)=\frac{1}{n}\sum \limits_{i=1}^{n}\left(\xi_{i}(t)-I(t)\right)$, and let $\xi_i(t)$ be independent copies of the random process $\xi(t)$.

Then for all $p\geq 1$ the inequality holds
\[
\left\|\left(\int\limits_{T}\left|\widehat{Z}_n(t)\right|^p d\mu(t)\right)^{1/p}\right\|\leq \frac{2 \sqrt{C_\psi}}{\sqrt{n}}\cdot \frac{\psi(p)}{\psi(1)}\left(\int\limits_{T}\left\|\xi(t)\right\|_\psi^p d\mu(t)\right)^{1/p},
\]
and $Z_{n}(t)$ approximates $I(t)$ with reliability $1-\delta$ and accuracy $\varepsilon$ in the space $L_p(T)$ for $n$ such that
\begin{equation}\label{eq:o7chapter-15}
\inf\limits_{u\geq1}\frac{\left(\frac{2 \sqrt{C_\psi}}{\sqrt{n}}\cdot\frac{\psi(p)}{\psi(1)}\right)^u\left(\int\limits_{T}\left\|\xi(t)\right\|^pd\mu(t)\right)^{u/p}(\psi(u))^u}{\varepsilon^u}\leq\delta.
\end{equation}
\end{theorem}

\begin{proof}
The theorem follows from Theorem \ref{th:o6chapter-2} if the inequality \eqref{eq:o7chapter-15} holds.
\end{proof}

\begin{example}
Consider the space $\mathbf{F}_\psi(\Omega)$, where $\psi(u)=u^\alpha$, $\alpha>\frac{1}{2}$. Then from Example \ref{ex:o1chapter_4} it follows that for this space the condition $\bf{H}$ with the constant $C_\psi=4 \cdot9^\alpha$ is satisfied, and from Example \ref{ex:o6chapter-4} and Theorem \ref{th:o7chapter-8} for
$$\varepsilon\geq\frac{4 (3pe)^\alpha\left(\int\limits_{T}\left\|\xi(t)\right\|_\psi^p d\mu(t)\right)^{1/p}}{\sqrt{n}}$$ it follows that
\begin{multline*}
\inf\limits_{u\geq1}\frac{\left(\frac{2 \sqrt{C_\psi}}{\sqrt{n}}\cdot\frac{\psi(p)}{\psi(1)}\right)^u\left(\int\limits_{T}\left\|\xi(t)\right\|^pd\mu(t)\right)^{u/p}(\psi(u))^u}{\varepsilon^u}\leq \\[1ex]
\leq\exp \left\{-\frac{\alpha}{e}\left(\frac{\sqrt{n}\varepsilon}{4 (3pe)^\alpha\left(\int\limits_{T}\left\|\xi(t)\right\|_\psi^p d\mu(t)\right)^{1/p}}\right)^{1/\alpha}\right\}.
\end{multline*}
Therefore, the inequality \eqref{eq:o7chapter-15} is satisfied when the following inequality is true
\[
\exp \left\{-\frac{\alpha}{e}\left(\frac{\sqrt{n}\varepsilon}{4 (3pe)^\alpha\left(\int\limits_{T}\left\|\xi(t)\right\|_\psi^p d\mu(t)\right)^{1/p}}\right)^{1/\alpha}\right\}\leq \delta,
\]
under condition
\[
n\geq\left(\frac{4 (3pe)^\alpha\left(\int\limits_{T}\left\|\xi(t)\right\|_\psi^p d\mu(t)\right)^{1/p}}{\varepsilon}\right)^2\left((-\ln \delta)\frac{e}{\alpha}\right)^{2\alpha}.
\]
Then
\[
n\geq\left(\frac{4 (3p)^\alpha\left(\int\limits_{T}\left\|\xi(t)\right\|_\psi^p d\mu(t)\right)^{1/p}}{\varepsilon}\right)^2 \max\left(1,\left(-\frac{\ln \delta}{\alpha}\right)^{2\alpha}\right).
\]
\end{example}

\begin{example}
Consider the space $\mathbf{F}_\psi(\Omega)$, where $\psi(u)=e^{au^\beta}$, where $a>0$, $0<\beta<1$. Then from Theorem \ref{th:o2chapter-14} it follows that two cases are possible.
In the first case, for $\frac{1}{(2a\beta)^{1/\beta}}=1$, the condition $\mathbf{H}$ with the constant $C_\psi=4e^{2^\beta a}$ is satisfied for the space $\mathbf{F}_\psi(\Omega)$.
Then from Example \ref{ex:o6chapter-5} and Theorem \ref{th:o7chapter-8} it follows that
\begin{multline*}
\inf\limits_{u\geq1}\frac{\left(\frac{2 \sqrt{C_\psi}}{\sqrt{n}}\cdot\frac{\psi(p)}{\psi(1)}\right)^u\left(\int\limits_{T}\left\|\xi(t)\right\|^pd\mu(t)\right)^{u/p}(\psi(u))^u}{\varepsilon^u}\leq \\[1ex]
\leq\exp \left\{-\frac{\beta}{a^{1/\beta}}\left(\frac{\ln \frac{\sqrt{n}\varepsilon}{4 e^{a\left(2^{\beta-1}+p^\beta-1\right)}\left(\int\limits_{T}\left\|\xi(t)\right\|_\psi^p d\mu(t)\right)^{1/p} }}{\beta+1}\right)^{\frac{\beta+1}{\beta}}\right\}.
\end{multline*}
Therefore, the inequality \eqref{eq:o7chapter-15} is satisfied when the following inequality is true
\[
\exp \left\{-\frac{\beta}{a^{1/\beta}}\left(\frac{\ln \frac{\sqrt{n}\varepsilon}{4 e^{a\left(2^{\beta-1}+p^\beta-1\right)}\left(\int\limits_{T}\left\|\xi(t)\right\|_\psi^p d\mu(t)\right)^{1/p} }}{\beta+1}\right)^{\frac{\beta+1}{\beta}}\right\}\leq \delta.
\]
It follows that
\[
n\geq\left(\frac{4 e^{a\left(2^{\beta-1}+p^\beta-1\right)}\left(\int\limits_{T}\left\|\xi(t)\right\|_\psi^p d\mu(t)\right)^{1/p}}{\varepsilon}\right)^2\times
\]
\[
\times\exp\left\{2(\beta+1)\left((-\ln\delta)\frac{a^{1/\beta}}{\beta}\right)^{\frac{\beta}{\beta+1}}\right\}.
\]
Then
\begin{multline*}
n\geq\left(\frac{4 e^{a\left(2^{\beta-1}+p^\beta-1\right)}\left(\int\limits_{T}\left\|\xi(t)\right\|_\psi^p d\mu(t)\right)^{1/p}}{\varepsilon}\right)^2 \times\\
\times\max\left(e^{a(\beta+1)}, \exp\left\{2(\beta+1)\left((-\ln\delta)\frac{a^{1/\beta}}{\beta}\right)^{\frac{\beta}{\beta+1}}\right\}\right).
\end{multline*}
In the second case, where $\frac{1}{(2a\beta)^{1/\beta}}>1$, we have
\begin{multline*}
\inf\limits_{u\geq1}\frac{\left(\frac{2 \sqrt{C_\psi}}{\sqrt{n}}\cdot\frac{\psi(p)}{\psi(1)}\right)^u\left(\int\limits_{T}\left\|\xi(t)\right\|^pd\mu(t)\right)^{u/p}(\psi(u))^u}{\varepsilon^u}\leq \\[1ex]
\leq\exp \left\{-\frac{\beta}{a^{1/\beta}}\left(\frac{\ln \frac{\sqrt{n}\varepsilon(2a\beta)^{\frac{1}{4\beta}}}{\frac{4e^{a\left(2^{\beta-1}+p^\beta-\frac{1}{2}\right)-\frac{1}{4\beta}-1}}{(2a\beta)^1/{4\beta}} \left(\int\limits_{T}\left\|\xi(t)\right\|_\psi^p d\mu(t)\right)^{1/p} }}{\beta+1}\right)^{\frac{\beta+1}{\beta}}\right\}
\end{multline*}
and the inequality \eqref{eq:o7chapter-15} holds when the following estimate is true
\[
\exp \left\{-\frac{\beta}{a^{1/\beta}}\left(\frac{\ln \frac{\sqrt{n}\varepsilon(2a\beta)^{\frac{1}{4\beta}}}{\frac{4e^{a\left(2^{\beta-1}+p^\beta-\frac{1}{2}\right)-\frac{1}{4\beta}}}{(2a\beta)^1/{4\beta}} \left(\int\limits_{T}\left\|\xi(t)\right\|_\psi^p d\mu(t)\right)^{1/p} }}{\beta+1}\right)^{\frac{\beta+1}{\beta}}\right\}\leq \delta
\]
under condition
\begin{multline*}
n\geq\left(\frac{\frac{4e^{a\left(2^{\beta-1}+p^\beta-\frac{1}{2}\right)-\frac{1}{4\beta}}}{(2a\beta)^1/{4\beta}} \left(\int\limits_{T}\left\|\xi(t)\right\|_\psi^p d\mu(t)\right)^{1/p}}{(2a\beta)^{\frac{1}{4\beta}}\varepsilon}\right)^2\times\\ \times\exp\left\{2(\beta+1)\left((-\ln\delta)\frac{a^{1/\beta}}{\beta}\right)^{\frac{\beta}{\beta+1}}\right\}.
\end{multline*}
Then
\begin{multline*}
n\geq\left(\frac{\frac{4e^{a\left(2^{\beta-1}+p^\beta-\frac{1}{2}\right)-\frac{1}{4\beta}}}{(2a\beta)^1/{4\beta}} \left(\int\limits_{T}\left\|\xi(t)\right\|_\psi^p d\mu(t)\right)^{1/p}}{(2a\beta)^{\frac{1}{4\beta}}\varepsilon}\right)^2\times\\ \times \max \left(e^{a(\beta+1)},\exp\left\{2(\beta+1)\left((-\ln\delta)\frac{a^{1/\beta}}{\beta}\right)^{\frac{\beta}{\beta+1}}\right\}\right).
\end{multline*}
\end{example}

\section*{Conclusions to chapter \ref{ch:o7series}}

In chapter \ref{ch:o7series}, the accuracy and reliability of the calculation of multiple integrals by the Monte Carlo method are found. Integrals that depend on the parameter are considered. The reliability and accuracy of the calculation of these integrals in $C(T)$ and $L_p(T)$ are found.

\chapter{Accuracy and reliability of modeling in $L_p(T)$ spaces of random processes that allow series expansions}
\label{omch:lpt1}

In this chapter we deal with the modeling of random processes from $Sub_\varphi(\Omega)$ spaces, which are a subclass of $K_\sigma$-spaces, with given reliability and accuracy in $L_p(T)$ spaces.

\section{Accuracy and reliability of modeling random processes in $L_p(T)$ spaces}

Let $\left(\Omega,{\cal F}, P \right)$ be a standard probability space, let $ L^{o}_2(\Omega)$ be the space of centered random variables with finite second moment,
$E\xi^2<\infty$, and let $\{\Lambda, {\cal U}, \mu\}$ be a measurable space, the measure $\mu$ is $\sigma$-finite.
Let $L_p(\Lambda,\mu)$ be the Banach space of functions integrable to the power of $p$ with measure $\mu$.

\begin{definition}
Let $X=\{X(t),t\in T\}$ be a random process from $Sub_\varphi(\Omega)$. A random process $X_N=\{X_N(t),t\in T\}$ from $Sub_\varphi(\Omega)$ will be called a model that approximates a random process $X$ with given reliability $1-\alpha$ and accuracy $\delta$ in the space $L_p(T)$ if
$$P\left\{\int_T( X(t)-X_N(t))^pdt)^{1/p}>\delta\right\}\leq\alpha.$$
\end{definition}

The following statement is proved in \cite{koz-kam}.

\begin{theorem}
\label{omkamn}
Let $\{T,\Lambda, M\}$ be a measurable space,  let $X=\{X(t),t\in T\}$ be a measurable random process from the space $Sub_\varphi (\Omega)$,
and let $\tau_\varphi(t)=\tau_\varphi(X(t))$.
Let the Lebesgue integral $\int_T(\tau_\varphi (t))^p d\mu (t)$ be well defined for $p\geq 1$.
Then with probability 1 there exists $\int_T |X(t)|^p d\mu(t)$ and  the following inequality holds true
$$P\left\{\int_T |X(t)|^p d\mu(t) >\delta\right\} \leq 2\exp\left\{-\varphi^* \left(\left(\frac{\delta}{c}\right)^{1/p}\right)\right\}$$
$\forall \delta:$
$\delta> c(f( \frac {c^{1/p}p}{\delta^{1/p}}))^p$, where $c=\int_0^T
(\tau_\varphi(t))^pd\mu(t)$, $f$ is a function such that  $\varphi(u)=\int_0^u f(v)dv$ for all $u>0$, and
where $\varphi^*$ is the Young-Fenchel transform of the function $\varphi$.
\end{theorem}

This theorem can be transformed to the following form.

\begin{theorem}
\label{omkamnd}
Let $\{T,\Lambda, M\}$ be a measurable space,  let $X=\{X(t),t\in T\}$ be a measurable random process from the space $Sub_\varphi (\Omega)$.
Let $f$ be a function such that  $\varphi(u)=\int_0^u f(v)dv$ for all $u>0$.
Let $$c_N=\int_T (\tau_\varphi(X(t)-X_N(t)))^pd\mu(t)<\infty.$$ The model $X_N(t)$
approximates the random process $X(t)$ with reliability $1-\alpha$ and accuracy $\delta$ in the space $L_p(T)$, when
\begin{equation}
\label{omcn-first}
c_N\leq \frac{\delta}{(\varphi^{*(-1)}(ln \frac{2}{\alpha}))^p}
\end{equation}
\noindent and
\begin{equation}
\label{omcn-second}
\delta> c_N\left(f\left( \frac {c_N^{1/p}p}{\delta^{1/p}}\right)\right)^p,
\end{equation}
\noindent where $\varphi^*$ is the Young-Fenchel transform of the function $\varphi$.
\end{theorem}

\begin{dov}

Since $|X(t)-X_N(t)|\in Sub_\varphi(\Omega)$, and by Theorem \ref{omkamn}
$$P\left\{\int_T |X(t)|^p d\mu(t) >\delta\right\} \leq 2\exp\left\{-\varphi^* \left(\left(\frac{\delta}{c}\right)^{1/p}\right)\right\},$$
\noindent then
$$P\left\{\int_T |X(t)-X_N(t)|^p d\mu(t) >\delta\right\} \leq 2\exp\left\{-\varphi^* \left(\left(\frac{\delta}{c_N}\right)^{1/p}\right)\right\},$$
\noindent therefore the inequality (\ref{omcn-second}) must hold and
$$2\exp\left\{-\varphi^* \left(\left(\frac{\delta}{c_N}\right)^{1/p}\right)\right\}\leq \alpha,$$
\noindent whence
$$\varphi^* \left(\left(\frac{\delta}{c_N}\right)^{1/p}\right)\geq ln \frac{2}{\alpha},$$
\noindent and, finally,
$$c_N^{1/p}\leq \frac{\delta^{1/p}}{\varphi^{*(-1)}(ln \frac{2}{\alpha})}.$$
\end{dov}

Based on this Theorem \ref{omkamnd}, we can prove the following statements.

\begin{theorem}
\label{lp_l2}
Let the random process $X=\{X(t),t\in [0,T]\}$ belong to $S u b_\varphi (\Omega)$,
$$\varphi(t)=\frac{t^\gamma}{\gamma},$$
\noindent for $1<\gamma\leq 2$. Let $$c_N=\int_0^T (\tau_\varphi(X(t)-X_N(t)))^pd\mu(t)<\infty.$$
The model $X_N(t)$ approximates the random process $X(t)$ with reliability $1-\alpha$ and accuracy $\delta$ in the space $L_p(T)$ when
$$\left\{\begin{array}{c}c_N\leq \delta/(\beta ln \frac{2}{\alpha})^{p/\beta}\\ c_N<\delta /p^{p\left(1-1/\gamma\right)}\end{array}\right..$$
Here $\beta$ is a number such that $\frac{1}{\beta}+\frac{1}{\gamma}=1$.

\end{theorem}

\begin{dov}
The first inequality follows immediately from Theorem \ref{omkamnd}. Indeed, since $\varphi(t)=t^\gamma/\gamma$, then $\varphi^*(t)=t^\beta/\beta$, which yields the desired result.

Now consider the second inequality. Again, since $\varphi(t)=\frac{t^\gamma}{\gamma}$, then $f(t)=t^{\gamma-1}$, $t>0$, and

$$\delta>c_N\left(\left(\frac{c_N^{1/p}p}{\delta^{1/p}}\right)^{\gamma-1}\right)^p = \frac{c_N^\gamma p^{p(\gamma-1)}}{\delta^{\gamma-1}},$$

\noindent therefore

$$c_N^\gamma<\frac{\delta^\gamma}{p^{p(\gamma-1)}}.$$

\end{dov}

\begin{theorem}
\label{lp_g2}
Let a random process $X=\{X(t),t\in [0,T]\}$ belong to $S u b_\varphi (\Omega)$,

$$\varphi(t)=\left\{\begin{array}{c}\frac{t^2}{\gamma}, t<1\\ \frac{t^\gamma}{\gamma},t\geq1\end{array}\right.,$$

\noindent where $\gamma>2$. Let $$c_N=\int_0^T (\tau_\varphi(X(t)-X_N(t)))^pd\mu(t)<\infty.$$
The model $X_N(t)$ approximates a random process $X(t)$ with reliability $1-\alpha$ and accuracy $\delta$ in the space $L_p[0,T]$, when

$$\left\{\begin{array}{c}c_N\leq \delta/(\beta ln \frac{2}{\alpha})^{p/\beta} \\ c_N<\delta/p^{p(1-1/\gamma)}\end{array}\right.,$$

\noindent where $\beta$ is such that $1/\beta+1/\gamma=1$.

\end{theorem}

\begin{dov}

Let's find the form of the function $f^{(-1)}(t)$. Since

$$\varphi(t)=\left\{\begin{array}{c}\frac{t^2}{\gamma}, t<1\\ \frac{t^\gamma}{\gamma},t\geq1\end{array}\right.,$$

\noindent then

$$f(t)=\left\{\begin{array}{c}\frac{2}{\gamma}t, t<1\\ t^{\gamma-1},t\geq1\end{array}\right.,$$

\noindent therefore

$$f^{(-1)}(t)=\left\{\begin{array}{c}\frac{\gamma}{2}t, t<\frac{2}{\gamma}\\ t^\frac{1}{\gamma-1},t\geq1\end{array}\right..$$

Now consider $\varphi^*(t)$. When $t>1$,

$$\varphi^*(t) = \int_0^{2/\gamma}\frac{\gamma}{2}udu+\int_{2/\gamma}^1du +\int_1^t u^\frac{1}{\gamma-1}du = \left.\frac{\gamma}{2}\frac{u^2}{2}\right|_0^{\gamma/2}+(1-\frac{2}{\gamma})+\left.u^{\frac{1}{\gamma-1}+1}\right|_1^t = $$

$$ = \frac{\gamma}{2}\frac{1}{2}\left(\frac{2}{\gamma}\right)^2+1-\frac{2}{\gamma}+\frac{t^\beta}{\beta}-\frac{1}{\beta} = \frac{t^\beta}{\beta}.$$

\noindent From this relation and from Theorem \ref{omkamnd} the first inequality of the theorem follows.

Now let us consider the second inequality. Let us first analyze the case where $(c_N^{1/p}p)/(\delta^{1/p})> 1$. In this case $f(x)=x^{\gamma-1}$, and

$$\delta>c_N\left(\left(\frac{c_N^{1/p}p}{\delta^{1/p}}\right)^{\gamma-1}\right)^p,$$

\noindent that is

$$c_N<\frac{\delta}{p^{p(1-1/\gamma)}},$$

\noindent therefore

$$\frac{\delta}{p^p}<c_N<\frac{\delta}{p^{p(1-1/\gamma)}}.$$

\noindent For the case $(c_N^{1/p}p)/(\delta^{1/p})< 1$ we have $f(x)=x\frac{2}{\gamma}$, and

$$\left\{\begin{array}{c} c_N<\frac{\delta}{p^p}\\ c_N<\frac{\delta}{p^p/2}\left(\frac{\gamma}{2}\right)^{p/2}. \end{array}\right.$$

\noindent since $\gamma>2$ and $p>1$, we get

$$c_N<\frac{\delta}{p^p}.$$

\noindent Finally, we have that

$$c_N<\frac{\delta}{p^{p(1-1/\gamma)}}.$$

\end{dov}

\section[Modeling of stochastic processes that admit series expansion] {Accuracy and reliability in $L_p(T)$
modeling of stochastic processes that admit series expansion with independent terms}

Let the stochastic process $X = \{X(t),t\in [0,T]\}$ be represented as a series
\begin{equation}
\label{omrozkld}
X(t) = \sum_{k=1}^\infty \xi_k a_k(t),
\end{equation}
\noindent where $\xi_k\in Sub_\varphi(\Omega)$, random variables $\xi_k$ are independent, and for this series the condition is satisfied
$$\sum_{k=1}^\infty \tau_\varphi^2(\xi_k)a^2_k(t)<\infty.$$
Usually, the sum of the first $N$ terms of such series is used as a model of such a random process.
However, it often happens that the functions $a_k(t)$ cannot be found explicitly.
In this case, as elements of the random process model, we can use $\hat{a}_k(t)$ which is an approximate values of the functions $a_k(t)$,
taking into account the influence of the error of such an approximation on the accuracy and reliability with which the model approximates the process.

\begin{definition}\label{ommodel-def} The model of the random process $X(t)$ will be called the expression
\begin{equation}
\label{omrozkld-xn}
X_N(t) = \sum_{k=1}^N \xi_k \hat{a}_k(t),
\end{equation}
\noindent where $\hat{a}_k(t)$ are approximate values of the functions $a_k(t)$, and $\xi_k\in Sub_\varphi(\Omega)$, $\xi_k$ are independent.
\end{definition}

Let us introduce the following notations:
$$\delta_k(t) = (a_k(t)-\hat{a}_k(t));$$
$$\Delta_N(t) = |X(t)-X_N(t)| = \left| \sum_{k=1}^N \xi_k \delta_k(t) + \sum_{k=N+1}^\infty \xi_k a_k(t) \right|.$$

We will say that the model $X_N$ approximates a random process $X$ with given reliability and accuracy in the space $L_p[0,T]$ if
$$P\left\{\left(\int_0^T(\Delta_N(t))^pdt\right)^{\frac{1}{p}}>\delta\right\} \leq \alpha.$$

Let us formulate a general theorem for assessing the reliability and accuracy of modeling such processes in $L_p[0,T]$.

\begin{theorem}
\label{omkamnd-lp} Let a random process $X=\{X(t),t\in [0,T]\}$ of the form \eqref{omrozkld} belong to $S u b_\varphi (\Omega)$, let $X_N$ of the form \eqref{omrozkld-xn} be a model of the process $X$.
Let
$$c_N=\int_0^T\left(\tau_\varphi\left(\sum_{k=1}^N \xi_k\delta_k(t) + \sum_{k=N+1}^\infty\xi_k a_k(t)  \right)\right)^{p} dt<\infty.$$
The model $X_N(t)$ approximates the random process $X(t)$ with reliability $1-\alpha$ and accuracy $\delta$ in the space $L_p(T)$ when

$$\left\{ \begin{array}{c}c_N\leq \delta/(\varphi^{*(-1)}(ln \frac{2}{\alpha}))^p\\ \delta> c_N(f( \frac {c_N^{1/p}p}{\delta^{1/p}}))^p\end{array}\right.,$$

\noindent where  $f$ is defined in Theorem $\ref{omkamn}$, and where $\varphi^*$ is the Young-Fenchel transformation of the function $\varphi$.
\end{theorem}

\begin{dov} The proof of this theorem is similar to the proof of Theorem \ref{omkamnd}.
\end{dov}

For further theorems we will need the following statement \cite{byld98}.

\begin{theorem}\label{omtau-square}
Let $\xi_1, \xi_1,\ldots,\xi_n$ $\in Sub_\varphi(\Omega)$ be independent random variables. If the function $\varphi(|x|^{1/s})$, $x\in \mathbb R$, is convex for $s\in(0,2]$, then
$$\tau_\varphi^s\left(\sum_{k=1}^n\xi_k\right)\leq \sum_{k=1}^n \tau_\varphi^s(\xi_k).$$
\end{theorem}

\begin{theorem}
\label{omlp_g2}
Let the random process $X=\{X(t),t\in [0,T]\}$ belong to $S u b_\varphi (\Omega)$,
$$\varphi(t)=\left\{\begin{array}{c}\frac{t^2}{\gamma}, t<1\\ \frac{t^\gamma}{\gamma},t\geq1\end{array}\right.,$$
\noindent where $\gamma>2$. Let $$c_N=\int_0^T\left(\sum_{k=1}^N \tau^2_\varphi(\xi_k)\delta^2_k(t) + \sum_{k=N+1}^\infty\tau^2_\varphi(\xi_k)a^2_k(t) \right)^{p/2} dt<\infty.$$
The model $X_N$ of the form \eqref{omrozkld-xn} approximates a random
process $X$ of the form \eqref{omrozkld} with reliability $1-\alpha$ and accuracy $\delta$ in the space
$L_p(0,T)$, when
$$\left\{\begin{array}{c}c_N\leq \delta/(\beta ln \frac{2}{\alpha})^{p/\beta}\\ c_N<\delta/p^{p(1-1/\gamma)}\end{array}\right.,$$

\noindent where $\beta$ is such that $1/\beta+1/\gamma=1$.
\end{theorem}

\begin{dov} The following inequalities follow from Theorem \ref{omtau-square} for $s=2$:
\begin{eqnarray*}
c_N&=&\int_0^T
(\tau_\varphi(\Delta_N(t)))^pd\mu(t) =
\\
&=&\int_0^T\left( \tau_\varphi\left(\sum_{k=1}^N\xi_k\delta_k(t)+\sum_{k=N+1}^\infty \xi_ka_k(t)\right) \right)^pdt
 \\
& \leq&\int_0^T\left(\sum_{k=1}^N \tau^2_\varphi(\xi_k)\delta^2_k(t) + \sum_{k=N+1}^\infty\tau^2_\varphi(\xi_k)a^2_k(t)  \right)^{p/2} dt.
 \end{eqnarray*}
\end{dov}

\begin{theorem}
\label{omlp_g1}
Let a random process $X=\{X(t),t\in [0,T]\}$ belong to $S u b_\varphi (\Omega)$,
$$\varphi(t)= \frac{t^\gamma}{\gamma},$$
\noindent where $1<\gamma<2$. Let $$c_N=\int_0^T\left(\sum_{k=1}^N \tau^{\gamma}_\varphi(\xi_k)\delta^{\gamma}_k(t) + \sum_{k=N+1}^\infty\tau^{\gamma}_\varphi(\xi_k)a^{\gamma}_k(t) \right)^{p/\gamma} dt<\infty.$$
The model $X_N$ of the form \eqref{omrozkld-xn} approximates a random
process $X$ of the form \eqref{omrozkld} with reliability $1-\alpha$ and accuracy $\delta$ in the space
$L_p(0,T)$, when
$$\left\{\begin{array}{c}c_N\leq \delta/(\beta ln \frac{2}{\alpha})^{p/\beta}\\ c_N<\delta/p^{p(1-1/\gamma)}\end{array}\right.,$$
\noindent where $\beta$ is such that $1/\beta+1/\gamma=1$.
\end{theorem}

\begin{dov}
The following inequalities follow from Theorem \ref{omtau-square} for $s=\gamma$:
\begin{multline*}
c_N=\int_0^T
(\tau_\varphi(\Delta_N(t)))^pd\mu(t) =
\\
=\int_0^T\left( \tau_\varphi\left(\sum_{k=1}^N\xi_k\delta_k(t)+\sum_{k=N+1}^\infty \xi_ka_k(t)\right) \right)^pdt \leq
 \\
 \leq\int_0^T\left(\sum_{k=1}^N \tau^\gamma_\varphi(\xi_k)\delta^\gamma_k(t) + \sum_{k=N+1}^\infty\tau^\gamma_\varphi(\xi_k)a^\gamma_k(t)  \right)^{p/\gamma} dt.
\end{multline*}
\end{dov}

\subsection{Modeling of stochastic processes using the Karhunen-Lo\`eve decomposition in $L_p[0,T]$}

In this section we will use the Karhunen-Lo\`eve decomposition of the second order centered stochastic processes (see \cite{koz-roz-turch} for more details).

\begin{theorem}\label{omkl-main}
Let $X(t)$, $t\in T$ be a second-order stochastic process, $EX(t)=0$, $ t\in T$, with the correlation function
$B(t,s)=EX(t)\overline{X(s)}$. Let $f(t,\lambda)$, $ t\in T$, be functions from the space $L_2(\Lambda,\mu)$.
The correlation function $B(t,s)$ admits the representation
$$B(t,s) = \int_\Lambda f(t,\lambda) f(s,\lambda) d\mu(\lambda)$$
\noindent if and only if the process $X(t)$,  $ t\in T$, is represented in the form
\begin{equation}
\label{omkl-decompos-x}
X(t)=\sum_{k=1}^\infty \frac{a_k(t)}{\sqrt{\lambda_k}}\xi_k,
\end{equation}
\noindent where $a_k(t)$, $k=1,2,\ldots$ are the eigenfunctions of the homogeneous Fredholm integral equation of the second kind
\begin{equation}
\label{omkl-eq}
a(t) = \lambda\int_D B(t,s) a(s) d\mu(s),
\end{equation}
\noindent $\xi_k$ are centered uncorrelated random variables, $E\xi_k=0$, $E\xi_n\xi_m=\delta_{nm}\lambda_n^{-1}$, and $\lambda_n$ are eigenvalues of the given equation \eqref{omkl-eq}, numbered in ascending order.
\end{theorem}

\begin{definition}
\label{ommodel-def-kl}
Let $X=\{X(t), t\in [0,T]\}$, $EX(t)=0$, $ t\in T$, be a random process from the space $Sub_\varphi(\Omega)$ with the correlation function $B(t,s)=EX(t)\overline{X(s)}$,
which can be represented using the Karhunen-Lo\`eve expansion \eqref{omkl-decompos-x}.
The Karhunen-Lo\`eve model of such a random process will be called the process $X_N=\{X_N(t),t\in B\}$ from the space $Sub_\varphi(\Omega)$, which has the representation
\begin{equation}
\label{omkl-decompos-x-n}
X_N(t)=\sum_{k=1}^N \frac{\hat{a}_k(t)}{\sqrt{\hat{\lambda}}_k}\xi_k,
\end{equation}
\noindent where $\hat{a}_k(t)$ is the approximation of the eigenfunctions of the homogeneous Fredholm integral equation of the second kind (\ref{omkl-eq}),
and where $\hat{\lambda}_k$ are approximations of the eigenvalues of the same equation (\ref{omkl-eq}).
\end{definition}

Let us formulate theorems that will allow us to build models of stochastic processes with such a decomposition in the space $L_p(T)$.

\begin{theorem}
\label{omkl-lp-square}
Let the random process $X=\{X(t),t\in [0,T]\} $ belong to $S u b_\varphi (\Omega)$, let
$$\varphi(t)=\left\{\begin{array}{c}\frac{t^2}{\gamma}, t<1\\ \frac{t^\gamma}{\gamma},t\geq1\end{array}\right.$$
\noindent for $\gamma>2$. Let the process $X(t)$ admits an representation of the form (\ref{omkl-decompos-x}). Let
$$c_N= \int_0^T\left(\sum_{k=1}^N \tau^2_\varphi(\xi_k)\left(\frac{\delta^2_k(t) }{\hat{\lambda}_k-\eta_k} +\hat{a}^2_k(t)\frac{(\sqrt{\hat{\lambda}_k}-\sqrt{\hat{\lambda}_k-\eta_k})^2}{\hat{\lambda}_k(\hat{\lambda}_k-\eta_k)}\right)+\right.$$
$$\left.+\sum_{k=N+1}^\infty\frac{\tau^2_\varphi(\xi_k)a^2_k(t)}{\lambda_k} \right)^{p/2}dt<\infty,$$
\noindent where $\delta_k(t)=|\varphi_k(t)-\hat{\varphi}_k(t)|$ is the error in finding the $k$-th eigenfunction of the equation (\ref{omkl-eq}),
 $\hat{\lambda}_k$ is the approximate value of the $k$-th eigenvalue, and $\eta_k=|\lambda_k-\hat{\lambda}_k|$ is the error in approximating the $k$-th eigenvalue.
 The $X_N$ model \eqref{omkl-decompos-x-n} approximates the random process $X(t)$ with reliability $1-\alpha$ and accuracy $\delta$ in the space $L_p[0,T]$, if the following conditions are fulfilled:
$$\left\{\begin{array}{c}c_N\leq \delta/(\beta ln \frac{2}{\alpha})^{p/\beta} \\ c_N<\delta/p^{p(1-1/\gamma)}
\end{array}\right..$$
\end{theorem}

\begin{dov}
The statement of this Theorem follows from Theorem \ref{omlp_g2} and Theorem \ref{omkl-main}. Indeed,
$$c_N=\int_0^T\left(\tau_\varphi\left(\sum_{k=1}^N \xi_k\left(\frac{a_k(t)}{\sqrt{\lambda_k}}-\frac{\hat{a}_k(t)}{\sqrt{\hat{\lambda}_k}}\right) + \sum_{k=N+1}^\infty\xi_k\frac{a_k(t)}{\sqrt{\lambda_k}}  \right)\right)^{p} dt = $$
$$ = \int_0^T\left(\tau_\varphi\left(\sum_{k=1}^N \xi_k\left(\frac{a_k(t)}{\sqrt{\lambda_k}}-\frac{\hat{a}_k(t)}{\sqrt{\lambda_k}} + \frac{\hat{a}_k(t)}{\sqrt{\lambda_k}} -\frac{\hat{a}_k(t)}{\sqrt{\hat{\lambda}_k}}\right) + \right.\right.$$
$$\left.\left.+ \sum_{k=N+1}^\infty\xi_k\frac{a_k(t)}{\sqrt{\lambda_k}}  \right)\right)^{p} dt= \int_0^T\left(\tau_\varphi\left(\sum_{k=1}^N \xi_k\frac{1}{\sqrt{\lambda_k}}\delta_k(t) +\right.\right.$$
$$\left.\left.+\sum_{k=1}^N\xi_k\hat{a}_k(t)\left(\frac{1}{\sqrt{\lambda_k}}-\frac{1}{\sqrt{\hat{\lambda}_k}}\right)  + \sum_{k=N+1}^\infty\tau_\varphi(\xi_k)\frac{a_k(t)}{\sqrt{\lambda_k}}  \right)\right)^{p} dt\leq$$
$$\leq \int_0^T\left(\sum_{k=1}^N \tau^2_\varphi(\xi_k)\left(\frac{\delta^2_k(t) }{\hat{\lambda}_k-\eta_k} +\hat{a}^2_k(t)\frac{(\sqrt{\hat{\lambda}_k}-\sqrt{\hat{\lambda}_k-\eta_k})^2}{\hat{\lambda}_k(\hat{\lambda}_k-\eta_k)}\right)+\right.$$
$$\left.+\sum_{k=N+1}^\infty\frac{\tau^2_\varphi(\xi_k)a^2_k(t)}{\lambda_k} \right)^{p/2}dt.$$

\end{dov}

In the case where $\tau_\varphi(\xi_k)\leq\tau$, $\forall k$, we obtain the following statement.

\begin{theorem}
\label{omlp-merser}
Let a random process $X=\{X(t),t\in [0,T]\} $ belong to $S u b_\varphi (\Omega)$,
$$\varphi(t)=\left\{\begin{array}{c}\frac{t^2}{\gamma}, t<1\\ \frac{t^\gamma}{\gamma},t\geq1\end{array}\right.$$
\noindent for $\gamma>2$, and let $\forall k:$ $\tau_\varphi(\xi_k)=\tau$ . Let the process $X$ admit a representation of the form (\ref{omkl-decompos-x}), and let
$$c_N=\tau^{p/2}\int_0^T\left(\left(B(t,t)-\sum_{k=1}^N \frac{(\hat{a}_k(t)-\delta_k(t))^2}{\hat{\lambda}_k+\eta_k}\right)+ \sum_{k=1}^N \left(\frac{\delta^2_k(t) }{\hat{\lambda}_k-\eta_k} + \right.\right.$$
$$\left.\left.+\hat{a}^2_k(t)\frac{(\sqrt{\hat{\lambda}_k}-\sqrt{\hat{\lambda}_k-\eta_k})^2}{\hat{\lambda}_k(\hat{\lambda}_k-\eta_k)}\right)\right)^{p/2} dt<\infty,$$
\noindent where $\delta_k(t)=|\varphi_k(t)-\hat{\varphi}_k(t)|$ is the error in finding the $k$-th eigenfunction of the equation (\ref{omkl-eq}), $\hat{\lambda}_k$ is the approximate value of the $k$-th eigenvalue, $\eta_k=|\lambda_k-\hat{\lambda}_k|$ is the error in approximating the $k$-th eigenvalue.
The  $X_N$ model \eqref{omkl-decompos-x-n} approximates the random process $X$ with reliability $1-\alpha$ and accuracy $\delta$ in the space
$L_p(0,T)$, if the following conditions are fulfilled:
$$\left\{\begin{array}{c}c_N\leq \delta/(\beta ln \frac{2}{\alpha})^{p/\beta} \\ c_N<\delta/p^{p(1-1/\gamma)}\end{array}\right..$$
\end{theorem}

\begin{dov}
The statement of the Theorem follows from Mercer's Theorem\cite{trikomi-1960} and the previous Theorem. Indeed, according to Mercer's Theorem

$$B(t,t)=\sum_{k=1}^\infty \frac{a^2_k(t)}{\lambda_k},$$

\noindent therefore

$$c_N= \int_0^T\left(\sum_{k=1}^N \tau^2\left(\frac{\delta^2_k(t) }{\hat{\lambda}_k-\eta_k} +\hat{a}^2_k(t)\frac{(\sqrt{\hat{\lambda}_k}-\sqrt{\hat{\lambda}_k-\eta_k})^2}{\hat{\lambda}_k(\hat{\lambda}_k-\eta_k)}\right)+\right.$$

$$\left.+\sum_{k=N+1}^\infty\frac{\tau^2a^2_k(t)}{\lambda_k} \right)^{p/2}dt= \tau^{p/2}\int_0^T\left(\sum_{k=1}^N \tau^2\left(\frac{\delta^2_k(t) }{\hat{\lambda}_k-\eta_k} +\right.\right.$$

$$\left.\left.+\hat{a}^2_k(t)\frac{\left(\sqrt{\hat{\lambda}_k}-\sqrt{\hat{\lambda}_k-\eta_k}\right)^2}{\hat{\lambda}_k(\hat{\lambda}_k-\eta_k)}\right) + B(t,t) - \sum_{k=1}^N\frac{a^2_k(t)}{\lambda_k} \right)^{p/2}dt\leq$$

$$\leq\tau^{p/2}\int_0^T\left(\left(B(t,t)-\sum_{k=1}^N \frac{(\hat{a}_k(t)-\delta_k(t))^2}{\hat{\lambda}_k+\eta_k}\right)+ \sum_{k=1}^N \left(\frac{\delta^2_k(t) }{\hat{\lambda}_k-\eta_k} + \right.\right.$$

$$\left.\left.+\hat{a}^2_k(t)\frac{(\sqrt{\hat{\lambda}_k}-\sqrt{\hat{\lambda}_k-\eta_k})^2}{\hat{\lambda}_k(\hat{\lambda}_k-\eta_k)}\right)\right)^{p/2} dt.$$
\end{dov}

\begin{theorem}
\label{omlp-gamma-2}
Let the random process $X=\{X(t),t\in [0,T]\} $ belong to $Sub_\varphi (\Omega)$,
$$\varphi(t)= \frac{t^\gamma}{\gamma}$$
\noindent for $1<\gamma<2$. Let
$$c_N= \int_0^T\left(\sum_{k=1}^N \tau^{\gamma}_\varphi(\xi_k)\left(\frac{\delta^{\gamma}_k(t) }{(\lambda_k-\eta_k)^\gamma} +\hat{a}^{\gamma}_k(t)\frac{(\sqrt{\hat{\lambda}_k}-\sqrt{\hat{\lambda}_k-\eta_k})^{\gamma}}{(\hat{\lambda}_k(\hat{\lambda}_k-\eta_k))^\gamma}\right)+\right.$$
$$+\left.\sum_{k=N+1}^\infty\frac{\tau^{\gamma}_\varphi(\xi_k)a^{2\gamma}_k(t)}{\lambda_k^\gamma} \right)^{p/\gamma}dt<\infty,$$
\noindent where $\delta_k(t)=|\varphi_k(t)-\hat{\varphi}_k(t)|$ is the error in finding the $k$-th eigenfunction of the equation (\ref{omkl-eq}), $\hat{\lambda}_k$ is the approximate value of the $k$-th eigenvalue, $\eta_k=|\lambda_k-\hat{\lambda}_k|$ is the error in approximating the $k$-th eigenvalue.
The  $X_N$ model \eqref{omkl-decompos-x-n} approximates the random process $X$ with reliability $1-\alpha$ and accuracy $\delta$ in the space $L_p[0,T]$
if the following conditions are satisfied
$$\left\{\begin{array}{c}c_N\leq \delta/(\beta ln \frac{2}{\alpha})^{p/\beta}\\ c_N<\delta/p^{p(1-1/\gamma)}\end{array}\right.$$
\end{theorem}

\begin{dov}The statement of the Theorem follows from the Theorem \ref{omlp_g1}, Theorem \ref{omkl-main} and the proof of Theorem \ref{omkl-lp-square}.
\end{dov}

\subsubsection{Modeling of random processes by means of their decomposition into certain bases in $L_p[0,T]$}

In this section we will use decompositions of the second order centered stochastic processes in series (see \cite{koz-roz-turch} for more details).

\begin{theorem} \label{ombasis-dec-main}
Let $X(t)$, $t\in T$ be a second-order stochastic process, $EX(t)=0$, $ t\in T$, with the correlation function
$B(t,s)=EX(t)\overline{X(s)}$.
Let $f(t,\lambda)$, $ t\in T$, be functions from the space $L_2(\Lambda,\mu)$.
Let $\{g_k(\lambda), k\in Z\}$ be an orthonormal basis in $L_2(\Lambda,\mu)$.
The correlation function $B(t,s)$ admits the representation
$$B(t,s) = \int_\Lambda f(t,\lambda)\overline{ f(s,\lambda)} d\mu(\lambda)$$
\noindent if and only if

\begin{equation}
\label{ombasis}
X(t)=\sum_{k=1}^\infty a_k(t)\xi_k,
\end{equation}
\noindent where

\begin{equation}
\label{ombasis-f}
a_k(t)=\int_\Lambda f(t,\lambda)\overline{g_k(\lambda)}d\mu(\lambda),
\end{equation}
\noindent $\xi_k$ are centered uncorrelated random variables such that $E\xi_k=0$, $E\xi_k\xi_l=\delta_{kl}$, $E\xi_k^2=1$.
\end{theorem}

\begin{definition}\label{ommodel-basis-ozn}
Let the random process $X=\{X(t),t\in T\}$ be represented in the form (\ref{ombasis}). A model of such a process will be called a process $X_N=\{X_N(t),t\in T\}$ of the form
\begin{equation}
\label{ommodel-basis}
X_N(t) = \sum_{k=1}^N \xi_k\hat{a}_k(t),
\end{equation}
\noindent where $\hat{a}_k(t)$ is an approximation of the functions $a_k(t)$ given in the form (\ref{ombasis-f}), $\xi_k$ are centered uncorrelated random variables: $E\xi_k=0$, $E\xi_k\xi_l=\delta_{kl}$, $E\xi_k^2=1$.
\end{definition}

Let us consider an example of a stochastic process on a finite interval with the representation
$$X(t)=\sum_{k=1}^{\infty} \xi_k a_k(t),$$
$$a_k(t) = \int_0^T f(t,\lambda)\cos \pi k\lambda d\lambda.$$

\begin{theorem}
\label{omcosdec}
Let the random process
$X=\{X(t),t\in [0,T]\}$ belong to the space $S u b_\varphi (\Omega)$,

$$\varphi(t)=\left\{\begin{array}{c}\frac{t^2}{\gamma}, t<1\\ \frac{t^\gamma}{\gamma},t\geq1\end{array}\right.$$

\noindent for $\gamma>2$. Let $f(t,s)$ be differentiable with respect to $s$, and the process $X(t)$ can be represented in the form (\ref{ombasis}) for $g_k(t) = cos(\pi k t)$.
Let in addition

$$c_N = \int_0^T \left(\delta^2_f(t)\sum_{k=N+1}^\infty \frac{4\tau^2_\varphi(\xi_k)}{\pi^2 k^2} + \sum_{k=1}^N \tau^2_\varphi(\xi_k)\delta^2_k(t)\right)^{p/2}dt<\infty,$$

\noindent $\delta_f(t) = f(t,T)-f(t,0).$ The model $X_N(t)$ \eqref{ommodel-basis} approximates the process $X(t)$ with given reliability $1-\alpha$ and accuracy $\delta$ in the spaces $L_p(0,T)$ when

$$\left\{\begin{array}{c}c_N\leq \delta/(\beta ln \frac{2}{\alpha})^{p/\beta} \\ c_N<\delta/p^{p(1-1/\gamma)}\end{array}\right.$$

\noindent for such $\beta$ that $1/\gamma+1/\beta=1$.

\end{theorem}

\begin{dov} From Theorem \ref{omlp_g2} we have that

$$c_N = \int_0^T\left(\tau_\varphi(\Delta_N(t))\right)^p d t = \int_0^T\left(\tau_\varphi\left(\sum_{k=1}^N\xi_k\delta_k(t) + \sum_{k=N+1}^\infty\xi_k a_k(t)\right)\right)^p dt \leq$$

$$\leq \int_0^T\left(\sum_{k=1}^N\tau^2_\varphi(\xi_k)\delta^2_k(t) + \sum_{k=N+1}^\infty\tau^2_\varphi(\xi_k) a^2_k(t)\right)^{p/2} dt .$$

\noindent Now the statement of Theorem \ref{omcosdec} follows from the fact that:

$$a_k(t) = \int_0^T f(t,\lambda) \cos(\pi k \lambda) d\lambda  = \left.\frac{1}{\pi k}f(t,\lambda)\sin(\pi k \lambda)\right|_0^T - $$

$$-\int_0^T \frac{\partial f(t,\lambda)}{\partial \lambda}\frac{\sin(\pi k \lambda)}{\pi k} d\lambda \leq (f(t,T)-f(t,0)) \frac{2}{\pi k}=\delta_f(t)\frac{2}{\pi k}.$$

\end{dov}

\begin{theorem}
\label{omcosdec-2}
Let the random process
$X=\{X(t),t\in [0,T]\}$ belong to the space $S u b_\varphi (\Omega)$,

$$\varphi(t)=\frac{t^\gamma}{\gamma},$$

\noindent for $1<\gamma<2$. Let $f(t,s)$ be differentiable with respect to $s$, $g_k(t) = cos(\pi k t)$, and the process $X(t)$ can be represented in the form (\ref{ombasis}).
Let in addition

$$c_N = \int_0^T \left(\delta^\gamma_f(t)\sum_{k=N+1}^\infty \frac{4\tau^\gamma_\varphi(\xi_k)}{\pi^\gamma k^\gamma} + \sum_{k=1}^N \tau^\gamma_\varphi(\xi_k)\delta^\gamma_k(t)\right)^{p/\gamma}dt<\infty,$$

\noindent $\delta_f(t)$ is given by Theorem \ref{omcosdec}. The model $X_N(t)$ \eqref{ommodel-basis} approximates the process $X(t)$ with given reliability $1-\alpha$ and accuracy $\delta$ in the spaces $L_p[0,T]$ when

$$\left\{\begin{array}{c}c_N\leq \delta/(\beta ln \frac{2}{\alpha})^{p/\beta}\\ c_N<\delta/p^{p(1-1/\gamma)}\end{array}\right.$$

\noindent for such $\beta$ that $1/\gamma+1/\beta=1$.
\end{theorem}

\begin{dov} The statement of Theorem  \ref{omcosdec-2} follows from Theorem \ref{omcosdec} and Theorem \ref{omlp_g1}.
\end{dov}

\subsubsection{Modeling of random processes using their Hermite basis decomposition in $L_p[0,T]$}

Let $X=\{X(t),t \in [0,T]\} \in Sub_\varphi(\Omega)$ be a second-order random process, $EX(t)=0$.
Let the covariance function of the process $B(t,s)=EX(t)\overline{X(s)}$ admit the representation

$$B(t,s)=\int_{-\infty}^\infty f(t,\lambda)f(s,\lambda)d\lambda,$$

\noindent where $f(t,\lambda)$, $t\in [0,T]$, $\lambda\in R$ is a family of functions from $L_2(R)$.
Since the Hermite functions \cite{hermite} form an orthonormal basis, then, according to the theorem \ref{ombasis-dec-main},
the process $X$ can be represented as

$$X(t) = \sum_{k=0}^{\infty}\xi_k \int_{-\infty}^{\infty} f(t,\lambda)g_k(\lambda)d\lambda,$$

\noindent where $\xi_k$ are $\varphi$-subGaussian independent uncorrelated random variables such that $E\xi_k^2=1$; $g_k(\lambda)$ are Hermite functions:

\begin{equation}
\label{omhermite_func}
g_k(\lambda)=\frac{H_k(\lambda)}{\sqrt{k!}}\frac{1}{\sqrt[4]{2\pi}}\exp\{-\frac{\lambda^2}{2}\},
\end{equation}

\noindent $H_k(\lambda)$ - Hermite polynomials, i.e. functions of the form

$$H_k(\lambda)=(-1)^ke^{\lambda^2/2}\frac{d^n}{d\lambda^n}e^{-\lambda^2/2}.$$

\begin{theorem}
\label{omhermdec}
Let the random process $X=\{X(t),t\in [0,T]\}$ belong to the space $ S u b_\varphi (\Omega)$,

$$\varphi(t)=\left\{\begin{array}{c}\frac{t^2}{\gamma}, t<1\\ \frac{t^\gamma}{\gamma},t\geq1\end{array}\right.$$

\noindent for $\gamma>2$, and let the process $X(t)$ be represented in the form (\ref{ombasis}), where $g_k(t)$ are Hermite functions.
Let

$$c_N = \int_0^T \left(K^2 \int_{-\infty}^\infty Z^2_f(t,\lambda)d\lambda \sum_{k=N+1}^\infty \frac{\tau^2_\varphi(\xi_k)}{k^2+3k+2} + \sum_{k=1}^N \tau^2_\varphi(\xi_k)\delta^2_k(t)\right)^{p/2}dt<\infty,$$
$$Z_f(t,\lambda) = \frac{\partial^2 f(t,\lambda)}{\partial \lambda^2} - \lambda \frac{\partial f(t,\lambda)}{\partial \lambda}+\frac{\lambda^2-2}{4}f(t,\lambda),$$

\noindent where $K\approx 1.086435,$ $f(t,s)$ is twice differentiable in $s$ and grows in $s$ no faster than $\exp^{s^2/4}$, $Z_f(\lambda)$ is integrable on $R$. 
The model $X_N(t)$, given in (\ref{ommodel-basis}), approximates the process $X(t)$ with given reliability $1-\alpha$ and accuracy $\delta$ in the spaces $L_p(0,T)$, when

$$\left\{\begin{array}{c}c_N\leq \delta/(\beta ln \frac{2}{\alpha})^{p/\beta} \\ c_N<\delta/p^{p(1-1/\gamma)}\end{array}\right.,$$

\noindent where $1/\gamma+1/\beta=1$.

\end{theorem}

\begin{dov}

Under the conditions of the theorem

$$a_k(t) = \int_{-\infty}^\infty f(t,\lambda)\hat{H}_k(\lambda)d\lambda = \frac{1}{\sqrt[4]{\pi}}\int_{-\infty}^\infty f(t,\lambda)\frac{H_k(\lambda)e^{-\lambda^2/4}}{\sqrt{k!}}d\lambda.$$
Using properties of Hermite polynomials\cite{hermite},  we can show that

$$\frac{\partial H_k(t)}{\partial t} = k H_{k-1}(t).$$
Using integration by parts, we get:

$$a_k(t) = \frac{1}{\sqrt[4]{\pi}}\int_{-\infty}^\infty f(t,\lambda)\frac{H_k(\lambda)e^{-\lambda^2/4}}{\sqrt{k!}}d\lambda = $$

$$=\frac{1}{\sqrt[4]{\pi}}\int_{-\infty}^\infty f(t,\lambda)\frac{e^{-\lambda^2/4}}{\sqrt{k+1}\sqrt{(k+1)!}}\frac{\partial H_{k+1}(\lambda)}{\partial \lambda}d\lambda = $$

$$ = \left.\frac{1}{\sqrt[4]{\pi}} f(t,\lambda)\frac{e^{-\lambda^2/4}}{\sqrt{k+1}\sqrt{(k+1)!}}H_{k+1}(\lambda)\right|_{\lambda=-\infty}^{\lambda=\infty} -$$

$$- \frac{1}{\sqrt[4]{\pi}}\int_{-\infty}^\infty \frac{\partial (f(t,\lambda) e^{-\lambda^2/4})}{\partial \lambda}\frac{H_{k+1}(\lambda)}{\sqrt{k+1}\sqrt{(k+1)!}}d\lambda.$$

Since $H_k(\lambda)\exp\{-\lambda^2/4\}$ tends to zero as $\lambda\to\pm\infty$, and $f(t,\lambda)$ is bounded, we get

$$a_k(t)=- \frac{1}{\sqrt[4]{\pi}}\int_{-\infty}^\infty \frac{\partial (f(t,\lambda) e^{-\lambda^2/4})}{\partial \lambda}\frac{H_{k+1}(\lambda)}{\sqrt{k+1}\sqrt{(k+1)!}}d\lambda.$$
Integration by parts one more time gives

$$a_k(t)=- \frac{1}{\sqrt[4]{\pi}}\int_{-\infty}^\infty \frac{\partial( f(t,\lambda) e^{-\lambda^2/4})}{\partial \lambda}\frac{H_{k+1}(\lambda)}{\sqrt{k+1}\sqrt{(k+1)!}}d\lambda=$$

$$=-\frac{1}{\sqrt[4]{\pi}}\int_{-\infty}^\infty \frac{\partial (f(t,\lambda) e^{-\lambda^2/4})}{\partial \lambda}\frac{1}{\sqrt{(k+1)(k+2)}\sqrt{(k+2)!}}\frac{\partial H_{k+2}(\lambda)}{\partial \lambda}d\lambda = $$

$$ = -\left.\frac{1}{\sqrt[4]{\pi}} \frac{\partial (f(t,\lambda) e^{-\lambda^2/4})}{\partial \lambda}\times\frac{1}{\sqrt{(k+1)(k+2)}\sqrt{(k+2)!}}H_{k+1}(\lambda)\right|_{\lambda=-\infty}^{\lambda=\infty} + $$

$$+ \frac{1}{\sqrt[4]{\pi}}\int_{-\infty}^\infty \frac{\partial^2 (f(t,\lambda) e^{-\lambda^2/4})}{\partial \lambda^2}\frac{H_{k+2}(\lambda)}{\sqrt{(k+1)(k+2)}\sqrt{(k+2)!}}d\lambda. $$

It is certain that $H_{k+1}(\lambda)\partial (f(t,\lambda)e^{-\lambda^2/4})/\partial \lambda$ tends to zero as $\lambda\to\pm\infty$, since $\partial (f(t,\lambda)e^{-\lambda^2/4})/\partial \lambda $ $=e^{-\lambda^2/4}(\partial f(t,\lambda)/\partial\lambda-\lambda f(t,\lambda))$, $H_{k+1}(\lambda)e^{-\lambda^2/4}\to 0$, and   $\partial f(t,\lambda)/\partial\lambda-\lambda f(t,\lambda)$ is bounded, because $f(t,\lambda)$ and $\lambda f(t,\lambda)$ are bounded due to the conditions of the theorem. Then,

$$a_k(t) = \frac{1}{\sqrt[4]{\pi}}\int_{-\infty}^\infty \frac{\partial^2 (f(t,\lambda) e^{-\lambda^2/4})}{\partial \lambda^2}\frac{H_{k+2}(\lambda)}{\sqrt{(k+1)(k+2)}\sqrt{(k+2)!}}d\lambda.$$
Because $\hat{H}_k$ is an orthonormal basis, $\int_{-\infty}^\infty \hat{H}^2_k(t)dt = 1$. That's why

$$\frac{1}{\sqrt[4]{\pi}}\int_{-\infty}^\infty \frac{\partial^2 (f(t,\lambda) e^{-\lambda^2/4})}{\partial \lambda^2}\frac{H_{k+2}(\lambda)}{\sqrt{(k+1)(k+2)}\sqrt{(k+2)!}}d\lambda = $$ $$=\frac{1}{\sqrt{(k+1)(k+2)}}\int_{-\infty}^\infty \frac{\partial^2 (f(t,\lambda) e^{-\lambda^2/4})}{\partial \lambda^2}\frac{H_{k+2}(\lambda)e^{-\lambda^2/4}}{\sqrt{(k+2)!}\sqrt[4]{2\pi}}e^{\lambda^2/4}d\lambda= $$$$= \frac{1}{\sqrt{(k+1)(k+2)}}\int_{-\infty}^\infty \frac{\partial^2 (f(t,\lambda) e^{-\lambda^2/4})}{\partial \lambda^2} g_{k+2}(\lambda)e^{\lambda^2/4}d\lambda\leq$$

$$\leq \frac{1}{\sqrt{(k+1)(k+2)}}\left(\int_{-\infty}^\infty \left(\frac{\partial^2 (f(t,\lambda) e^{-\lambda^2/4})}{\partial \lambda^2}\times e^{\lambda^2/4}\right)^2d\lambda\right)^\frac{1}{2}\left(\int_{-\infty}^\infty g^2_{k+2}(\lambda) d\lambda\right)^\frac{1}{2}=$$

$$= \frac{1}{\sqrt{(k+1)(k+2)}}\left(\int_{-\infty}^\infty \left(\frac{\partial^2 (f(t,\lambda) e^{-\lambda^2/4})}{\partial \lambda^2} e^{\lambda^2/4}\right)^2d\lambda\right)^\frac{1}{2}.$$
Besides,

$$\frac{\partial^2 (f(t,\lambda) e^{-\lambda^2/4})}{\partial \lambda^2}e^{\lambda^2/4} = $$$$=\frac{\partial}{\partial \lambda}\left(\frac{\partial f(t,\lambda)}{\partial \lambda}e^{-\lambda^2/4}-\frac{1}{2}e^{-\lambda^2/4}\lambda f(t,\lambda)\right)e^{\lambda^2/4}=$$$$=\left(\frac{\partial}{\partial \lambda}\left(\frac{\partial f(t,\lambda)}{\partial \lambda}e^{-\lambda^2/4}\right)-\frac{\partial}{\partial \lambda}\left(\frac{1}{2}e^{-\lambda^2/4}\lambda f(t,\lambda)\right)\right)e^{\lambda^2/4}=$$$$=\left(\frac{\partial^2 f(t,\lambda)}{\partial \lambda^2}e^{-\lambda^2/4}-\frac{1}{2}e^{-\lambda^2/4}\lambda \frac{\partial f(t,\lambda)}{\partial \lambda} - \frac{1}{2}e^{-\lambda^2/4}\lambda \frac{\partial f(t,\lambda)}{\partial \lambda} -\right.$$

$$\left.-f(t,\lambda)\frac{\partial}{\partial \lambda}\left(\frac{1}{2}e^{-\lambda^2/4}\lambda\right)\right)e^{\lambda^2/4} = $$$$=\left(\frac{\partial^2 f(t,\lambda)}{\partial \lambda^2}e^{-\lambda^2/4}-e^{-\lambda^2/4}\lambda \frac{\partial f(t,\lambda)}{\partial \lambda}\right.-$$

$$\left. - f(t,\lambda)\left(\frac{1}{2}e^{-\lambda^2/4}-\frac{1}{4}e^{-\lambda^2/4}\lambda^2\right) \right)e^{\lambda^2/4} = $$

$$=\frac{\partial^2 f(t,\lambda)}{\partial \lambda^2} - \lambda \frac{\partial f(t,\lambda)}{\partial \lambda}+\frac{\lambda^2-2}{4}f(t,\lambda).$$

Finally, we obtain

$$a_k(t)\leq \left(\frac{1}{(k+1)(k+2)}\int_{-\infty}^\infty Z_f^2(t,\lambda) d\lambda\right)^\frac{1}{2},$$

\noindent where

$$Z_f(t,\lambda) = \left|\frac{\partial^2 f(t,\lambda)}{\partial \lambda^2} - \lambda \frac{\partial f(t,\lambda)}{\partial \lambda}+\frac{\lambda^2-2}{4}f(t,\lambda)\right|.$$

\noindent Therefore,

$$c_N = \int_0^T\left(\sum_{k=1}^N \tau^2_\varphi(\xi_k)\delta^2_k(t) +\right.$$$$\left.+ \sum_{k=N+1}^\infty\tau^2_\varphi(\xi_k)a^2_k(t)  \right)^{p/2} dt \leq $$

$$\leq\int_0^T \left( \int_{-\infty}^\infty Z^2_f(t,\lambda)d\lambda \sum_{k=N+1}^\infty \frac{\tau^2_\varphi(\xi_k)}{k^2+3k+2} + \sum_{k=1}^N \tau^2_\varphi(\xi_k)\delta^2_k(t)\right)^{p/2}dt.$$

\noindent Finally, the statement of Theorem  \ref{omhermdec} follows from Theorem \ref{omlp_g2}.

\end{dov}

\begin{theorem}
\label{omhermdec-2}

Let the random process $X=\{X(t),t\in [0,T]\}$ belong to the space $ S u b_\varphi (\Omega)$,

$$\varphi(t)=\left\{\begin{array}{c}\frac{t^2}{\gamma}, t<1\\ \frac{t^\gamma}{\gamma},t\geq1\end{array}\right.$$

\noindent for $\gamma>2$, and let the process $X(t)$ be represented in the form (\ref{ombasis})
and the series $\hat{H}_k(t)$ of Hermite functions is the basis. Let $\tau_\varphi(\xi_k)=\tau \omega^k$, $|\omega|<1$, and

$$c_N = \int_0^T \left( \frac{\tau}{\sqrt{(1-\omega^2)}}\left(\int_{-\infty}^\infty f^2(t,\lambda)d\lambda\right)^{1/2}-\sum_{k=0}^N \tau \omega^k\hat{a}_k(t)\right)^{p}dt<\infty.$$

Model $X_N(t)$, provided in (\ref{ommodel-basis}), approximates $X(t)$ with given reliability $1-\alpha$ and accuracy $\delta$ in the space $L_p(0,T)$, if

$$\left\{\begin{array}{c}c_N\leq \delta/(\beta ln \frac{2}{\alpha})^{p/\beta} \\ c_N<\delta/p^{p(1-1/\gamma)}\end{array}\right.,$$

\noindent where $1/\gamma+1/\beta=1$.

\end{theorem}

\begin{proof}

According to the statement of the theorem,

$$a_k(t) = \int_{-\infty}^\infty f(t,\lambda)\hat{H}_k(\lambda)d\lambda = $$ $$=\frac{1}{\sqrt[4]{\pi}}\int_{-\infty}^\infty f(t,\lambda)\frac{H_k(\lambda)e^{-\lambda^2/2}}{\sqrt{k!}}d\lambda.$$

Using the properties of Hermite polynomials $\{H_k(\lambda)\}$, we get

$$GF_{H^2}(\lambda,\omega) = \sum_{k=1}^\infty \frac{H^2_k(\lambda)}{k!2^k}\omega^k = \frac{1}{\sqrt{1-\omega^2}}\exp\left\{\frac{2\lambda^2\omega}{1+\omega}\right\}.$$

\noindent Moreover, such series converges for $|\omega|<1$. Under the conditions of the theorem, $\tau_\varphi(\xi_k)=\tau \omega^k$. That's why, for the process $X(t)$ the next condition holds true:

$$\tau_\varphi(X(t)) = \tau_\varphi\left(\sum_{k=1}^\infty \xi_k\int_{-\infty}^\infty f(t,\lambda)\hat{H}_k(\lambda)d\lambda\right) \leq  \tau\int_{-\infty}^\infty f(t,\lambda) \sum_{k=0}^\infty \omega^k\hat{H}_k(\lambda)d\lambda\leq$$

 $$\leq \tau \left(\int_{-\infty}^\infty f^2(t,\lambda)d\lambda\right)^{1/2}  \left(\int_{-\infty}^\infty \left(\sum_{k=0}^\infty\frac{H_k(\lambda)\exp\{-\lambda^2/2\}}{\sqrt{k!2^k\sqrt{\pi}}}\omega^k\right)^2d\lambda\right)^{1/2} \leq $$

 $$\leq \tau \left(\int_{-\infty}^\infty f^2(t,\lambda)d\lambda\right)^{1/2} \left(\int_{-\infty}^\infty \sum_{k=0}^\infty\frac{H^2_k(\lambda)\exp\{-\lambda^2\}}{k!2^k\sqrt{\pi}}\omega^{2k}d\lambda\right)^{1/2} =$$

 $$= \tau \left(\int_{-\infty}^\infty f^2(t,\lambda)d\lambda\right)^{1/2} \left(\frac{1}{\sqrt{\pi}}\int_{-\infty}^\infty GF_{H^2}(\lambda,\omega^2)\exp\{-\lambda^2\}d\lambda\right)^{1/2} =$$

 $$= \tau \left(\int_{-\infty}^\infty f^2(t,\lambda)d\lambda\right)^{1/2} \left(\frac{1}{\sqrt{\pi(1-\omega^4)}}\int_{-\infty}^\infty\exp\left\{\frac{2\lambda^2\omega}{1+\omega}\right\}\exp\{-\lambda^2\}d\lambda\right)^{1/2} = $$

 $$ = \tau \left(\int_{-\infty}^\infty f^2(t,\lambda)d\lambda\right)^{1/2} \left(\frac{1}{\sqrt{\pi(1-\omega^4)}}\int_{-\infty}^\infty\exp\left\{\lambda^2\frac{\omega^2-1}{\omega^2+1}\right\}d\lambda\right)^{1/2} = $$

 $$ = \tau \left(\int_{-\infty}^\infty f^2(t,\lambda)d\lambda\right)^{1/2} \left(\frac{1}{\sqrt{\pi(1-\omega^4)}} \sqrt{\pi\frac{1+\omega^2}{1-\omega^2}}\right)^{1/2} =$$

 $$ = \frac{\tau}{\sqrt{(1-\omega^2)}}\left(\int_{-\infty}^\infty f^2(t,\lambda)d\lambda\right)^{1/2}.$$

 Given these considerations, the estimator of the model of the process will take the next form:

 $$\tau_\varphi(X(t)-X_N(t)) = \tau_\varphi\left(\sum_{k=0}^\infty \xi_ka_k(t)-\sum_{k=0}^N \xi_k\hat{a}_k(t)\right)= $$$$=\tau_\varphi\left(\sum_{k=0}^N \xi_k\delta_k(t) + \sum_{k=N+1}^\infty \xi_ka_k(t)\right)\leq \sum_{k=0}^N \tau_\varphi(\xi_k)\delta_k(t) + \sum_{k=N+1}^\infty \tau_\varphi(\xi_k) a_k(t)  =$$$$= \sum_{k=0}^\infty \tau \omega^k a_k(t)-\sum_{k=0}^N \tau \omega^k\hat{a}_k(t)\leq \frac{\tau}{\sqrt{(1-\omega^2)}}\left(\int_{-\infty}^\infty f^2(t,\lambda)d\lambda\right)^{1/2}-\sum_{k=0}^N \tau \omega^k\hat{a}_k(t).$$

The statement of Theorem \ref{omhermdec-2} follows from  Theorem \ref{omlp_g1} and Theorem \ref{omhermdec}.

\end{proof}

\subsubsection{Modeling of Stochastic Processes in $L_p(0,T)$ Using the Chebyshev Polynomials}

Let the process $X(t)$ has the same properties as the process from the previous section. Let orthonotmal Chebyshev polynomials be used as the basis:

$$\hat{T}_n(\lambda) = \sqrt{\frac{2}{\pi}}T_n(\lambda),$$ where

$$T_n(\lambda) = \cos(n \arccos \lambda).$$

In such a case we can proof the next theorem.

\begin{theorem}

 Let a stochastic process $X=\{X(t),t\in [0,T]\}$ belong to the space $ S u b_\varphi (\Omega)$ with

$$\varphi(t)=\left\{\begin{array}{c}\frac{t^2}{\gamma}, t<1\\ \frac{t^\gamma}{\gamma},t\geq1\end{array}\right.$$

\noindent for  $\gamma>2$, let $X(t)$ is representated in the form (\ref{ombasis}), and let the system of orthonormal Chebyshev polynomials $\{\hat{T}_k(t)\}$ be used as the basis. Let $\tau_\varphi(\xi_k)=\tau \omega^k$, $0<\omega<1$, and

$$c_N = \int_0^T \left(  \sqrt{\frac{2}{\pi}}\tau\left(\int_{-1}^1 f^2(t,\lambda)d\lambda\right)^{1/2} \sqrt{D_T(\omega)}-\sum_{k=0}^N \tau \omega^k\hat{a}_k(t)\right)^{p}dt<\infty,$$

$$D_T(\omega)=2\frac{1}{{\omega(4+3\omega^2+\omega^4)}}\left(\omega(5+5\omega^2+2\omega^2)+\right.$$

$$+\left.(4+7\omega^2+4\omega^4+\omega^6)\ln\{(\omega^2-\omega+2)/(\omega^2+\omega+2)\}\right).$$

Model $X_N(t)$, determined in (\ref{ommodel-basis}), approximates the process $X(t)$ with given reliability $1-\alpha$ and accuracy $\delta$ in the space $L_p(0,T)$, if

$$\left\{\begin{array}{c}c_N\leq \delta/(\beta ln \frac{2}{\alpha})^{p/\beta} \\ c_N<\delta/p^{p(1-1/\gamma)}\end{array}\right.,$$

\noindent where $1/\gamma+1/\beta=1$.
\end{theorem}

\begin{proof} Under the conditions of the theorem,

$$a_k(t) = \int_{-1}^1 f(t,\lambda)\hat{T}_n(\lambda)d\lambda = \int_{-1}^1 f(t,\lambda)\sqrt{\frac{2}{\pi}}T_n(\lambda)d\lambda.$$

The generating function of orthogonal Chebyshev polynomials of the first kind $\{T_k(\lambda)\}$ has the next form:

$$GF_T(\lambda,\omega) = \sum_{k=0}^\infty T_k(\lambda)\omega^k = \frac{1-\omega \lambda}{2-\omega\lambda+\omega^2}$$ for $0<\omega<1$. Under the conditions of the theorem, $\tau_\varphi(\xi_k)=\tau \omega^k$. Thats why for the process $X(t)$ the next  condition is true:

$$\tau_\varphi(X(t)) = \tau_\varphi\left(\sum_{k=1}^\infty \xi_k\int_{-1}^1 f(t,\lambda)\hat{T}_k(\lambda)d\lambda\right) \leq  \tau\int_{-1}^1 f(t,\lambda) \sum_{k=0}^\infty \omega^k\hat{T}_k(\lambda)d\lambda\leq$$

 $$\leq \tau \left(\int_{-1}^1 f^2(t,\lambda)d\lambda\right)^{1/2}   \left(\int_{-1}^1 \left(\sum_{k=0}^\infty\sqrt{\frac{2}{\pi}}T_k(\lambda)\right)^2d\lambda\right)^{1/2} = $$

$$ =  \sqrt{\frac{2}{\pi}}\tau\left(\int_{-1}^1 f^2(t,\lambda)d\lambda\right)^{1/2}   \left(\int_{-1}^1\left(\frac{1-\omega \lambda}{2-\omega\lambda+\omega^2} \right)^2d\lambda\right)^{1/2}.$$

Lets calculate the second integral of the last expression separetely.

$$\int_{-1}^1\left(\frac{1-\omega \lambda}{2-\omega\lambda+\omega^2} \right)^2d\lambda = $$ $$=2\frac{1}{{\omega(4+3\omega^2+\omega^4)}}\left(\omega(5+5\omega^2+2\omega^2)+\right.$$

$$+\left.(4+7\omega^2+4\omega^4+\omega^6)\ln\left\{\frac{\omega^2-\omega+2}{\omega^2+\omega+2)}\right\}\right):=D_T(\omega).$$

\noindent Then the estimator of $\tau_\varphi(X(t))$ will take the form:

$$\tau_\varphi(X(t))\leq\sqrt{\frac{2}{\pi}}\tau\left(\int_{-1}^1 f^2(t,\lambda)d\lambda\right)^{1/2} \sqrt{D_T(\omega)}.$$

Given these considerations, the estimator of the model of the process will take the following form:

 $$\tau_\varphi(X(t)-X_N(t)) = \tau_\varphi\left(\sum_{k=0}^\infty \xi_ka_k(t)-\sum_{k=0}^N \xi_k\hat{a}_k(t)\right)=$$$$= \tau_\varphi\left(\sum_{k=0}^N \xi_k\delta_k(t) + \sum_{k=N+1}^\infty \xi_ka_k(t)\right)\leq \sum_{k=0}^N \tau_\varphi(\xi_k)\delta_k(t) + \sum_{k=N+1}^\infty \tau_\varphi(\xi_k) a_k(t)  = $$$$=\sum_{k=0}^\infty \tau \omega^k a_k(t)-\sum_{k=0}^N \tau \omega^k\hat{a}_k(t)\leq$$

 $$\leq\sqrt{\frac{2}{\pi}}\tau\left(\int_{-1}^1 f^2(t,\lambda)d\lambda\right)^{1/2} \sqrt{D_T(\omega)}-\sum_{k=0}^N \tau \omega^k\hat{a}_k(t).$$

 Statement of this theorem follows from  Theorem \ref{omlp_g1} and the last inequality.

\end{proof}

We can also proof the similar theorem in the case of Chebyshev polynomials of the second kind:

$$U_n(\lambda)=\frac{\sin((n+1)\arccos \lambda)}{\sqrt{1-\lambda^2}}$$

$$\hat{U}_n(\lambda) = \sqrt{\frac{2}{\pi}}U_n(\lambda).$$ Generating function of the series $\{U_n(\lambda)\}$ has the form

$$GF_U(\lambda,\omega) = \sum_{k=0}^\infty \omega^k U_n(\lambda).$$

Let the process $X(t)$ has the same properties as in the previous theorem.
Then the next proposition holds true.

\begin{theorem}  Let stochastic process $X=\{X(t),t\in [0,T]\}$ belong to the space $ S u b_\varphi (\Omega)$,

$$\varphi(t)=\left\{\begin{array}{c}\frac{t^2}{\gamma}, t<1\\ \frac{t^\gamma}{\gamma},t\geq1\end{array}\right.$$

\noindent for $\gamma>2$, let process $X(t)$ can be represented in the form (\ref{ombasis}), and let series $\{\hat{U}_k(t)\}$ of orthonormal Chebyshev polynomials be the basis. Let also $\tau_\varphi(\xi_k)=\tau \omega^k$, $0<\omega<1$, and

$$c_N = \int_0^T \left(  \frac{2\tau}{\sqrt{\pi}(\omega^2-1)}\left(\int_{-1}^1 f^2(t,\lambda)d\lambda\right)^{1/2}-\sum_{k=0}^N \tau \omega^k\hat{a}_k(t)\right)^{p}dt<\infty.$$

Model $X_N(t)$, determined in (\ref{ommodel-basis}), approximates the process $X(t)$ with given reliability $1-\alpha$ and accuracy $\delta$ in the space $L_p(0,T)$, if

$$\left\{\begin{array}{c}c_N\leq \delta/(\beta ln \frac{2}{\alpha})^{p/\beta} \\ c_N<\delta/p^{p(1-1/\gamma)}\end{array}\right.,$$

\noindent where $1/\gamma+1/\beta=1$.
\end{theorem}

\begin{proof} Under the conditions of the theorem, we have:

$$a_k(t) = \int_{-1}^1 f(t,\lambda)\hat{U}_n(\lambda)d\lambda = \int_{-1}^1 f(t,\lambda)\sqrt{\frac{2}{\pi}}U_n(\lambda)d\lambda.$$

Generation function of the series of orthogonal Chebyshev polynomials of the second kind $\{U_k(\lambda)\}$ has the following form:

$$GF_T(\lambda,\omega) = \sum_{k=0}^\infty U_k(\lambda)\omega^k = \frac{1}{1-2\omega\lambda+\omega^2}$$ for $0<\omega<1$. As stated in the theorem's conditions, $\tau_\varphi(\xi_k)=\tau \omega^k$. That's why for stochastic process $X(t)$ the next statement holds true:

$$\tau_\varphi(X(t)) = \tau_\varphi\left(\sum_{k=1}^\infty \xi_k\int_{-1}^1 f(t,\lambda)\hat{U}_k(\lambda)d\lambda\right) \leq$$$$\leq  \tau\int_{-1}^1 f(t,\lambda) \sum_{k=0}^\infty \omega^k\hat{U}_k(\lambda)d\lambda\leq$$

 $$\leq \tau \left(\int_{-1}^1 f^2(t,\lambda)d\lambda\right)^{1/2} \left(\int_{-1}^1 \left(\sum_{k=0}^\infty\sqrt{\frac{2}{\pi}}U_k(\lambda)\right)^2d\lambda\right)^{1/2} = $$

$$ =  \sqrt{\frac{2}{\pi}}\tau\left(\int_{-1}^1 f^2(t,\lambda)d\lambda\right)^{1/2} \left(\int_{-1}^1\left(\frac{1}{1-2\omega\lambda+\omega^2} \right)^2d\lambda\right)^{1/2}=$$

$$=\sqrt{\frac{2}{\pi}}\tau\left(\int_{-1}^1 f^2(t,\lambda)d\lambda\right)^{1/2} \frac{\sqrt{2}}{(\omega^2-1)} = $$$$=\frac{2\tau}{\sqrt{\pi}(\omega^2-1)}\left(\int_{-1}^1 f^2(t,\lambda)d\lambda\right)^{1/2}.$$

Taking into account these considerations, the estimator of the model of the stochastic process will fit the next condition:

 $$\tau_\varphi(X(t)-X_N(t)) = \tau_\varphi\left(\sum_{k=0}^\infty \xi_ka_k(t)-\sum_{k=0}^N \xi_k\hat{a}_k(t)\right)=\tau_\varphi\left(\sum_{k=0}^N \xi_k\delta_k(t) + \sum_{k=N+1}^\infty \xi_ka_k(t)\right)\leq$$

 $$\leq \sum_{k=0}^N \tau_\varphi(\xi_k)\delta_k(t) + \sum_{k=N+1}^\infty \tau_\varphi(\xi_k) a_k(t)  = \sum_{k=0}^\infty \tau \omega^k a_k(t)-\sum_{k=0}^N \tau \omega^k\hat{a}_k(t)\leq$$

 $$\leq \frac{2\tau}{\sqrt{\pi}(\omega^2-1)}\left(\int_{-1}^1 f^2(t,\lambda)d\lambda\right)^{1/2}-\sum_{k=0}^N \tau \omega^k\hat{a}_k(t).$$

Statement of this theorem follows from Theorem \ref{omlp_g1} and the last inequality.
\end{proof}

\subsubsection{Legendre polynomials}

The Legendre polynomials $P_k(t)$ are defined for $t\in[-1;1]$ as solutions to the Legendre differential equation

$$(1-t^2)y''-2yy'+k(k+1)y=0,\ y=y(t).$$

\noindent The Rodriguez formula for Legendre polynomials if the following:

$$P_k(t)=\frac{1}{2^k k!} \frac{d^k}{dt^k}(t^2-1)^k.$$

\noindent  The generation function for the Legendre polynomials is the following:

$$GF_{Legendre}(t,w)=\sum_{k=0}^\infty P_k(t)w^k = \frac{1}{\sqrt{1-2 t w + w^2}}.$$

\noindent  These polynomials form the orthogonal system:

$$\int_{-1}^1 P_n(t)P_m(t)dt = \frac{2}{2n+1} \delta_{nm},$$
where $\delta_{nm}$ is the Kronecker delta.

\noindent Therefore, this polynomial system can generate the orthonormal basis $\{\hat{P}_k(t),\,t\in[-1;1]\}$, where
\begin{equation}
\label{P-legendre}
\hat{P}_k(t)=\sqrt{\frac{2k+1}{2}}P_k(t).
\end{equation}
Unlike the polynomial systems considered in \cite{mokl2022}, the generating function of the system $\{\hat{P}_k(t)\}$ is hard to derive. It may even not be possible to express this generating function in elementary functions. Therefore, to compensate for elements that obstruct calculation of this function, we will introduce the upper bound on the sub-Gaussian norm in the following form:

\begin{equation}
\label{tau-legendre}
\tau_\varphi(\xi_k)\leq \tau_{Legendre}(w,k)=\sqrt{\frac{2}{2k+1}}\tau w^k,
\end{equation}
$\tau$ is a constant.
From this moment and further in the paper, we will also assume that
$\tau_\varphi (a_k(u))\leq \tau_\varphi (\hat{a_k}(u))$, $\forall u\in [-1,1]$.

To estimate $C_N$, we will apply the following estimates:
$$\tau_\varphi(\Delta_N(t))\leq  \sum_{k=0}^\infty \tau_\varphi(\xi_k) a_k(t) - \sum_{k=0}^N \tau_\varphi(\xi_k) \hat{a}_k(t)  $$
$$=\sum_{k=0}^\infty \tau_\varphi(\xi_k) \int_{-1}^1 f(t,\lambda) \hat{P}_k(\lambda) d\lambda - \sum_{k=0}^N \tau_\varphi(\xi_k) \hat{a}_k(t) $$
$$= \int_{-1}^1 \left(
 {\sum_{k=0}^{\infty}}\tau_\varphi(\xi_k)f(t,\lambda) \hat{P}_k(\lambda) \right)d\lambda- \sum_{k=0}^N \tau_\varphi(\xi_k) \hat{a}_k(t) $$
$$\leq \int_{-1}^1 \left( {\sum_{k=0}^{\infty}}\tau_\varphi(\xi_k)f(t,\lambda) \hat{P}_k(\lambda) \right)d\lambda- \sum_{k=0}^N \tau_\varphi(\xi_k) \hat{a}_k(t) $$
$$\leq \left( \int_{-1}^1 |f(t,\lambda)|^2d\lambda\right)^{1/2}\left(\int_{-1}^1 \left( {\sum_{k=0}^{\infty} }\tau_\varphi(\xi_k) \hat{P}_k(\lambda)\right)^2 d\lambda\right)^{1/2}- \sum_{k=0}^N \tau_\varphi(\xi_k) \hat{a}_k(t). $$

Applying inequality (\ref{tau-legendre}), we will have
$$\left( \int_{-1}^1 |f(t,\lambda)|^2d\lambda\right)^{1/2}\left(\int_{-1}^1 \left( {\sum_{k=0}^{\infty}} \tau_\varphi(\xi_k) \hat{P}_k(\lambda)\right)^2 d\lambda\right)^{1/2}- \sum_{k=0}^N \tau_\varphi(\xi_k) \hat{a}_k(t)  $$
$$\leq \left( \int_{-1}^1 |f(t,\lambda)|^2d\lambda\right)^{1/2}\left(\int_{-1}^1 \left( {\sum_{k=0}^{\infty}} \sqrt{\frac{2}{2k+1}}\tau w^k \hat{P}_k(\lambda)\right)^2 d\lambda\right)^{1/2}- \sum_{k=0}^N \sqrt{\frac{2}{2k+1}}\tau w^k \hat{a}_k(t) $$
$$= \left( \int_{-1}^1 |f(t,\lambda)|^2d\lambda\right)^{1/2}\left(\int_{-1}^1 \left( {\sum_{k=0}^{\infty}} \sqrt{\frac{2}{2k+1}}\tau w^k \sqrt{\frac{2k+1}{2}}P_k(\lambda)\right)^2 d\lambda\right)^{1/2}
$$
$$- \sum_{k=0}^N \sqrt{\frac{2}{2k+1}}\tau w^k \hat{a}_k(t) $$

Canceling the square roots and substituting the first series with the corresponding generating function, we obtain
$$ {\tau} \left( \int_{-1}^1 |f(t,\lambda)|^2d\lambda\right)^{1/2}\left(\int_{-1}^1 \left( {\sum_{k=0}^{\infty}} w^k P_k(\lambda)\right)^2 d\lambda\right)^{1/2}- \sum_{k=0}^N \sqrt{\frac{2}{2k+1}}\tau w^k \hat{a}_k(t)
$$
$$= \tau \left( \int_{-1}^1 |f(t,\lambda)|^2d\lambda\right)^{1/2}\left(\int_{-1}^1 \left(\frac{1}{\sqrt{1-2 \lambda w + w^2}}\right)^2 d\lambda\right)^{1/2}- \sum_{k=0}^N \sqrt{\frac{2}{2k+1}}\tau w^k \hat{a}_k(t)
$$
$$= \tau \left( \int_{-1}^1 |f(t,\lambda)|^2d\lambda\right)^{1/2}\left(\int_{-1}^1 \frac{1}{1-2 \lambda w + w^2}d\lambda\right)^{1/2}  - \sum_{k=0}^N \sqrt{\frac{2}{2k+1}}\tau w^k \hat{a}_k(t)
$$
$$= \tau \left( \int_{-1}^1 |f(t,\lambda)|^2d\lambda\right)^{1/2}\left(\left.\frac{\ln(1-2\lambda w + w^2)}{2w}\right|_{-1}^1\right)^{1/2}  - \sum_{k=0}^N \sqrt{\frac{2}{2k+1}}\tau w^k \hat{a}_k(t)
$$
$$= \tau \left( \int_{-1}^1 |f(t,\lambda)|^2d\lambda\right)^{1/2}\frac{1}{\sqrt{w}}\sqrt{\ln\left(\frac{w+1}{w-1}\right)}  - \sum_{k=0}^N \sqrt{\frac{2}{2k+1}}\tau w^k \hat{a}_k(t). $$

Therefore, we have the following estimate for the $C_N$:

$$C_N\leq {\int_0^T }\left(\tau \left( \int_{-1}^1 |f(t,\lambda)|^2d\lambda\right)^{1/2}\frac{1}{\sqrt{w}}\sqrt{\ln\left(\frac{w+1}{w-1}\right)}  - \sum_{k=0}^N \sqrt{\frac{2}{2k+1}}\tau w^k \hat{a}_k(t)\right)^p d\mu(t).$$

Let us denote the right side of the last inequality by
$C_{N,Legendre}$.

It follows from  Theorem \ref{lp_l2} and  Theorem \ref{lp_g2}, that the following statements are true.

\begin{theorem}
\label{Legendre1}
 Let a stochastic process $X=\{X(t),t\in [0,T]\}$ belong to the space ${Sub}_\varphi (\Omega)$ with the Orlicz $N$-function
  $$\varphi(t)=\frac{t^\gamma}{\gamma},\,\,1<\gamma\leq 2,$$
  and let the process $X(t)$ admits the orthogonal decomposition  (\ref{ombasis})
 based on Legendre orthonormal polynomial families (\ref{P-legendre}).
Assume that $C_{N,Legendre}<\infty$,
$\tau_\varphi (a_k(u))\leq \tau_\varphi (\hat{a}_k(u))$, $\forall u\in [0,\infty)$ and condition \eqref{tau-legendre} holds true.
The model (\ref{ommodel-basis}) $X_N(t)=\sum_{k=0}^N \xi_k \hat{a}_k(t)$ approximates the stochastic process $X(t)$
with given reliability $1-\alpha$ and accuracy $\delta$ in the space $L_p[0,T]$, if
$$\left\{\begin{array}{c}C_{N,Legendre}\leq \delta/(\beta \ln \frac{2}{\alpha})^{p/\beta},\\  C_{N,Legendre}<\delta /p^{p\left(1-1/\gamma\right)},\end{array}\right.$$
where $\beta$ if a number that fits the condition $\frac{1}{\beta}+\frac{1}{\gamma}=1$.
\end{theorem}

\begin{theorem}
\label{Legendre2}
 Let a stochastic process $X=\{X(t),t\in [0,T]\}$ belong to the space ${Sub}_\varphi (\Omega)$ with the Orlicz $N$-function
   $$\varphi(t)=\left\{\begin{array}{c}\frac{t^2}{\gamma}, t<1,\\  \frac{t^\gamma}{\gamma},t\geq1,\end{array}\right.,$$
 where $\gamma>2$,
  and let the process $X(t)$ admits the orthogonal decomposition  (\ref{ombasis})
 based on Legendre orthonormal polynomial families (\ref{P-legendre}).
Assume that $C_{N,Legendre}<\infty$,
$\tau_\varphi (a_k(u))\leq \tau_\varphi (\hat{a}_k(u))$, $\forall u\in [-1,1]$ and condition \eqref{tau-legendre} holds true.
 The model (\ref{ommodel-basis}) $X_N(t)=\sum_{k=0}^N \xi_k \hat{a}_k(t)$ approximates the stochastic process $X(t)$
with given reliability $1-\alpha$ and accuracy $\delta$ in the space $L_p[0,T]$, if
$$\left\{\begin{array}{c}C_{N,Legendre}\leq \delta/(\beta \ln \frac{2}{\alpha})^{p/\beta},\\  C_{N,Legendre}<\delta /p^{p\left(1-1/\gamma\right)},\end{array}\right.$$
where $\beta$ if a number that fits the condition $\frac{1}{\beta}+\frac{1}{\gamma}=1$.
\end{theorem}

\subsubsection{Generalized Laguerre polynomials}

The Laguerre polynomials $L_k(t)$ are defined for $t\in[0;\infty]$ as solutions of the Laguerre differential equation

$$ty''+(1-t)y'+ky=0, y=y(t),$$
where $k$ is a non-negative integer. The generalized Laguerre polynomials are solutions of the following equation:

$$ty''+(\alpha + 1-t)y'+ky=0,\ y=y(t).$$
Here, an arbitrary real $\alpha$ is introduced in the equation. Hence, the generalized Laguerre polynomials are denoted as $L_k^{(\alpha)}(t)$.

\noindent The Rodriguez formula for generalized Laguerre polynomials has the form:

$$L_k^{(\alpha)}(t)=\frac{t^{-\alpha}}{k!}\frac{d^k}{dt^k}(e^{-x}x^{n+\alpha})$$

\noindent The generation function for generalized Laguerre polynomials is the following:

$$GF_{Laguerre}(t,w)=\sum_{k=0}^\infty L^{(\alpha)}_k(t)w^k = \frac{1}{(1-w)^{\alpha+1}}e^{-wt/(1-w)}$$

\noindent The generalized Laguerre polynomials are orthogonal over $[0,\infty]$ with respect to the measure with weighting function $t^\alpha e^{-t}:$

$$\int_0^\infty x^\alpha e^{-x}L^{(\alpha)}_n(t)L^{(\alpha)}_m(t)dt=\frac{\Gamma(n+\alpha+1)}{n!}\delta_{nm},$$
therefore, the orthonormal version of the generalized Laguerre polynomials will take the form
\begin{equation}
\label{P-laguerre}
\hat{L}^{(\alpha)}_k(t)=\sqrt{\frac{k!}{\Gamma(k+\alpha+1)}}x^{\alpha/2}e^{-x/2}L^{(\alpha)}_k(t).
\end{equation}
\noindent The upper bound for the sub-Gaussian norm of generalized Laguerre polynomials will take the form

\begin{equation}
\label{tau-laguerre}
\tau_\varphi(\xi_k)\leq\tau_{Laguerre} (w,k)= \sqrt{\frac{\Gamma(k+\alpha+1)}{k!}}\tau w^k
\end{equation}

\noindent Let us estimate the value of $C_N$  in this case. Following the same procedure as provided above for Legendre polynomials, with inequality (\ref{tau-laguerre}) in mind,

$$\tau_\varphi(\Delta_N(t))\leq  \sum_{k=0}^\infty \tau_\varphi(\xi_k) a_k(t) - \sum_{k=0}^N \tau_\varphi(\xi_k) \hat{a}_k(t)  $$
$$=\sum_{k=0}^\infty \tau_\varphi(\xi_k) \int_0^\infty f(t,\lambda) \hat{L}^{(\alpha)}_k(\lambda) d\lambda - \sum_{k=0}^N \tau_\varphi(\xi_k) \hat{a}_k(t) $$
$$= \int_0^\infty
\left(
\sum_{k=0}^{\infty}
\tau_\varphi(\xi_k)f(t,\lambda) \hat{L}^{(\alpha)}_k(\lambda) \right)d\lambda- \sum_{k=0}^N \tau_\varphi(\xi_k) \hat{a}_k(t) $$
$$\leq \left( \int_0^\infty |f(t,\lambda)|^2d\lambda\right)^{1/2}\left(\int_0^\infty \left(
\sum_{k=0}^{\infty}
\tau_\varphi(\xi_k) \hat{L}^{(\alpha)}_k(\lambda)\right)^2 d\lambda\right)^{1/2}- \sum_{k=0}^N \tau_\varphi(\xi_k) \hat{a}_k(t) $$
$$\leq  \left( \int_0^\infty|f(t,\lambda)|^2d\lambda\right)^{1/2}
$$
$$\times
\left(\int_0^\infty \left(
\sum_{k=0}^{\infty}
\sqrt{\frac{\Gamma(k+\alpha+1)}{k!}}\tau w^k \sqrt{\frac{k!}{\Gamma(k+\alpha+1)}}\lambda^{\alpha/2}e^{-\lambda/2}L^{(\alpha)}_k(\lambda)\right)^2 d\lambda\right)^{1/2}$$
$$- \sum_{k=0}^N \sqrt{\frac{\Gamma(k+\alpha+1)}{k!}}\tau w^k \hat{a}_k(t)  $$
$$= \tau \left( \int_0^\infty |f(t,\lambda)|^2d\lambda\right)^{1/2}\left(\int_0^\infty \lambda^{\alpha}e^{-\lambda}\left(
\sum_{k=0}^{\infty}
 w^k L^{(\alpha)}_k(\lambda)\right)^2 d\lambda\right)^{1/2}
 $$
 $$
 - \sum_{k=0}^N \sqrt{\frac{\Gamma(k+\alpha+1)}{k!}}\tau w^k \hat{a}_k(t)
 $$
$$= \tau \left( \int_0^\infty |f(t,\lambda)|^2d\lambda\right)^{1/2}\left(\int_0^\infty \lambda^{\alpha}e^{-\lambda}\frac{1}{(1-w)^{2\alpha+2}}e^{-2w\lambda/(1-w)} d\lambda\right)^{1/2}
$$
$$- \sum_{k=0}^N \sqrt{\frac{\Gamma(k+\alpha+1)}{k!}}\tau w^k \hat{a}_k(t)
$$
$$= \tau \left( \int_0^\infty |f(t,\lambda)|^2d\lambda\right)^{1/2}\left(\left.-\lambda^{\alpha+1} \left(\frac{1-w}{\lambda(1+w)}\right)^{\alpha+1}\Gamma\left(\alpha+1,\frac{\lambda(1+w)}{1-w}\right)\right|_0^\infty\right)^{1/2}
 $$
 $$- \sum_{k=0}^N \sqrt{\frac{\Gamma(k+\alpha+1)}{k!}}\tau w^k \hat{a}_k(t),$$

\noindent where $\Gamma(s,x)$ is the upper incomplete gamma function. Applying the integration boundaries, we obtain

$$\tau \left( \int_0^\infty |f(t,\lambda)|^2d\lambda\right)^{1/2}\left(\left.-\lambda^{\alpha+1} \left(\frac{1-w}{\lambda(1+w)}\right)^{\alpha+1}\Gamma\left(\alpha+1,\frac{\lambda(1+w)}{1-w}\right)\right|_0^\infty\right)^{1/2} $$
 $$- \sum_{k=0}^N \sqrt{\frac{\Gamma(k+\alpha+1)}{k!}}\tau w^k \hat{a}_k(t)
 $$
$$= \tau \left( \int_0^\infty |f(t,\lambda)|^2d\lambda\right)^{1/2}\left(\frac{1-w}{1+w}\right)^{\alpha+1}\Gamma(\alpha+1) - \sum_{k=0}^N \sqrt{\frac{\Gamma(k+\alpha+1)}{k!}}\tau w^k \hat{a}_k(t), $$
and this expression is defined for $\alpha>-1$.

\noindent As a result, we have

$$C_N\leq
{\int_0^T}
 \left(\tau \biggl( \int_0^\infty |f(t,\lambda)|^2d\lambda\right)^{1/2}\left(\frac{1-w}{1+w}\right)^{\alpha+1}\Gamma(\alpha+1)$$
  $$- \sum_{k=0}^N \sqrt{\frac{\Gamma(k+\alpha+1)}{k!}}\tau w^k \hat{a}_k(t)\biggr)^p d\mu(t).$$

Let us denote the right side of the last inequality by
$C_{N,Laguerre}$.

It follows from  Theorem \ref{lp_l2} and  Theorem \ref{lp_g2}, that the following statements are true.

\begin{theorem}
\label{Laguerre1}
 Let a stochastic process $X=\{X(t),t\in[0,T]\}$ belong to the space ${Sub}_\varphi (\Omega)$ with the Orlicz $N$-function
  $$\varphi(t)=\frac{t^\gamma}{\gamma},\,\,1<\gamma\leq 2,$$
  and let the process $X(t)$ admits the orthogonal decomposition  (\ref{ombasis})
 based on Laguerre orthonormal polynomial families (\ref{P-laguerre}).
Assume that $C_{N,Laguerre}<\infty$,
$\tau_\varphi (a_k(u))\leq \tau_\varphi (\hat{a}_k(u))$, $\forall u\in [-1,1]$ and condition \eqref{tau-laguerre} holds true.

\noindent The model (\ref{ommodel-basis}) $X_N(t)=\sum_{k=0}^N \xi_k \hat{a}_k(t)$ approximates the stochastic process $X(t)$
with given reliability $1-\alpha$ and accuracy $\delta$ in the space $L_p[0,T]$, if
$$\left\{\begin{array}{c}
C_{N,Laguerre}\leq \delta/(\beta \ln \frac{2}{\alpha})^{p/\beta},\\
 C_{N,Laguerre}<\delta /p^{p\left(1-1/\gamma\right)},
\end{array}\right.$$
where $\beta$ if a number that fits the condition $\frac{1}{\beta}+\frac{1}{\gamma}=1$.
\end{theorem}

\begin{theorem}
\label{Laguerre2}
 Let a stochastic process $X=\{X(t),t\in[0,T]\}$ belong to the space ${Sub}_\varphi (\Omega)$ with the Orlicz $N$-function
   $$\varphi(t)=\left\{\begin{array}{c}\frac{t^2}{\gamma}, t<1,\\
    \frac{t^\gamma}{\gamma},t\geq1,
   \end{array}\right.,$$
 where $\gamma>2$,
  and let the process $X(t)$ admits the orthogonal decomposition  (\ref{ombasis})
 based on Laguerre orthonormal polynomial families (\ref{P-laguerre}).
Assume that $C_{N,Laguerre}<\infty$,
$\tau_\varphi (a_k(u))\leq \tau_\varphi (\hat{a}_k(u))$, $\forall u\in[0;\infty)$ and condition \eqref{tau-laguerre} holds true.
 The model (\ref{ommodel-basis}) $X_N(t)=\sum_{k=0}^N \xi_k \hat{a}_k(t)$ approximates the stochastic process $X(t)$
with given reliability $1-\alpha$ and accuracy $\delta$ in the space $L_p[0,T]$, if
$$\left\{\begin{array}{c}
C_{N,Laguerre}\leq \delta/(\beta \ln \frac{2}{\alpha})^{p/\beta},\\
 C_{N,Laguerre}<\delta /p^{p\left(1-1/\gamma\right)},
\end{array}\right.$$
where $\beta$ if a number that fits the condition $\frac{1}{\beta}+\frac{1}{\gamma}=1$.
\end{theorem}

\subsection{Gegenbauer polynomials}

The Gegenbauer polynomials $C^{(\alpha)}_k(t)$ are defined for $t\in[-1; 1]$ as solutions of the Gegenbauer differential equation:

$$(1-t^2)y''-(2\alpha+1)ty'+k(k+2\alpha)y=0.$$
When $\alpha=1/2$, this equation reduces to the Legendre equation. Therefore, for the mentioned $\alpha=1/2$, the Gegenbauer polynomials are reduced to the Legendre polynomials.

\noindent The Gegenbauer polynomials can be represented as Gaussian hypergeometric series when this series is finite:

$$C_k^{(\alpha)}(t) = \frac{(2\alpha)_k}{k!}{}_2F_1\left(-k,2\alpha+k;\alpha+\frac{1}{2};\frac{1-t}{2}\right),$$
where ${}_2F_1(a,b;c;z)$ is a special function which is a solution of the Euler hypergeometric differential equation:

$$t(1-t)\frac{d^2 y }{dt^2} + (c-(a+b+1)t)\frac{dy}{dt} - a b y =0.$$

\noindent This function can also be represented in the form of a series

$${}_2F_1(a,b;c;z) = \frac{\Gamma(c)}{\Gamma(a)\Gamma(b)}\sum_{n=0}^\infty \frac{\Gamma(a+n)\Gamma(b+n)}{\Gamma(c+n)}\frac{z^n}{n!}.$$

\noindent The Rodrigues formula for the Gegenbauer polynomials is the following:

$$C_k^{(\alpha)}(t)=\frac{(-1)^k}{2^k k!} \frac{\Gamma(\alpha+\frac{1}{2})\Gamma(k+2\alpha)}{\Gamma(2\alpha)\Gamma(\alpha+k+\frac{1}{2})}(1-t^2)^{1/2-\alpha} \frac{d^k}{dt^k}(1-t^2)^{k+\alpha-1/2}.$$

\noindent The generating function is given as

$$\sum_{k=0}^\infty C_k^{(\alpha)} (t)w^k = \frac{1}{(1-2tw+w^2)^\alpha}.$$

\noindent  The Gegenbauer polynomials are orthogonal:

$$\int_{-1}^1 (1-t^2)^{\alpha-1/2}C_n^{(\alpha)} (t)C_m^{(\alpha)} (t)dt=\frac{\pi 2^{1-2\alpha}\Gamma(k+2\alpha)}{k! (k+\alpha) \Gamma^2 (\alpha)}\delta_{nm}.$$

\noindent Therefore, the orthonormal version of the generalized Gegenbauer polynomials is
\begin{equation}
\label{P-gegenbauer}
\hat{C}^{(\alpha)}_k(t)=\frac{\Gamma (\alpha)\sqrt{k! (k+\alpha)}}{\sqrt{\pi 2^{1-2\alpha}}\sqrt{\Gamma(k+2\alpha)}}{\sqrt{(1-t^2)^{\alpha-1/2}}}C^{(\alpha)}_k(t).
\end{equation}
\noindent We can assume that
\begin{equation}
\label{tau-gegenbauer}
\tau_\varphi(\xi_k)\leq\tau_{Gegenbauer} (w,k)= \sqrt{\frac{k! (k+\alpha)}{\Gamma(k+2\alpha)}}\tau w^k
\end{equation}

\noindent Let us estimate the value of $C_N$. We have

$$\tau_\varphi(\Delta_N(t))\leq  \sum_{k=0}^\infty \tau_\varphi(\xi_k) a_k(t) - \sum_{k=0}^N \tau_\varphi(\xi_k) \hat{a}_k(t)
$$
$$=\sum_{k=0}^\infty \tau_\varphi(\xi_k) \int_{-1}^1 f(t,\lambda) \hat{C}^{(\alpha)}_k(\lambda) d\lambda - \sum_{k=0}^N \tau_\varphi(\xi_k) \hat{a}_k(t)
$$
$$= \int_{-1}^1 \left(
\sum_{k=0}^{\infty}
\tau_\varphi(\xi_k)f(t,\lambda) \hat{C}^{(\alpha)}_k(\lambda) \right)d\lambda- \sum_{k=0}^N \tau_\varphi(\xi_k) \hat{a}_k(t)
$$
$$\leq \left( \int_{-1}^1 |f(t,\lambda)|^2d\lambda\right)^{1/2}\left(\int_{-1}^1 \left(
\sum_{k=0}^{\infty}
\tau_\varphi(\xi_k) \hat{C}^{(\alpha)}_k(\lambda)\right)^2 d\lambda\right)^{1/2}- \sum_{k=0}^N \tau_\varphi(\xi_k) \hat{a}_k(t)
$$
$$\leq  \left( \int_{-1}^1 |f(t,\lambda)|^2d\lambda\right)^{1/2}\left(\int_{-1}^1 \left(
\sum_{k=0}^{\infty}
\frac{\Gamma (\alpha)}{\sqrt{\pi 2^{1-2\alpha}}}\tau w^k{\sqrt{(1-\lambda ^2)^{\alpha-1/2}}}C^{(\alpha)}_k(\lambda)\right)^2 d\lambda\right)^{1/2}
$$
$$- \sum_{k=0}^N \sqrt{\frac{k! (k+\alpha)}{\Gamma(k+2\alpha)}}\tau w^k \hat{a}_k(t)
$$
$$=   \frac{\tau\Gamma (\alpha)}{\sqrt{\pi 2^{1-2\alpha}}}
\left( \int_{-1}^1 |f(t,\lambda)|^2d\lambda\right)^{1/2}
\left(\int_{-1}^1 (1-\lambda ^2)^{\alpha-1/2}\left(\sum_{k=0}^N  w^k C^{(\alpha)}_k(\lambda)\right)^2 d\lambda\right)^{1/2}
$$
$$- \sum_{k=0}^N \sqrt{\frac{k! (k+\alpha)}{\Gamma(k+2\alpha)}}\tau w^k \hat{a}_k(t)
$$
$$=  \frac{\tau\Gamma (\alpha)}{\sqrt{\pi 2^{1-2\alpha}}}\left( \int_{-1}^1 |f(t,\lambda)|^2d\lambda\right)^{1/2}\left(\int_{-1}^1 \frac{(1-\lambda ^2)^{\alpha-1/2}}{(1-2\lambda w+w^2)^{2\alpha}} d\lambda\right)^{1/2}-$$ $$- \sum_{k=0}^N \sqrt{\frac{k! (k+\alpha)}{\Gamma(k+2\alpha)}}\tau w^k \hat{a}_k(t)
$$
$$=  \frac{\tau\Gamma (\alpha)}{\sqrt{\pi 2^{1-2\alpha}}}
\left( \frac{1}{(2 \alpha -1) w} 4^{-\alpha } \left(1-\lambda ^2\right)^{\alpha - 1/2}    \left(\frac{(\lambda -1) (\lambda +1) w^2}{(w+1)^2(w-1)^2}\right)^{1/2-\alpha }  \left(w^2-2    \lambda  w+1\right)^{1-2 \alpha }\right.$$
$$\times \left.\left. F_1\left(1-2 \alpha ;\frac{1}{2}-\alpha    ,\frac{1}{2}-\alpha ;2-2 \alpha ;\frac{w^2-2 \lambda  w+1}{(w+1)^2},\frac{w^2-2    \lambda  w+1}{(w-1)^2}\right)  \right|_{-1}^1 \right)^{1/2}
   $$
   $$
   \times
   \left( \int_{-1}^1 |f(t,\lambda)|^2d\lambda\right)^{1/2}- \sum_{k=0}^N \sqrt{\frac{k! (k+\alpha)}{\Gamma(k+2\alpha)}}\tau w^k \hat{a}_k(t)
   $$
$$=  \frac{\tau\Gamma (\alpha)}{\sqrt{\pi 2^{1-2\alpha}}}
\left( \sqrt{\pi } \Gamma \left(\alpha +\frac{1}{2}\right) \left(w^2+1\right)^{-2
   \alpha } \, _2\tilde{F}_1\left(\alpha ,\alpha +\frac{1}{2};\alpha +1;\frac{4
   w^2}{\left(w^2+1\right)^2}\right) \right)^{1/2}
   $$
   $$
   \times
   \left( \int_{-1}^1 |f(t,\lambda)|^2d\lambda\right)^{1/2} - \sum_{k=0}^N \sqrt{\frac{k! (k+\alpha)}{\Gamma(k+2\alpha)}}\tau w^k \hat{a}_k(t)
   $$
$$=  \frac{\tau\Gamma (\alpha)}{\left(w^2+1\right)^{
   \alpha }}\sqrt{\frac{\Gamma \left(\alpha +\frac{1}{2}\right)}{\sqrt{\pi} 2^{1-2\alpha}}}    \, _2\tilde{F}_1^{1/2}\left(\alpha ,\alpha +\frac{1}{2};\alpha +1;\frac{4
   w^2}{\left(w^2+1\right)^2}\right)\left( \int_{-1}^1 |f(t,\lambda)|^2d\lambda\right)^{1/2} $$
   $$- \sum_{k=0}^N \sqrt{\frac{k! (k+\alpha)}{\Gamma(k+2\alpha)}}\tau w^k \hat{a}_k(t), $$
where $\, _2\tilde{F}_1(a,b;c;z)$ is the regularized hypergeometric function $\, _2F_1(a,b;c;z)$

As a result, we have

$$C_N\leq
{\int_0^T}
\left(\frac{\tau\Gamma (\alpha)}{\left(w^2+1\right)^{
   \alpha }}\sqrt{\frac{\Gamma \left(\alpha +\frac{1}{2}\right)}{\sqrt{\pi} 2^{1-2\alpha}}}    \, _2\tilde{F}_1^{1/2}\left(\alpha ,\alpha +\frac{1}{2};\alpha +1;\frac{4
   w^2}{\left(w^2+1\right)^2}\right)\right.
   $$
   $$\left.
   \times\left( \int_{-1}^1 |f(t,\lambda)|^2d\lambda\right)^{1/2}-
    \sum_{k=0}^N \sqrt{\frac{k! (k+\alpha)}{\Gamma(k+2\alpha)}}\tau w^k \hat{a}_k(t)\right)^p d\mu(t).$$

Let us denote the right side of the last inequality by
$C_{N,Gegenbauer}$.

It follows from  Theorem \ref{lp_l2} and  Theorem \ref{lp_g2}, that the following statements are true.

\begin{theorem}
\label{Gegenbauer1}
 Let a stochastic process $X=\{X(t),t\in [0,T]\}$ belong to the space ${Sub}_\varphi (\Omega)$ with the Orlicz $N$-function
  $$\varphi(t)=\frac{t^\gamma}{\gamma},\,\,1<\gamma\leq 2,$$
  and let the process $X(t)$ admits the orthogonal decomposition  (\ref{ombasis})
 based on Gegenbauer orthonormal polynomial families (\ref{P-gegenbauer}).
Assume that $C_{N,Gegenbauer}<\infty$,
$\tau_\varphi (a_k(u))\leq \tau_\varphi (\hat{a}_k(u))$, $\forall u\in [-1,1]$ and condition \eqref{tau-gegenbauer} holds true.
 The model (\ref{ommodel-basis}) $X_N(t)=\sum_{k=0}^N \xi_k \hat{a}_k(t)$ approximates the stochastic process $X(t)$
with given reliability $1-\alpha$ and accuracy $\delta$ in the space $L_p[0,T]$, if
$$\left\{\begin{array}{c}C_{N,Gegenbauer}\leq \delta/(\beta \ln \frac{2}{\alpha})^{p/\beta},\\
 C_{N,Gegenbauere}<\delta /p^{p\left(1-1/\gamma\right)},\end{array}\right.$$
where $\beta$ if a number that fits the condition $\frac{1}{\beta}+\frac{1}{\gamma}=1$.
\end{theorem}

\begin{theorem}
\label{Legendre2}
 Let a stochastic process $X=\{X(t),t\in [0,T]\}$ belong to the space ${Sub}_\varphi (\Omega)$ with the Orlicz $N$-function
   $$\varphi(t)=\left\{\begin{array}{c}\frac{t^2}{\gamma}, t<1,\\  \frac{t^\gamma}{\gamma},t\geq1,\end{array}\right.,$$
 where $\gamma>2$,
  and let the process $X(t)$ admits the orthogonal decomposition  (\ref{ombasis})
 based on Gegenbauer orthonormal polynomial families (\ref{P-gegenbauer}).
Assume that $C_{N,Gegenbauer}<\infty$,
$\tau_\varphi (a_k(u))\leq \tau_\varphi (\hat{a}_k(u))$, $\forall u\in [-1,1]$ and condition \eqref{tau-gegenbauer} holds true.
 The model (\ref{ommodel-basis}) $X_N(t)=\sum_{k=0}^N \xi_k \hat{a}_k(t)$ approximates the stochastic process $X(t)$
with given reliability $1-\alpha$ and accuracy $\delta$ in the space $L_p[0,T]$, if
$$\left\{\begin{array}{c}C_{N,Gegenbauer}\leq \delta/(\beta \ln \frac{2}{\alpha})^{p/\beta},\\ C_{N,Gegenbauer}<\delta /p^{p\left(1-1/\gamma\right)},\end{array}\right.$$
where $\beta$ if a number that fits the condition $\frac{1}{\beta}+\frac{1}{\gamma}=1$.
\end{theorem}

\section*{Conclusions to chapter \ref{omch:lpt1}}

In chapter \ref{omch:lpt1}, we consider random processes from the spaces $Sub_\varphi(\Omega)$ that admit a series expansion whose elements cannot be found explicitly.
A mechanism for assessing the accuracy and reliability of constructing models of such processes in the spaces $L_p(T)$ is developed.
The cases $\varphi(t)=|t|^\gamma/\gamma$, $1<\gamma<2$, and $\varphi(t)=|t|^\gamma/\gamma$, when $|t|>1$, and $\varphi(t)=|t|^2/\gamma$, $|t|<1$, when $\gamma>2$ are considered.
The above calculations were used to assess the accuracy and reliability of modeling in $L_p[0,T]$ of random processes using the
Karhunen-Lo\`eve expansion in the case where the eigenvalues and eigenvectors cannot be found explicitly.

\chapter{Accuracy and reliability of modeling in $C(T)$ random processes that allow series expansions with independent terms.}
\label{omch:lpt2}

In this chapter we deal with the modeling of random processes from $Sub_\varphi(\Omega)$ spaces, which are a subclass of $K_\sigma$-spaces, with given reliability and accuracy in $C(T)$ spaces.

\section[{Accuracy and reliability of modeling of random processes in $C(T)$ spaces}]{Assessment of accuracy and reliability of modeling of random processes in $C(T)$ spaces}

Let $X=\{X(t),t\in B\}$ be a random process from $Sub_\varphi(\Omega)$. Let $(B,\rho)$ be a separable compact set and the process $X$ be separable on $(B,\rho)$.
Suppose that there exists a continuously monotonically increasing function $\sigma=\{\sigma(h),h>0\}$ such that $\sigma(h)\to0$, $h\to0$, and the inequality holds:
\begin{equation}
\label{ome2.1}
\sup_{\rho(t,s)\leq h}\tau_\varphi(X(t)-X(s))\leq \sigma(h),
\end{equation}
\noindent and the process $X$ is continuous. Let in addition $\beta>0$ be such that
$$\beta=\sigma\left(\inf_{s\in B}\sup_{t\in B}\rho(t,s)\right),$$
\noindent and let $N_B(u)$ be the minimal number of closed balls of radius $u$ that cover $(B,\rho)$.

The following theorem is a modification of the theorem from \cite{vasyl2008}.

\begin{theorem}
\label{omtk.1}
Let a random process $X=\{X(t),t\in B\}$ belong to $Sub_\varphi(\Omega)$, let $(B,\rho)$ be a separable compact set, and let $X$ be separable on $(B,\rho)$.
Let the process $X$ satisfy the condition \eqref{ome2.1}, and let $r_1=\{r_1(u):u\geq1\}$ be a continuous function such that $r_1(u)>0$ for $u>1$, and let the function $s(t)=r_1(\exp\{t\})$, $t>0$ be convex. Then, under the condition

$$\int_0^\beta r_1(N_B(\sigma^{(-1)}(u)))du<\infty,$$

\noindent the process $X(t)$ is bounded with probability one, and for all $p\in (0,1)$ and $x>0$ the inequalities hold:

$$P\{\sup_{t\in (0,T)}X(t)>x\}\leq Z_{r_1}(p,\beta,x),$$
$$P\{\inf_{t\in (0,T)}X(t)<-x\}\leq Z_{r_1}(p,\beta,x),$$
$$P\{\sup_{t\in (0,T)}|X(t)|>x\}\leq 2Z_{r_1}(p,\beta,x),$$

$$Z_{r_1}(p,\beta,x) = \inf_{\lambda>0}\exp\left\{ \theta_\varphi(\lambda,p)+p\varphi(\frac{\lambda \beta}{1-p}) - \lambda x \right\}\times$$

$$\times r_1^{(-1)}\left(\frac{1}{\beta p}\int_0^{\beta p}r_1(N_B(\sigma^{(-1)}(u)))du\right),$$

\noindent where $\theta_\varphi(\lambda,p) = \sup_{u\in B}\left((1-p)\varphi(\frac{\gamma(u)\lambda}{1-p})\right)$, $\gamma(u)=\tau_\varphi(X(u)).$

\end{theorem}
\begin{theorem}
\label{omctf-square}
Let the random process $X=\{X(t),t\in B\}$ belong to $Sub_\varphi(\Omega)$, let $(B,\rho)$ be a separable compact set, and let $X$ be separable on $(B,\rho)$.
Let $\varphi(t)=\frac{t^\zeta}{\zeta}$, $\zeta\geq2$, $t>1$, and for the process $X$ conditions of the previous theorem are fulfilled.
Then the process $X(t)$ is bounded with probability one, and for all $p\in (0,1)$ the inequalities hold:

$$P\{\sup_{t\in (0,T)}X(t)>x\}\leq Z_{r_1}(p,\beta,x),$$
$$P\{\inf_{t\in (0,T)}X(t)<-x\}\leq Z_{r_1}(p,\beta,x),$$
$$P\{\sup_{t\in (0,T)}|X(t)|>x\}\leq 2Z_{r_1}(p,\beta,x),$$

\noindent when

$$Z_{r_1}(p,\beta,x)= \exp\left\{\frac{(1-\zeta)(x(1-p))^\frac{\zeta}{\zeta-1}}{\zeta(\gamma^\zeta(1-p)+p\beta^\zeta)^\frac{1}{\zeta-1}}\right\}\times$$

$$\times r_1^{(-1)}\left(\frac{1}{\beta p}\int_0^{\beta p}r_1(N_B(\sigma^{(-1)}(u)))du\right),$$ $$x>\frac{\gamma^\zeta(1-p)+p\beta^\zeta}{(1-p)v^{\zeta-1}},$$ where $v=\min\{\beta,\gamma\}$, $\gamma^\zeta=\sup_{u\in B} \tau^\zeta_\varphi(X(u))$.
\end{theorem}

\begin{dov}
Let
\begin{equation}
\label{omneq-gb}
(1-p)\leq \gamma\lambda \mbox{  and  } (1-p)\leq \beta\lambda.
 \end{equation}
\noindent Consider the power of the exponent in the expression for $Z_{r_1}(p,x,\beta)$ from the previous theorem. Substitute the values for $\varphi(t)$ and $\theta_\varphi(\lambda,p):$
\begin{equation}
\label{omtodif}
\sup_{u\in (0,T)}\left((1-p)\varphi\left(\frac{\gamma(u)\lambda}{1-p}\right)\right)+p\varphi\left(\frac{\lambda \beta}{1-p}\right) - \lambda x  =
\end{equation}
$$=\sup_{u\in B}\left(\frac{\gamma^\zeta(u)\lambda^\zeta}{\zeta(1-p)^{\zeta-1}}\right)+p\frac{\lambda^\zeta \beta^\zeta}{\zeta(1-p)^\zeta} - \lambda x   = $$
$$ = \frac{\lambda^\zeta}{\zeta}\left( \frac{\sup_{u\in B}\gamma^\zeta(u)(1-p)+p\beta^\zeta}{(1-p)^\zeta} \right)-\lambda x.$$
\noindent Let us denote $\gamma^\zeta:=\sup_{u\in B}\gamma^\zeta(u) = \sup_{u\in B}\tau_\varphi^\zeta(X(u))$. Then the minimum of the expression (\ref{omtodif}) is reached at the point
$$\lambda = \left( \frac{x(1-p)^\zeta}{\gamma^\zeta(1-p)+p\beta^\zeta} \right)^\frac{1}{\zeta-1}.$$
\noindent Substituting this value into (\ref{omtodif}), we obtain that the minimum of this expression is equal to
$$\frac{(1-\zeta)(x(1-p))^\frac{\zeta}{\zeta-1}}{\zeta(\gamma^\zeta(1-p)+p\beta^\zeta)^\frac{1}{\zeta-1}}.$$
In this case, for
$$\lambda = \left( \frac{x(1-p)^\zeta}{\gamma^\zeta(1-p)+p\beta^\zeta} \right)^\frac{1}{\zeta-1}$$
\noindent inequalities (\ref{omneq-gb}) are satisfied when the expression (\ref{omtodif}) reaches a minimum, that is, when
$$x>\frac{\gamma^\zeta(1-p)+p\beta^\zeta}{(1-p)v^{\zeta-1}},$$
where $v=\min\{\beta,\gamma\}$.
\end{dov}

\begin{definition}
\label{omct-model}
The model $X_N(t)$ approximates the process $X(t)$ with given reliability $1-\nu$ and accuracy $\delta$ in the space $C(B)$ if
$$P\{\sup_{t\in (B)}|\Delta_N(t)|>\delta\}\leq \nu,$$ where $$\Delta_N(t) = X(t)-X_N(t).$$
\end{definition}

Let there exist a continuously monotonically increasing function $\sigma_N=\{\sigma_N(h),h>0\}$ such that $\sigma_N(h)\to0$, $h\to0$, and the inequality holds:

\begin{equation}
\label{omsigma-N}
 \sup_{\rho(t,s)\leq h}\tau_\varphi(X_N(t)-X_N(s))\leq \sigma_N(h).
\end{equation}

\begin{corollary}
Let the random process $X=\{X(t),t\in B\}$ belong to $Sub_\varphi(\Omega)$, let $(B,\rho)$ be a separable compact set, and let $X$ be separable on $(B,\rho)$.
Let $\varphi(t)=\frac{t^\zeta}{\zeta}$, $\zeta\geq2$, $t>1$, and for the process $X$ the conditions of Theorem \ref{omtk.1} are fulfilled.
The model $X_N(t)$ approximates the process $X(t)$ with given reliability $1-\nu$ and accuracy $\delta$ in the space $C(B)$, if
$$\nu\leq 2\exp\left\{-\frac{(\zeta-1)(\delta(1-p))^\frac{\zeta}{\zeta-1}}{\zeta(\gamma_N^\zeta(1-p)+p\beta^\zeta)^\frac{1}{\zeta-1}}\right\}r_1^{(-1)}\left(\frac{1}{\beta p}\int_0^{\beta p}r_1(N_B(\sigma_N^{(-1)}(u)))du\right),$$
$$\delta>\frac{\gamma_N^\zeta(1-p)+p\beta^\zeta}{(1-p)v^{\zeta-1}},$$
where $v=\min\{\beta,\gamma_N \}$, $\gamma_N^\zeta=\sup_{u\in [0,T]} \tau^\zeta_\varphi(\Delta_N(u))$.
\end{corollary}

\begin{theorem}

Let the random process $X=\{X(t),t\in B\}$ belong to $Sub_\varphi(\Omega)$, let $(B,\rho)$ be a separable compact set, and let $X$ is separable on $(B,\rho)$.
Let $\varphi(t)=\frac{t^\zeta}{\zeta}$, $1<\zeta\leq2$, and for the process $X$ the conditions of Theorem \ref{omtk.1} are fulfilled.
The model $X_N(t)$ approximates the process $X(t)$ with given reliability $1-\nu$ and accuracy $\delta$ in the space $C(B)$, if
$$\nu\leq 2\exp\left\{-\frac{(\zeta-1)(\delta(1-p))^\frac{\zeta}{\zeta-1}}{\zeta(\gamma_N^\zeta(1-p)+p\beta^\zeta)^\frac{1}{\zeta-1}}\right\}r_1^{(-1)}\left(\frac{1}{\beta p}\int_0^{\beta p}r_1(N_B(\sigma_N^{(-1)}(u)))du\right),$$
\noindent where $\gamma_N^\zeta=\sup_{u\in [0,T]} \tau^\zeta_\varphi(\Delta_N(u)).$
\end{theorem}

\begin{dov} This theorem follows from Theorem \ref{omctf-square}, but in this case the restriction on $\delta$ is not necessary.
\end{dov}

\begin{theorem}
\label{omctfz-square}
Let a random process $X=\{X(t),t\in [0,T]\}$ belong to the space $Sub_\varphi(\Omega)$, $\varphi(t)=\frac{t^\zeta}{\zeta}$, $\zeta\geq2$, $t>1$.
Let $X$ be separable, and the condition (\ref{ome2.1}) holds for it, where $\sigma(h) = Ch^{\ae}$.
Then the process $X(t)$ is bounded with probability one, and the inequalities hold:

$$P\{\sup_{t\in (0,T)}X(t)>x\}\leq Z_{r_1}(x),$$
$$P\{\inf_{t\in (0,T)}X(t)<-x\}\leq Z_{r_1}(x),$$
$$P\{\sup_{t\in (0,T)}|X(t)|>x\}\leq 2Z_{r_1}(x),$$

\noindent where

$$Z_{r_1}(x) = \exp\left\{-(x-\gamma)^\frac{\zeta}{\zeta-1}\frac{(\zeta-1)x^\frac{1}{\zeta-1}}{\zeta(\gamma^\zeta(x-\gamma)+\beta^\zeta\gamma)^\frac{1}{\zeta-1}}\right\}2\left(ex\right)^{1/\ae},$$ moreover,
$$x>\frac{\gamma(v^{\zeta-1}+1)+\sqrt{\gamma^2(v^{\zeta-1}+1)^2+4v^{\zeta-1}(C^\zeta(T/2)^{\ae\zeta}-\gamma^2)}}{2v^{\zeta-1}},$$
when $v=\min\{C(T/2)^{\ae},\gamma\}$, $\gamma^\zeta=\sup_{u\in [0,T]} \tau^\zeta_\varphi(X(u))$.

\end{theorem}
\begin{dov}
The theorem follows from Theorem \ref{omctf-square}. As $r_1$ we choose the function $r_1(t)=t^\alpha$, $0<\alpha<\ae$.
Under the conditions of the theorem, the second component of the function $Z_{r_1}(p,\beta,x)$ of Theorem \ref{omctf-square} can be transformed as follows:

$$r_1^{(-1)}\left( \frac{1}{\beta p}     \int_0^{\beta p} r_1(N_B(\sigma^{(-1)}(u)))du      \right) \leq $$

$$\leq r_1^{(-1)}\left( \frac{1}{\beta p}   \int_0^{\beta p} r_1\left(\frac{T}{2\sigma^{(-1)}(u)}+1\right)du        \right) , $$

\noindent because the metric massiveness  $N_B(u) \leq \frac{T}{2u}+1$ on the interval $[0,T]$.
Since $u\leq\beta p \leq \beta = \sigma(T/2)$, then $\sigma^{(-1)}(u)\leq T/2$ for such $u$.
Therefore, $T/2\sigma^{(-1)}\geq 1,$ therefore

$$r_1^{(-1)}\left( \frac{1}{\beta p}   \int_0^{\beta p} r_1\left(\frac{T}{2\sigma^{(-1)}(u)}+1\right)du        \right) \leq \left( \frac{1}{\beta p}   \int_0^{\beta p} \left(\frac{T}{\sigma^{(-1)}(u)}\right)^\alpha du \right)^{1/\alpha}.$$
Substituting $\sigma(u)$ into this expression, we obtain
$$\left( \frac{1}{\beta p}   \int_0^{\beta p} \left(\frac{T}{\sigma^{(-1)}(u)}\right)^\alpha du        \right)^{1/\alpha} = \left( \frac{1}{\beta p}   \int_0^{\beta p} \frac{T^\alpha C^{\alpha/\ae}}{u^{\alpha/\ae}}du      \right)^{1/\alpha} =$$
 $$= \left( \frac{1}{\beta p}   \frac{T^\alpha C^{\alpha/\ae}(\beta p)^{1-\frac{\alpha}{\ae}}}{1-\alpha/\ae}       \right)^{1/\alpha} \leq \frac{TC^{1/\ae}}{(\beta p)^{1/\ae}(1-\frac{\alpha}{\ae})^{1/\alpha}}.$$ For $\alpha\to 0$ we will have
$$\frac{TC^{1/\ae}}{(\beta p)^{1/\ae}(1-\frac{\alpha}{\ae})^{1/\alpha}}\to \frac{TC^{1/\ae}}{(\beta p)^{1/\ae}}e^{1/\ae},$$

\noindent and, since $\beta=\sigma(\inf_{t\in (0,T)}\sup_{t\in (0,T)}\rho(t,s)) = C (T/2)^{\ae}$, we obtain

$$\frac{TC^{1/\ae}}{(\beta p)^{1/\ae}}e^{1/\ae} = 2\left(\frac{e}{p}\right)^{1/\ae}.$$
Therefore,

$$Z_{r_1}(p,\beta,x) =\exp\left\{\frac{(1-\zeta)(x(1-p))^\frac{\zeta}{\zeta-1}}{\zeta(\gamma^\zeta(1-p)+p\beta^\zeta)^\frac{1}{\zeta-1}}\right\}2\left(\frac{e}{p}\right)^{1/\ae}.$$

The statement of the theorem follows from the last equality if we put $p=\gamma/x$.
In this case, the condition for $x$ of the previous theorem is transformed into a condition of the form

 $$x>\frac{\gamma(1-1/x)+(1/x)C^\zeta (T/2)^{\ae\zeta}}{(1-1/x)v^{\zeta-1}},$$
therefore
 $$x>\frac{\gamma(v^{\zeta-1}+1)+\sqrt{\gamma^2(v^{\zeta-1}+1)^2+4v^{\zeta-1}(C^\zeta(T/2)^{\ae\zeta}-\gamma^2)}}{2v^{\zeta-1}},$$
where $v=\min\{\beta,\gamma\}.$

\end{dov}

\begin{corollary}

Let a random process $X=\{X(t),t\in [0,T]\}$ belong to the space $Sub_\varphi(\Omega)$, $\varphi(t)=\frac{t^\zeta}{\zeta}$, $\zeta\geq2$, $t>1$.
Let $X$ be separable, and the condition (\ref{ome2.1}) holds for it, where $\sigma(h) = Ch^{\ae}$.
The model $X_N(t)$ approximates the process $X(t)$ with given reliability $1-\nu$ and accuracy $\delta$ in the space $C(0,T)$, if
$$\nu\leq 2\exp\left\{-\frac{(\zeta-1)\delta^\frac{1}{\zeta-1}(\delta-\gamma_N)^\frac{\zeta}{\zeta-1}}{\zeta(\gamma_N^\zeta(\delta-\gamma_N)+
\beta^\zeta\gamma_N)^\frac{1}{\zeta-1}}\right\}2\left(e\delta\right)^{1/\ae},$$
moreover
$$\delta>\frac{\gamma_N(v^{\zeta-1}+1)+\sqrt{\gamma_N^2(v^{\zeta-1}+1)^2+4v^{\zeta-1}(C^\zeta(T/2)^{\ae\zeta}-\gamma_N^2)}}{2v^{\zeta-1}},$$
where $v=\min\{\beta,\gamma_N\}$, $\gamma_N^\zeta=\sup_{u\in [0,T]} \tau^\zeta_\varphi(\Delta_N(u))$.
\end{corollary}

\begin{theorem}
Let a random process $X=\{X(t),t\in [0,T]\}$ belong to the space $Sub_\varphi(\Omega)$, $\varphi(t)=\frac{t^\zeta}{\zeta}$, $1<\zeta\leq2$.
Let $X$ be separable, and the condition (\ref{ome2.1}) holds for it, where $\sigma(h) = Ch^{\ae}$.
The model $X_N(t)$ approximates the process $X(t)$ with given reliability $1-\nu$ and accuracy $\delta$ in the space $C(0,T)$, if

$$\nu\leq 2\exp\left\{-\frac{(\zeta-1)\delta^\frac{1}{\zeta-1}(\delta-\gamma_N)^\frac{\zeta}{\zeta-1}}{\zeta(\gamma_N^\zeta(\delta-\gamma_N)+
\beta^\zeta\gamma_N)^\frac{1}{\zeta-1}}\right\}2\left(e\delta\right)^{1/\ae},$$
where $\gamma_N^\zeta=\sup_{u\in [0,T]} \tau^\zeta_\varphi(\Delta_N(u))$.

\end{theorem}

\begin{dov} This theorem follows from Theorem \ref{omctfz-square}, but in this case the restriction on $\delta$ is not required.
\end{dov}

\section[Models of random processes in $C(T)$ that admit series expansions]{Construction of models of random processes with $Sub_\varphi(\Omega)$ in $C(T)$
that admit series expansions with independent terms}

Let a random process $X=\{X(t),t\in B\}$ admit the representation
\begin{equation}
\label{omct-decomp-2}
X(t)=\sum_{k=1}^{\infty} a_k(t)\xi_k,
\end{equation}
where $\xi_k\in Sub_\varphi(\Omega)$.
Let $\delta_k(t) = |a_k(t)-\hat{a}_k(t)|$, $\hat{a}_k(t)$ is an approximation of $a_k(t)$, $\sigma_N$ are given in the definition (\ref{omsigma-N}).
We will consider $X_N$ \eqref{ommodel-basis}, given in Definition \ref{ommodel-def}, as a model of such a process.
In this case, we can prove the following theorems.

\begin{theorem}
\label{omcdelta-1}
Let the random process $X=\{X(t),t\in B\}$ belong to $Sub_\varphi(\Omega)$, let $(B,\rho)$ be a separable compact set, and let $X$ be separable on $(B,\rho)$.
Let $\varphi(t)=\frac{t^\zeta}{\zeta}$, $\zeta\geq2$, $t>1$, and for the process $X$ the conditions of Theorem \ref{omtk.1} are satisfied.
The model $X_N(t)$ \eqref{ommodel-basis} approximates the process $X(t)$ with given reliability $1-\nu$ and accuracy $\delta$ in the space $C(B)$, if
$$\nu\leq 2\exp\left\{-\frac{(\zeta-1)(\delta(1-p))^\frac{\zeta}{\zeta-1}}{\zeta(\gamma_N^\zeta(1-p)+p\beta^\zeta)^\frac{1}{\zeta-1}}\right\}r_1^{(-1)}\left(\frac{1}{\beta p}\int_0^{\beta p}r_1(N_B(\sigma_N^{(-1)}(u)))du\right),$$
$$\delta>\frac{\gamma_N^\zeta(1-p)+p\beta^\zeta}{(1-p)v^{\zeta-1}},$$
where $v=\min\{\beta,\gamma_N \}$, $\gamma_N^\zeta=\sup_{u\in B} \tau^\zeta_\varphi(\Delta_N(u))$,
 $$\gamma_N^\zeta\leq\left( \sum_{k=1}^N \tau_\varphi^2(\xi_k)\sup_{u\in B}\delta_k^2(u) + \sum_{k=N+1}^\infty \tau_\varphi^2(\xi_k)\sup_{u\in B}a^2_k(u) \right)^{\zeta/2}.$$
\end{theorem}

\begin{dov}
$$\gamma_N^\zeta=\sup_{u\in B}\tau_\varphi^\zeta(\Delta_N(u)) = \sup_{u\in B}\tau_\varphi^\zeta\left(\sum_{k=1}^N \xi_k\delta_k(u) + \sum_{k=N+1}^\infty \xi_ka_k(u)\right) \leq $$

$$ \leq \left( \sum_{k=1}^N \tau_\varphi^2(\xi_k)\sup_{u\in B}\delta_k^2(u) + \sum_{k=N+1}^\infty \tau_\varphi^2(\xi_k)\sup_{u\in B}a^2_k(u) \right)^{\zeta/2}.$$

\noindent The last inequality follows from the properties of the function $\tau_\varphi$ and Theorem~\ref{omtau-square}.
\end{dov}

\begin{theorem}
\label{omcdelta-2}
Let the random process $X=\{X(t),t\in B\}$ belong to $Sub_\varphi(\Omega)$, let $(B,\rho)$ be a separable compact set, and let $X$ be separable on $(B,\rho)$.
Let $\varphi(t)=\frac{t^\zeta}{\zeta}$, $1<\zeta\leq 2$, and for the process $X$ the conditions of Theorem \ref{omtk.1} are fulfilled.
The model $X_N(t)$ \eqref{ommodel-basis} approximates the process $X(t)$ with given reliability $1-\nu$ and accuracy $\delta$ in the space $C(B)$, if
$$\nu\leq 2\exp\left\{-\frac{(\zeta-1)(\delta(1-p))^\frac{\zeta}{\zeta-1}}{\zeta(\gamma_N^\zeta(1-p)+p\beta^\zeta)^\frac{1}{\zeta-1}}\right\}r_1^{(-1)}\left(\frac{1}{\beta p}\int_0^{\beta p}r_1(N_B(\sigma_N^{(-1)}(u)))du\right),$$
where $v=\min\{\beta,\gamma_N\}$,
$$\gamma_N^\zeta\leq\sum_{k=1}^N \tau_\varphi^\zeta(\xi_k)\sup_{u\in B}\delta_k^\zeta(u) + \sum_{k=N+1}^\infty \tau_\varphi^\zeta(\xi_k)\sup_{u\in B}a^\zeta_k(u) .$$

\end{theorem}

\begin{dov}
$$\gamma_N^\zeta=\sup_{u\in B}\tau_\varphi^\zeta(\Delta_N(u)) = \sup_{u\in B}\tau_\varphi^\zeta\left(\sum_{k=1}^N \xi_k\delta_k(u) + \sum_{k=N+1}^\infty \xi_ka_k(u)\right) \leq $$

$$ \leq \sum_{k=1}^N \tau_\varphi^\zeta(\xi_k)\sup_{u\in B}\delta_k^\zeta(u) + \sum_{k=N+1}^\infty \tau_\varphi^\zeta(\xi_k)\sup_{u\in B}a^\zeta_k(u).$$

\noindent The last inequality follows from the properties of the function $\tau_\varphi$ and Theorem~\ref{omtau-square}.
\end{dov}

Let the process $X(t)$ be given on $[0,T]$, and the following condition holds:

\bigskip
{\bf(C1)} $\sigma(h)=Ch^{\ae}$, $\delta_k(t)-\delta_k(s)\leq\hat{C}_kh^{\ae}$, $a_k(t)-a_k(s)\leq \tilde{C}_kh^{\ae}$.

\bigskip
\noindent Then we can prove the following statements.

\begin{theorem}
\label{omcdelta-3}
Let a random process $X=\{X(t),t\in [0,T]\}$ belong to the space $Sub_\varphi(\Omega)$, $\varphi(t)=\frac{t^\zeta}{\zeta}$, $\zeta\geq2$, $t>1$.
Let $X$ be separable, and the conditions (\ref{ome2.1}) and {\bf (C1)} are satisfied for it.
The model $X_N(t)$ \eqref{ommodel-basis} approximates the process $X(t)$ with given reliability $1-\nu$ and accuracy $\delta$ in the space $C(0,T)$, if
$$\nu\leq 2\exp\left\{-\frac{(\zeta-1)\delta^\frac{1}{\zeta-1}(\delta-\gamma_N)^\frac{\zeta}{\zeta-1}}{\zeta(\gamma_N^\zeta(\delta-\gamma_N)+\beta^\zeta\gamma_N)^\frac{1}{\zeta-1}}\right\}2\left(e\delta\right)^{1/\ae},$$ moreover
$$\delta>\frac{\gamma_N(v^{\zeta-1}+1)+\sqrt{\gamma_N^2(v^{\zeta-1}+1)^2+4v^{\zeta-1}(C^\zeta(T/2)^{\ae\zeta}-\gamma_N^2)}}{2v^{\zeta-1}},$$
where $v=\min\{\beta,\gamma_N\}$,
$$\gamma_N^\zeta\leq\left( \sum_{k=1}^N \tau_\varphi^2(\xi_k)\sup_{u\in[0,T]}\delta_k^2(u) + \sum_{k=N+1}^\infty \tau_\varphi^2(\xi_k)\sup_{u\in[0,T]}a^2_k(u) \right)^{\zeta/2}.$$
\end{theorem}

\begin{dov}
This statement follows from Theorem \ref{omctfz-square} and Theorem \ref{omcdelta-1}.
\end{dov}

\begin{theorem}
\label{omcdelta-4}
Let the random process $X=\{X(t),t\in [0,T]\}$ belong to the space $Sub_\varphi(\Omega)$, $\varphi(t)=\frac{t^\zeta}{\zeta}$, $1<\zeta\leq2$.
Let $X$ be separable, and the conditions (\ref{ome2.1}) and {\bf (C1)} are satisfied for it.
The model $X_N(t)$ \eqref{ommodel-basis} approximates the process $X(t)$ with given reliability $1-\nu$ and accuracy $\delta$ in the space $C(0,T)$, if
$$\nu\leq 2\exp\left\{-\frac{(\zeta-1)\delta^\frac{1}{\zeta-1}(\delta-\gamma_N)^\frac{\zeta}{\zeta-1}}{\zeta(\gamma^\zeta(\delta-\gamma)+\beta^\zeta\gamma)^\frac{1}{\zeta-1}}\right\}2\left(e\delta\right)^{1/\ae},$$
where
$$\gamma^\zeta\leq\sum_{k=1}^N \tau_\varphi^\zeta(\xi_k)\sup_{u\in[0,T]}\delta_k^\zeta(u) + \sum_{k=N+1}^\infty \tau_\varphi^\zeta(\xi_k)\sup_{u\in[0,T]}a^\zeta_k(u) .$$

\end{theorem}

\begin{dov}
This statement follows from the previous Theorem \ref{omcdelta-3} and Theorem \ref{omcdelta-2}.
\end{dov}

\subsection{Modeling of stochastic processes using the Karhunen-Lo\`eve decomposition in $C(B)$}

Let $X=\{X(t), t\in B\}$ be a centered second-order stochastic process with correlation function
$K(t,s)=EX(t)\overline{X(s)}$, let $f(t,\lambda)\in L_2(\Lambda,\mu)$, and let $K(t,s)$ admit the representation
$$K(t,s) = \int_\Lambda f(t,\lambda) f(s,\lambda)d\mu(\lambda).$$
Then the process $X$ admits a representation in the form of the Karhunen-Lo\`eve decomposition given in Theorem \ref{omkl-main}.
As a model of such a process, we will use $X_N$ process, given in Definition \ref{ommodel-def-kl}.

We will formulate theorems that will allow us to build models of stochastic processes with the Karhunen-Lo\`eve decomposition in $C(B)$.

\begin{theorem}
Let the random process $X=\{X(t),t\in B\}$ belong to $Sub_\varphi(\Omega)$, let $(B,\rho)$ be a separable compact set, and let $X$ be separable on
$(B,\rho)$. Let the process $X$ can be represented as a Karhunen-Loeve decomposition. Let $\varphi(t)=\frac{t^\zeta}{\zeta}$,
$\zeta\geq2$, $t>1$, and for the process $X$ the conditions of Theorem \ref{omtk.1} are satisfied.
The model $X_N(t)$ \eqref{ommodel-basis} approximates the process $X(t)$ with given reliability $1-\nu$ and accuracy $\delta$ in the space $C(B)$, if
$$\nu\leq 2\exp\left\{-\frac{(\zeta-1)(\delta(1-p))^\frac{\zeta}{\zeta-1}}{\zeta(\gamma_N^\zeta(1-p)+p\beta^\zeta)^\frac{1}{\zeta-1}}\right\}r_1^{(-1)}\left(\frac{1}{\beta p}\int_0^{\beta p}r_1(N_B(\sigma_N^{(-1)}(u)))du\right),$$
$$\delta>\frac{\gamma_N^\zeta(1-p)+p\beta^\zeta}{(1-p)v^{\zeta-1}},$$ where $v=\min\{\beta,\gamma_N\}$, $$\gamma_N^\zeta\leq \left(\sum_{k=1}^N \tau^2_\varphi(\xi_k)\left(\frac{\sup_{u\in B}\delta^2_k(u) }{\hat{\lambda}_k-\eta_k} +\sup_{u\in B}\hat{a}^2_k(u)\frac{(\sqrt{\hat{\lambda}_k}-\sqrt{\hat{\lambda}_k-\eta_k})^2}{\hat{\lambda}_k(\hat{\lambda}_k-\eta_k)}\right)+\right.$$
$$+\left.\sum_{k=N+1}^\infty\frac{\tau^2_\varphi(\xi_k)\sup_{u\in B}a^2_k(u)}{\lambda_k} \right)^{\zeta/2},$$

\noindent where $\delta_k(t)$ is the error in finding the $k$-th eigenfunction, $\hat{\lambda}_k$ is the approximate value of the $k$-th eigenvalue, $\eta_k$ is the error in approximating the $k$-th eigenvalue.

\end{theorem}

\begin{dov}
The theorem follows from Theorem \ref{omkl-lp-square} and Theorem \ref{omcdelta-1}.
\end{dov}

In the case where $\forall k:$ $\tau_\varphi(\xi_k)=\tau$, we obtain the following statement.

\begin{theorem}
\label{omct-merser}
Let for a random process $X=\{X(t),t\in B\}$ with $Sub_\varphi(\Omega)$, $\varphi(t)=\frac{t^\zeta}{\zeta}$, $\zeta\geq2$, $t>1$ the conditions of theorem \ref{omtk.1} be satisfied.
The model $X_N(t)$ \eqref{ommodel-basis} approximates the process $X(t)$ with given reliability $1-\nu$ and accuracy $\delta$ in the space $C(B)$, when

$$\nu\leq 2\exp\left\{-\frac{(\zeta-1)(\delta(1-p))^\frac{\zeta}{\zeta-1}}{\zeta(\gamma_N^\zeta(1-p)+p\beta^\zeta)^\frac{1}{\zeta-1}}\right\}r_1^{(-1)}\left(\frac{1}{\beta p}\int_0^{\beta p}r_1(N_B(\sigma_N^{(-1)}(u)))du\right),$$

$$\delta>\frac{\gamma_N^\zeta(1-p)+p\beta^\zeta}{(1-p)v^{\zeta-1}},$$ where $v=\min\{\beta,\gamma_N\}$, $$\gamma_N^\zeta\leq \left(\sup_{u\in B}K(u,u)+ \sum_{k=1}^N \left(\frac{\sup_{u\in B}\delta^2_k(u) }{\hat{\lambda}_k-\eta_k} +\right.\right.$$

$$\left.\left.+\sup_{u\in B}\hat{a}^2_k(u)\frac{(\sqrt{\hat{\lambda}_k}-\sqrt{\hat{\lambda}_k-\eta_k})^2}{\hat{\lambda}_k(\hat{\lambda}_k-\eta_k)}-\frac{\sup_{u\in B}(\hat{a}_k(u)-\delta_k(u))^2}{\hat{\lambda}_k+\eta_k}\right)\right)^{\zeta/2} ,$$

\noindent where $\delta_k(t)$ is the error in finding the $k$-th eigenfunction, $\lambda_k$ is the approximate value of the $k$-th eigenvalue, $\eta_k$ is the error in approximating the $k$-th eigenvalue.
\end{theorem}

\begin{theorem}

Let for a random process $X=\{X(t),t\in B\}$ with $Sub_\varphi(\Omega)$, $\varphi(t)=\frac{t^\zeta}{\zeta}$, $1<\zeta\leq2$ the conditions of theorem \ref{omtk.1} be satisfied.
The model $X_N(t)$ \eqref{ommodel-basis} approximates the process $X(t)$ with given reliability $1-\nu$ and accuracy $\delta$ in the space $C(B)$, when

$$\nu\leq 2\exp\left\{-\frac{(\zeta-1)(\delta(1-p))^\frac{\zeta}{\zeta-1}}{\zeta(\gamma_N^\zeta(1-p)+p\beta^\zeta)^\frac{1}{\zeta-1}}\right\}r_1^{(-1)}\left(\frac{1}{\beta p}\int_0^{\beta p}r_1(N_B(\sigma_N^{(-1)}(u)))du\right),$$

\noindent where

$$\gamma_N^\zeta\leq \sum_{k=1}^N \tau^{\gamma}_\varphi(\xi_k)\left(\frac{\sup_{u\in B}\delta^{\gamma}_k(u) }{(\lambda_k-\eta_k)^\gamma} +\sup_{u\in B}\hat{a}^{\gamma}_k(u)\frac{(\sqrt{\hat{\lambda}_k}-\sqrt{\hat{\lambda}_k-\eta_k})^{\gamma}}{(\hat{\lambda}_k(\hat{\lambda}_k-\eta_k))^\gamma}\right)+$$

$$+\sum_{k=N+1}^\infty\frac{\tau^{\gamma}_\varphi(\xi_k)\sup_{u\in B}a^{2\gamma}_k(u)}{\lambda_k^\gamma}$$

\noindent $\delta_k(t)$ is the error in finding the $k$-th eigenfunction, $\lambda_k$ is the approximate value of the $k$-th eigenvalue, $\eta_k$ is the error in approximating the $k$-th eigenvalue.
\end{theorem}

\begin{dov} The statement follows from Theorem \ref{omlp-merser} and Theorem \ref{omcdelta-2}.
\end{dov}

Now let $B=[0,T]$. Then we can prove the following theorems.

\begin{theorem}
Let the random process $X=\{X(t),t\in [0,T]\}$ belong to the space $Sub_\varphi(\Omega)$, $\varphi(t)=\frac{t^\zeta}{\zeta}$, $\zeta\geq2$, $t>1$.
Let $X$ be represented in the form of a Karhunen-Loeve decomposition (\ref{omkl-decompos-x}), and the conditions (\ref{ome2.1}) and {\bf (C1)} are satisfied for it.
The model $X_N(t)$ \eqref{ommodel-basis} approximates the process $X(t)$ with given reliability $1-\nu$ and accuracy $\delta$ in the space $C(0,T)$, when

$$\nu\leq 2\exp\left\{-\frac{(\zeta-1)\delta^\frac{1}{\zeta-1}(\delta-\gamma_N)^\frac{\zeta}{\zeta-1}}{\zeta(\gamma_N^\zeta(\delta-\gamma_N)+\beta^\zeta\gamma_N)^\frac{1}{\zeta-1}}\right\}2\left(e\delta\right)^{1/\ae},$$ moreover

$$\delta>\frac{\gamma_N(v^{\zeta-1}+1)+\sqrt{\gamma_N^2(v^{\zeta-1}+1)^2+4v^{\zeta-1}(C^\zeta(T/2)^{\ae\zeta}-\gamma_N^2)}}{2v^{\zeta-1}},$$
where $v=\min\{\beta,\gamma_N\}$,

$$\gamma_N^\zeta\leq \left(\sum_{k=1}^N \tau^2_\varphi(\xi_k)\left(\frac{\sup_{u\in[0,T]}\delta^2_k(u) }{\hat{\lambda}_k-\eta_k} +\sup_{u\in[0,T]}\hat{a}^2_k(u)\frac{(\sqrt{\hat{\lambda}_k}-\sqrt{\hat{\lambda}_k-\eta_k})^2}{\hat{\lambda}_k(\hat{\lambda}_k-\eta_k)}\right)+\right.$$

$$\left.+\sum_{k=N+1}^\infty\frac{\tau^2_\varphi(\xi_k)\sup_{u\in[0,T]}a^2_k(u)}{\lambda_k} \right)^{\zeta/2},$$

$$C = \sqrt{\sum_{k=1}^N \tau^2_\varphi(\xi_k)\left(\frac{\hat{\lambda}_k\hat{C}_k+(\sqrt{\hat{\lambda}_k}-\sqrt{\hat{\lambda}_k-\eta_k})\tilde{C}_k}{\hat{\lambda}_k(\hat{\lambda}_k-\eta_k)}\right)^2+\sum_{k=N+1}^\infty \tau^2_\varphi(\xi_k)\frac{\tilde{C}^2_k}{\lambda^2_k}},$$

\noindent $\delta_k(t)$ is the error in finding the $k$-th eigenfunction, $\lambda_k$ is the approximate value of the $k$-th eigenvalue, $\eta_k$ is the error in approximating the $k$-th eigenvalue.
\end{theorem}

\begin{dov}
Consider $\sigma(h)$ for this case. Indeed,

$$\sup_{|t-s|<h}\tau_\varphi(\Delta_N(t)-\Delta_N(s)) = \sup_{|t-s|<h}\tau_\varphi\left( \sum_{k=1}^N \xi_k \left(\frac{a_k(t)-a_k(s)}{\lambda_k}-\right.\right.$$

$$-\left.\left.\frac{\hat{a}_k(t)-\hat{a}_k(s)}{\hat{\lambda}_k}\right) + \sum_{k=N+1}^\infty\xi_k\left(\frac{ a_k(t)-a_k(s)}{\lambda_k}\right) \right)\leq $$

$$\leq h^{\ae}\left(\sum_{k=1}^N \tau^2_\varphi(\xi_k)\frac{\hat{\lambda}_k\hat{C}_k+(\sqrt{\hat{\lambda}_k}-\sqrt{\hat{\lambda}_k-\eta_k})\tilde{C}_k}{\hat{\lambda}_k(\hat{\lambda}_k-\eta_k)}+\right.$$

$$\left.+\sum_{k=N+1}^\infty \tau^2_\varphi(\xi_k)\frac{\tilde{C}_k}{\lambda_k}\right)^{1/2} = Ch^{\ae} = \sigma(h).$$

\noindent The last inequality uses the proof of the theorem
\ref{omkl-lp-square}. Now the statement of the theorem follows from Theorem \ref{omkl-lp-square} and Theorem \ref{omcdelta-3}.
\end{dov}

\begin{theorem}
\label{omct-merser-2}
Let for a random process $X=\{X(t),t\in (0,T)\}$ from $Sub_\varphi(\Omega)$, $\varphi(t)=\frac{t^\zeta}{\zeta}$, $\zeta\geq2$, $t>1$,
the conditions of Theorem \ref{omtk.1} are satisfied, and this process can be represented in the form (\ref{omkl-decompos-x}).
Let for $\Delta_N$ of this process the condition {\bf (C1)} is satisfied, and also $\tau_\varphi(\xi_k)=\tau$.
The model $X_N(t)$ \eqref{ommodel-basis} approximates the process $X(t)$ with given reliability $1-\nu$ and accuracy $\delta$ in the space $C(0,T)$ when

$$\nu\leq 2\exp\left\{-\frac{(\zeta-1)\delta^\frac{1}{\zeta-1}(\delta-\gamma_N)^\frac{\zeta}{\zeta-1}}{\zeta(\gamma_N^\zeta(\delta-\gamma_N)+
\beta^\zeta\gamma_N)^\frac{1}{\zeta-1}}\right\}2\left(e\delta\right)^{1/\ae},$$
moreover $$\delta>\frac{\gamma_N(v^{\zeta-1}+1)+\sqrt{\gamma_N^2(v^{\zeta-1}+1)^2+4v^{\zeta-1}(C^\zeta(T/2)^{\ae\zeta}-\gamma_N^2)}}{2v^{\zeta-1}},$$
where $v=\min\{\beta,\gamma_N\}$,

$$\gamma_N^\zeta\leq \tau^\zeta\left(\sup_{u\in[0,T]}K(u,u)+ \sum_{k=1}^N \left(\frac{\sup_{u\in[0,T]}\delta^2_k(u) }{\hat{\lambda}_k-\eta_k} +\right.\right.$$

$$\left.\left.+\sup_{u\in[0,T]}\hat{a}^2_k(u)\frac{(\sqrt{\hat{\lambda}_k}-\sqrt{\hat{\lambda}_k-\eta_k})^2}{\hat{\lambda}_k(\hat{\lambda}_k-\eta_k)}-\frac{\sup_{u\in[0,T]}(\hat{a}_k(u)-\delta_k(u))^2}{\hat{\lambda}_k+\eta_k}\right)\right)^{\zeta/2} ,$$

$$C = \tau\left(\sum_{k=1}^N \left(\frac{\hat{\lambda}_k\hat{C}_k+(\sqrt{\hat{\lambda}_k}-\sqrt{\hat{\lambda}_k-\eta_k})\tilde{C}_k}{\hat{\lambda}_k(\hat{\lambda}_k-\eta_k)}\right)^2+C_B-\sum_{k=1}^N \frac{\tilde{C_k}^2}{\hat{\lambda}_k+\eta_k}\right)^{1/2},$$

$$\sup_{|t-s|<h}|K(t,u)-K(s,u)|\leq h^{\ae}\tilde{K}(u),$$

$$\sup_{|t-s|<h}|\tilde{K}(t)-\tilde{K}(s)|\leq h^{\ae} C_K,$$

\noindent $K(t,s)$ is the correlation function of the process, $\delta_k(t)$ is the error in finding the $k$-th eigenfunction, $\lambda_k$ is the approximate value of the $k$-th eigenvalue, $\eta_k$ is the error in approximating the $k$-th eigenvalue.
\end{theorem}

\begin{dov}
Consider $\sigma(h)$ for this case. Indeed,
$$\sup_{|t-s|<h}\tau_\varphi(\Delta_N(t)-\Delta_N(s)) = \sup_{|t-s|<h}\tau_\varphi\left( \sum_{k=1}^N \xi_k \left(\frac{a_k(t)-a_k(s)}{\lambda_k}-\right.\right.$$
$$\left.\left.-\frac{\hat{a}_k(t)-\hat{a}_k(s)}{\hat{\lambda}_k}\right) + \sum_{k=N+1}^\infty\xi_k\left(\frac{ a_k(t)-a_k(s)}{\lambda_k}\right) \right)\leq $$
$$\leq h^{\ae}\tau\left(\sum_{k=1}^N \left(\frac{\hat{\lambda}_k\hat{C}_k+(\sqrt{\hat{\lambda}_k}-\sqrt{\hat{\lambda}_k-\eta_k})\tilde{C}_k}{\hat{\lambda}_k(\hat{\lambda}_k-\eta_k)}\right)^2+\right.$$
$$\left.+\sum_{k=N+1}^\infty \frac{(a_k(t)-a_k(s))^2}{\lambda_k}\right)^{1/2} = h^{\ae}\tau\left(\sum_{k=1}^N \left(\frac{\hat{C}_k}{(\hat{\lambda}_k-\eta_k)}+\right.\right.$$
$$\left.\left.+\frac{(\sqrt{\hat{\lambda}_k}-\sqrt{\hat{\lambda}_k-\eta_k})\tilde{C}_k}{\hat{\lambda}_k(\hat{\lambda}_k-\eta_k)}\right)^2+\sum_{k=N+1}^\infty \frac{1}{\lambda_k}(a_k^2(t)-2a_k(t)a_k(s)+a^2_k(s))\right)^{1/2} = $$
$$=h^{\ae}\tau\left(\sum_{k=1}^N \left(\frac{\hat{\lambda}_k\hat{C}_k+(\sqrt{\hat{\lambda}_k}-\sqrt{\hat{\lambda}_k-\eta_k})\tilde{C}_k}{\hat{\lambda}_k(\hat{\lambda}_k-\eta_k)}\right)^2+\right.$$
$$+\left.K(t,t)+K(s,s)-2K(t,s)-\sum_{k=1}^N \frac{1}{\lambda_k}(a_k^2(t)-2a_k(t)a_k(s)+a^2_k(s))\right)^{1/2} \leq $$
$$\leq h^{\ae}\tau\left(\sum_{k=1}^N \left(\frac{\hat{\lambda}_k\hat{C}_k+(\sqrt{\hat{\lambda}_k}-\sqrt{\hat{\lambda}_k-\eta_k})\tilde{C}_k}{\hat{\lambda}_k(\hat{\lambda}_k-\eta_k)}\right)^2+C_K-\sum_{k=1}^N \frac{\tilde{C_k}^2}{\hat{\lambda}_k+\eta_k}\right)^{1/2}=$$
$$=Ch^{\ae} = \sigma(h).$$

\noindent Now the statement of the Theorem follows from Theorem \ref{omkl-lp-square} and Theorem \ref{omcdelta-3}.
\end{dov}

The following theorem is a modification of the theorem from \cite{ostrovsky1973}.

\begin{theorem}
\label{omphik-supremum}
Let $\varphi_k(t)$ be the eigenfunctions of the homogeneous Fredholm integral equation of the second kind

$$\phi(t)=\lambda\int_0^TK(t,s)\phi(s)ds,$$

\noindent $\lambda_k$ are the eigenvalues of this equation, numbered in ascending order. Then the inequality holds

$$\sup_{|t-s|<h}|\phi_k(t)-\phi_k(s)|\leq \lambda_k \omega_K(h),$$

$$\omega_K(h) = \sup_{|t-s|<h}\left(\int_0^T (K(t,u)-K(s,u))^2du\right)^{1/2}$$
\end{theorem}

\begin{dov}
Let $|t-s|<h$. Then the relations hold

$$|\phi_k(s)-\phi_k(t)|=\left|\int_0^T \lambda_k\phi_k(u)(K(t,u)-K(s,u))du\right|\leq$$

$$\leq \left(\int_0^T \lambda^2_k\phi^2_k(u)du\right)^{1/2}\times\left(\int_0^T (K(t,u)-K(s,u))^2du\right)^{1/2}\leq \lambda_k\omega_K(h).$$
\end{dov}

Based on this statement, we can slightly change the conditions of Theorem \ref{omct-merser-2}.

\begin{theorem}
\label{omct-merser-3}
Let for a random process $X=\{X(t),t\in (0,T)\}$ from $Sub_\varphi(\Omega)$, $\varphi(t)=\frac{t^\zeta}{\zeta}$, $\zeta\geq2$, $t>1$,
the conditions of Theorem \ref{omtk.1} are satisfied, and this process can be represented in the form (\ref{omkl-decompos-x}).
Let for $\Delta_N$ of this process the condition {\bf (C1)} is satisfied, and also $\tau_\varphi(\xi_k)=\tau$.
The model $X_N(t)$ \eqref{ommodel-basis} approximates the process $X(t)$ with given reliability $1-\nu$ and accuracy $\delta$ in the space $C(0,T)$ if

$$\nu\leq 2\exp\left\{-\frac{(\zeta-1)\delta^\frac{1}{\zeta-1}(\delta-\gamma_N)^\frac{\zeta}{\zeta-1}}{\zeta(\gamma_N^\zeta(\delta-\gamma_N)+\beta^\zeta\gamma_N)^\frac{1}{\zeta-1}}\right\}2\left(e\delta\right)^{1/\ae},$$ moreover $$\delta>\frac{\gamma_N(v^{\zeta-1}+1)+\sqrt{\gamma_N^2(v^{\zeta-1}+1)^2+4v^{\zeta-1}(C^\zeta(T/2)^{\ae\zeta}-\gamma_N^2)}}{2v^{\zeta-1}},$$
where $v=\min\{\beta,\gamma_N\}$,

$$\gamma_N^\zeta\leq \tau^\zeta\left(\sup_{u\in[0,T]}K(u,u)+ \sum_{k=1}^N \left(\frac{\sup_{u\in[0,T]}\delta^2_k(u) }{\hat{\lambda}_k-\eta_k} +\right.\right.$$

$$\left.\left.+\sup_{u\in[0,T]}\hat{a}^2_k(u)\frac{(\sqrt{\hat{\lambda}_k}-\sqrt{\hat{\lambda}_k-\eta_k})^2}{\hat{\lambda}_k(\hat{\lambda}_k-\eta_k)}-\frac{\sup_{u\in[0,T]}(\hat{a}_k(u)-\delta_k(u))^2}{\hat{\lambda}_k+\eta_k}\right)\right)^{\zeta/2} ,$$

$$C = \tau\left(\sum_{k=1}^N \left(\frac{\hat{C}_k}{\hat{\lambda}_k-\eta_k}+\frac{(\sqrt{\hat{\lambda}_k}-\sqrt{\hat{\lambda}_k-\eta_k})}{\hat{\lambda}_k}\tilde{C}_K\right)^2+C_K-\right.$$

$$\left.-\tilde{C}_K\sum_{k=1}^N (\hat{\lambda}_k-\eta_k)\right)^{1/2},$$

$$\sup_{|t-s|<h}|K(t,u)-K(s,u)|\leq h^{\ae}\tilde{K}(u),$$

$$\sup_{|t-s|<h}|\tilde{K}(t)-\tilde{K}(s)|\leq h^{\ae} C_K,$$

$$\left(\int_0^T \tilde{K}^2(u)du\right)^{1/2} =  \tilde{C}_K,$$

\noindent where $K(t,s)$ is the correlation function of the process, $\delta_k(t)$ is the error in finding the $k$-th eigenfunction, $\lambda_k$ is the approximate value of the
$k$-th eigenvalue, $\eta_k$ is the error in approximating the $k$-th eigenvalue.
\end{theorem}

\begin{dov}

Consider $\sigma(h)$ for this case. From Theorem \ref{omct-merser-2} we have that

$$C = \tau\left(\sum_{k=1}^N \left(\frac{\hat{\lambda}_k\hat{C}_k+(\sqrt{\hat{\lambda}_k}-\sqrt{\hat{\lambda}_k-\eta_k})\tilde{C}_k}{\hat{\lambda}_k(\hat{\lambda}_k-\eta_k)}\right)^2+C_K-\sum_{k=1}^N \frac{\tilde{C_k}^2}{\hat{\lambda}_k+\eta_k}\right)^{1/2}.$$

\noindent Applying Theorem \ref{omphik-supremum}, we obtain

$$\sup_{|t-s|<h}|a_k(t)-a_k(s)| \leq \lambda_k\left(\int_0^T(K(t,u)-K(s,u))^2du\right)^{1/2}\leq$$

$$\leq\lambda_k h^{\ae}\left(\int_0^T(K(u))^2du\right)^{1/2},$$

\noindent i.e. $\tilde{C}_k = \lambda_k\tilde{C}_B$. Substituting this value into the conditions of the theorem \ref{omct-merser-2}, we obtain

$$C = \tau\left(\sum_{k=1}^N \left(\frac{\hat{\lambda}_k\hat{C}_k+(\sqrt{\hat{\lambda}_k}-\sqrt{\hat{\lambda}_k-\eta_k})(\hat{\lambda}_k+\eta_k)\tilde{C}_K}{\hat{\lambda}_k(\hat{\lambda}_k-\eta_k)}\right)^2+C_K-\right.$$

$$\left.-\tilde{C}_K\sum_{k=1}^N (\hat{\lambda}_k-\eta_k)\right)^{1/2}.$$

\end{dov}

\begin{theorem}
Let for a random process $X=\{X(t),t\in (0,T)\}$ from $Sub_\varphi(\Omega)$, $\varphi(t)=\frac{t^\zeta}{\zeta}$, $1<\zeta\leq2$,
the conditions of Theorem \ref{omtk.1} are satisfied, and this process can be represented in the form (\ref{omkl-decompos-x}).
Let for $\Delta_N$ of this process the condition {\bf (C1)} is satisfied.
The model $X_N(t)$ \eqref{ommodel-basis} approximates the process $X(t)$ with given reliability $1-\nu$ and accuracy $\delta$ in the space $C(0,T)$, if

$$\nu\leq 2\exp\left\{-\frac{(\zeta-1)\delta^\frac{1}{\zeta-1}(\delta-\gamma)^\frac{\zeta}{\zeta-1}}{\zeta(\gamma^\zeta(\delta-\gamma)+\beta^\zeta\gamma)^\frac{1}{\zeta-1}}\right\}2\left(e\delta\right)^{1/\ae},$$
where

$$\gamma^\zeta\leq \sum_{k=1}^N \tau^{\gamma}_\varphi(\xi_k)\left(\frac{\sup_{u\in[0,T]}\delta^{\gamma}_k(u) }{(\lambda_k-\eta_k)^\gamma} +\sup_{u\in[0,T]}\hat{a}^{\gamma}_k(u)\frac{(\sqrt{\hat{\lambda}_k}-\sqrt{\hat{\lambda}_k-\eta_k})^{\gamma}}{(\hat{\lambda}_k(\hat{\lambda}_k-\eta_k))^\gamma}\right)+$$

$$+\sum_{k=N+1}^\infty\frac{\tau^{\gamma}_\varphi(\xi_k)\sup_{u\in[0,T]}a^{\gamma}_k(u)}{\lambda_k^\gamma},$$

$$C = \sqrt{\sum_{k=1}^N \tau^2_\varphi(\xi_k)\left(\frac{\hat{\lambda}_k\hat{C}_k+(\sqrt{\hat{\lambda}_k}-\sqrt{\hat{\lambda}_k-\eta_k})\tilde{C}_k}{\hat{\lambda}_k(\hat{\lambda}_k-\eta_k)}\right)^2+\sum_{k=N+1}^\infty \tau^2_\varphi(\xi_k)\frac{\tilde{C}^2_k}{\lambda^2_k}}.$$

\end{theorem}

\begin{dov}
This statement follows from the previous Theorem and Theorem \ref{omcdelta-4}.
\end{dov}

\subsection{Modeling of stochastic processes using their decomposition into certain bases in $C(0,T)$}

Let $X=\{X(t),t \in [0,T]\} \in Sub_\varphi(\Omega)$ be a second-order stochastic process, $EX(t)=0$.
Let the covariance function of the process $B(t,s)=EX(t)\overline{X(s)}$ admit the representation
$$B(t,s)=\int_{-\infty}^\infty f(t,\lambda)f(s,\lambda)d\lambda,$$
where $f(t,\lambda)$, $t\in [0,T]$, $\lambda\in R$ are a family of functions from $L_2(R)$.
According to the theorem \ref{ombasis-dec-main}, the process $X$ can be represented as

$$X(t) = \sum_{k=0}^{\infty}\xi_k a_k(t),$$

$$a_k(t)=\int_{0}^{T} f(t,\lambda)\cos \pi k\lambda d\lambda,$$

\noindent where $\xi_k$ are $\varphi$-sub-Gaussian independent random variables such that $E\xi_k^2=1$.
As a model of such a process, we use $X_N$ \eqref{ommodel-basis},
given in Definition \ref{ommodel-basis-ozn}.

\begin{theorem}
\label{omcosdecC-1}
Let for a random process $X=\{X(t),t\in (0,T)\}$ from $Sub_\varphi(\Omega)$, $\varphi(t)=\frac{t^\zeta}{\zeta}$, $\zeta\geq2$, $t>1$,
the conditions of Theorem \ref{omtk.1} are satisfied, and this process can be represented in the form (\ref{ombasis}), $g_k(s)=\cos \pi k s$, $f(t,s)$
is differentiable with respect to $s$.
Let for $\Delta_N$ of this process the condition {\bf (C1)} is satisfied.
The model $X_N(t)$ \eqref{ommodel-basis} approximates the process $X(t)$ with given reliability $1-\nu$ and accuracy $\delta$ in the space $C(0,T)$, if

$$\nu\leq 2\exp\left\{-\frac{(\zeta-1)\delta^\frac{1}{\zeta-1}(\delta-\gamma)^\frac{\zeta}{\zeta-1}}{\zeta(\gamma^\zeta(\delta-\gamma)+\beta^\zeta\gamma)^\frac{1}{\zeta-1}}\right\}2\left(e\delta\right)^{1/\ae},$$ moreover ,$$\delta>\frac{\gamma(v^{\zeta-1}+1)+\sqrt{\gamma^2(v^{\zeta-1}+1)^2+4v^{\zeta-1}(C^\zeta(T/2)^{\ae\zeta}-\gamma^2)}}{2v^{\zeta-1}},$$
where $v=\min\{\beta,\gamma\}$,

$$\gamma^\zeta\leq\left(\delta^2_f(t)\sum_{k=N+1}^\infty \frac{4\tau^2_\varphi(\xi_k)}{\pi^2 k^2} + \sum_{k=1}^N \tau^2_\varphi(\xi_k)\delta^2_k(t)\right)^{\zeta/2},$$
$$\delta_f(t) = f(t,T)-f(t,0),$$

$$C = \sqrt{\sum_{k=1}^N \tau_\varphi^2(\xi_k)\hat{C}_k^2 +
C_f^2\sum_{k=N+1}^\infty
\frac{4\tau_\varphi^2(\xi_k)}{\pi^2k^2}},$$

$$\sup_{|t-s|<h}|f(s,u)-f(t,u)|\leq h^{\ae}\tilde{f}(u),$$

$$C_f=\tilde{f}(T)-\tilde{f}(0).$$

\end{theorem}

\begin{dov}

Consider $\sigma(h)$ for this case. Indeed,

$$\sup_{|t-s|<h}\tau_\varphi(\Delta_N(t)-\Delta_N(s)) = \sup_{|t-s|<h}\tau_\varphi\left( \sum_{k=1}^N \xi_k (\delta_k(t)-\delta_k(s))+\right.$$

$$\left.+\sum_{K=N+1}^\infty \xi_k (a_k(t)-a_k(s))\right)\leq $$

$$\leq h^{\ae}\left(\sum_{k=1}^N \tau_\varphi^2(\xi_k)\hat{C}_k^2 +
C_f\sum_{k=N+1}^\infty
\frac{4\tau_\varphi^2(\xi_k)}{\pi^2k^2}\right)^{1/2} = Ch^{\ae} =\sigma(h).$$

The last inequality uses the proof of Theorem \ref{omcosdec}.
Now the necessary statement follows from Theorem \ref{omcosdec} and Theorem \eqref{ommodel-basis}.
\ref{omcdelta-3}.
\end{dov}

\begin{theorem}
\label{omcosdecC-2}
Let for a random process $X=\{X(t),t\in (0,T)\}$ from $Sub_\varphi(\Omega)$, $\varphi(t)=\frac{t^\zeta}{\zeta}$, $1<\zeta\leq2$,
the conditions of Theorem \ref{omtk.1} be satisfied, and this process can be represented in the form (\ref{ombasis}), $g_k(s)=\cos \pi k s$, $f(t,s)$ is differentiable with respect to $s$.
Let for $\Delta_N$ of this process the condition {\bf (C1)} be satisfied. The model $X_N(t)$ \eqref{ommodel-basis} approximates the process $X(t)$ with given reliability $1-\nu$
and accuracy $\delta$ in the space $C(0,T)$ if

$$\nu\leq 2\exp\left\{-\frac{(\zeta-1)\delta^\frac{1}{\zeta-1}(\delta-\gamma)^\frac{\zeta}{\zeta-1}}{\zeta(\gamma^\zeta(\delta-\gamma)+\beta^\zeta\gamma)^\frac{1}{\zeta-1}}\right\}2\left(e\delta\right)^{1/\ae},$$
where

$$\gamma^\zeta\leq\left(\delta^2_f(t)\sum_{k=N+1}^\infty \frac{4\tau^2_\varphi(\xi_k)}{\pi^2 k^2} + \sum_{k=1}^N \tau^2_\varphi(\xi_k)\delta^2_k(t)\right)^{\zeta/2},$$

$$\delta_f(t) = f(t,T)-f(t,0).$$

\end{theorem}

\begin{dov}
This statement follows from the previous theorem and the theorem
\ref{omcdelta-4}
\end{dov}

\subsection{Modeling of stochastic processes using their Hermite basis decomposition in $C(0,T)$}

Let $X=\{X(t),t \in [0,T]\} \in Sub_\varphi(\Omega)$ be a second-order stochastic process, $EX(t)=0$.
Let the covariance function of the process $B(t,s)=EX(t)\overline{X(s)}$ admit the representation

$$B(t,s)=\int_{-\infty}^\infty f(t,\lambda)f(s,\lambda)d\lambda,$$
where $f(t,\lambda)$, $t\in [0,T]$,
$\lambda\in R$ is a family of functions from $L_2(R)$. According to Theorem \ref{ombasis-dec-main}, the process $X$ can be represented as

$$X(t) = \sum_{k=0}^{\infty}\xi_k a_k(t),$$

$$a_k(t)=\int_{-\infty}^{\infty} f(t,\lambda)g_k(\lambda)d\lambda,$$

\noindent where $\xi_k$ are $\varphi$-sub-Gaussian independent random variables such that $E\xi_k^2=1$, $g_k(\lambda)$ are Hermite functions,
given in (\ref{omhermite_func}). As a model of such a process, we use $X_N$ \eqref{ommodel-basis}, given in the definition of \ref{ommodel-basis-ozn}.

\begin{theorem}
\label{omct-hermite-1}
Let for a random process $X=\{X(t),t\in (0,T)\}$ from $Sub_\varphi(\Omega)$, $\varphi(t)=\frac{t^\zeta}{\zeta}$, $\zeta\geq2$, $t>1$,
the conditions of Theorem \ref{omtk.1} be satisfied, and let this process can be represented in the form (\ref{ombasis}).
Let for $\Delta_N$ of this process the condition {\bf (C1)} be satisfied. The model $X_N(t)$ \eqref{ommodel-basis} approximates the process $X(t)$ with given reliability $1-\nu$ and accuracy $\delta$ in the space $C(0,T)$, if

$$\nu\leq 2\exp\left\{-\frac{(\zeta-1)\delta^\frac{1}{\zeta-1}(\delta-\gamma)^\frac{\zeta}{\zeta-1}}{\zeta(\gamma^\zeta(\delta-\gamma)+\beta^\zeta\gamma)^\frac{1}{\zeta-1}}\right\}2\left(e\delta\right)^{1/\ae},$$ moreover
$$\delta>\frac{\gamma(v^{\zeta-1}+1)+\sqrt{\gamma^2(v^{\zeta-1}+1)^2+4v^{\zeta-1}(C^\zeta(T/2)^{\ae\zeta}-\gamma^2)}}{2v^{\zeta-1}},$$
where $v=\max\{\beta,\gamma\}$,

$$\gamma^\zeta\leq\left(K^2 \sup_{u\in[0,T]}\int_{-\infty}^\infty Z^2_f(u,\lambda)d\lambda \sum_{k=N+1}^\infty \frac{\tau^2_\varphi(\xi_k)}{k^2+3k+2} + \right.$$

$$\left.+\sum_{k=1}^N \tau^2_\varphi(\xi_k)\sup_{u\in[0,T]}\delta^2_k(u)\right)^{\zeta/2},$$

$$Z_f(t,\lambda) = \frac{\partial^2 f(t,\lambda)}{\partial \lambda^2} - \lambda \frac{\partial f(t,\lambda)}{\partial \lambda}+\frac{\lambda^2-2}{4}f(t,\lambda),$$

$$C = \sqrt{K^2\tilde{C}^2\sum_{k=N+1}^\infty \frac{\tau_\varphi^2(\xi_k)}{k^2+3k+2}+\sum_{k=1}^N \tau_\varphi^2(\xi_k)\hat{C}^2_k},$$

$$\tilde{C}=\int_{-\infty}^\infty \left(Z_f(\lambda)\right)d\lambda, $$

$$\sup_{|t-s|<h}|f(t,u)-f(s,u)|\leq h^{\ae}\tilde{f}(u),$$

\noindent when $K\approx 1.086435,$ $f(t,s)$ is twice differentiable in $s$ and grows in $s$ no faster than $\exp^{s^2/4}$, $Z_f(\lambda)$ is integrable on $R$.

\end{theorem}

\begin{dov}

Consider $\sigma(h)$ for this case. Indeed,

$$\sup_{|t-s|<h}\tau_\varphi(\Delta_N(t)-\Delta_N(s)) = \sup_{|t-s|<h}\tau_\varphi\left( \sum_{k=1}^N \xi_k (\delta_k(t)-\delta_k(s))+\right.$$

$$\left.+\sum_{K=N+1}^\infty \xi_k (a_k(t)-a_k(s))\right)\leq h^{\ae}\left(\sum_{k=1}^N \tau_\varphi^2(\xi_k)\hat{C}_k^2 +\right.$$

$$\left.+K^2\int_{-\infty}^\infty \left(\frac{\partial^2 \tilde{f}(\lambda)}{\partial \lambda^2} - \lambda \frac{\partial \tilde{f}(\lambda)}{\partial \lambda}+\frac{\lambda^2-2}{4}\tilde{f}(\lambda)\right)^2d\lambda\sum_{k=N+1}^\infty
\frac{4\tau_\varphi^2(\xi_k)}{\pi^2k^2}\right)^{1/2} = $$

$$=Ch^{\ae} =\sigma(h).$$

The last inequality uses the proof of Theorem \ref{omhermdec}.
Now the necessary statement follows from Theorem \ref{omhermdec} and Theorem \ref{omcdelta-3}.

\end{dov}

\begin{theorem}
\label{omct-hermite-2}
Let for a random process $X=\{X(t),t\in (0,T)\}$ from $Sub_\varphi(\Omega)$, $\varphi(t)=\frac{t^\zeta}{\zeta}$, $1<\zeta\leq2$,
the conditions of Theorem \ref{omtk.1} be satisfied, and let this process can be represented in the form (\ref{ombasis}).
Let for $\Delta_N$ of this process the condition {\bf (C1)} is satisfied.
The model $X_N(t)$ \eqref{ommodel-basis}  approximates the process $X(t)$ with given reliability $1-\nu$ and accuracy $\delta$ in the space $C(0,T)$ if

$$\nu\leq 2\exp\left\{-\frac{(\zeta-1)\delta^\frac{1}{\zeta-1}(\delta-\gamma)^\frac{\zeta}{\zeta-1}}{\zeta(\gamma^\zeta(\delta-\gamma)+\beta^\zeta\gamma)^\frac{1}{\zeta-1}}\right\}2\left(e\delta\right)^{1/\ae},$$
where

$$\gamma^\zeta\leq K^\zeta \sup_{u\in[0,T]}\int_{-\infty}^\infty Z^\zeta_f(u,\lambda)d\lambda \sum_{k=N+1}^\infty \frac{\tau^\zeta_\varphi(\xi_k)}{(k^2+3k+2)^{\zeta/2}} + $$

$$+\sum_{k=1}^N \tau^\zeta_\varphi(\xi_k)\sup_{u\in[0,T]}\delta^\zeta_k(u),$$

$$Z_f(t,\lambda) = \frac{\partial^2 f(t,\lambda)}{\partial \lambda^2} - \lambda \frac{\partial f(t,\lambda)}{\partial \lambda}+\frac{\lambda^2-2}{4}f(t,\lambda),$$

\noindent when $K\approx 1.086435,$ $f(t,s)$ is twice differentiable in $s$ and grows in $s$ no faster than $\exp\{s^2/4\}$, $Z_f(\lambda)$ is integrable on $R$.

\end{theorem}

\begin{dov}
This statement follows from the previous Theorem and Theorem \ref{omcdelta-4}.
\end{dov}

\section*{Conclusions to chapter \ref{omch:lpt2}}

In this chapter \ref{omch:lpt2}, we consider stochastic processes from the spaces $Sub_\varphi(\Omega)$ that admit a series expansion whose elements cannot be found explicitly.
We evaluate the construction of models of such processes with specified accuracy and reliability in the spaces $C(T)$.
Theorems are proved that allow us to evaluate the accuracy and reliability of the constructed model in cases where
$\varphi(t)=|t|^\gamma/\gamma$, $1<\gamma<2$, and $\varphi(t)=|t|^\gamma/\gamma$, for $|t|>1$, and $\varphi(t)=|t|^2/\gamma$, $|t|<1$, for $\gamma>2$.
The developed mechanisms are applied to evaluate the accuracy and reliability of modeling in $C(T)$ of random processes using the Karhunen-Lo\`eve decomposition in the case where the eigenvalues and eigenvectors cannot be found explicitly.


\addcontentsline{toc}{chapter}{References}

\begin{thebibliography}{99}

\bibitem{abza85}
Abzhanov E. A., Kozachenko Yu.V. Some properties of random processes in Banach $K_\sigma$-spaces.
Ukr. Mat. Zh. -- 1985. -- Vol. 37, No.~3. -- P. 275-280.


\bibitem{abza86}
Abzhanov E. A., Kozachenko Yu.V. Stochastic processes in quasi-{Banach} $K_{\sigma}$-spaces of random variables.
Probabilistic methods for the investigation of systems with an infinite number of degrees of freedom. {Collect}. sci. {Works},
{Ky\"\i v}, Institute of mathematics NAS USSR, -- 1986. p. 4-11.

\bibitem{adle90}
Adler~R.~J.
An introduction to continuity, extrema and related topics  for general Gaussian processes. -- Hayward, 1990. -- (Institute of Mathematical Statistics.
Lecture Notes - Monograph Series Vol. 12).

\bibitem{albi90}
Albin~J. M.~P.
On extremal theory for stationary processes.
Ann. Probab.  -- 1990. -- Vol.~15. -- p. 92-128.

\bibitem{albi98}
Albin~J. M.~P. On extremal theory for self-similar processes.
Ann. Probab.  -- 1998. -- Vol.~26. -- Pp. 723-793.

\bibitem{askey2010}
Askey, R.A., Roy, R.: Gamma function. In: Olver, F., Lozier, D., Boisvert, R., Clark, Ch. (eds.) NIST Handbook of Mathematical Functions, pp. 135--148. Cambridge University Press, Cambridge (2010)


\bibitem{bagr92}
Bagro~S.~V.
Some probabilistic inequalities and the central limit   theorem in function spaces.
Theory of Probability and Math. Statist. -- 1992. -- Vol.~44. -- p. 7-13.



\bibitem{Barrasa:1995}
Barrasa de la Krus, E., Kozachenko, Yu. V.: Boundary-value problems for equations of mathematical physics with strictly Orlicz random initial conditions. Random Operators and Stochastic Equations. 3(3), 201--220 (1995). https://doi.org/10.1515/rose.1995.3.3.201




\bibitem{bely60}
Belyaev ~Yu.~K.
Local properties of the sample functions of stationary Gaussian processes.
Teor. Veroyatn. Primen. -- 1960. - V. 5, No. 1. -- p. 128-131.

\bibitem{bely61}
Belyaev~Yu.~K.
Continuity and Holder's conditions for sample  functions of stationary Gaussian processes.
Proc. Fourth Berkeley Symp. on Math. and Probability.
   -- 1961. -- Vol.~2. -- Pp. 23-33.

\bibitem{bely68}
Belyaev ~Yu.~K.
On the number of exits across the boundary of a region by a vector stochastic process.
Teor. Veroyatn. Primen. -- 1968. -- V. 13, No.~2. -- p. 333-337.


\bibitem{berm73}
Berman~S.~M.
Excursions of stationary Gaussian processes above high moving barriers.
The Annals of Probability.   -- 1973. -- Vol.~1, no.~3. -- p. 133-184.


\bibitem{bond92}
Bondarenko~I.~V., Ivanov A. V.
On sample properties of random fields with stable increments.
Theory Probab. Math. Stat. -- 1992. - v. 47. -- p. 11-15.


\bibitem{bore78}
Borell~C.
Tail probabilities in Gaussian space. Vector space  measures. Appl. I.
Lect. Notes Math.
   -- 1978. -- No.~644. -- Pp. 73-82. (Proc. Conf. Dublin 1977).

\bibitem{byld77}
Buldygin~V.~V.
Sub-Gaussian processes and convergence of random series in functional spaces.
Ukr. Mat. Zh. -- 1977. - V. 29, No.~4. -- p. 443-454.


\bibitem{byld801}
Buldygin~V.~V.
Convergence of random elements in topological spaces.--
Ky\"\i v:  ''{Naukova} {Dumka}'', 1980. -- 239 p.

\bibitem{byld74}
Buldygin~V.~V., Kozachenko Yu.V.
On local properties of realizations of some stochastic processes and fields.
Theory Probab. Math. Stat. -- 1974. -- v. 10. -- p. 39-47.


\bibitem{byld80}
Buldygin~V.~V., Kozachenko Yu.V.
Sub-Gaussian random variables.
Ukr. Mat. Zh. -- 1980. -- V. 32, No.~6. -- p. 723-730.

\bibitem{byld87}
Buldygin~V.~V., Kozachenko Yu.V.
Sub-Gaussian random vectors and processes.
Theory Probab. Math. Stat. -- 1987. - V.. 36. -- p. 10-22.

\bibitem{byld92}
Buldygin~V.~V., Kozachenko Yu.V.
Exponential estimates of the distribution of the supremum of one class of random processes.
  Dopov. Akad. Nauk Ukr.  -- 1992. - No. 2. -- p. 31-34.

\bibitem{byld93}
Buldygin~V.~V., Kozachenko Yu.V.
Estimates of the supremum distribution for a certain class of random processes.
Ukr. Mat. Zh. -- 1993. - V. 45, No.~5. -- p. 596-608.

\bibitem{byld95}
Buldygin~V.~V., Kozachenko Yu. V.
Bernstein functionals and exponential inequalities for distributions of sums of random variables.
Mathematics today.   -- 1994. - vol. 9. - p. 55-79.

\bibitem{byld98}
Buldygin~V.~V., Kozachenko Yu.V.
Metric characteristics of random variables and random processes.
Transl. Math. Monogr. vol. 188, Providence, RI: AMS -- 2000. -- 257p.


\bibitem{byld97}
Buldygin~V.~V., Solntsev S.~A.
Asymptotic behaviour of linearly transformed sums of  random variables.
-- Dordrecht: Kluwer Academic Press, 1997.-- 516p.

\bibitem{busenko1978}
Buslenko~N.~P.
Modeling of complex systems.
-- M.: "Nauka", 1978. -- 356 p.

\bibitem{bykov1978}
Bykov~V.~V.
Digital modeling in statistical radiotechnics.
 -- M.: Sov. Radio, 1971. -- 328 p.


\bibitem{chui92}
Chui~C.~K.
An introduction to wavelets. -- New York: Academic Press, 1992. -- 266 p.

\bibitem{cram67}
Cram\'er~H., Leadbetter M.~R.
Stationary and related stochastic processes. Sample function properties and their applications.
 -- New York-London-Sydney: Wiley, 1967. -- 348 p.


\bibitem{daub92}
Daubechies~I.
Ten lecture on wavelets. -- Philadelphia:
Soc. Industrial and Appl. Math., 1992. -- 324 p.

\bibitem{dari11}
Darijchuk~I.~V., Kozachenko Yu.~V., Perestyuk M.~M.
Random processes from Orlicz spaces.
Chernivtsi: "Zoloti Lytavry", 2011. -- 212 p.


\bibitem{delp64}
Delporte~J.
Fonctions al\'eatoires presque s\^urement continues sur un intervalle ferm\'e.
Annales de l'I H. P.  -- 1964. -- Vol.~1, no.~2. -- p. 111-215.



\bibitem{dudl65}
Dudley~R.~M. Gaussian processes on several parameters.
Ann. Math. Statist.
   -- 1965. -- Vol.~36, no.~3. -- p. 771-788.

\bibitem{dudl73}
Dudley~R.~M. Sample functions of the Gaussian processes.
The Annals of Probability.
   -- 1973. -- Vol.~1, no.~1. -- p. 3-68.


\bibitem{hermite2}
{Erd{\'e}lyi, Arthur and Magnus, W. and Oberhettinger, F. and Tricomi, Francesco G.}
Higher transcendental functions, Vol. II. -- New-York: McGraw-Hill, 1953. -- 396 p.

\bibitem{erma75}
Ermakov~S.~M.
Monte Carlo method and related topics. - M.: Nayka, 1975.-- 471 p.

\bibitem{erma76}
Ermakov~S.~M., Mikhajlov G.~A.
A course in statistical modeling.
-- M. ``{Nauka}''. 1976. -- 319 p.

\bibitem{ermmikh1982}
Ermakov~S.~M., Mikhajlov G.~A.
Statistical simulation. Textbook.
 -- M. ``{Nauka}''. 1982. -- 296p.

\bibitem{erma86}
Ermakov~S.~M., Ostrovskij, E.I.
Continuity conditions, exponential estimates and the central limit theorem for random fields /
VINITI Dep. -- 1986. --  No.~3752-B.86.0. -- 42 p.

\bibitem{fern70}
Fernique~X.
Int\'egrabilite des vecteurs gaussiens.
C. R. Acad. Sci.
   -- 1970. -- Vol.~270, no.~7. -- p. 1698-1699.

\bibitem{fern75}
Fernique~X.
R\'{e}gularit\'{e} des trajectoires des fonctions  al\'{e}atoires gaussiennes.
Lecture Notes in   Mathematics. -- 1975. -- Vol.~480. -- p. 1-96.

\bibitem{fern83}
Fernique~X.
R\'{e}gularit\'{e} de fonctions al\'{e}atoires non  gaussiennes.
Lecture Notes in  Mathematics. -- 1983. -- Vol.~976. -- p. 1-74.


\bibitem{fuku90}
Fukuda~R.
Exponential integrability of sub-Gaussian vectors.
Probab. Theory Related Fields.
 -- 1990. -- Vol.~85, no. 4. -- p. 505-521.

\bibitem{dzyl02}
 Giuliano~A.~R., Kozachenko Yu.~V., Tegza A.~ M.
Inequalities for the norms of sub-Gaussian vectors and the accuracy of the modelling of random processes.
Teor. Jmovirn. Mat. Stat. -- 2002. -- V. 66.  -  p. 58-66.

\bibitem{julantonini2002-2}
Giuliano~A.~R., Kozachenko Yu.~V., Tegza A.~ M.
Accuracy of simulation in $L_p$ of Gaussian random processes.
Visn., Mat. Mekh., Ky\"\i v Univ. Im. Tarasa Shevchenka. -- 2002.  -  No. 5 -- p. 7-14


\bibitem{giul03}
Giuliani A. R., Kozachenko Yu.~V., Nikitina T.
Spaces of $\varphi$-subgaussian random variables.
Rediconti Accademia Nazionale dell Scienze detta dei XL, Memorie di Matimatica e Applicazioni. -- 2003. -- No. 121. -- p. 95 - 124.


\bibitem{giul2013}
Giuliano Antonini Rita, Hu Tien-Chung, Kozachenko Yuriy, Volodin, Andrei.
An application of $(\varphi)$-subGaussian technique to {Fourier} analysis.
J. Math. Anal. Appl. - 2013. -  v. 408,  n. 1. -  p. 114- 124.

\bibitem{Hardle98}
 H\"ardle~W., Kerkyacharian~G., Picard~D., Tsybakov~A. Wavelets, approximation and statistical applications.
 New York: Springer, 1998. 265~p. DOI:~10.1007/978-1-4612-2222-4.


\bibitem{Ian2022}
Ianevych Tetiana, Rozora Iryna, Pashko, Anatolii.
On one way of modeling a stochastic process with given accuracy and reliability.
Monte Carlo Methods Appl. -  2022. -  v. 28,  n. 2. -  p. 135- 147.

\bibitem{Ian2024}
Ianevych Tetiana, Vasylyk Olga, Doshchuk, Julia.
On modeling {Gaussian} stationary {Ornstein}-{Uhlenbeck} processes with given reliability and accuracy in $(L_p)$-spaces.
Visn., Ser. Fiz.-Mat. Nauky, Ky{\"{\i}}v. Univ. Im. Tarasa Shevchenka. -  2024. -  v. 78, n. 1. -  p. 51- 56.



\bibitem{jain75}
Jain~N.~C., Marcus M.~B.
Integrability of infinite sums of independent  vector-valued random variables.
Trans.  Amer. Math. Soc. -- 1975. -- Vol.~212, no. 1. -- p. 1-36.

\bibitem{jain78}
Jain~N.~C., Marcus M.~B.
Continuity of sub-Gaussian processes.
Adv. Probab. -- 1978. -- Vol.~4. -- p. 81-196.

\bibitem{kaha60}
Kahane~J.~P. Propri\'et\'es locales des fonctions \`a series de {F}ouries al\'eatories.
Studia Math.  -- 1960. -- Vol.~19, no. 1. -- p. 1-25.

\bibitem{kaha68}
Kahane~J.~P.
Some random series of functions.
Heath Mathematical Monographs. -- Lexington, Mass: D.C. Heath and Company, 1968.
-- 184 p.


\bibitem{hermite}
Kamp{\'e} de F{\'e}riet, J. and Campbell, R. and Petiau, G. and Vogel, T.
Functions of mathematical physics.
-- M.: ``State Publishing House of Physical and Mathematical Literature'', 1963. -- 102 p.


\bibitem{kame07}
Kamenschykova~O.
Approximation of random processes in the space  ${L}_{2}({T})$.
Theory of Stochastic  Processes. -- 2007. -- Vol.~13, no. 29. -- p. 64-68.

\bibitem{kame071}
Kamenshchikova~O.~E.
Approximation of random processes by Bernstein polynomials.
Prykl. Stat., Aktuarna Finans. Mat. 2007, No. 2, p. 112- 117



\bibitem{kant84}
 Kantorovich~L.~V., Akilov G.~P.
 Functional analysis.
 -- M. ``{Nauka}''. 1984. -- 752 p.

\bibitem{Kolm}
Kolmogorov A. N.
Selected works II. Probability theory and mathematical statistics. Ed. by A. N. Shiryayev.  Dordrecht: Springer. - 2019. - 597 p.

\bibitem{kono80}
K\^{o}no~N.
Sample path properties of stochastic processes.
J. Math. Kyoto Univ. -- 1980. -- Vol.~20, no. 2. -- p. 295-313.

\bibitem{koz68}
Kozachenko~Yu.~V.
Sufficient conditions of continuity with unity probability of sub-Gaussian random processes.
Dopov. Akad. Nauk Ukr. RSR, Ser. A, 1968.  -  No.~2. -- p. 113-115.

\bibitem{koz83}
Kozachenko~Yu.~V.
On the uniform convergence of stochastic integrals with respect to the norm of the {Orlicz} space.
Teor. {Veroyatn}. {Mat}. {Stat}. -- 1983.  -  V.~29. -- p. 52-64.

\bibitem{koz84}
Kozachenko~Yu.~V.
Stochastic processes in {Orlicz} spaces.
Teor. {Veroyatn}. {Mat}. {Stat}.-- 1984.  -  No.~30. -- p. 92-107.

\bibitem{koz841}
Kozachenko~Yu.~V.
Stochastic processes in {Orlicz} spaces. II.
Teor. {Veroyatn}. {Mat}. {Stat}.-- 1984.  -  No.~31. -- p. 44-50.

\bibitem{koz851}
Kozachenko~Yu.~V.
Stochastic processes in Orlicz spaces.
Trajectory property, convergence of series and integrals.
Dissertation for the degree of Doctor of Physical and Mathematical Sciences: 01.01.05. -- Kyiv, 1985. -- 296 p.


\bibitem{koz04}
Kozachenko~Yu. ~V.
Lectures on wavelet analysis.
 -- Ky\"\i v: TViMS, 2004. -- 147 p.



\bibitem{kozkoz1991}
Kozachenko~Yu.~V., Kozachenko L.~F.
Accuracy of modeling stationary Gaussian stochastic processes in $L^2(0,T)$.
Vychisl. Prikl. Mat. -- 1991.  -  No. 75. -- p. 108-115.

\bibitem{kozkoz1992}
Kozachenko~Yu.~V., Kozachenko L.~F.
On accuracy of modeling of Gaussian stochastic processes in $L^2(0,T)$.
Vychisl. Prikl. Mat.-- 1992.  -  No. 74. -- p. 88-93.


\bibitem{kozDar2009}
Kozachenko Yu. V., Darijchuk I. V.
The distribution of the supremum of $(\Theta)$-pre-{Gaussian} shot noise processes.
Teor. {\u{I}}movirn. Mat. Stat. -  2009. -  v. 80. -  p. 76--90.


\bibitem{koz-kam}
Kozachenko~Yu.~V., Kamenshchikova O.~E.
Approximation of $SSub_\varphi(\Omega)$ stochastic processes in the space  $L_p(T)$.
Theory Probab. Math. Stat. -- 2009.  --  Iss. 79.  --  p. 83-88.


\bibitem{kozKam-2014}
Kozachenko Yu. V., Kamenshchikova O. E.
An approximation of stochastic processes belonging to the {Orlicz} space in the norm of the space $(C[0,\infty))$.
Theory Probab. Math. Stat. -  2014. -- v. 88. --- p. 123-138.


\bibitem{mlav111}
Kozachenko~Yu.~V., Mlavets' Yu. Yu.
The accuracy and reliability of the calculation of integrals by the Monte Carlo method.
Dopov. Nats. Akad. Nauk Ukr., Mat. Pryr. Tekh. Nauky.-- 2011.  -  No.~8. -- p. 18-20.


\bibitem{mlav11}
Kozachenko~Yu.~V., Mlavets Yu.~Yu.
Probability of large deviations of sums of random processes from Orlicz space.
Monte Carlo Methods Appl.   -- 2011. -- Vol.~17. -- p. 155-168.

\bibitem{mlav12}
Kozachenko~Yu.~V., Mlavets' Yu. Yu.
The Banach spaces $\mathbf{F}_\psi({\Omega})$ of random variables.
Theory Probab. Math. Stat. -  2013. --v. 86. -  p. 105- 121.

\bibitem{mlav143}
Kozachenko~Yu.~V., Mlavets Yu.~Yu.
Stochastic processes from $\mathbf{F}_\psi(\Omega)$ spaces.
Contemporary Mathematics and Statistics. -- 2014. -- Vol.~2, no.~1. -- p. 55-75.

\bibitem{mlav2015}
Kozachenko~Yu.~V., Mlavets Yu.~Yu.
Reliability and accuracy in the space $(L_{p}(T))$ for the calculation of integrals depending on a parameter by the {Monte} {Carlo} method.
Monte Carlo Methods Appl. -  2015. -  vol. 21, no. 3. -  p. 233- 244.


\bibitem{mlav2016}
Kozachenko Yu. V.,  Mlavets, Yu. Yu.
An application of the theory of spaces $(\mathbf{F}_\psi(\Omega))$ for evaluating multiple integrals by using the {Monte} {Carlo} method.
Theory Probab. Math. Stat. -  2016. -- Vol. 92. -  p. 59- 69.


\bibitem{mlav2018}
Kozachenko Yu.V., Mlavets Yu.Yu., Yurchenko N.V.
Weak convergence of random processes from space $\mathbf{F}_\psi(\Omega)$.
Statistics, Optimization and Information Computing. -- 2018, -- Vol.6, Iss.2 -- p. 266 -  277.


\bibitem{kozMokl2011a}
Kozachenko Yu. V., Mokliachuk O. M.
Random processes in the spaces $(D_{V,W})$.
Theory Probab. Math. Stat. -  2011. -  v. 82. --  p. 43 --56.


\bibitem{kozMokl2011}
Kozachenko Yu. V., Mokliachuk O. M.
Sample continuity and modeling of stochastic processes from the spaces $(D_{V,W})$.
Theory Probab. Math. Stat. -  2011. -  v. 83. --  p. 95 --110.

\bibitem{Kozachenko:2011}
Kozachenko, Y., Olenko, A., Polosmak, O.: Uniform Convergence of Wavelet Expansions of Gaussian Random Processes. Stochastic Analysis and Applications. 29(2), 169--184 (2011). https://doi.org/10.1080/07362994.2011.532034


\bibitem{kozOlenko2015}
Kozachenko Yuriy, Olenko Andriy, Polosmak, Olga.
Convergence in $(L_{p}([0, {T}]))$ of wavelet expansions of  $(\varphi)$-sub-{Gaussian} random processes.
Methodol. Comput. Appl. Probab. -  2015. -  vol. 17, no. 1. -  p. 139- 153.


\bibitem{kozOlenko-2016}
Kozachenko Yu., Olenko A.
Whitaker-Kotelnikov-Shanon approximation of $\varphi$-sub-Gaussian random processes.
Journal of Mathematical Analysis and Applications. -  2016. -- Vol. 442, Iss.2. -- p. 924 -  946.


\bibitem{koz85}
Kozachenko~Yu.~V., Ostrovskij E.~I.
Banach spaces of random variables of sub-Gaussian type.
Teor. Veroyatn. Mat. Stat. -- 1985. -- v.~32. -- p. 42-53.

\bibitem{kozpash1998}
Kozachenko~Yu.~V., Pashko A.~O.
Accuracy of simulation of stochastic processes in norms of Orlicz spaces. I.
Teor. Jmovirn. Mat. Stat.. -- 1998.  -  No. 58. -- p. 45--60.


\bibitem{kozpash1999}
Kozachenko~Yu.~V., Pashko A.~O.
Accuracy of simulation of stochastic processes in norms of Orlicz spaces. II.
Teor. Jmovirn. Mat. Stat.. -- 1999.  -  No. 59. -- p. 77--92.

\bibitem{koz-pash}
Kozachenko~Yu.~V., Pashko A.~O.
Modelling of random processes.
-- Ky\"\i v: Vydavnychy\u\i\, Tsentr ``Ky\"\i vs'ky\u\i\, Universytet''.-- 1999. -- 223 p.

\bibitem{kozpash1999-2}
Kozachenko~Yu.~V., Pashko A.~O.
On the simulation of random fields. I.
Teor. Jmovirn. Mat. Stat.-- 1999.  -  No. 61. -- p. 59--71.

\bibitem{kozpash2000}
Kozachenko~Yu.~V., Pashko A.~O.
On the simulation of random fields. II.
Teor. Jmovirn. Mat. Stat.-- 2000.  -  No. 62.  -  p.. 48--60.

\bibitem{koz-pash-roz2007}
Kozachenko, Yu. V., Pashko, A. O., Rozora, I. V.
Modelling of stochastic processes and fields.
Ky\"\i v: Zadruga, 2007. -- 230 p.


\bibitem{kozPash2018}
Kozachenko Yu. V., Pashko A. O., Vasylyk O. I.
 Simulation of a fractional {Brownian} motion in the space $(L_p([0,T]))$.
Theory Probab. Math. Stat.- 2018. - v. 97 -- p. 99- 111.

\bibitem{kozPash2018b}
Kozachenko Yu. V., Pashko A. M., Vasylyk, O. I.
Simulation of generalized fractional {Brownian} motion in $(C([0,T]))$.
Monte Carlo Methods Appl. -- 2018. -  v. 24, n. 3 -- p. 179 -192.


\bibitem{pere07}
Kozachenko~Yu.~V., Perestyuk M.~M.
On the uniform convergence of wavelet expansions of random processes from Orlicz spaces of random variables. I.
  Ukr. Mat. Zh. -- 2007.  -  V.~59, No.~12. -- p. 1647-1660.

\bibitem{pere08}
Kozachenko~Yu.~V., Perestyuk M.~M.
On the uniform convergence of wavelet expansions of random processes from Orlicz spaces of random variables. II.
  Ukr. Mat. Zh.-- 2008.  -  V.~60, No.~6. -- p. 759-775.

\bibitem{pere06}
Kozachenko Yu. V., Perestyuk M. M., Vasylyk O. I.
On uniform convergence of wavelet expansions of $(\varphi)$-sub-{Gaussian} random processes.
Random Oper. Stoch. Equ. -  2006. -- v. 14, n. 3,--  p. 209- 232.


\bibitem{kozPetr2017}
Kozachenko Yu.V., Petranova M. Yu.
Simulation of Gaussian stationary Ornstein-Uhlenbeck process with given reliability and accuracy in space $C([0,T])$.
Monte Carlo Methods and Applications.-- 2017. -- Vol. 23, Iss. 4 -- p. 277 - 286


\bibitem{kozPogor2007}
Kozachenko Yu. V., Pogorilyak, O. O.
Modelling log {Gaussian} {Cox} processes with a given reliability and accuracy.
Teor. {\u{I}}movirn. Mat. Stat. -  2007. -  v. 76. -  p. 70- 83.


\bibitem{kozPogor2007b}
Kozachenko Yu. V., Pogorilyak, O. O.
A method of modelling log {Gaussian} {Cox} processes.
Teor. {\u{I}}movirn. Mat. Stat. -  2007. -  v. 77. -  p. 82--95.

\bibitem{kozBook-2016}
Kozachenko Yu. V., Pogorilyak, O.O., Tegza, A.M.
Modelling of Gaussian stochastic processes and Cox processes.
Uzhgorod: "Karpaty", 2012. -- 194p.

\bibitem{kozBook-2016}
Kozachenko, Yuriy; Pogorilyak, Oleksandr; Rozora, Iryna; Tegza, Antonina.
Simulation of stochastic processes with given accuracy and reliability. Amsterdam: Elsevier/ISTE Press, 2016. -- 333 p.


\bibitem{kozPolos2009}
Kozachenko Yu. V., Polos'mak O. V.
Conditions for the uniform convergence in probability of wavelet decompositions for stochastic processes from the space $(\text{Exp}_{\varphi}(\Omega))$.
Teor. {\u{I}}movirn. Mat. Stat. -  2009. -  v. 81. -  p. 76--87.



\bibitem{koz-roz2019}
Kozachenko, Yu. V., Rozora, I. V.
 Construction of the {Karhunen}-{Lo{\`e}ve} model for an input {Gaussian} process in a linear system by using the output process,  Theory Probab. Math. Stat. -- 2019.--
v. 99, -- p. 113- 124.

\bibitem{Kozachenko:2016-1}
Kozacenko, Y.V., Rozora, I.V.: A Criterion for Testing Hypothesis about Impulse Response Function. Statistics, Optimization \& Information Computing. 4(3), 214--232 (2016). https://doi.org/10.19139/soic.v4i3.222

\bibitem{koz-roz-turch}
Kozachenko~Yu.~V., Rozora I.~V., Turchyn Ye.~V.
On expansion of random process in series.
Random Operators and Stochastic Equations. -- 2007. -- Vol. 15, no 1. -- p. 15--35.

\bibitem{Kozachenko:2011-1}
Kozachenko, Y.V., Rozora, I.V., Turchyn, Y.V.:
Properties of Some Random Series. Communications in Statistics -- Theory and Methods. 40(19-20), 3672--3683 (2011). https://doi.org/10.1080/03610926.2011.581188


\bibitem{koz-ryaza1992}
Kozachenko, Yu. V., Ryazantseva, V. V.
Conditions for boundedness and continuity in terms of majorizing measures of random processes in certain {Orlicz} spaces.
Theory Probab. Math. Stat. --  1992. -- n.44, --  p. 69-75.

\bibitem{koztegza2002}
Kozachenko~Yu.~V., Tegza A.~M.
Applications of the theory of ${\bf Sub}_\varphi(\Omega)$ spaces of random variables for modelling stationary Gaussian processes.
Teor. Jmovirn. Mat. Stat.-- 2002.  -  No. 67.  -  p. 71-78.

\bibitem{kozTegza2017}
Kozachenko Yu.V., Tegza A.M.
Construction of a model of Gaussian stationary random process in some Orlicz spaces with given accuracy and reliability.
Journal of Applied Mathematics and Statistics. -  2017. -- Vol. 4, Iss. 2 --p. 70 -  77.

\bibitem{koz2020}
Kozachenko Yu. V., Tegza A.M., Troshki N.V.
The accuracy of modeling of Gaussian stochastic processes in some Orlicz space.
Statistics Optimization and information computing.- 2020. -- Vol. 8, Iss.7, -- p. 127- 135.

\bibitem{koz03}
Kozachenko~Y.~V.,Vasylyk O. I., Yamnenko R.Ye.
Upper estimate of overrunning by  $sub_\varphi(\omega)$ random process the level specified by continuous   function.
Random  Oper. and Stochastic Equations. -- 2003. -- Vol.~11, no. 1. -- p. 1-20.


\bibitem{yako06}
Kozachenko~Yu.~V., Yakovenko T.~O.
Stochastic processes in Sobolev-Orlicz spaces.
  Ukr. Mat. Zh.-- 2006.  -  V.~58, No.~10. -- p. 1517-1537.


\bibitem{kozYamn-2014}
Kozachenko Yu. V., Yamnenko R. E.
Application of $(\varphi)$-sub-{Gaussian} random processes in queueing theory.
Springer Optimization and Its Applications. -  2014. -  Vol.90. -  p. 21-38.


\bibitem{koz76}
Kozachenko~Yu.~V., Yadrenko M.~Y.
Local properties of sample functions of random fields. {I}.
Teor. {Veroyatn}. {Mat}. {Stat}.-- 1976.  -  No.~14. -- p. 53-66.

\bibitem{koz761}
Kozachenko~Yu.~V., Yadrenko M.~Y.
Local properties of sample functions of random fields. {II.}
Teor. {Veroyatn}. {Mat}. {Stat}.-- 1976.  -  No.~15. -- p. 82-98.


\bibitem{kozZatula2019}
Kozachenko Yu. V., Zatula D.
Estimates for distributions of H\"older semi-norms of random processes from $\mathbf{F}_\psi({\Omega})$ spaces, defined on the interval [0, R+].
Statistics, optimization and information computing. -- 2019. -- Vol.7, Iss.1, -- p. 198 -- 210.


\bibitem{kras58}
 Krasnosel'ski{\u{\i}}~M.~A., Rutitski{\u{\i}} Ya.~B.
Convex functions and Orlicz spaces.
-- M.: Fizmatgiz, 1958. -- 271 p.


\bibitem{land70}
Landau~H.~J., Shepp L.~A.
On the supremum of Gaussian processes.
Sankhya. -- 1970. -- Vol.~32, no. 4. -- p. 369-378.

\bibitem{lead83}
Leadbetter~M.~R., Lindgren G., Rootzen H.
 Extremes and related properties of random sequences  and processes.
-- Berlin: Springer, 1983.-- 336 p.

\bibitem{ledo90}
Ledoux~M.
A note on large deviations for Wiener chaos .
Lecture Notes in  Mathematics. -- 1990. -- Vol.~1426. -- p. 1-14.  -
 (Seminaire de probebilities XXIV 1988/89).

\bibitem{ledo91}
Ledoux~M., Talagrand M.
Probability in Banach space. -- Berlin-New York: Springer-Verlag, 1991. -- 480 p.

\bibitem{lifs95}
Lifshits~M.~A.,
Gaussian random functions.
 -  Ky\"\i v: TViMS, 1995. -- 246 p.


\bibitem{lind71}
Lindgren~G. Extreme values of stationary normal processes.
Z. Wahrscheinlichkeitstheor. Verw. Geb. -- 1971. -- No. 17. -- p. 39-47.


\bibitem{lyka79}
Lukacs~E.
Characteristics functions.
-- M. ``{Nauka}''. 1979. -- 423 p.


\bibitem{mall98}
Mallat~S.~G.
A wavelet tour of signal processing.
-- San Diego: Academic Press, 1998. -- 637 p.


\bibitem{marc88}
Marcus~M.~B.
Continuity in $l^p$ of certain Ornstein-Uhlenbeck processes.
Probability in Banach space 7, Proc. 7th Int. Conf., Oberwolfach/FRG. -- Oberwolfach, 1988. -- p. 139-145.

\bibitem{marc81}
Marcus~M.~B., Pisier G.
Random Fourier series with applications to harmonic analysis.
-- Prinston: Prinston University Press, 1981. -- 152 p.  -
 (Annals of Mathematics Studies Vol. 101).

\bibitem{mats88}
Matsak~I.~K., Plichko A.~N.
Some inequalities for sums of independent random variables in Banach spaces.
Theory of Probability and Math. Statist. -- 1988. -- Vol.~38. -- p. 81-88.

\bibitem{metr49}
Metropolis~N., Ulam S.
The Monte Carlo method.
Journal of the American Statistical Association.
-- 1949. -- Vol.~44, no.~247. -- p. 335-341.

\bibitem{meye90}
Meyer~Y. Ondelettes et Op\'erateurs.
-- Paris: ``Hermann'', 1990. -- Vol. I, II.


\bibitem{mlav121}
Mlavets'~Yu.~Yu.
On the distribution of suprema of the increments of stochastic processes from the spaces $\mathbf{F}_\psi({\Omega})$.
Nauk. Visn. Uzhgorod. Univ., Ser. Mat. Inform.  -- 2012.
  -- Vol.~23, No.~1. -- P. 79-88.

\bibitem{mlav122}
Mlavets'~Yu.~Yu.
$\mathbf{F}_\psi({\Omega})$-spaces of random variables with exponential function $\psi$.
Visn., Ser. Fiz.-Mat. Nauky, Ky\"\i v. Univ. Im. Tarasa Shevchenka -- 2012.
  -- No.~2 -- P. 19-22.

\bibitem{mlav141}
Mlavets'~Yu.~Yu.
The connection of Orlicz spaces of random variables with spaces $\mathbf{F}_\psi({\Omega})$.
Nauk. Visn. Uzhgorod. Univ., Ser. Mat. Inform. -- 2014.
  -- Vol.~25, No.~1. -- P. 77-84.

\bibitem{mlav142}
Mlavets'~Yu.~Yu.
 Conditions for uniform convergence of random functional series from spaces
 $\mathbf{F}_\psi({\Omega})$.
Prykl. Stat., Aktuarna Finans. Mat.
 -- 2014. -- No.~1. -- P. 97-103.

\bibitem{mlav144}
Mlavets'~Yu.~Yu.
Condition "H" for Orlicz spaces of exponential type.
Nauk. Visn. Uzhgorod. Univ., Ser. Mat. Inform. -- 2014.
  -- Vol.~26, No.~2. -- P. 118-122.








\bibitem{mokl2007}
Moklyachuk O.M.
Simulation of random processes with known correlation function with the help of {Karhunen}-{Lo{\`e}ve} decomposition.
Theory Stoch. Process. -  2007. -  v. 13,  n. 4. --  p.163 -169.

\bibitem{mokl2009}
Moklyachuk, O. M.
Simulation of {Gaussian} random sequences in $L_p(T)$.
Visn., Ser. Fiz.-Mat. Nauky, Ky{\"{\i}}v. Univ. Im. Tarasa Shevchenka. -- 2009. --  n.~1. --
p. 23-26.


\bibitem{mokl2012}
Moklyachuk, O. M.
Modelling with given reliability and accuracy in the space ${L_p}({T})$ of stochastic processes from $Sub_{{\varphi}}({{\Omega}})$ decomposable in series with independent elements.
Applied statistics. Actuarial and financial mathematics. -  2012. -- No. 2. -  p. 13- 23.

\bibitem{mokl2014}
Moklyachuk, O. M.
Modeling of stochastic processes in ${L_p}({T})$  using orthogonal polynomials.
Universal Journal of Applied Mathematics. -  2014. - Vol. 2(3). -- p. 141-147.

\bibitem{mokl2018}
Moklyachuk, O. M.
Estimation of accuracy and reliability of models of $\varphi$-sub-Gaussian stochastic processes in $C(T)$ spaces.
Research Bulletin of the National Technical University of Ukraine "Kyiv Polytechnics Institute". - 2017 - No. 4 -- p. 17-24.

\bibitem{mokl2022}
Moklyachuk, O. M.
Quasi-Banach spaces of random variables and modeling of stochastic processes.
In "Stochastic processes. Fundamentals and emerging applications",
{New York, NY: Nova Science Publishers}, -2022. -  p. 351- 414.



\bibitem{nano78}
Nanopoulos~C., Nobelis P.
Regularite et proprietes limites des fonctions aleatoires
Semin. Probab. XII, Univ.  Strasbourg 1976/77. -- Lect. Notes Math.
-- 1978. -- Vol.~649. -- p. 567-690.

\bibitem{ogorodnik1996}
Ogorodnikov~V.~A., Prigarin  S.~M.
Numerical modeling of random processes and fields: algorithms and applications.
  -- Utrecht: VSP, -- 1996. -- 240 p.

\bibitem{ostrovsky1973}
Ostrovskij~E.~I.
Convergence of the canonical expansion for normal fields.
Mat. Zamet. -- 1973. -- Vol. 14, No. 4. -- P. 565-572.

\bibitem{ostr82}
Ostrovskij~E.~I.
Generalization of the Buldygin-Kozachenko norm and the central limit theorem in Banach spaces.
Theory Probab. Appl. -- 1982.  -- Vol.~27, No.~3. -- P. 618-623.

\bibitem{ostr90}
Ostrovskij~E.~I.
Exponential bounds for the distribution of the maximum of a non-Gaussian random field.
Theory Probab. Appl. -- 1990.  -- Vol.~35, No.~3. -- P. 482-493.


\bibitem{pash98}
Pashko~A.~O.
Uniform convergence of sub-Gaussian integrals.
Theory Probab. Appl. -- 1998.  -- Iss.~43, No.~4. -- P. 650-655.


\bibitem{pashko2000}
Pashko~A.~O.
Assessment of simulation accuracy of Gaussian isotropic random fields on a sphere in the uniform metric.
Theory Probab. Math. Stat. -- 2000. -- Vol. 60. -- P. 149--157.

\bibitem{pashko2001}
 Pashko~A.~O.
The estimation of accuracy of simulation of sub-Gaussian random processes.
Visn., Mat. Mekh., Ky\"\i v Univ. Im. Tarasa Shevchenka. -- 2001. -- No. 6. -- P. 42-47.

\bibitem{pashko2001-2}
Pashko~A.~O.
An estimation of accuracy of the simulation of sub-Gaussian random fields in a uniform metric.
 Dopov. Nats. Akad. Nauk Ukr. -- 2001. -- No. 2. -- P. 30-36.

\bibitem{pashko2020}
Pashko A. O., Rozora I. V., Ianevych T. O.
 On modelling of {Gaussian} process with accuracy and reliability in the space $L_p[0,T]$.
Nauk. Visn. Uzhgorod. Univ., Ser. Mat. Inform. -  2020. -  v. 37, n. 2, -- p. 91- 100.


\bibitem{pearson1901}
Pearson~K.
On lines and planes of closest fit to systems of points in space.
Philosophical Magazine. -- 1901 -- No. 2. -- p. 559-572.

\bibitem{petranova2020}
Petranova, M. Yu.
Simulation of a {Gaussian} stationary process with a stable correlation function with a given reliability and accuracy.
Visn., Ser. Fiz.-Mat. Nauky, Ky{\"{\i}}v. Univ. Im. Tarasa Shevchenka. --
2020. -  n. 3, p. 89-95.

\bibitem{pick69}
Pickands~J.~III.
Upcrossing probabilities for stationary Gaussian processes.
Trans. Am. Math. Soc. -- 1969. -- Vol.~145. -- p. 51-73.

\bibitem{pisi79}
Pisier~G.
Conditions d'entropie assurant la continuite de certains  processus et applications
 \`a l'analyse harmonique.
Semin. Anal. Fonct.  1979-1980. -- No.~13-14. -- p. 43.

\bibitem{pisi83}
Pisier~G.
Some applications of the metric entropy condition to harmonic analysis.
Lect. Notes Math.  -- 1983. -- Vol.~995. -- p. 123-154.

\bibitem{pite96}
Piterbarg~V.~I.
 Asymptotic methods in the theory of Gaussian processes and fields.
-- Providence: AMS, 1996. -- 205 p. -- (Transl. of Math. Monographs, Vol. 148).

\bibitem{pite82}
Piterbarg~V.~I.
 Large deviations of random processes close to Gaussian ones.
Theory Probab. Appl.   -- 1982. -- Iss.~27, No.~3. -- P. 474-491.



\bibitem{prodayvoda1999}
Prodayvoda~G.~T., Vyzhva S. A.
 Mathematical modeling of geophysical parameters.
Ky\"\i v: Vydavnycho-Poligrafichny\"\i\ Tsentr, Ky\"\i vsky\"\i \ Universytet, 1999. -- 112 p.

\bibitem{prodayvoda2009}
Prodayvoda~G.~T., Vyzhva S. A., Bezrodnyj D. A., Bezrodna I.M.
Mathematical modeling of tectonofacies of metamorphic rocks of the Kryvbas epizone.
Geoinformatics. -- 2009. -- No. 3. -- P. 68-73.

\bibitem{rakhimovyadremko1993}
Rakhimov, G., Yadrenko M.~I.
Statistical simulation of a homogeneous isotropic random field on the plane.
Theory Probab. Math. Stat.  -- 1993. -- Iss. 49. -- P. 245-251.


\bibitem{raom91}
Rao~M.~M., Ren Z.~D.
Theory of Orlicz spaces.
-- New York -- Basel -- Hong Kong: Marsel Dekker, 1991. -- 445 p.

\bibitem{raom01}
Rao~M.~M., Ren Z.~D.
Applications of Orlicz spaces
-- New York -- Basel: Marsel Dekker, 2001. -- 464 p.

\bibitem{ripley1987}
Ripley~B.~D.
Stochastic simulation.  -- New York: Wiley, 1987. -- 237 p.


\bibitem{rozora2004}
Rozora~I.~V.
Accuracy and reliability of models of stochastic processes of the space $Sub_\varphi(\Omega)$.
Teor. Jmovirn. Mat. Stat. -- 2004. -- Iss. 71. -- P. 93-105.

\bibitem{rozora2009}
Rozora~I.~V.
Conditions for finding a model of a stochastic process with given accuracy and reliability in the space  $C[0,T]$.
Visn., Ser. Fiz.-Mat. Nauky, Ky{\"{\i}}v. Univ. Im. Tarasa Shevchenka. --
2009.  --  n. 1, -- p. 139-144.

\bibitem{rozora2015}
Rozora~I.~V.
Model accuracy and reliability of stochastic processes with discrete spectrum with respect to linear filter.
Visn., Ser. Fiz.-Mat. Nauky, Ky{\"{\i}}v. Univ. Im. Tarasa Shevchenka. --
2015. --  n. 1, -- p. 163-167.

\bibitem{rozora2018}
Rozora~I.~V.
 On simulation accuracy and reliability in the space ${L_p}[0,{T}]$ for the input {Gaussian} process served by the linear system taking into account the output.
Visn., Ser. Fiz.-Mat. Nauky, Ky{\"{\i}}v. Univ. Im. Tarasa Shevchenka. --
2018. - n. 2, -- p. 75-80.

\bibitem{rozora2018b}
Rozora~I.~V.
 Convergence rate for the estimator of impulse response function in the space of continuous functions.
Visn., Ser. Fiz.-Mat. Nauky, Ky{\"{\i}}v. Univ. Im. Tarasa Shevchenka. --
2018. --  n. 3, --  p. 30-36.

\bibitem{rozora2018c}
Rozora~I.~V.
On the convergence rate for the estimation of impulse response function in the space $L_p(T)$.
Visn., Ser. Fiz.-Mat. Nauky, Ky{\"{\i}}v. Univ. Im. Tarasa Shevchenka. --
2018. - n. 4, - p. 36-41.

\bibitem{rozora2020}
Rozora I. V., Pereyaslov O. S.
The development of software for simulation of random processes with a given accuracy and reliability.
Visn., Ser. Fiz.-Mat. Nauky, Ky{\"{\i}}v. Univ. Im. Tarasa Shevchenka. --
2020. -- n. 1-2, -- p. 83-87.


\bibitem{rozora2022}
Rozora Iryna, Ianevych Tetiana, Pashko Anatoliy, Zatula Dmytro.
Simulation of stochastic processes with given reliability and accuracy.
in: Stochastic processes. Fundamentals and emerging applications. New York, NY: Nova Science Publishers. -  2022. -  p. 415- 452.

\bibitem{rozora2024a}
Rozora Iryna, Sheptukha, Yevhenii.
Simulation of the fractional {Brownian} process with given accuracy and reliability.
Visn., Ser. Fiz.-Mat. Nauky, Ky{\"{\i}}v. Univ. Im. Tarasa Shevchenka. -  2024. --
v. 78,  n. 1, -  p. 147- 153.

\bibitem{rozora2024b}
Rozora Iryna, Mlavets Yurii, Vasylyk Olga, Polishchuk Volodymyr.
 On convergence of the uniform norm and approximation for stochastic processes from the space $\mathbf{F}_{{\psi}} ({{\Omega}})$.
J. Theor. Probab. - 2024. -  v. 37,  n. 2, -- p. 1627- 1653.



\bibitem{shaligin1985}
Shalygin, A. S., Palagin  Yu. I.
Applied methods of statistical modelling.
  -- Leningrad: "Mashinostroenie" , 1986. -- 320 p.


\bibitem{skor70}
Skorokhod~A.~V.
A note on Gaussian measures in a Banach space. Theor. Probab. Appl.
  -- 1970. -- Vol.~15, No.~3. -- P. 519-520.

\bibitem{skor73}
Skorokhod~A.~V.
The theorem on the continuity of a random function on a compact set in a Hilbert space.
Theor. Probab. Appl.  -- 1973. -- Vol.~18, No.~4. -- P. 809-811.


\bibitem{Slut37}
Slutsky~E.
Alcuno proposizioni sulla teoria delle funzioni aleatorie.
Giorn. Inst. Italiano degli Attuari.
-- 1937. -- Vol.~8. -- Pp. 183-199.


\bibitem{sun1997}
Sun~T.~C., Chaika~Milton.
On simultaion of Gaussian stationary process.
Journal of Time Series Analysis. -- 1997. -- Vol. 18, No. 1. -- p. 79--93.

\bibitem{tala87}
Talagrand~M.
Regularity of Gaussian processes
Acta Math.-- 1987. -- No.~159. -- Pp. 99-149.

\bibitem{Tegza2021}
Tegza A. M.
Simulation of {Gaussian} stationary random process with bounded spectrum with a given accuracy and reliability in the uniform metric.
Nauk. Visn. Uzhgorod. Univ., Ser. Mat. Inform. -  2021. --  v. 39, n. 2, -- p. 91- 99.



\bibitem{Tegza2022}
Tegza A. M., Slyvka-Tylyshchak G. I., Gerych M. S., Pogorilyak O. O., Boyaryshcheva T. V.
 Modeling of {Gaussian} stationary random process with unbounded spectrum using the theory of $L_2(\Omega)$-processes.
Nauk. Visn. Uzhgorod. Univ., Ser. Mat. Inform. -  2022. --  v. 40, n. 1, -- p. 75- 81.


\bibitem{trikomi-1960}
Tricomi~F.~G.
Integral Equations.
Pure and Applied Mathematics, Vol. 5. New York: Interscience Publishers.- 1957.-- 238 p.


\bibitem{vaha85}
Vakhaniya~N.~N., Tarieladze V.~I., Chobanyan S.~A.
Probability distributions in Banach spaces.
-- M. "Nauka", 1985. -- 368 p.


\bibitem{vasyl2005}
Vasylyk~O.~I., Kozachenko Yu.~V., Yamnenko R.~E.
Upper estimate of overrunning by $\text{Sub}_{\varphi}(\Omega)$ random process the level specified by continuous function.
Random Oper. Stoch. Equ. -  2005. -  v. 13,  n. 2, -- p. 111- 128.



\bibitem{vasyl2008}
Vasylyk~O.~I., Kozachenko Yu.~V., Yamnenko R.~E.
$\varphi$- sub-Gaussian random processes.
 - Ky\"\i v: Vydavnycho-Poligrafichny\u\i\ Tsentr, Ky\"\i vsky\u\i\  Universytet, 2008. -- 231 p.

\bibitem{vasyl2019}
Vasylyk O. I., Hopkalo O. M., Kozachenko Yu. V., Sakhno L. M.
Some properties and estimates for $(\varphi)$-sub-{Gaussian} stochastic processes.
Visn., Ser. Fiz.-Mat. Nauky, Ky{\"{\i}}v. Univ. Im. Tarasa Shevchenka. -- 2019, n. 4, -- p. 18--22.

\bibitem{vasyl2021}
Vasylyk O. I., Rozora I. V., Ianevych T. O., Lovyts'ka, I. I.
On some method on model construction for strictly $(\varphi)$-sub-{Gaussian} generalized fractional {Brownian} motion.
Visn., Ser. Fiz.-Mat. Nauky, Ky{\"{\i}}v. Univ. Im. Tarasa Shevchenka. 2021. --
n. 2. --  p. 18- 25.


\bibitem{vyzhvayadrenko2000}
Vyzhva~Z.~O., Yadrenko M.~I.
Statistical simulation of isotropic random fields on a sphere.
Visn., Mat. Mekh., Ky\"\i v. Univ. Im. Tarasa Shevchenka, 2000, No. 5, -- p. 5-11.

\bibitem{vyzhva2003}
Vyzhva~Z.~O.
About approximation of 3-D random fields and statistical simulation.
Random Operators and Stochastic Equations. -- 2003. -- Vol. 11, No. 3. -- p. 255-266.

\bibitem{vyzhva2014}
Vyzhva Z. O., Fedorenko K. V.
 Statistical modelling of a 3D random field by using the {Kotelnikov}-{Shannon} decomposition.
Theory Probab. Math. Stat. -2014. - v. 88, -- p. 19-34.

\bibitem{vyzhva2019}
Vyzhva Z.O.
Statistical modeling in geology.
Ky\"\i v: Vydavnycho-Poligrafichny\u\i\ Tsentr, Ky\"\i vsky\u\i\  Universytet, -- 2019. -- 395 p.


\bibitem{webe81}
Weber~M.
Analyse infinitesimale de fonctions aleatorie.
Ecole d'Ete de Probabilites de St-Flour, 1981,
Lecture Notes in  Mathematics. -- Berlin -- Heidelberg: Springer Verlag, 1983. -- Vol.~976. -- p. 381-465.


\bibitem{jadr67}
Yadrenko~M.~I.
Local properties of sample functions of random fields.
Visn., Mat. Mekh., Ky\"\i v. Univ.  -- 1967. -- No.~9. -- P. 103-112.


\bibitem{jadr68}
Yadrenko~M.~I.
On continuity of the sample functions of the Gaussian random fields on the Hilbert space.
  Dopov. Akad. Nauk Ukr. RSR -- 1968. -- P. 734-737.


\bibitem{jadr80}
Yadrenko~M.~I.
Spectral theory of stochastic fields.
 -- Ky\"\i v: ``Vyshcha Shkola'', 1980. -- 208 p.

\bibitem{jadr83}
Yadrenko~M.~I.
Spectral theory of random fields.
Translation {Series} in {Mathematics} and {Engineering}. {New} {York}: {Optimization} {Software}, {Inc}., {Publications} {Division}; {New} {York}-{Heidelberg}-{Berlin}: {Springer}-{Verlag}. {III}, 1983. -- 259 p.


\bibitem{yamn2014}
Yamnenko R., Kozachenko Yu., Bushmitch D.
Generalized sub-{Gaussian} fractional {Brownian} motion queueing model.
Queueing Syst. -  2014. -  v. 77,  n. 1. -  p. 75- 96.

\bibitem{yamn2020}
Yamnenko R., Yurchenko N.
On an estimate of probability of exceeding a line by weighted aggregate of sub-{Gaussian} random process.
Nauk. Visn. Uzhgorod. Univ., Ser. Mat. Inform. - 2020. -- v. 37,  n. 2. -  p. 122- 129.

\bibitem{yako04}
Yakovenko~T. O.
Condition under which processes belong to Orlicz
spaces in case of noncompact parametric set.
Theory of Stochastic Processes.
-- 2004. -- Vol. 10, no.~26, issue 1-2. -- p. 178-183.

\bibitem{yuri95}
Yurinsky~V.
Sums and Gaussian vectors.
Lecture  Notes in Mathematics, Vol. 1617 -- Berlin: Springer, 1995. -- 315 p.

\bibitem{Zatula2013}
Zatula, D. V.
Modules of continuity of random processes from {Orlicz} spaces of random variables, defined on the interval.
Visn., Ser. Fiz.-Mat. Nauky, Ky{\"{\i}}v. Univ. Im. Tarasa Shevchenka. - 2013. --
n. 2, -- p. 23- 28.

\bibitem{Zatula2014}
Zatula D. V., Kozachenko, Yu. V.
Lipschitz conditions for random processes from {Banach} spaces $\mathbb{F}_{\psi}(\Omega)$ of random variables.
Dopov. Nats. Akad. Nauk Ukr., Mat. Pryr. Tekh. Nauky. -  2014. --
n. 9. -- p. 19- 24.


\bibitem{Zatula2015}
Zatula D. V., Kozachenko, Yu. V.
Lipschitz conditions for stochastic processes in the {Banach} spaces $\mathbb{F}_{{\psi}} ({{\Omega}})$ of random variables.
Theory Probab. Math. Stat. -  2015. -  vol. 91. -  p.  43- 60.

\bibitem{zelepugkoz1982}
Zelepugina~I.~P., Kozachenko Yu. V.
On the question of the simulation of {Gaussian} stochastic processes.
Some questions of the theory of stochastic processes. {Collect}. sci. {Works} { Ky\"\i v}. 1982. -- P. 47-56.


\bibitem{zelepugkoz1988}
Zelepugina~I.~P., Kozachenko Yu. V.
On accuracy estimations in modelling random fields in spaces $L_p$.
Issled. Oper. ASU. -- 1988.  -  No.32.  -  p. 10-14.

\bibitem{zygm65}
Zygmund~A.
Trigonometric series. Volumes I and II combined. 3rd ed. Cambridge: Cambridge University Press -- 2002. -- 364p.





\end{thebibliography}
\end{document}